\documentclass[12pt]{amsart}
\usepackage{amssymb,amscd,amsthm,amsmath,color}
\usepackage{fullpage}
\usepackage{epsfig}
\usepackage{graphicx}
\usepackage{xypic}
\usepackage{url}
\usepackage{color}
\usepackage{hyperref}
\usepackage{enumerate}
\usepackage{cjhebrew}

\numberwithin{equation}{section}

\newcommand{\OO}{\mathcal{O}}

\newcommand{\FF}{\mathcal{F}}

\newcommand{\ev}{\operatorname{ev}}

\newcommand{\Spec}{\operatorname{Spec}}

\newcommand{\rk}{\operatorname{rk}}
\newcommand{\cok}{\operatorname{cok}}

\newcommand{\Eff}{\overline{\operatorname{Eff}}}
\newcommand{\Nef}{\operatorname{Nef}}
\newcommand{\Supp}{\operatorname{Supp}}
\newcommand{\Mor}{\operatorname{Mor}}

\newcommand{\Aut}{\operatorname{Aut}}
\newcommand{\Gal}{\operatorname{Gal}}

\newcommand{\Sym}{\operatorname{Sym}}

\newcommand{\SP}{\operatorname{SP}}

\newcommand{\cU}{\mathcal{U}}

\newcommand{\Hilb}{\operatorname{Hilb}}
\newcommand{\Sec}{\operatorname{Sec}}
\newcommand{\vol}{\operatorname{vol}}

\newcommand{\Br}{\mathrm{Br}}

\newtheorem{theorem}{Theorem}[section]
\newtheorem{lemma}[theorem]{Lemma}
\newtheorem{proposition}[theorem]{Proposition}
\newtheorem{corollary}[theorem]{Corollary}

\theoremstyle{definition}
\newtheorem{definition}[theorem]{Definition}
\newtheorem{conjecture}[theorem]{Conjecture}
\newtheorem{remark}[theorem]{Remark}

\newtheorem{example}[theorem]{Example}
\newtheorem{construction}[theorem]{Construction}

\newtheorem{claim}[theorem]{Claim}
\newtheorem{simp}{Batyrev's heuristic}
\newtheorem{guid}{Geometric Manin's Conjecture}

\xyoption{all}

\begin{document}

\title{Non-free sections of Fano fibrations}
\author{Brian Lehmann}
\address{Department of Mathematics \\
Boston College  \\
Chestnut Hill, MA \, \, 02467}
\email{lehmannb@bc.edu}

\author{Eric Riedl}
\address{Department of Mathematics \\
University of Notre Dame  \\
255 Hurley Hall \\
Notre Dame, IN 46556}
\email{eriedl@nd.edu}

\author{Sho Tanimoto}
\address{Graduate School of Mathematics, Nagoya University, Furocho Chikusa-ku, Nagoya, 464-8602, Japan}
\email{sho.tanimoto@math.nagoya-u.ac.jp}


\begin{abstract}
Let $B$ be a smooth projective curve and let $\pi: \mathcal{X} \to B$ be a smooth integral model of a geometrically integral Fano variety over $K(B)$.  Geometric Manin's Conjecture predicts the structure of the irreducible components $M \subset \Sec(\mathcal{X}/B)$ which parametrize non-relatively free sections of sufficiently large anticanonical degree.  Over the complex numbers, we prove that for any such component $M$ the sections come from morphisms $f: \mathcal{Y} \to \mathcal{X}$ such that the generic fiber of $\mathcal{Y}$ has Fujita invariant $\geq 1$.  Furthermore, we prove that there is a bounded family of morphisms $f$ which together account for all such components $M$.  These results verify the first part of Batyrev's heuristics for Geometric Manin's Conjecture over $\mathbb{C}$.  Our result has ramifications for Manin's Conjecture over global function fields: if we start with a Fano fibration over a number field and reduce mod $p$, we obtain upper bounds of the desired form by first letting the prime go to infinity, then the height.
\end{abstract}

\maketitle

\tableofcontents

%\sho{We should not use full page as that's what Memoirs does.}  \brian{I agree, but this should be a last step since we will need to do some typesetting.  I prefer the full page version for our personal use.}

\section{Introduction}

Since Mori's groundbreaking work in \cite{Mori79} and \cite{Mori82} the moduli space of curves has played a central role in the analysis of Fano varieties.  The irreducible components of the moduli space that parametrize free curves -- that is, curves for which the restriction of the tangent bundle is sufficiently positive -- are generically smooth and have other desirable geometric properties.  By contrast, the irreducible components that only parametrize non-free curves frequently exhibit pathological behavior.
Our goal is to classify these ``non-free components'' for Fano varieties over $\mathbb{C}$.

We show that the non-free components that parametrize curves of sufficiently large degree must come from morphisms which increase the Fujita invariant.  We call such morphisms ``accumulating maps''; their geometry is strongly constrained by the Minimal Model Program.  Furthermore, we show that all non-free components can be accounted for by a bounded family of accumulating maps.  The analogous statements are still true for the moduli space of sections of a $\mathbb{C}$-Fano fibration over a curve and we will work in this more general setting for the rest of the paper.

Our approach to this problem is motivated by arithmetic geometry.  In \cite{Bat88} Batyrev developed a heuristic argument for Manin's Conjecture over a global function field based on some assumptions about the geometry of the space of curves on an $\mathbb{F}_{q}$-Fano variety. Experts anticipated that some aspects of these heuristics should hold for curves over geometric fields (\cite{Manin95}, \cite{Bourqui09}, \cite{Bourqui11}) and due to many groundbreaking advances this expectation has been validated for some particular examples  (\cite{Peyre04}, \cite{Bourqui09}, \cite{Bourqui12}, \cite{CLL16}, \cite{Bilu}). 

Geometric Manin's Conjecture adapts Batyrev's assumptions into a precise set of conjectures about the structure of the moduli space of sections on a $k$-Fano fibration over a curve for arbitrary fields $k$.  Our work completely resolves the first prediction of Geometric Manin's Conjecture for a $\mathbb{C}$-Fano fibration over a curve: non-free components come from accumulating maps.  Our results provide evidence for Batyrev's heuristics in characteristic $p$, and in some circumstances we can deduce an arithmetic statement: if we start with a Fano fibration over a number field and reduce mod $p$, we obtain upper bounds on the counting function in Manin's Conjecture by first letting the prime go to infinity, then the height.

\subsection{Fano fibrations}
Let $B$ be a smooth irreducible projective curve over a field $k$ and let $\eta$ denote its generic point.  A Fano fibration over $B$ is a flat $k$-morphism $\pi: \mathcal{X} \to B$ from a smooth projective variety $\mathcal{X}$ whose generic fiber $\mathcal{X}_{\eta}$ is a geometrically integral Fano variety over $K(B)$.  We will denote by $\Sec(\mathcal{X}/B)$ the moduli space of sections of $\pi$.

The following definition identifies the analogue of a free curve in the setting of fibrations.

\begin{definition}
A section $C$ of $\pi: \mathcal{X} \to B$ is relatively free if $T_{\mathcal{X}/B}|_{C}$ is globally generated and $H^{1}(C,T_{\mathcal{X}/B}|_{C}) = 0$.
\end{definition}

\subsection{Geometric Manin's Conjecture} \label{sect:introgmc}
Geometric Manin's Conjecture is based on an influential heuristic for Manin's Conjecture developed by Baytrev (\cite{Bat88}).  The main invariant in Batyrev's heuristic is the Fujita invariant.

\begin{definition}
\label{defi:a-invariant}
Let $X$ be a smooth projective variety over a field of characteristic $0$ and let $L$ be a big and nef $\mathbb{Q}$-Cartier divisor on $X$.  The Fujita invariant of $(X,L)$ is
\begin{equation*}
a(X,L) = \min \{ t\in \mathbb{R} \mid  K_X + tL \textrm{ is pseudo-effective }\}.
\end{equation*}
If $L$ is nef but not big, we formally set $a(X,L) = \infty$.  If $X$ is singular, choose a resolution of singularities $\phi: X' \to X$ and define $a(X,L)$ to be $a(X',\phi^{*}L)$.  (The choice of resolution does not affect the value by \cite[Proposition 2.7]{HTT15}.)
\end{definition}

Geometric Manin's Conjecture (based on \cite{Bourqui09}, \cite{EVW16}, \cite{LT19}, and \cite{LST18}) predicts that sections of Fano fibrations over any ground field are governed by two key principles:
\begin{itemize}
\item (Exceptional set)  ``Pathological'' families of sections are controlled by the Fujita invariant.
\item (Stability) ``Non-pathological'' families of sections exhibit homological or motivic stability as the degree increases.
\end{itemize}
Our main theorems establish a precise version of the first principle over the complex numbers: the Fujita invariant controls the failure of relative freeness.  In accordance with the asymptotic nature of Manin's Conjecture, our results apply to any family of sufficiently large degree.

\subsection{Main results}
Suppose $\pi: \mathcal{X} \to B$ is a Fano fibration over $\mathbb{C}$.  Our first main result, Theorem \ref{theo:maintheorem1}, shows that every non-relatively free section of sufficiently large degree comes from an accumulating map which does not decrease the Fujita invariant along the generic fiber.  This verifies a conjecture of \cite{LST18} and generalizes earlier results for del Pezzo surface fibrations (\cite{LT21a}, \cite{LT21b}) to fibrations of arbitrary dimension using completely different techniques.

Theorem \ref{theo:maintheorem1} has two cases which correspond to the two ways in which a section $C$ could fail to be relatively free.  First, the deformations of $C$ could fail to dominate $\mathcal{X}$, in which case the sections sweep out a subvariety $\mathcal{Y} \subsetneq \mathcal{X}$.  Second, the deformations of $C$ could dominate $\mathcal{X}$ but $T_{\mathcal{X}/B}|_{C}$ could have a low degree quotient.  In the latter case we expect $C$ to deform more in some directions than in others so that there is a subvariety of $\mathcal{X}$ swept out by the ``most positive'' deformations of $C$.  In fact, there is a rational map on a generically finite cover of $\mathcal{X}$ such that most deformations of $C$ are tangent to the fibers of the map; our subvariety is swept out by the fibers meeting $C$. In both cases Theorem \ref{theo:maintheorem1} shows that the relevant subvariety must have a large Fujita invariant along the generic fiber.

\begin{theorem} \label{theo:maintheorem1}
Let $\pi: \mathcal{X} \to B$ be a Fano fibration.  There is a constant $\xi = \xi(\pi)$ with the following properties.  Let $M$ be an irreducible component of $\Sec(\mathcal{X}/B)$ parametrizing a family of non-relatively free sections $C$ which satisfy $-K_{\mathcal{X}/B} \cdot C \geq \xi$.  Let $\mathcal{U}^{\nu}$ denote the normalization of the universal family over $M$ and let $ev: \mathcal{U}^{\nu} \to \mathcal{X}$ denote the evaluation map.  Then either: 
\begin{enumerate}
\item $ev$ is not dominant.  Then the subvariety $\mathcal{Y}$ swept out by the sections parametrized by $M$ satisfies
\begin{equation*}
a(\mathcal{Y}_{\eta},-K_{\mathcal{X}/B}|_{\mathcal{Y}_{\eta}}) \geq a(\mathcal{X}_{\eta},-K_{\mathcal{X}/B}|_{\mathcal X_\eta}).
\end{equation*}
\item $ev$ is dominant.  Letting $f: \mathcal{Y} \to \mathcal{X}$ denote the finite part of the Stein factorization of $ev$, we have
\begin{equation*}
a(\mathcal{Y}_{\eta},-f^{*}K_{\mathcal{X}/B}|_{\mathcal{Y}_{\eta}}) = a(\mathcal{X}_{\eta},-K_{\mathcal{X}/B}|_{\mathcal X_\eta}).
\end{equation*}
Furthermore, there is a dominant rational map $\phi: \mathcal{Y} \dashrightarrow \mathcal{Z}$ over $B$ with connected fibers such that the dimension of $\mathcal Z$ is at least $2$ and  the following properties hold. Let $C'$ denote a general section of $\mathcal{Y} \to B$ parametrized by $M$ and let $\mathcal{W}' \subset \mathcal{Y}$ denote the unique irreducible component of the closure of $\phi^{-1}(\phi(C'))$ which maps dominantly to $\phi(C')$.  There is a resolution $\psi: \mathcal{W} \to \mathcal{W}'$ such that the locus where $\psi^{-1}$ is well-defined intersects $C'$ and $\psi$ has the following properties.  
\begin{enumerate}
\item  We have $a(\mathcal{W}_{\eta},-\psi^{*}f^{*}K_{\mathcal{X}/B}|_{\mathcal{W}_{\eta}}) = a(\mathcal{X}_{\eta},-K_{\mathcal{X}/B}|_{\mathcal X_\eta})$.
\item The Iitaka dimension of $K_{\mathcal{W}_{\eta}} - a(\mathcal{W}_{\eta},-\psi^{*}f^{*}K_{\mathcal{X}/B}|_{\mathcal{W}_{\eta}})\psi^{*}f^{*}K_{\mathcal{X}/B}|_{\mathcal{W}_{\eta}}$ is $0$.
\item The general deformation of the strict transform of $C'$ in $\mathcal W$ is relatively free in $\mathcal W$.
\item There is a constant $T = T(\pi)$ depending only on $\pi$, but not $M$, such that the sublocus of $M$ parametrizing deformations of the strict transform of $C'$ in $\mathcal{W}$ has codimension at most $T$ in $M$.
\end{enumerate}
\end{enumerate}
\end{theorem}

\begin{remark} \label{rema:introexplain}
Since Theorem \ref{theo:maintheorem1}.(2) is somewhat technical, we give a brief explanation here.  First, since we would like to remain in the setting of varieties with generically finite morphisms to $\mathcal{X}$, we always insist on passing to the finite part $\mathcal{Y}$ of the Stein factorization of the evaluation map.  In fact, since Stein factorizations are only birationally invariant for normal projective varieties, we also include a normalization of the universal family when defining $\mathcal{Y}$.

Conceptually, the sections $C'$ ``move a lot'' in the direction of the fibers of the rational map $\phi: \mathcal{Y} \dashrightarrow \mathcal{Z}$.  The cleanest way to express this fact is to look at the variety $\mathcal{W}$ which resolves the union of the fibers meeting $C'$.  Claim (d) expresses in what sense the sections ``move a lot'' in $\mathcal{W}$: the space of deformations in $\mathcal{W}$ have bounded codimension in the space of deformations in $\mathcal{Y}$.  Claim (a) verifies that this property is reflected by the Fujita invariant of $W_{\eta}$.  The birational geometry is further clarified in claim (b) -- loosely speaking, we expect curves to ``move a lot'' along the fibers of the Iitaka fibration, and thus the fact that the Iitaka dimension is $0$ verifies that we have indeed identified the correct construction.  Finally claim (c) is similar to claim (d); the sections move so much in $\mathcal{W}$ that they are relatively free.  
\end{remark}

\begin{remark}
Suppose $X$ is a Fano variety and we are studying irreducible components of $\Mor(B,X)$.  Any family of non-free curves on $X$ leads to a family of non-relatively free sections of $\pi: X \times B \to B$ and thus Theorem \ref{theo:maintheorem1} gives a classification result for such families.  However, it is natural to ask whether non-free curves can be described using generically finite maps $f: Y \to X$ instead of generically finite maps $f: \mathcal{Y} \to X \times B$.  The answer is ``yes,'' but due to length constraints we will give the details in a supplementary paper (\cite{LRT22}).
\end{remark}

Theorem \ref{theo:maintheorem1} has two key consequences.  First, the Fujita invariant can be computed using tools from the Minimal Model Program and thus Theorem \ref{theo:maintheorem1} gives a practical way to classify families of non-relatively free sections.

\begin{example}
In Example \ref{exam:cubichyp} we analyze the moduli spaces $\Sec(\mathcal{X}/B)$ when $\pi: \mathcal{X} \to B$ is a cubic hypersurface fibration and $B$ is a curve of arbitrary genus.  When $\dim(\mathcal{X}_{\eta}) \geq 5$, we show that the ``exceptional set'' is empty so that every component of $\Sec(\mathcal{X}/B)$ of sufficiently high degree will generically parametrize relatively free sections.  For dimension $4$ a similar analysis allows us to describe the families of non-relatively free sections of large degree.
\end{example}

Second, \cite{Birkar21} imposes strong finiteness constraints on Fujita invariants and Theorem \ref{theo:maintheorem1} allows us to deduce finiteness results for families of non-relatively free sections.  Recall that Theorem \ref{theo:maintheorem1} shows that families of non-relatively free sections come from maps $f: \mathcal{Y} \to \mathcal{X}$ such that the Fujita invariant of $\mathcal{Y}_{\eta}$ is at least $a(\mathcal{X}_{\eta},-K_{\mathcal{X}/B}|_{\mathcal{X}_{\eta}}) = 1$.  If we pass to an algebraic closure $\overline{K(B)}$ then the set of maps $f_{\overline{\eta}}: \mathcal{Y}_{\overline{\eta}} \to \mathcal{X}_{\overline{\eta}}$ such that $a(\mathcal{Y}_{\overline{\eta}},-f_{\overline{\eta}}^{*}K_{\mathcal{X}_{\overline{\eta}}}) \geq 1$ satisfy certain types of boundedness (see e.g.~\cite[Theorem 1.7]{LST18}).  However, the analogous boundedness statements over $K(B)$ are no longer true since a map over $\overline{K(B)}$ can correspond to infinite families of twists over $K(B)$.

In Section \ref{sect:twists} we systematically study the set of twists of the map $f_{\eta}: \mathcal{Y}_{\eta} \to \mathcal{X}_{\eta}$.  Theorem \ref{theo:introtwists} shows that amongst all the twists only a bounded subfamily carry a family of sections which is dense in an irreducible component of $\Sec(\mathcal{X}/B)$.  In this way, we conclude that all non-relatively free sections will come from a bounded family of accumulating maps $f: \mathcal{Y} \to \mathcal{X}$.

\begin{theorem} \label{theo:maintheorem2}
Let $\pi: \mathcal{X} \to B$ be a Fano fibration.
\begin{enumerate}
\item There is a proper closed subset $\mathcal{R} \subsetneq \mathcal{X}$ such that if $M \subset \Sec(\mathcal{X}/B)$ is an  irreducible component parametrizing a non-dominant family of sections then the sections parametrized by $M$ are contained in $\mathcal{R}$.
\item There is a constant $\xi = \xi(\pi)$, a proper closed subset $\mathcal{V} \subset \mathcal{X}$, and a bounded family of smooth projective $B$-varieties $\mathcal{Y}_{s}$, $s \in \mathfrak{S}$ equipped with $B$-morphisms $f_{s}: \mathcal{Y}_{s} \to \mathcal{X}$ satisfying:
\begin{enumerate}
\item $\dim(\mathcal{Y}_{s}) < \dim(\mathcal{X})$ and $f_{s}$ is generically finite onto its image; 
\item $a(\mathcal Y_{s,\eta}, -f_{s}^*K_{\mathcal{X}/B}|_{\mathcal{Y}_{s,\eta}}) = a(\mathcal{X}_{\eta},-K_{\mathcal{X}/B}|_{\mathcal X_\eta})$ and the Iitaka dimension of $K_{\mathcal{Y}_{s,\eta}} -f_{s}^*K_{\mathcal X/B}|_{\mathcal{Y}_{s,\eta}}$ is $0$;
\item if $M \subset \Sec(\mathcal{X}/B)$ is an irreducible component  that generically parametrizes non-relatively free sections $C$ with $-K_{\mathcal{X}/B} \cdot C \geq \xi$ then for a general section $C$ parametrized by $M$ we have either
\begin{enumerate}
\item $C \subset \mathcal{V}$, or
\item for some $f_{s}: \mathcal{Y}_{s} \to \mathcal{X}$ in our family there is a relatively free section $C'$ of $\mathcal{Y}_{s}/B$ such that $C = f_{s}(C')$.
\end{enumerate}
\end{enumerate}
\end{enumerate}
\end{theorem}

\begin{remark}
Theorem \ref{theo:maintheorem2} establishes geometric analogues of various conjectures about the exceptional set in Manin's Conjecture.  For example, suppose that $B$ is a smooth projective $\mathbb{F}_{q}$-curve and $\pi: \mathcal{X} \to B$ is a Fano fibration equipped with an adelic metrization on the relative canonical bundle.  Weak Manin's Conjecture predicts that  there exist a constant $C>0$ and a closed subset $\mathcal{R} \subset \mathcal{X}$ such that for any $\epsilon > 0$ the number of sections meeting $\mathcal{X} \backslash \mathcal{R}$ of height at most $d$ is bounded above by $Cq^{d(1+\epsilon)}$.  Using the heuristic estimate $\# M(\mathbb{F}_{q}) \approx q^{\dim M}$, this means that $\mathcal{R}$ should contain all sections parametrized by a family $M$ such that $\dim(M)/\mathrm{expdim}(M) \geq 1 + \epsilon$ (with perhaps finitely many exceptions).  Theorem \ref{theo:maintheorem2}.(1) shows that over $\mathbb{C}$ there exists a closed set $\mathcal{R}$ with this property.   
\end{remark}

\subsection{A geometric application}
Suppose $M$ is an irreducible component of $\Sec(\mathcal{X}/B)$ and let $N \subset M$ denote the sublocus parametrizing sections $C$ such that $T_{\mathcal{X}/B}|_{C}$ is not generically globally generated.  One would like to find a lower bound on the codimension of $N$.  This problem has been previously studied when $X$ is a smooth Fano variety and $M \subset \Mor(\mathbb{P}^{1},X)$, in which case $N \subset M$ is simply the non-free locus.  For example, \cite[Theorem 1.2]{BS18} shows that when $X$ is a smooth hypersurface whose dimension is much larger than the degree and $M$ is an irreducible component of $\Mor(\mathbb{P}^{1},X)$ then the codimension of $N \subset M$ grows linearly in the anticanonical degree of the curves parametrized by $M$.

We prove the first general statement for arbitrary Fano fibrations: the codimension of the non-generically-globally-generated locus grows linearly in the degree unless there is a clear geometric reason why it cannot.

\begin{theorem} \label{theo:maintheorem3}
Let $\pi: \mathcal{X} \to B$ be a Fano fibration.  There is a linear function $Q(d)$ whose leading coefficient is a positive number depending only on $\dim(\mathcal{X})$ such that the following property holds.

Suppose that $M$ is an irreducible component of $\Sec(\mathcal{X}/B)$ parametrizing a family of sections $C$ which satisfy $-K_{\mathcal{X}/B} \cdot C = d$.  Let $N \subset M$ be a subvariety parametrizing sections $C$ such that $T_{\mathcal{X}/B}|_{C}$ is not generically globally generated.  Then either
\begin{enumerate}
\item the codimension of $N$ in $M$ is at least $\sup \{ \lfloor Q(d) \rfloor, 0 \}$, or
\item the sections parametrized by $N$ sweep out a subvariety $\mathcal{Y} \subsetneq \mathcal{X}$ satisfying
\begin{equation*}
a(\mathcal{Y}_{\eta},-K_{\mathcal{X}/B}|_{\mathcal{Y}_{\eta}}) \geq a(\mathcal{X}_{\eta},-K_{\mathcal{X}/B}|_{\mathcal X_\eta}).
\end{equation*}
\end{enumerate}
\end{theorem}

In Example \ref{exam:largenonfreelocus} we will show that the codimension of the non-generically-globally-generated locus can be constant as the degree increases, demonstrating that case (2) of Theorem \ref{theo:maintheorem3} must be included.

\begin{remark}
Over $\mathbb{C}$, Theorem \ref{theo:maintheorem3} shows that families of sections such that $T_{\mathcal{X}/B}|_{C}$ is not generically globally generated are either contained in the exceptional locus or they have large codimension in a component of $\Sec(\mathcal{X}/B)$.  Assuming the analogous statement over a global function field, we can expect such sections to make a negligible contribution to the counting function for Manin's Conjecture.  Thus Theorem \ref{theo:maintheorem3} supports the novel formulation of Manin's Conjecture due to \cite{Peyre17} which only counts rational points which are ``free'' in a suitable sense. 
\end{remark}

\subsection{An arithmetic application}
\label{subsec:arithmeticapp}
Our results can be applied to prove an upper bound of Manin type over a global function field.  Let $F$ be a number field and let $B$ be a smooth projective curve over $F$. Let $S$ be a finite set of places of $F$ including all archimedean places and let $\mathfrak o_{F, S}$ be the ring of $S$-integers in $F$.

Let $\pi : X \to B$ be a Fano fibration defined over $F$.  After perhaps enlarging $S$, we can find an integral model $\widetilde{\pi} : \mathcal X \to \mathcal B$ of $\pi$ over $\mathfrak o_{F, S}$ such that $\mathcal{X}$ and $\mathcal{B}$ are smooth over $\mathfrak o_{F, S}$. 
Let $R \subset X$ be the Zariski closure of the union of the loci swept out by non-dominant families of sections in $\mathrm{Sec}(X/B)$. Theorem~\ref{theo:maintheorem2}.(1) implies that the base change of $R$ to $\mathbb C$ is contained in a proper closed subset, so in particular $R$ itself is a proper closed subset. We consider the flat closure $\mathcal R \subset \mathcal X$ of $R$.

Let $v$ be a non-archimedean place of $F$ not contained in $S$ and consider the reduction $\pi_v : X_v \to B_v$ at $v$ which is defined over a finite field $k_v$. Let $R_{v}$ be the reduction of $\mathcal R$ at $v$ and let $\mathrm{Sec}(X_v/B_v, R_v)_{\leq d}$ be the open subset of $\mathrm{Sec}(X_v/B_v)$ parametrizing sections $C \not\subset R_v$ of anticanonical height $\leq d$. Then we consider the following counting function:
\[
N(X_v\setminus R_v, -K_{X_v/B_v}, d) = \# \mathrm{Sec}(X_v/B_v, R_v)_{\leq d}(k_v).
\]
Weak Manin's Conjecture over $K(B_v)$ predicts that for any $\epsilon > 0$ we have
\[
N(X_v\setminus R_v, -K_{X_v/B_v}, d) = o(q_v^{d(1+\epsilon)}),
\]
as $d \to \infty$ where $q_v = \#k_v$. 
The following approximation of this conjecture was suggested to us by Jordan Ellenberg and Melanie Matchett Wood:

\begin{theorem}
\label{theorem:arithmeticapp}
Let $F, S, \widetilde{\pi} : \mathcal X \to \mathcal B$ be as above.  Then assuming $d\epsilon > \dim X_\eta$, we have
 \[
 \frac{ N(X_v\setminus R_v, -K_{X_v/B_v}, d)}{q_v^{d(1+\epsilon)}} \to 0
 \]
 as $v \to \infty$. 
\end{theorem}

This result fits into the recent trend of taking hard arithmetic questions that are asymptotic in a different parameter and making them more accessible by first letting $v$ go to $\infty$.  This technique has been explored in the contexts of Malle's Conjecture and Cohen-Lenstra heuristics over global function fields; see e.g.~\cite{Ach06, EVW16, FLR22, PW21, LWZB19, ETW17}.

\subsection{Stability and restriction}
The proof of our main results requires a number of statements which are interesting in their own right.
We highlight one such result here concerning the behavior of vector bundles under restrictions.  We will need the notions of slope stability and Harder-Narasimhan filtrations of torsion-free sheaves with respect to nef curve classes as developed by \cite{CP11}.

Suppose $\mathcal{E}$ is a torsion-free sheaf and $M$ is an irreducible component of $\Sec(\mathcal{X}/B)$ parametrizing a dominant family of sections.  We would like to understand the Harder-Narasimhan filtration of $\mathcal{E}|_{C}$ for a general section $C$ parametrized by $M$.  Intuitively, we can hope that the stability of the restricted bundle depends on the global stability properties of $\mathcal{E}$ with respect to the class $[C]$.

In Section \ref{sect:gm} we prove that when $\mathcal{E}$ is a torsion-free sheaf on $\mathcal{X}$ that is semistable with respect to the numerical class $[C]$ of a flat family of curves then $\mathcal{E}|_C$ is ``almost'' semistable. This result should be compared against the Mehta-Ramanthan theorem for complete intersection curves: although the conclusion is a little weaker, the result applies to a much broader class of curves than just complete intersection curves. The following statement is a special case designed for families of non-relatively free sections.

\begin{theorem} \label{theo:introgm}
Let $\pi: \mathcal{X} \to B$ be a flat morphism with connected fibers from a smooth projective variety $\mathcal{X}$ to a smooth projective curve $B$ and let $\mathcal{E}$ be a torsion-free sheaf on $\mathcal X$.  Let $M$ be an irreducible component of $\Sec(\mathcal{X}/B)$ and let $\mathcal{U}^{\nu}$ denote the normalization of the universal family over $M$.  Suppose that
\begin{enumerate}
\item the evaluation map $ev: \mathcal{U}^{\nu} \to X$ is dominant with connected fibers and
\item for some open subset $M_{red}^{\circ} \subset M_{red}$ the restriction of $ev$ to the preimage of $M_{red}^{\circ}$ is flat.
\end{enumerate}
For a general curve $C$ parametrized by $M$, write
\begin{equation*}
0 = \mathcal{F}_{0} \subset \mathcal{F}_{1} \subset \ldots \subset \mathcal{F}_{r} = \mathcal{E}|_{C}
\end{equation*}
for the Harder-Narasimhan filtration of $\mathcal{E}|_{C}$.

Suppose that $\mathcal{E}$ is $[C]$-semistable.  Then for every index $i$ we have
\begin{equation*}
\left| \mu(\mathcal{E}|_{C}) - \mu(\mathcal{F}_{i}/\mathcal{F}_{i-1}) \right| \leq  (g(B)\dim(\mathcal{X}) - g(B) + 1)^{2} (\rk(\mathcal{E})-1).
\end{equation*}
\end{theorem}

Note that no Fano assumption is necessary in Theorem \ref{theo:introgm}.  Under some weak positivity assumptions on the normal bundle of $C$, the constant occurring in Theorem \ref{theo:introgm} can be dramatically improved.

\begin{remark}
We emphasize that the conditions (1) and (2) in Theorem \ref{theo:introgm} are both necessary and exhaustive.  If these conditions do not hold, then the evaluation map over $M$ factors rationally through a generically finite morphism $f: \mathcal{Y} \to \mathcal{X}$ such that the analogous conditions do hold on $\mathcal{Y}$.  (For example, we can let $f$ be a flattening of the finite part of the Stein factorization of the evaluation map.)  Then Theorem \ref{theo:introgm} shows that the restricted bundle is controlled by the stability of $f^{*}\mathcal{E}$ on $\mathcal{Y}$ and not the stability of $\mathcal{E}$ on $\mathcal{X}$.
\end{remark}

\subsection{Strategy}  
Our strategy for the proof of Theorem \ref{theo:maintheorem1} relies on several tools: the theory of foliations and slope stability, the MMP, and local-to-global principles over function fields of complex curves.  We outline the proof of Theorem \ref{theo:maintheorem1}.(2) where $M$ parametrizes a dominant family of non-relatively free sections on $\mathcal{X}$.  For simplicity we assume that $ev : \mathcal U^\nu \to \mathcal X$ has connected fibers so that the Stein factorization $\mathcal{Y}$ of $ev$ is equal to $\mathcal{X}$.  %We must construct a rational map $\phi$ on $\mathcal{X}$ that captures the geometry of this dominant family of sections.   

Since by hypothesis $C$ is not relatively free, we see that $T_{\mathcal{X}/B}|_{C}$ must have a low slope quotient.  Applying Theorem \ref{theo:introgm} to a birational model flattening the family, we can ``lift'' this quotient to all of $\mathcal{X}$.  The result is a foliation $\mathcal{F} \subset T_{\mathcal{X}}$ of large slope and the pioneering results of \cite{CP19} show that $\mathcal{F}$ is induced by a rational map $\phi$.  We carry out this construction in Section \ref{sect:genpoints}.

It then remains to verify the desired properties of $\phi$.  The most difficult is the computation of the Fujita invariant as in Theorem \ref{theo:maintheorem1}.(2).(b).  By appealing to Birkar's recent boundedness results in the Minimal Model Program,  we prove the following general criterion for computing the Fujita invariant in Section \ref{sec:fujinv}.

\begin{theorem} \label{theo:introainvariant}
Let $\pi: \mathcal{X} \to B$ be a Fano fibration.  Fix a positive rational number $a$ and a positive integer $T$.  There is some constant $\xi = \xi(\pi,a,T)$ with the following property.

Suppose that $\psi: \mathcal{Y} \to B$ is a flat morphism with connected fibers from a smooth projective variety $\mathcal{Y}$ and $f: \mathcal{Y} \to \mathcal{X}$ is a $B$-morphism that is generically finite onto its image.   Suppose that $N$ is an irreducible component of $\Sec(\mathcal{Y}/B)$ parametrizing a dominant family of sections on $\mathcal{Y}$ and let $M$ denote the irreducible component of $\Sec(\mathcal{X}/B)$ containing $f_{*}N$.

Assume that the sections $C$ parametrized by $N$ satisfy $-f^{*}K_{\mathcal{X}/B} \cdot C \geq \xi$ and that
\begin{equation*}
\dim(N) \geq a \cdot \dim(M) - T.
\end{equation*}
Then
\begin{equation*}
a(\mathcal{Y}_{\eta},-f^{*}K_{\mathcal{X}/B}|_{\mathcal Y_\eta}) \geq a.
\end{equation*}
\end{theorem}

\begin{remark}
We prove statements analogous to Theorem \ref{theo:introainvariant} in the more general setting of pairs $(\mathcal{X},L)$ where $\mathcal{X}$ is a smooth projective variety admitting a flat morphism with connected fibers $\pi: \mathcal{X} \to B$ and $L$ is a generically relatively big and semiample Cartier divisor on $\mathcal{X}$.
\end{remark}

Returning to the setting of Theorem \ref{theo:maintheorem1}.(2), we prove Theorem \ref{theo:maintheorem1}.(2).(b) by combining Theorem \ref{theo:introainvariant} with Theorem \ref{theo:maintheorem1}.(2).(d).

We next outline the proof of Theorem \ref{theo:maintheorem2}.(2).  Suppose we have a dominant generically finite morphism $f_{\eta}: \mathcal{Y}_{\eta} \to \mathcal{X}_{\eta}$.  As discussed earlier, the key is to understand how the set of twists $f_{\eta}$ interacts with the behavior of sections on an integral model.  We systematically analyze this relationship in Section \ref{sect:twists}; in particular, we prove the following statement that undergirds Theorem \ref{theo:maintheorem2}.

\begin{theorem} \label{theo:introtwists}
Let $\pi: \mathcal{X} \to B$ be a Fano fibration.  Let $\mathcal{Y}$ be a normal projective variety equipped with a flat morphism with connected fibers $\psi: \mathcal{Y} \to B$ and with a dominant generically finite $B$-morphism $f: \mathcal{Y} \to \mathcal{X}$.

Suppose that $\widetilde{\mathcal{Y}}$ is a $B$-variety which is smooth and projective and $\widetilde{f}: \widetilde{\mathcal{Y}} \to \mathcal{X}$ is a dominant generically finite morphism such that $\widetilde{f}_{\eta}: \widetilde{\mathcal{Y}}_{\eta} \to \mathcal{X}_{\eta}$ is birational to a twist of $f_{\eta}$.  Fix a positive integer $T$ and suppose that there exists an irreducible component $\widetilde{N} \subset \Sec(\widetilde{\mathcal{Y}}/B)$ parametrizing a dominant family of sections on $\widetilde{\mathcal{Y}}$ such that the pushforward of $\widetilde{N}$ has codimension at most $T$ in a component of $\Sec(\mathcal{X}/B)$.

Then there exist constants $d = d(\mathcal Y/\mathcal X)$ and $n = n(\mathcal Y/\mathcal X, T)$ and a finite Galois morphism $B'\to B$ of degree at most $d$ with at most $n$ branch points such that the base changes of $f_{\eta}$ and $\widetilde{f}_{\eta}$ to $K(B')$ are birationally equivalent.
\end{theorem}

This result implies that the set of $\widetilde{\mathcal Y}$ satisfying the conditions of Theorem \ref{theo:introtwists} is birationally bounded. We prove this boundedness by constructing a parameter space of twists which is of finite type over the Hurwitz stack parametrizing finite covers $B' \to B$.

\subsection{History}

Ever since the seminal results due to Mori and his coauthors (\cite{Mori79}, \cite{Mori82}, \cite{MM86}) the moduli space of curves has played a prominent role in the study of Fano varieties.  The notion of free rational curves goes back to pioneering work by Koll\'ar--Miyaoka--Mori (\cite{KMM92}, \cite{Kollar}) on rational connectedness of Fano varieties.  Since then there have been many breakthroughs in the description of the moduli spaces $\Mor(B,X)$ for Fano varieties $X$, most notably when $B = \mathbb P^1$.  One particularly influential example is the analysis of rational curves on Fano hypersurfaces pioneered by \cite{HRS04} and subsequently developed by \cite{CS09}, \cite{BK13}, and \cite{RY19}. (\cite{BV17}, \cite{BS18} provide a different approach to this problem using an idea from analytic number theory.) 
Another important class of examples is the moduli spaces of curves on various homogeneous spaces (\cite{Thomsen98}, \cite{KP01}, \cite{Bourqui16}). However for a long time it was unclear what structure to expect for arbitrary Fano varieties. 

The situation was clarified by the introduction of ideas from arithmetic geometry.  Manin's Conjecture is a conjectural asymptotic formula for the counting function of rational points on Fano varieties formulated and refined in \cite{FMT89}, \cite{BM90}, \cite{Peyre95}, \cite{BT98}, and \cite{LST18}.   In \cite{Peyre17} Peyre proposed another version of Manin's Conjecture using the notion of freeness of rational points which is inspired by the concept of free rational curves.   A motivic version of Manin's Conjecture has been established for toric varieties in \cite{Bourqui09} and for equivariant compactifications of vector groups in \cite{CLL16} and \cite{Bilu}.  In Manin's Conjecture, it is important to exclude the contribution to the counting function from ``exceptional sets'' where rational points accumulate too quickly.  The relationship between exceptional sets and Fujita invariants was developed in \cite{HTT15}, \cite{LTT18}, \cite{HJ16}, \cite{LT17}, \cite{Sengupta21}, \cite{LST18}, and \cite{LTRMS}. These developments culminated in the main theorem of \cite{LST18} proving that the contribution to the exceptional set coming from maps $f: Y \to X$ such that $Y$ has larger $a$ and $b$ invariants will be contained in a thin set of rational points.  This result is a source of our main theorem (Theorem~\ref{theo:maintheorem2}) showing that pathological components come from a bounded family.

In his influential notes (\cite{Bat88}) Batyrev gave a heuristic for the global function field version of Manin's Conjecture.  Over time the principles underlying Batyrev's heuristic were made into precise conjectures and extended to arbitrary ground fields.  First, building upon ideas of Peyre \cite{Bourqui09} proposed a motivic version of Manin's Conjecture over arbitrary ground fields.  Building upon earlier work on homological stability by \cite{Segal79}, \cite{CJS00}, and many others, \cite{EVW16} highlighted the connection between homological stability and rational point counts via the Grothendieck-Lefschetz 
trace formula.  Second, based on the analysis of the exceptional set described earlier, \cite{LT19} predicted the geometry underlying ``pathological'' families of rational curves on Fano varieties and obtained a first prototype result of Theorem~\ref{theo:maintheorem1} in this setting. Further works leveraged this intuition to study rational curves for Fano varieties of dimension $\leq 3$ and for sections of del Pezzo fibrations (\cite{Castravet}, \cite{Testa09}, \cite{LT19}, \cite{LTJAG}, \cite{LT21a}, \cite{LT21b}, \cite{BLRT20}, \cite{MTiB20}, \cite{ST22}, and \cite{BJ22}). Higher dimensional cases also have been explored previously in the case of hypersurfaces and homogeneous spaces as mentioned above and more recently in \cite{Okamura}. This series of work is one of instances for the geometrization of arithmetic conjectures which has been pioneered by many authors, e.g., \cite{GHS03}, \cite{deJongStarr}, \cite{GHMS}, \cite{HT06}, \cite{Xu12a}, \cite{Xu12b}, \cite{Tian15}, \cite{Tian17}, \cite{TZ18}, \cite{TZ19}, \cite{SX20}, and \cite{STZ22}. 

Together, the two principles in Batyrev's heuristic (stated in Section \ref{sect:introgmc}) are known as Geometric Manin's Conjecture.  Geometric Manin's Conjecture unifies many disparate examples and clarifies the conjectural structure of $\Mor(B,X)$ for arbitrary Fano varieties. Theorems \ref{theo:maintheorem1} and \ref{theo:maintheorem2} are the first statements in Geometric Manin's Conjecture which have been proved for arbitrary Fano fibrations over curves of arbitrary genus. This answers a question implicitly raised by Batyrev in his heuristics (\cite{Bat88}).

\

\noindent
{\bf Acknowledgments:}
The authors thank Shintarou Yanagida for a helpful conversation about stacks. The authors thank Jordan Ellenberg and Melanie Matchett Wood for suggesting an arithmetic application of our work and Lars Hesselholt for recommending the reference \cite{BS15}. The authors also thank Izzet Coskun and Zhiyu Tian for comments on this paper.
The authors thank referees for detailed suggestions which significantly improved the exposition of the paper.

Part of this  project was conducted at the SQuaRE workshop ``Geometric Manin's Conjecture in characteristic $p$'' at the American Institute of Mathematics. The authors would like to thank AIM for the generous support. 

Brian Lehmann was supported by Simons Foundation grant Award Number 851129.
Eric Riedl was supported by NSF CAREER grant DMS-1945944.  Sho Tanimoto was partially supported by JST FOREST program Grant number JPMJFR212Z, by JSPS KAKENHI Grand-in-Aid (B) 23K25764, by JSPS Bilateral Joint Research Projects Grant number JPJSBP120219935, and by JSPS KAKENHI Early-Career Scientists Grant number 19K14512.

\section{Background} \label{sect:background}

Throughout all our schemes will be assumed to be separated and every connected component will have finite type over the base ring (which is usually $\mathbb{C}$ or the function field of a complex curve). Recall that in this situation the normalization of a scheme $X$ is isomorphic to the normalization of $X_{red}$.  A variety is a separated integral scheme of finite type over the base field.  Given a coherent sheaf $\mathcal{F}$ on a variety $V$, we denote by $\mathcal{F}_{tors}$ the torsion subsheaf and by $\mathcal{F}_{tf}$ the quotient of $\mathcal{F}$ by its torsion subsheaf.

Given a dominant generically finite morphism of projective varieties $f: Y \to X$, we denote by $\Aut(Y/X)$ the automorphism group of $Y$ over $X$ and by $\mathrm{Bir}(Y/X)$ the birational automorphism group of $Y$ over $X$. 

\begin{definition} 
Suppose we have a dominant morphism of varieties $f: U \to V$ such that the general fiber of $f$ is geometrically irreducible.  Suppose that $T \subset V$ is a subvariety that meets the open locus over which $f$ has geometrically irreducible fibers.  Then $f^{-1}(T)$ has a unique irreducible component which dominates $T$ under $f$.  We call this the ``main component'' of $f^{-1}(T)$.
\end{definition}

When $X$ is a projective variety, we will let $N^{1}(X)_{\mathbb{R}}$ denote the space of $\mathbb{R}$-Cartier divisors up to numerical equivalence.  In this finite-dimensional vector space we have the pseudo-effective cone $\Eff^{1}(X)$ and the nef cone $\Nef^{1}(X)$.  
Dually, we will let $N_{1}(X)_{\mathbb{R}}$ denote the space of $\mathbb{R}$-$1$-cycles up to numerical equivalence.  Inside $N_{1}(X)_{\mathbb{R}}$ we have  the pseudo-effective cone $\Eff_{1}(X)$ and the nef cone $\Nef_{1}(X)$.  Given a curve $C$, we will denote its numerical class by $[C]$.

We also use the standard definitions and techniques from the Minimal Model Program. See \cite{KM98} and \cite{BCHM} for more details.

\begin{definition}
Let $f: X \to Y$ be a projective morphism of varieties and let $L$ be a $\mathbb{Q}$-Cartier divisor on $X$.  For a property $P$ of $\mathbb{Q}$-Cartier divisors (such as ample, big, nef, semiample, etc.), we say that $L$ is generically relatively $P$ if the restriction of $L$ to the generic fiber of $f$ satisfies $P$. 
\end{definition}

\begin{definition}
Let $X$ be a normal projective variety and let $L$ be a pseudo-effective $\mathbb{Q}$-Cartier divisor on $X$.  The stable base locus of $L$ is
\begin{equation*}
\mathbf{B}(L) := \bigcap_{m \in \mathbb{Z}_{>0}, mL\textrm{ Cartier}} \mathrm{Bs}(mL).
\end{equation*}
The augmented base locus is
\begin{equation*}
\mathbf{B}_{+}(L) := \bigcap_{A \textrm{ ample }\mathbb{Q}\textrm{-Cartier}} \mathbf{B}(L-A).
\end{equation*}
\cite[Proposition 1.5]{ELMNP06} verifies that there is some ample $\mathbb{Q}$-Cartier divisor $A$ such that $\mathbf{B}_{+}(L) = \mathbf{B}(L-A)$.  In particular $\mathbf{B}_{+}(L)$ is a closed subset.
\end{definition}

\subsection{Convex geometry}

If $V$ is a finite dimensional $\mathbb{R}$-vector space, we will say that a subset $\mathcal{C} \subset V$ is a (closed convex) polyhedron if it is a finite intersection of closed affine half spaces, i.e.
\begin{equation*}
\mathcal{C} = \bigcap_{i=1}^{t} \{ v \in V | \ell_{i}(v) \leq c_{i} \}
\end{equation*}
where the $\ell_{i}$ are non-zero linear functions on $V$ and $c_{i} \in \mathbb{R}$.  If $V$ carries a rational structure -- that is, we identify $V = V_{\mathbb{Q}} \otimes_{\mathbb{Q}} \mathbb{R}$ for some $\mathbb{Q}$-vector space $V_{\mathbb{Q}}$ -- then we say that $\mathcal{C}$ is rational polyhedral if every $c_{i}$ is in $\mathbb{Q}$ and every $\ell_{i}$ is induced by a $\mathbb{Q}$-linear function on $V_{\mathbb{Q}}$.

Suppose that $\mathcal{C}$ is a (closed, convex, salient, full-dimensional) cone in a finite dimension vector space $V$.  Given an element $\ell \in V^{\vee}$, we denote by $C_{\ell \geq 0}$ the intersection of $\mathcal{C}$ with the half-space $\{v \in V | \ell(v) \geq 0 \}$.  Similarly, $C_{\ell = 0}$ denotes the intersection with the hyperplane $\ell^{\vee}$.

We will need the following result of \cite{HW07}:

\begin{theorem}[{\cite[Theorem 1.1.(a)]{HW07}}]
Let $V$ be a finite dimensional $\mathbb{R}$-vector space equipped with a full-rank discrete lattice $\Lambda$.  Let $S$ be a subset of $\Lambda$.  Suppose $\mathcal{C}$ is a rational polyhedral cone in $V$ and $Q$ is a compact subset of $V$ (in the metric topology) such that $S \subset Q + \mathcal{C}$.  Then there is a finite set of elements $\{s_{i}\}_{i=1}^{r}$ in $S$ such that every element of $S$ can be written as $s_{i} + c$ for some $c \in \mathcal{C} \cap \Lambda$. 
\end{theorem}

\subsection{Vector bundles on curves}

Let $B$ be a smooth projective curve and let $\mathcal{E}$ be a vector bundle of rank $r$ on $B$.  Write the Harder-Narasimhan filtration of $\mathcal{E}$ as
\begin{equation*}
0 = \mathcal{F}_{0} \subset \mathcal{F}_{1} \subset \mathcal{F}_{2} \subset \ldots \subset \mathcal{F}_{k} = \mathcal{E}.
\end{equation*}
We denote by $\mu^{max}(\mathcal{E})$ the maximal slope of any torsion-free subsheaf, i.e.,~$\mu^{max}(\mathcal{E}) = \mu(\mathcal{F}_{1})$.  We denote by $\mu^{min}(\mathcal{E})$ the minimal slope of any torsion-free quotient, i.e.,~$\mu^{min}(\mathcal{E}) = \mu(\mathcal{E}/\mathcal{F}_{k-1})$.  Note that by the mediant inequality for every index $1 < i \leq k$ we have
\begin{equation}
\label{eq:mediant}
\mu(\mathcal{F}_{i}) = \frac{c_{1}(\mathcal{F}_{i-1}) + c_{1}(\mathcal{F}_{i}/\mathcal{F}_{i-1})}{\rk(\mathcal{F}_{i-1}) + \rk(\mathcal{F}_{i}/\mathcal{F}_{i-1})} < \mu(\mathcal{F}_{i-1}).
\end{equation}

\begin{lemma} \label{lemm:minslopelower}
Let $f: Y \to S$ be a smooth projective morphism of varieties with relative dimension $1$.  Suppose that $\mathcal{E}$ is a locally free sheaf on $Y$.  Then $\mu^{min}(\mathcal{E}|_{Y_{s}})$ is a lower-semicontinuous function on $S$.
\end{lemma}

\begin{proof}
By \cite[Theorem 2.3.2]{HL97} there is a dense open set $U \subset S$ and a torsion-free sheaf $\mathcal{F}$ on $Y_{U}$ such that $(\mathcal{E}/\mathcal{F})|_{Y_{t}}$ is the minimal slope quotient of $\mathcal{E}|_{Y_{t}}$ for the fiber $Y_{t}$ over any point $t \in U$.  Arguing by Noetherian induction, it suffices to show that if $s$ denotes an arbitrary point of $S$ then $\mu^{min}(\mathcal{E}|_{Y_{s}}) \leq \mu^{min}(\mathcal{E}|_{Y_{t}})$.  By projectivity of the Quot scheme, for any point $s \in S$ there is a surjection $\mathcal{E}|_{Y_{s}} \to \mathcal{Q}_{s}$ where $\mathcal{Q}_{s}$ has the same degree and rank as $(\mathcal{E}/\mathcal{F})|_{Y_{t}}$.  In particular $\mu(\mathcal{Q}_{s,tf}) \leq \mu((\mathcal{E}/\mathcal{F})|_{Y_{t}})$ finishing the proof.
\end{proof}

\begin{definition}
We say that a coherent sheaf $\mathcal{E}$ on a smooth projective curve $B$ is generically globally generated if the evaluation map
\begin{equation*}
H^{0}(B,\mathcal{E}) \otimes \mathcal{O}_{B} \to \mathcal{E}
\end{equation*}
is surjective at the generic point of $B$.
\end{definition}

\begin{lemma} \label{lemm:genericallygloballygenerated}
Let $B$ be a smooth projective curve.  Suppose that $\mathcal{E}$ is a generically globally generated vector bundle on $B$.  Then $\mu^{min}(\mathcal{E}) \geq 0$.
\end{lemma}

\begin{proof}
Since the evaluation map on global sections has torsion cokernel, we have
\begin{equation*}
\mu(\mathcal{E}) \geq \mu(H^{0}(B,\mathcal{E}) \otimes \mathcal{O}_{B}) = 0.
\end{equation*}
Denote the Harder-Narasimhan filtration of $\mathcal{E}$ by
\begin{equation*}
0 = \mathcal{F}_{0} \subset \mathcal{F}_{1} \subset \ldots \subset \mathcal{F}_{k} = \mathcal{E}.
\end{equation*}
Since $\mathcal{F}_{1}$ is the maximal destabilizing subsheaf, we see that
\begin{equation*}
\mu(\mathcal{F}_{1}) \geq \mu(\mathcal{E}) \geq 0.
\end{equation*}
Since $\mathcal{E}$ is generically globally generated, its quotient $\mathcal{E}/\mathcal{F}_{1}$ is also generically globally generated, and we conclude by induction on the length $k$ of the Harder-Narasimhan filtration.
\end{proof}

\subsubsection{Cohomology bounds} We next recall some bounds on the cohomology groups of semistable and generically globally generated vector bundles.

\begin{lemma} \label{lemm:semistablevanishing}
Let $\mathcal{E}$ denote a semistable vector bundle on a smooth projective curve $B$.  Suppose that $\mu(\mathcal{E}) > (2g(B)-2)$. Then $H^{1}(B,\mathcal{E}) = 0$.
\end{lemma}

\begin{proof}
By Serre duality it suffices to show that $H^{0}(B,\mathcal{E}^{\vee} \otimes \omega_{B}) = 0$.  Since $\mathcal{E}$ is semistable, $\mathcal{E}^{\vee} \otimes \omega_{B}$ is as well.  Since $\mu(\mathcal{E}^{\vee} \otimes \omega_{B}) < 0$ there are no non-zero morphisms $\mathcal{O}_{B} \to \mathcal{E}^{\vee} \otimes \omega_{B}$.
\end{proof}

\begin{corollary} \label{coro:checkinggg}
Let $\mathcal{E}$ denote a vector bundle on the smooth projective curve $B$.
Define $d = (2g(B)-2) - \mu^{min}(\mathcal{E})$.  Then:
\begin{enumerate}
\item For any line bundle $\mathcal{L}$ of degree $> d$ we have that $H^{1}(B,\mathcal{E} \otimes \mathcal{L}) = 0$.  
\item For any line bundle $\mathcal{T}$ of degree $> d+1$ we have that $\mathcal{E} \otimes \mathcal{T}$ is globally generated.  
\end{enumerate}
\end{corollary}

\begin{proof}
Write the Harder-Narasimhan filtration of $\mathcal{E}$ as
\begin{equation*}
0 = \mathcal{F}_{0} \subset \mathcal{F}_{1} \subset \mathcal{F}_{2} \subset \ldots \subset \mathcal{F}_{k} = \mathcal{E}.
\end{equation*}

(1)  Since the slopes $\mu(\mathcal{F}_{i}/\mathcal{F}_{i-1})$ are strictly decreasing in $i$ we have $\mu(\mathcal{F}_{i}/\mathcal{F}_{i-1} \otimes \mathcal{L}) > 2g(B)-2$ for every index $i=1,2,\ldots,k$. Thus for $i$ in this range we have $H^{1}(B,\mathcal{F}_{i}/\mathcal{F}_{i-1} \otimes \mathcal{L}) = 0$ by Lemma \ref{lemm:semistablevanishing}.
Using the exact sequences
\begin{equation*}
H^{1}(B, \mathcal{F}_{i-1}  \otimes \mathcal{L}) \to H^{1}(B, \mathcal{F}_{i}  \otimes \mathcal{L}) \to H^{1}(B, \mathcal{F}_{i}/\mathcal{F}_{i-1} \otimes \mathcal{L}) \to 0
\end{equation*}
and arguing by induction on $i$ we see that $H^{1}(B,\mathcal{E} \otimes \mathcal{L}) = 0$.

(2) follows immediately from (1) and the long exact sequence of cohomology associated to the inclusion
\begin{equation*}
\mathcal{E} \otimes \mathcal{T} \otimes \mathcal{O}_{B}(-p) \hookrightarrow \mathcal{E} \otimes \mathcal{T}
\end{equation*}
where $p$ is any closed point of $B$.
\end{proof}

\begin{lemma} \label{lemm:ggh1bound}
Let $B$ be a smooth projective curve of genus $g(B)$.  Suppose that $\mathcal{E}$ is a generically globally generated bundle on $B$.  Then
\begin{enumerate}
\item $h^{0}(B,\mathcal{E}) \leq \deg(\mathcal{E}) + \rk(\mathcal{E})$.
\item $h^{1}(B,\mathcal{E}) \leq g(B) \rk(\mathcal{E})$.
\end{enumerate}
\end{lemma}

\begin{proof}
Write $0 = \mathcal{F}_{0} \subset \mathcal{F}_{1} \subset \ldots \subset \mathcal{F}_{k} = \mathcal{E}$ for the Harder-Narasimhan filtration of $\mathcal{E}$.  Since $\mathcal{E}$ is generically globally generated, Lemma \ref{lemm:genericallygloballygenerated} shows that every successive quotient $\mathcal{F}_{i}/\mathcal{F}_{i-1}$ has degree $\geq 0$.

If $0 \leq \mu(\mathcal{F}_{i}/\mathcal{F}_{i-1}) \leq 2g(B)-2$, Clifford's Theorem for semistable bundles as in \cite[Theorem 2.1]{BGN97} shows that
\begin{equation*}
h^{0}(B,\mathcal{F}_{i}/\mathcal{F}_{i-1}) \leq \frac{1}{2} \deg(\mathcal{F}_{i}/\mathcal{F}_{i-1}) + \rk(\mathcal{F}_{i}/\mathcal{F}_{i-1})
\end{equation*}
and that
\begin{equation*}
h^{1}(B,\mathcal{F}_{i}/\mathcal{F}_{i-1}) = h^{0}(B,\mathcal{F}_{i}/\mathcal{F}_{i-1}) - \chi(\mathcal{F}_{i}/\mathcal{F}_{i-1}) \leq \frac{-1}{2} \deg(\mathcal{F}_{i}/\mathcal{F}_{i-1}) + g(B)\rk(\mathcal{F}_{i}/\mathcal{F}_{i-1}).
\end{equation*}
On the other hand, if $2g(B)-2 < \mu(\mathcal{F}_{i}/\mathcal{F}_{i-1})$ then $h^{1}(B,\mathcal{F}_{i}/\mathcal{F}_{i-1}) = 0$ and
\begin{equation*}
h^{0}(B,\mathcal{F}_{i}/\mathcal{F}_{i-1}) = \deg(\mathcal{F}_{i}/\mathcal{F}_{i-1}) + \rk(\mathcal{F}_{i}/\mathcal{F}_{i-1})(1-g(B))
\end{equation*}
by Lemma \ref{lemm:semistablevanishing} and Riemann-Roch.  Since
\begin{equation*}
h^{0}(B,\mathcal{E}) \leq \sum_{i=1}^{s} h^{0}(B,\mathcal{F}_{i}/\mathcal{F}_{i-1})  \qquad \textrm{and} \qquad h^{1}(B,\mathcal{E})  \leq \sum_{i=1}^{s} h^{1}(B,\mathcal{F}_{i}/\mathcal{F}_{i-1})
\end{equation*}
we obtain the desired statement using the additivity of $\deg$ and $\rk$ in exact sequences.
\end{proof}

\subsection{Slope stability for smooth projective varieties}  
The notion of slope stability with respect to movable curve classes was developed by \cite{CP11}, \cite{GKP14}, \cite{GKP16}.

\begin{definition}
Let $X$ be a smooth projective variety and let $\alpha \in \Nef_{1}(\mathcal{X})$.  For any torsion-free sheaf $\mathcal{E}$ on $X$, we define
\begin{equation*}
\mu_{\alpha}(\mathcal{E}) = \frac{c_{1}(\mathcal{E}) \cdot \alpha}{\rk(\mathcal{E})}.
\end{equation*}
We say that $\mathcal{E}$ is $\alpha$-semistable if for every non-zero torsion-free subsheaf $\mathcal{F} \subset \mathcal{E}$ we have $\mu_{\alpha}(\mathcal{F}) \leq \mu_{\alpha}(\mathcal{E})$.
\end{definition}

Every torsion free sheaf admits a maximal destabilizing subsheaf with respect to this slope function.  Thus we get a theory of $\alpha$-Harder-Narasimhan filtrations for torsion-free sheaves on $X$.  The following definition captures the slopes of the pieces of the Harder-Narasimhan filtration.

\begin{definition}
Let $X$ be a smooth projective variety and let $\alpha \in \Nef_{1}(\mathcal{X})$.  Suppose that $\mathcal{E}$ is a torsion-free sheaf of rank $r$.  Write
\begin{equation*}
0 = \mathcal{F}_{0} \subset \mathcal{F}_{1} \subset \mathcal{F}_{2} \subset \ldots \subset \mathcal{F}_{k} = \mathcal{E}
\end{equation*}
for the $\alpha$-Harder-Narasimhan filtration of $\mathcal{E}$.
The slope panel $\SP_{X,\alpha}(\mathcal{E})$ is the $r$-tuple of rational numbers obtained by combining for every index $i$ the list of $\rk(\mathcal{F}_{i}/\mathcal{F}_{i-1})$ copies of $\mu_{\alpha}(\mathcal{F}_{i}/\mathcal{F}_{i-1})$ (arranged in non-increasing order):
\begin{equation*}
\SP_{X, \alpha}(\mathcal{E}) = ( \underbrace{\mu_{\alpha}(\mathcal{F}_{1}/\mathcal{F}_{0}), \ldots}_{\rk(\mathcal{F}_{1}/\mathcal{F}_{0}) \textrm{ copies}}, \underbrace{\mu_{\alpha}(\mathcal{F}_{2}/\mathcal{F}_{1}), \ldots}_{\rk(\mathcal{F}_{2}/\mathcal{F}_{1}) \textrm{ copies}},\ldots,\underbrace{\mu_{\alpha}(\mathcal{F}_{k}/\mathcal{F}_{k-1}), \ldots}_{\rk(\mathcal{F}_{k}/\mathcal{F}_{k-1}) \textrm{ copies}})
\end{equation*}
We denote by $\mu^{max}_{\alpha}(\mathcal{E})$ the maximal slope of any torsion-free subsheaf, i.e.,~$\mu^{max}_{\alpha}(\mathcal{E}) = \mu_{\alpha}(\mathcal{F}_{1})$.  We denote by $\mu^{min}_{\alpha}(\mathcal{E})$ the minimal slope of any torsion-free quotient, i.e.,~$\mu^{min}_{\alpha}(\mathcal{E}) = \mu_{\alpha}(\mathcal{E}/\mathcal{F}_{k-1})$.

When discussing slope panels in the case when $X$ is a curve, we will always let $\alpha$ be an ample class of degree $1$ and thus we will simply write $\SP_{X}(\mathcal{E})$.
\end{definition}

Variations of the next result have been proved many times in the literature.

\begin{theorem}[{\cite[Proposition 1.3.32]{Pang}}] \label{theo:HNisfoliation}
Let $X$ be a smooth projective variety and let $\alpha \in \Nef_{1}(X)$ be a nef curve class.  Denote the $\alpha$-Harder-Narasimhan filtration of $T_{X}$ by
\begin{equation*}
0 = \mathcal{F}_{0} \subset \mathcal{F}_{1} \subset \ldots \subset \mathcal{F}_{k} = T_{X}.
\end{equation*}
Then every term $\mathcal{F}_{i}$ such that $\mu_{\alpha}^{min}(\mathcal{F}_{i}) > 0$ defines a foliation on $X$.
\end{theorem}

Suppose that $f: X \dashrightarrow Y$ is a rational map from a smooth projective variety $X$ to a normal projective variety $Y$.  Let $U$ be the open locus where $f$ is defined.  There is a unique foliation on $X$ whose restriction to $U$ is the saturation in $T_{U}$ of the kernel of $T_{U} \to f^{*}T_{Y}$.  We call this the foliation induced by $f$.  Note that if $f': X \dashrightarrow Y'$ is a rational map birationally equivalent to $f$ then $f$ and $f'$ induce the same foliation.

\subsection{Fujita invariant}
Recall from Definition~\ref{defi:a-invariant} that if $X$ is a smooth projective variety over a field of characteristic $0$ and $L$ is a big and nef $\mathbb{Q}$-Cartier divisor on $X$ then
\begin{equation*}
a(X,L) = \min \{ t\in \mathbb{R} \mid  K_X + tL \in \Eff^{1}(X) \}.
\end{equation*}

By \cite{BDPP13} the Fujita invariant will be positive if and only if $X$ is geometrically uniruled.  We will rely on the following boundedness result.

\begin{theorem}[{\cite[Theorem 1.2]{Dicerbo17}}, {\cite[Theorem 1.3]{HL20}}]
\label{theo:Dicerbo}
Fix a positive integer $n$ and fix $\epsilon > 0$.  As we vary $X$ over all smooth projective varieties of dimension $n$ defined over a field of characteristic $0$ and vary $L$ over all big and nef Cartier divisors on $X$, there are only finitely many possible values of $a(X,L)$ in the range $(\epsilon,\infty)$.
\end{theorem}

The Fujita invariant is most useful for analyzing pairs satisfying an additional assumption.

\begin{definition}
Let $X$ be a smooth projective variety and let $L$ be a big and nef $\mathbb{Q}$-divisor on $X$.  We say that $(X,L)$ is adjoint rigid if $K_{X} + a(X,L)L$ has Iitaka dimension $0$.  If $X$ is singular and $L$ is a big and nef $\mathbb{Q}$-Cartier divisor, we say that $(X,L)$ is adjoint rigid if $(X',\phi^{*}L)$ is adjoint rigid for some resolution of singularities $\phi: X' \to X$. This definition does not depend on the choice of resolution.
\end{definition}

\subsection{Boundedness and the Fujita invariant} 
In this section we recall some results of \cite{LST18}.  Our first construction shows that the family of subvarieties of $X$ which are adjoint rigid and have the same Fujita invariant as $X$ is bounded.

\begin{construction} \label{cons:rigidsubvarieties}
Let $k$ be a field of characteristic $0$.  Let $X$ be a geometrically uniruled geometrically integral smooth projective $k$-variety and let $L$ be a big and nef $\mathbb{Q}$-Cartier divisor on $X$.  By \cite[Theorem 4.19]{LST18} there exist a proper closed subset $V$, finitely many projective varieties $W_i \subset \mathrm{Hilb}(X)$, proper families $p_{i}: U_{i} \to W_{i}$ where $U_i$ is a smooth birational model of the universal family $U_i' \to W_i$, and dominant generically finite morphisms $s_i :  U_i \to X$
such that
\begin{itemize}
\item over $\overline{k}$, a general fiber of $p_{i,\overline{k}} : U_{i,\overline{k}} \rightarrow W_{i,\overline{k}}$ is an integral uniruled projective variety which is mapped birationally by $s_{i,\overline{k}}$ onto the subvariety of $X_{\overline{k}}$ parametrized by the corresponding point of $\Hilb(X_{\overline{k}})$;
\item a general fiber $Z$ of $p_{i}$ is a smooth projective variety satisfying $a(Z,s_{i}^{*}L|_{Z}) = a(X,L)$ and is adjoint rigid with respect to $s_{i}^{*}L|_{Z}$;
and
\item for every subvariety $Y \subset X$ not contained in the augmented base locus $\mathbf{B}_{+}(L)$ of $L$ which satisfies $a(Y,L|_{Y}) \geq a(X,L)$ and which is adjoint rigid with respect to $L$, either $Y$ is contained in $V$ or there is some index $i$ and a smooth fiber of $p_{i}$ that is mapped birationally to $Y$ under the map $s_{i}$.
\end{itemize}
\end{construction}

In fact more is true: the next result shows that the morphisms $f: Y \to X$ such that $Y$ is adjoint rigid and has the same Fujita invariant as $X$ also form a bounded family up to twisting.

\begin{theorem} \label{theo:ainvboundedandtwists}
Let $k$ be a field of characteristic $0$.  Let $X$ be a geometrically uniruled geometrically integral smooth projective $k$-variety and let $L$ be a big and nef $\mathbb{Q}$-Cartier divisor on $X$.  Denote by $\{ p_{i}: {U}_{i} \to W_{i}\}$ the finite set of families equipped with maps $s_{i}: U_{i} \to X$ and by $V$ the closed subset of Construction \ref{cons:rigidsubvarieties}.  There is a closed set $R \subset X$ and a finite set of smooth projective varieties $Y_{i, j}$ equipped with dominant morphisms $r_{i, j}: Y_{i, j} \to T_{i, j}$ with connected fibers and dominant morphisms $h_{i, j}: Y_{i, j} \to U_i$ forming commuting diagrams
\begin{equation*}
\xymatrix{ {Y}_{i, j} \ar[r]^{h_{i, j}} \ar[d]_{r_{i, j}} &  {U}_i \ar[d]_{p_{i}} \\
T_{i, j} \ar[r]_{t_{i, j}} & W_i}
\end{equation*}
that satisfy the following properties:
\begin{enumerate}
\item each map $h_{i, j}$ is generically finite and $f_{i, j} = s_i \circ h_{i, j}$ is not birational; 
\item $t_{i, j}$ is a finite Galois cover and $T_{i, j}$ is normal;
\item $\mathrm{Bir}(Y_{i, j}/U_i) = \mathrm{Aut}(Y_{i, j}/U_i)$;
\item every twist $Y_{i, j}^{\sigma}$ of $Y_{i, j}$ over $U_i$ admits a morphism $r_{i, j}^{\sigma}: Y_{i, j}^{\sigma} \to T_{i,j}^{\sigma}$ which is a twist of $r_{i, j}$;
\item we have $a(Y_{i, j},f_{i, j}^{*}L) = a(X,L)$;
\item suppose that $Y$ is a geometrically integral smooth projective variety and that $f: Y \to X$ is a morphism that is generically finite onto its image but not birational such that $a(Y,f^{*}L) \geq a(X,L)$.  Suppose furthermore that $y \in Y(k)$ satisfies $f(y) \not \subset R$.  Then:
\begin{enumerate}
\item there are indices $i,j$ and a twist $h_{i, j}^{\sigma}: Y_{i, j}^{\sigma} \to U_i$ of $h_{i, j}$ such that $f(y) \in s_{i}(h_{i, j}^{\sigma}(Y_{i, j}^{\sigma}(k)))$, and
\item if $(Y,f^{*}L)$ is adjoint rigid then furthermore $f$ factors rationally through $h_{i, j}^{\sigma}$ and $f$ maps $Y$ birationally to a fiber of $r_{i, j}^{\sigma}$.
\end{enumerate}
\end{enumerate}
\end{theorem}

\begin{proof}
Consider the families $p_{i}: {U}_{i} \to W_{i}$.  By applying \cite[Lemma 7.3]{LST18} there exists a Zariski open subset $W^\circ_{i}$ such that each map ${U}^\circ_{i} \to W^\circ_{i}$ is a good family in the sense of \cite[Definition 8.2]{LST18}.

We may then apply \cite[Lemma 8.3]{LST18} to each ${U}_{i}$ equipped with $s_{i}^{*}L$.  The result is a closed set $D_{i} \subset U_{i}$ and a finite set of smooth projective varieties $Y_{i,j}$ equipped with morphisms $r_{i,j}: Y_{i,j} \to T_{i,j}$, $h_{i,j}: Y_{i,j} \to U_{i}$, and $t_{i, j} : T_{i, j} \to W_i$ that have the following properties.  
First, since $r_{i,j}$ is constructed as the Stein factorization of $Y_{i, j} \to U_i \to W_i$ we see that every twist $Y_{i, j}^{\sigma}$ of $Y_{i,j}$ over ${U}_{i}$ admits a morphism $r_{i, j}^{\sigma}$ that is a twist of $r_{i, j}$.  Second, suppose that $q: Y \to U_{i}$ is a generically finite morphism such that $a(Y,q^{*}s_{i}^{*}L) = a(X,L)$ and a general fiber of the Iitaka fibration for $K_{Y} + a(Y,q^{*}s_{i}^{*}L)q^{*}s_{i}^{*}L$ maps generically finitely onto a general fiber of ${U}_{i} \to W_{i}$. 
Suppose furthermore that $y \in Y(k)$ is a rational point such that $q(y)$ is not contained in $D_{i}$.  Then there is some index $j$ and a twist $h_{i,j}^{\sigma}: Y_{i,j}^{\sigma} \to {U}_{i}$ such that $q(y)$ is in $h_{i,j}^{\sigma}(Y_{i,j}^{\sigma}(k))$.  Furthermore every general fiber of the canonical fibration for $K_{Y} + a(Y,q^{*}s_{i}^{*}L)q^{*}s_{i}^{*}L$ is birational to a fiber of $r^\sigma_{i,j}$.

By \cite[Theorem 4.18]{LST18} there is a proper closed subset $V' \subset X$ which is the union of all subvarieties $Y$ satisfying $a(Y,L|_{Y}) > a(X,L)$.  We let $R$ be the union of $V$ and $V'$ with $\cup_{i}s_{i}(D_{i})$. 
We also enlarge $R$ by adding the images of singular fibers of $r_{i, j}$.

We verify each property. (1), (2), (3) follow from \cite[Lemma 8.3]{LST18}.  We have already verified (4). (5) follows from \cite[Lemma 8.3 (i) and (ii)]{LST18}.

Now we verify (6). Suppose $f: Y \to X$ is as in the statement.  By assumption $f(Y) \not \subset V'$.  In particular this implies that $a(Y,f^{*}L) = a(f(Y),L|_{f(Y)}) = a(X,L)$.  Let $F$ be the closure of a general fiber of the canonical fibration for $K_{Y} + a(Y,f^{*}L)f^{*}L$ so that $(F,f^{*}L|_{F})$ is adjoint rigid.  Then by \cite[Lemma 4.9]{LST18} we see that $(f(F),L|_{f(F)})$ is also adjoint rigid and thus is birational to a fiber of some map $p_{i}: {U}_{i} \to W_{i}$.  This induces a rational map $T \dashrightarrow \mathrm{Hilb}(X)$ where $T$ is the base of the canonical fibration of $K_{Y} + a(Y,f^{*}L)f^{*}L$.  Since $U_{i}$ is birational to the universal family over $W_{i}$, we also obtain a rational map $Y \dashrightarrow U_{i}$. 
Since the desired statement only depends on the birational equivalence class of $f: Y \to X$ (and not the choice of birational model of $Y$), after blowing up $Y$ we may suppose that $Y$ admits a morphism to ${U}_{i}$ such that the general fiber of the canonical fibration on $Y$ maps generically finitely onto a fiber of the map ${U}_{i} \to W_{i}$.  Then the desired containment of rational points follows from \cite[Lemma 8.3 (vi)]{LST18} as described above.  When $(Y,f^{*}L)$ is adjoint rigid, the factoring statement also follows from \cite[Lemma 8.3 (vi)]{LST18}.
\end{proof}

\section{Sections of good fibrations}  \label{sect:goodfibration}

\begin{definition}
We say that a morphism $\pi: \mathcal{Z} \to B$ is a good fibration if:
\begin{enumerate}
\item $\mathcal{Z}$ is a smooth projective variety,
\item $B$ is a smooth projective curve, and
\item $\pi$ is flat and has connected fibers.
\end{enumerate}
\end{definition}

Suppose that $\pi: \mathcal{Z} \to B$ is a good fibration.  We let $\Sec(\mathcal{Z}/B)$ denote the open subset of the Hilbert scheme parametrizing sections of $\pi$.  If $M \subset \Sec(\mathcal{Z}/B)$ is an irreducible component, the expected dimension of $M$ is
\begin{equation*}
\chi(T_{\mathcal{Z}/B}|_{C}) = -K_{\mathcal{Z}/B} \cdot C + (\dim \mathcal{Z}-1)(1-g(B))
\end{equation*}
where $C$ is any section parametrized by $M$.  The expected dimension is a lower bound for the dimension of $M$.  An upper bound is
\begin{equation*}
\dim H^{0}(B,T_{\mathcal{Z}/B}|_{C}) = -K_{\mathcal{Z}/B} \cdot C + (\dim \mathcal{Z}-1)(1-g(B)) + \dim H^{1}(B,T_{\mathcal{Z}/B}|_{C}).
\end{equation*}

One of the basic facts about sections of a good fibration is the Northcott property, which in our setting should be interpreted in the following way.

\begin{corollary}[{\cite[Lemma 2.2]{LT21b}}] \label{coro:boundednegativity}
Let $\pi: \mathcal{Z} \to B$ be a good fibration and let $L$ be a generically relatively ample $\mathbb{Q}$-Cartier divisor on $\mathcal{Z}$.  If we fix a constant $Q$, then there are only finitely many components of $\Sec(\mathcal{Z}/B)$ parametrizing sections $C$ satisfying $L \cdot C \leq Q$. 
\end{corollary}

\subsection{Relatively free sections and general points}

Suppose $\pi: \mathcal{Z} \to B$ is a good fibration.  Fix points $q_{1},\ldots,q_{m} \in \mathcal{Z}$ which are contained in different fibers of $\pi$.  We let $\Sec(\mathcal{Z}/B,q_{1},\ldots,q_{m})$ denote the sublocus of $\Sec(\mathcal{Z}/B)$ parametrizing sections containing the points $q_{1},\ldots,q_{m}$. In particular, if $M \subset \Sec(\mathcal{Z}/B,q_{1},\ldots,q_{m})$ is an irreducible component then the expected dimension of $M$ is
\begin{equation*}
\chi(T_{\mathcal{Z}/B}|_{C}(-q_{1}-\ldots-q_{m})) = -K_{\mathcal{Z}/B} \cdot C + (\dim \mathcal{Z}-1)(1-g(B)) - m (\dim(\mathcal{Z})-1)
\end{equation*}
where $C$ is any section parametrized by $M$.  The expected dimension is a lower bound for the dimension of $M$.  An upper bound is
\begin{align*}
\dim H^{0}(B,T_{\mathcal{Z}/B}|_{C}(-q_{1}-\ldots-q_{m})) = -K&_{\mathcal{Z}/B} \cdot C + (\dim \mathcal{Z}-1)(1-g(B))  - m (\dim(\mathcal{Z})-1) \\
& + \dim H^{1}(C,T_{\mathcal{Z}/B}|_{C}(-q_{1}-\ldots-q_{m})).
\end{align*}
The following result describes how the normal bundle of a section $C$ controls the number of general points contained in deformations of $C$.

\begin{proposition}[\cite{LT21b} Proposition 3.3] \label{prop:deffixpoints} 
Let $\pi : \mathcal Z \rightarrow B$ be a good fibration.  
Fix points $q_{1},\ldots,q_{m}$ of $\mathcal Z$ contained in different fibers of $\pi$.  Let $M$ denote an irreducible component of $\Sec(\mathcal{Z}/B,q_{1},\ldots,q_{m})$ and suppose that the sections parametrized by $M$ dominate $\mathcal{Z}$.  Then for a general section $C$ parametrized by $M$ and for a general point $p \in B$ we have that $H^{0}(C, T_{\mathcal{Z}/B}|_{C}(-q_{1}-\ldots-q_{m})) \to T_{\mathcal{Z}/B}|_{C}|_{p}$ is surjective.

Conversely, suppose we fix a section $C$.  Suppose that $q_{1},\ldots,q_{m}$ are distinct points of $C$ such that $H^1(C, T_{\mathcal{Z}/B}|_{C}(-q_{1}-\ldots-q_{m}))= 0$.  Let $M \subset \Sec(\mathcal{Z}/B,q_{1},\ldots,q_{m})$ denote the unique irreducible component containing $C$.  If for a general point $p \in C$ we have that $H^{0}(C, T_{\mathcal{Z}/B}|_{C}(-q_{1}-\ldots-q_{m})) \to T_{\mathcal{Z}/B}|_{C}|_{p}$ is surjective, then $M$ parametrizes a dominant family of sections on $\mathcal{Z}$.
\end{proposition}

\begin{corollary} \label{coro:domfamilyexpdim}
Let $\pi: \mathcal{Z} \to B$ be a good fibration.  Suppose that $M$ is an irreducible component of $\Sec(\mathcal{Z}/B)$ parametrizing a dominant family of sections.  Letting $C$ denote a general section parametrized by $M$, we have
\begin{equation*}
-K_{\mathcal{Z}/B} \cdot C + (\dim \mathcal{Z}-1)(1-g(B)) \leq \dim(M) \leq -K_{\mathcal{Z}/B} \cdot C + \dim \mathcal{Z}-1.
\end{equation*}
\end{corollary}

\begin{proof}
By Proposition \ref{prop:deffixpoints} the bundle $T_{\mathcal{Z}/B}|_{C}$ is generically globally generated.  Thus we have $h^{1}(C,T_{\mathcal{Z}/B}|_{C}) \leq g(B)(\dim(\mathcal{Z})-1)$ by Lemma \ref{lemm:ggh1bound}.  The desired statement follows.
\end{proof}

Recall that a section $C$ is relatively free if $H^{1}(C,T_{\mathcal{Z}/B}|_{C}) = 0$ and $T_{\mathcal{Z}/B}|_{C}$ is globally generated.  Proposition \ref{prop:deffixpoints} shows that any relatively free section deforms in a dominant family on $\mathcal{Z}$.  It is easiest to work with relatively free sections when we impose further conditions on the positivity of the terms of the Harder-Narasimhan filtration of $T_{\mathcal{Z}/B}|_{C}$.

\begin{definition}
Let $\pi: \mathcal{Z} \to B$ be a good fibration.  We say that a section $C$ is HN-free if
\begin{equation*}
\mu^{min}(T_{\mathcal{Z}/B}|_{C}) \geq 2g(B).
\end{equation*}
\end{definition}

The following result summarizes the key properties of HN-free sections.

\begin{lemma} \label{lemma:hnfreecurves}
Let $\pi: \mathcal{Z} \to B$ be a good fibration.  Suppose that $C$ is a HN-free section of $\pi$.  Then:
\begin{enumerate}
\item $H^{1}(C,T_{\mathcal{Z}/B}|_{C}) = 0$ and for any closed point $p \in B$ we have $H^{1}(C,T_{\mathcal{Z}/B}|_{C}(-p)) = 0$.
\item $T_{\mathcal{Z}/B}|_{C}$ is globally generated.
\item $C$ is relatively free.
\item Let $b = \mu^{min}(T_{\mathcal{Z}/B}|_{C})$.  Then deformations of $C$ can pass through at least $\lfloor b \rfloor - 2g(B) + 1$ general points of $\mathcal{Z}$.
\end{enumerate}
\end{lemma}

\begin{proof}
(1) and (2) follow from Corollary \ref{coro:checkinggg} and (3) follows from (1) and (2).  To see (4) we apply Corollary \ref{coro:checkinggg} to see that for any points $q_{1},\ldots,q_{m}$ on $C$ the twist $T_{\mathcal{Z}/B}|_{C}(-q_{1}-\ldots-q_{m})$ is globally generated and has vanishing $H^{1}$ so long as $m \leq b - 2g(B)$. The desired statement follows from Proposition \ref{prop:deffixpoints}.
\end{proof}

The next proposition shows that sections through sufficiently many general points must be HN-free.

\begin{proposition} \label{prop:generalimplieshnfree}
Let $\pi: \mathcal{Z} \to B$ be a good fibration.  Let $M$ be an irreducible component of $\Sec(\mathcal{Z}/B)$. Suppose that the sections parametrized by $M$ pass through $\geq 2g(B)+1$ general points of $\mathcal{Z}$.  Then the general section parametrized by $M$ is HN-free.
\end{proposition}

\begin{proof}
If we fix a general section $C$ parametrized by $M$ and a set of $2g(B)$ general points $\{q_{i}\}_{i=1}^{2g(B)}$ on $C$ then Proposition \ref{prop:deffixpoints} shows that $T_{\mathcal{Z}/B}|_{C}(-q_{1}-\ldots-q_{2g(B)})$ is generically globally generated. Lemma \ref{lemm:genericallygloballygenerated} shows that
\begin{equation*}
\mu^{min}(T_{\mathcal{Z}/B}|_{C}(-q_{1} - \ldots - q_{2g(B)})) \geq 0
\end{equation*}
and we conclude that $\mu^{min}(T_{\mathcal{Z}/B}|_{C}) \geq 2g(B)$.
\end{proof}

We will also need to know the following avoidance property of HN-free sections.

\begin{lemma} \label{lemm:HNavoidscodim2}
Let $\pi: \mathcal{Z} \to B$ be a good fibration.  Suppose that $C$ is a HN-free section of $\pi$.  Then for any codimension $2$ closed subset $\mathcal{W} \subset \mathcal{Z}$ there is a deformation of $C$ which is HN-free and avoids $\mathcal{W}$.
\end{lemma}

\begin{proof}
Assume for a contradiction that every deformation of $C$ meets with $\mathcal{W}$.
Then there exists $p \in \mathcal W$ such that a general deformation of $C$ containing $p$ is HN-free and the dimension of the family parametrizing such deformations is greater than or equal to $-K_{\mathcal Z/B} \cdot C + (\dim \mathcal Z -1)(1-g(B)) - \dim \mathcal W$.  Note that this is larger than the expected dimension $-K_{\mathcal Z/B} \cdot C -(\dim \mathcal Z -1)g(B)$ for the parameter space of sections through $p$.  But this contradicts with Lemma \ref{lemma:hnfreecurves} which shows that $H^{1}(C,T_{\mathcal{Z}/B}|_{C}(-p)) = 0$.
\end{proof}

\section{Batyrev's heuristic} \label{sect:bh}

This section is devoted to an in-depth examination of Batyrev's heuristic and its reformulation into a set of precise conjectures which we call Geometric Manin's Conjecture.  In Section~\ref{subsection:Batyrev} we outline Batyrev's heuristic as described in \cite{Bat88}, \cite{Tsc09}, and \cite{Bourqui11}.  We have tried to capture the spirit of the heuristic rather than exactly replicating earlier work.

Section~\ref{subsection:overnumberfields} reviews the statement of Manin's Conjecture over a number field.  It details some proposals to ``correct'' the original formulation of the conjecture.

Section~\ref{subsec:GMC} explains how the various assumptions in Batyrev's heuristic can be replaced by precise mathematical statements.  As described in the introduction, there are roughly two parts to Geometric Manin's Conjecture for a Fano fibration $\pi: \mathcal{X} \to B$: the classification of components of $\Sec(\mathcal{X}/B)$ and a homological or motivic stability for the ``good'' components.  For classification results, we present a fairly complete outline of what to expect.  For stability results, we currently do not have enough examples to present a precise conjecture.

\subsection{Batyrev's heuristic}
\label{subsection:Batyrev}

Throughout this section we fix a smooth projective curve  $B$ over $\mathbb{F}_{q}$ and a Fano fibration $\pi: \mathcal{X} \to B$.  Our goal is to describe a heuristic for counting the number of $K(B)$-points on the generic fiber $\mathcal{X}_{\eta}$ of bounded anticanonical height.  Recall that a $K(B)$-point of $\mathcal{X}_{\eta}$ is the same as a section of $\pi$, so we can equivalently count sections of $\pi$ of bounded anticanonical degree.  Our strategy is to first identify the possible numerical classes $\alpha \in N_{1}(\mathcal{X})_{\mathbb{Z}}$ representing sections of $\pi$ and then to count $\mathbb{F}_{q}$-points on the irreducible components of $\Sec(\mathcal{X}/B)$ which have class $\alpha$.

\begin{definition}
Let $d$ be an integer.
We define the ``naive'' counting function
\begin{equation*}
N_{\text{naive}}(\mathcal{X},-K_{\mathcal{X}/B}, d) = \# \{ \textrm{sections }C \textrm{ of }\pi \textrm{ such that }-K_{\mathcal{X}/B} \cdot C \leq d \}.
\end{equation*}
\end{definition}

\subsubsection{Identifying numerical classes of sections}

Our first task is to identify the possible numerical classes of sections.  In certain advantageous geometric situations -- for example, when $\mathcal{X}$ admits a transitive group action -- every section will deform to dominate $\mathcal{X}$ and thus will have nef numerical class.  This property gives our first heuristic assumption:

\begin{simp} \label{simp:one}
Every section of $\pi: \mathcal{X} \to B$ is nef.
\end{simp}

To give structure to the set of sections, we start with a few reminders about the behavior of sections over an algebraically closed ground field.  If we start from a section, we can glue on some free rational curves in fibers of $\pi$ and smooth the resulting comb to obtain a new section of higher degree.  Conversely, starting from a high degree section, we can break it into the union of a lower degree section and some $\pi$-vertical trees of rational curves.  Combining these operations, we can hope that the numerical classes of sections have a structure which is ``captured'' by free $\pi$-vertical curves.

Recall that we have an inclusion $N_{1}(\mathcal{X}_{\eta}) \hookrightarrow N_{1}(\mathcal{X})$ dual to the restriction map on Cartier divisors.  This map takes $\Nef_{1}(\mathcal{X}_{\eta})$ into $\Nef_{1}(\mathcal{X})$, and in fact one can show (see Lemma \ref{lemm:relativenef}) that $\Nef_{1}(\mathcal{X}_{\eta})$ is exactly the same as the intersection of $\Nef_{1}(\mathcal{X})$ with the subspace
\begin{equation*}
V = \{ \alpha \in N_{1}(\mathcal{X}) | F \cdot \alpha = 0\},
\end{equation*}
where $F$ is a general fiber of $\pi$.
On the other hand, all classes of sections are contained in the affine translate
\begin{equation*}
A = \{ \alpha \in N_{1}(\mathcal{X}) | F \cdot \alpha = 1\}.
\end{equation*}
According to the discussion above, we might expect to understand the numerical classes of sections via the ``gluing and breaking'' structure by translating $\Nef_{1}(\mathcal{X}_{\eta})$ into $A$.  For example, we might hope that a translate of $\Nef_{1}(\mathcal{X}_{\eta})$ into $A$ coincides with the set of classes of nef sections.  Our next assumption posits that this simplest possible relationship holds:

\begin{simp} \label{simp:two}
There is a translation from $V$ to $A$ that defines a bijection between $\Nef_{1}(\mathcal{X}_{\eta})_{\mathbb{Z}}$ and all numerical classes of sections.
We let $\alpha_0$ denote the numerical class corresponding to the translation of the origin.
\end{simp}

In summary, our assumptions lead to a convenient way of identifying all numerical classes of sections: they are in bijection with the lattice points contained in the polyhedral cone $\Nef_{1}(\mathcal{X}_{\eta})$.

\subsubsection{Counting points}
Our second task is to count the number of $\mathbb{F}_{q}$-points on the irreducible components $M$ of $\Sec(\mathcal{X}/B)$ which parametrize sections of a given numerical class.
We will organize this count using the Grothendieck-Lefschetz formula.  Recall that this formula expresses the number of $\mathbb{F}_{q}$-points of $M$ as an alternating sum of the traces of the Frobenius action on the \'etale cohomology groups with compact supports for $M$.

The (dominant) exponential term in Manin's Conjecture arises from the dimension of $M$.  Indeed, the trace of the Frobenius action on the top cohomology is always given by $q^{\dim M}$ and all other terms will have lower order in $q$.  To ensure the top order term has the expected behavior, we need to assume that every irreducible component $M$ has the expected dimension:

\begin{simp} \label{simp:three}
Every irreducible component $M$ of $\Sec(\mathcal{X}/B)$ has the expected dimension:
\begin{equation*}
\dim(M) = -K_{\mathcal{X}/B} \cdot C + (\dim \mathcal{X}-1)(1-g(B)),
\end{equation*}
where $C$ is a section parametrized by $M$.
\end{simp}

The (subdominant) polynomial term in Manin's Conjecture arises from the number of irreducible components of $\Sec(\mathcal{X}/B)$ representing each fixed numerical class.  Of course, there might be many such irreducible components representing a given numerical class, or none at all.  Again we will assume that this relationship is as simple as possible.

\begin{simp} \label{simp:four}
Each numerical class in $A$ corresponding to $\Nef_{1}(\mathcal{X}_{\eta})_{\mathbb{Z}}$ (according to the bijection in Batyrev's heuristic \ref{simp:two}) is represented by a unique irreducible component of $\Sec(\mathcal{X}/B)$.
\end{simp}

Finally, we assume that the lower degree cohomology groups do not affect the asymptotic growth rate:

\begin{simp} \label{simp:five}
To compute the asymptotics of $N_{\text{naive}}(\mathcal{X},-K_{\mathcal{X}/B},d)$ we may use the approximation $|M(\mathbb{F}_{q})| \approx q^{\dim M}$.
\end{simp}

Together our assumptions imply that the asymptotics of $N_{\text{naive}}(\mathcal{X},-K_{\mathcal{X}/B},d)$ are identical to the asymptotics of the function
\begin{equation}
\label{equation:countingfunction}
Q(\mathcal{X},-K_{\mathcal{X}/B},d) = q^{(\dim \mathcal{X}-1)(1-g(B))} \sum_{\substack{\alpha \in \Nef_{1}(\mathcal{X}_{\eta})_{\mathbb{Z}} \\ -K_{\mathcal{X}/B} \cdot \alpha \leq d}} q^{-K_{\mathcal{X}/B} \cdot \alpha}.
\end{equation}

This sum of exponentials over lattice points can be analyzed directly using a combinatorial argument, or using zeta function techniques.  To understand the asymptotics of this counting function, we introduce the following constant:

\begin{definition}
Let $\pi : \mathcal X \to B$ be a Fano fibration and let $L$ be a $\mathbb Q$-divisor on $\mathcal{X}$ such that $L|_{\mathcal X_\eta}$ is big and nef. Let $N^1(\mathcal X_\eta)_{\mathbb Z}$ be the N\'eron-Severi lattice. 
We define the face associated to $L|_{\mathcal X_\eta}$ as
\[
\mathcal F(\mathcal X_\eta, L|_{\mathcal X_\eta}) := \Nef_1(\mathcal X_\eta) \cap\{\alpha \in N_{1}(\mathcal{X}_{\eta}) \, | \, (a(\mathcal X_\eta, L|_{\mathcal X_\eta})L|_{\mathcal X_\eta} + K_{\mathcal X_\eta}) \cdot \alpha = 0\}.
\]
This is an extremal face of the rational polyhedral cone $\Nef_1(\mathcal X_\eta)$.
Let $V_{L|_{\mathcal X_\eta}}$ be the subspace of $N_1(\mathcal X_\eta)$ spanned by $\mathcal F(\mathcal X_\eta, L|_{\mathcal X_\eta})$.
We consider the lattice $\Lambda \subset V_{L|_{\mathcal X_\eta}}$ generated by integral curve classes in $V_{L|_{\mathcal X_\eta}}$
and we assign $V_{L|_{\mathcal X_\eta}}$ the Lebesgue measure such that the fundamental domain for $\Lambda$ has volume $1$.  We define the $\alpha$-constant as
\[
\alpha(\mathcal X_\eta, L|_{\mathcal X_\eta}) = \dim V_{L|_{\mathcal X_\eta}} \cdot \vol (\{ \alpha \in V_{L|_{\mathcal X_\eta}} \, |\,  L|_{\mathcal X_\eta}. \alpha \leq 1 \}).
\]
We also define the index of $(\mathcal X_\eta, L|_{\mathcal X_\eta})$ to be
\[
r(\mathcal X_\eta, L|_{\mathcal X_\eta}) = \min\{ L|_{\mathcal X_\eta}.\alpha \, | \, \alpha \in N_1(\mathcal X_\eta)_{\mathbb Z}\cap V_{L|_{\mathcal X_\eta}}, L|_{\mathcal X_\eta}.\alpha >0\}.
\]
\end{definition}

With these constants, we have the following proposition:

\begin{proposition}
We have
\begin{align*}
&Q(\mathcal{X},-K_{\mathcal{X}/B},-K_{\mathcal X/B}.\alpha_0 + dr(\mathcal X_\eta)) \\ &\sim_{d \to \infty} \frac{\alpha(\mathcal X_\eta, -K_{\mathcal X_\eta}) r(\mathcal X_\eta, -K_{\mathcal X_\eta})q^{-K_{\mathcal X/B}.\alpha_0 + (\dim \mathcal X-1)(1-g(B))}}{1-q^{-r(\mathcal X_\eta, -K_{\mathcal X_\eta})}}q^{dr(\mathcal X_\eta, -K_{\mathcal X_\eta})}(dr(\mathcal X_\eta, -K_{\mathcal X_\eta}))^{\rho(\mathcal{X}_{\eta}) -1}.
\end{align*}
\end{proposition}

\begin{proof}
Let $P(d)$ be the Ehrhart quasi-polynomial for the compact rational convex set
\begin{equation*}
\Nef_1(\mathcal X_\eta)\cap \{\alpha \in N_{1}(\mathcal{X}_{\eta}) \, | \, -K_{\mathcal X_\eta} \cdot \alpha = 1\}.
\end{equation*} 
This is a quasi-polynomial of degree $\rho(\mathcal X_\eta)-1$ and the top coefficient is given by 
\[
\alpha(\mathcal X_\eta, -K_{\mathcal X_\eta})r(\mathcal X_\eta, -K_{\mathcal X_\eta}).
\]  
(See for example \cite[Exercise 3.34]{Erhart}.)
Then the function $Q(\mathcal{X},-K_{\mathcal{X}/B},-K_{\mathcal X/B}.\alpha_0 + dr(\mathcal X_\eta,-K_{\mathcal{X}_{\eta}}))$ can be written as
\[
\sum_{e = 0}^d q^{-K_{\mathcal X/B}.\alpha_0 + er(\mathcal X_\eta, -K_{\mathcal X_\eta}) + (\dim \mathcal X -1)(1-g(B))}P(er(\mathcal X_\eta, -K_{\mathcal X_\eta})).
\]
After dividing this function by $q^{dr(\mathcal X_\eta,-K_{\mathcal{X}_{\eta}})}(dr(\mathcal X_\eta, -K_{\mathcal X_\eta}))^{\rho(\mathcal{X}_{\eta}) -1}$ and letting $d \to \infty$, the above expression converges to
\[
\alpha(\mathcal X_\eta, -K_{\mathcal X_\eta})r(\mathcal X_\eta, -K_{\mathcal X_\eta})q^{-K_{\mathcal X/B}.\alpha_0+ (\dim \mathcal X -1)(1-g(B))} \left(\sum_{e = 0}^\infty q^{-er(\mathcal X_\eta, -K_{\mathcal X_\eta})} \right).
\]
Thus our assertion follows.
\end{proof}

\begin{remark} \label{remark:arbitrarypolarization}
One can formulate a similar counting question with respect to an arbitrary polarization.  Let $\pi: \mathcal{X} \to B$ be a Fano fibration and let $L$ be a divisor on $\mathcal{X}$ whose restriction to $\mathcal{X}_{\eta}$ is big and nef.  We then define the ``naive'' counting function
\begin{equation*}
N_{\text{naive}}(\mathcal{X},L,d) = \# \{ \textrm{sections }C \textrm{ of }\pi \textrm{ such that }L \cdot C \leq d \}.
\end{equation*}
Repeating the argument above leads to an asymptotic heuristic
\begin{align*} \label{eq:conjecturalasymptotics}
&N(\mathcal{X},L,-K_{\mathcal X/B}.\alpha_0 + dr(\mathcal X_\eta), L|_{\mathcal X_\eta})) \\ 
&\sim_{d \to \infty} \frac{\alpha(\mathcal X_\eta, L|_{\mathcal X_\eta}) r(\mathcal X_\eta, L|_{\mathcal X_\eta}) q^{-K_{\mathcal X/B}.\alpha_0 + (\dim \mathcal X-1)(1-g(B))}}{1-q^{-a(\mathcal{X}_{\eta},L|_{\mathcal X_\eta})r(\mathcal X_\eta, L|_{\mathcal X_\eta})}} q^{da(\mathcal{X}_{\eta},L|_{\mathcal X_\eta})r(\mathcal X_\eta, L|_{\mathcal X_\eta})} (dr(\mathcal X_\eta, L|_{\mathcal X_\eta}))^{b(K(B),\mathcal{X}_{\eta},L|_{\mathcal X_\eta}) - 1}
\end{align*}
where 
\begin{itemize}
\item $a(\mathcal{X}_{\eta},L|_{\mathcal X_\eta})$ denotes the Fujita invariant, and
\item $b(K(B),\mathcal{X}_{\eta},L|_{\mathcal X_\eta})$ is the $b$-invariant defined by 
\[
b(K(B),\mathcal{X}_{\eta},L|_{\mathcal X_\eta}) = \dim V_{L|_{\mathcal X_\eta}}.
\]
\end{itemize}
The geometric version of Manin's Conjecture for non-canonically-polarized varieties has a number of subtleties about how best to formulate it; we will not discuss this case further here.
\end{remark}

Unfortunately none of the simplifying assumptions used in Batyrev's heuristic are valid in general.  In Section~\ref{subsec:GMC} we will address these assumptions one-by-one and explain conjectural replacements. An additional advantage is that these replacements make sense over other ground fields besides $\mathbb{F}_{q}$.

\subsection{Manin's Conjecture over number fields}
\label{subsection:overnumberfields}

The original version of Manin's Conjecture predicts an asymptotic formula for the counting function of rational points of bounded height on a smooth geometrically integral Fano variety defined over a number field $\kappa$. Let $X$ be a smooth projective geometrically integral Fano variety over $\kappa$ with an adelically metrized $\mathbb Q$-divisor $\mathcal L = (L, \{\|\cdot \|_v\})$. (See \cite{CLT10} for the definition of adelic metrics.) Then one can associate the height function 
\[
\mathsf H_{\mathcal L} : X(\kappa) \to \mathbb R_{\geq 0}
\]
to a triple $(\kappa, X, \mathcal L)$. Moreover when $L$ is ample, the set of rational points of height $\leq T$ is finite, so for any subset $Q \subset X(\kappa)$, we define the counting function
\[
N(Q, \mathcal L, T) = \{ x \in Q \, |\, \mathsf H_{\mathcal L}(x) \leq T\}. 
\]
Manin's Conjecture predicts the asymptotic growth rate of $N(Q, \mathcal L, T)$ for an appropriate choice of $Q$. 

From the beginning of the study of Manin's Conjecture, it has been recognized that one cannot take $Q =X(F)$ because rational points can be accumulating along subvarieties. Thus it is important to introduce an exceptional set and exclude the contribution of this exceptional set from the counting function. Originally it was expected that one can choose an exceptional set that is contained in a proper closed subset (\cite{BM90}).  However there are several counterexamples to this expectation; the first was found in \cite{BT96} and more recent examples have been identified in \cite{LeRu19, BHB20}. Peyre was the first to suggest that the exceptional set in Manin's Conjecture should be a thin set (\cite{Peyre03}). Here is its definition:

\begin{definition} 
Let $\kappa$ be a field of characteristic $0$ and $X$ be a variety defined over $\kappa$. We say a morphism $f: Y \to X$ from a variety $Y$ is a thin map if it is generically finite to the image and there is no rational section of $f$ from $X$.

A thin set is any subset of a finite union 
\[
\bigcup_i f_i(Y_i(\kappa)) \subset X(\kappa)
\]
where $f_i : Y_i \to X$ is a thin map.
\end{definition}

The thin set version of Manin's Conjecture has been formulated in a series of works by many authors (\cite{FMT89}, \cite{BM90}, \cite{Peyre95}, \cite{BT98}, \cite{Peyre03}, \cite{Peyre17}, and \cite{LST18}).

\begin{conjecture}[Batyrev--Manin--Peyre--Tschinkel]
\label{conjecture:Manin}
Let $\kappa$ be a number field and $X$ be a smooth Fano variety defined over $\kappa$. Let $\mathcal L = (L, \{\|\cdot\|_v\})$ be an adelically metrized big and nef $\mathbb Q$-divisor on $X$.
Suppose that $X(\kappa)$ is not a thin set. Then there exists a thin set $Z \subset X(\kappa)$ such that we have
\[
N(X(\kappa) \setminus Z, \mathcal L, T) \sim c(\kappa, \mathcal L, Z)T^{a(X, L)}(\log T)^{b(\kappa, X, L)-1},
\]
as $T \to \infty$. Here $b(\kappa, X, L)$ is the $b$-invariant defined in Remark~\ref{remark:arbitrarypolarization}. The leading constant $c(\kappa, \mathcal L, Z)$ is Peyre's constant introduced in \cite{Peyre95} and \cite{BT98}.
\end{conjecture}

Here the set $Z$ is called the exceptional set.  There have been several proposals about how to identify the exceptional set (\cite{Peyre03}, \cite{Peyre17}, and \cite{LST18}).  We will focus on the approach developed in the series of works \cite{HTT15, LTT18, HJ16, LT17, LT19, Sengupta21, LST18, LTRMS} using birational geometry.

In \cite{LST18}, Akash Kumar Sengupta and the second and third authors propose a conjectural description of the exceptional set in Manin's Conjecture. 
Suppose there is a generically finite morphism $f: Y \to X$ such that the $(a,b)$ invariants of $Y$ are larger in the lexicographic order than the corresponding invariants on $X$. 
According to Remark \ref{remark:arbitrarypolarization}, the expected growth rate of points on $Y$ is higher than that on $X$.  Thus the image $f(Y(\kappa))$ must be removed from $X(\kappa)$ if we are to have any hope of achieving the expected growth rate.  By taking into account all such morphisms $f: Y \to X$, we obtain the exceptional set of \cite{LST18}.  Using results from higher dimensional algebraic geometry such as the minimal model program (\cite{BCHM}) and the boundedness of singular Fano varieties (\cite{Birkar21}), \cite{LST18} shows that there are a finite set of such maps $f$ whose twists account for all such contributions, and in particular, that the proposed set is always a thin set.

To give a precise description of the proposed exceptional set, we first must introduce several definitions:

\begin{definition}
\label{definition:breakingthinmap}
Let $\kappa$ be a field of characteristic $0$.
Let $X$ be a smooth Fano variety defined over $\kappa$.
Let $f : Y \to X$ be a thin map from a smooth projective variety $Y$ defined over $\kappa$. 

We say $f$ is an accumulating map if $a(Y, -f^*K_X) \geq a(X, -K_X) = 1$.

We say $f$ is a breaking thin map if we have an inequality
\[
(a(X, -K_X), b(\kappa, X, -K_X)) \leq (a(Y, -f^*K_X), b(\kappa, Y, -f^*K_X)),
\]
in the lexicographic order and furthermore if the equality holds then either
\begin{itemize}
\item $\dim Y < \dim X$;
\item $\dim Y = \dim X$ and $\kappa(-f^*K_X + K_Y) >0$;
\item $\dim Y = \dim X$, $\kappa(-f^*K_X + K_Y) =0$, and $f$ is geometrically non-Galois, or;
\item $\dim Y = \dim X$, $\kappa(-f^*K_X + K_Y) =0$, and $f$ is geometrically Galois and face contracting in the sense of \cite[Definition 4.26]{LST18}.
\end{itemize}
\end{definition}

A version of the conjectural exceptional set $Z$ constructed in \cite[Section 5]{LST18} is the following: in the settings of Conjecture~\ref{conjecture:Manin}, let $f : Y \to X$ run over all breaking thin maps and we define
\[
Z = \bigcup_{f} f(Y(\kappa)).
\]
As explained above, \cite[Theorem 5.7]{LST18} showed that this is indeed a thin set.

\begin{remark}
The above exceptional set is potentially slightly bigger than the one constructed in \cite[Section 5]{LST18}. The difference is that when $f$ is a geometrically non-Galois cover, here we do not impose the face-contracting condition. The above exceptional set is still a thin set because of \cite[Proposition 8.2]{LT17}. We do not know whether the two sets coincide in general.
\end{remark}

\subsection{Geometric Manin's Conjecture}
\label{subsec:GMC}

Let us come back to the original situation. Let $\pi : \mathcal X \to B$ be a Fano fibration over a smooth projective curve $B$ defined over the ground field $\kappa$. In this section we discuss Geometric Manin's Conjecture, which is given by using refinements of the assumptions of Batyrev's heuristic to predict various properties of components of $\Sec(\mathcal X/B)$.

For the sake of consistency, we will work exclusively over an algebraically closed field $\kappa$ of characteristic $0$ in this section. While one can hope for a similar structural framework in characteristic $p$, we currently do not have enough evidence to support the extension of these conjectures and results.

\subsubsection{The exceptional set in Geometric Manin's Conjecture}

We first discuss the exceptional set in the context of Geometric Manin's Conjecture.  We will give two slightly different (but closely related) definitions of what it means for an irreducible component $M \subset \Sec(\mathcal{X}/B)$ to be ``exceptional''.  %The first definition is simpler 
%and relies only on the Fujita invariant.  The second definition comes from the proposed exceptional set in \cite{LST18}.

We first define the notion of an accumulating component.  Loosely speaking, we say that $M$ is an accumulating component if there exists a family of accumulating maps that accounts for most of the sections parametrized by $M$. 

\begin{definition}
\label{definition:accumulating}
%\sho{I made changes to this definition.}
Let $\kappa$ be an algebraically closed field of characteristic $0$ and let $\pi : \mathcal X \to B$ be a Fano fibration over a smooth projective curve $B$ defined over $\kappa$ and $T > 0$ be a positive real number. Let $M \subset \Sec(\mathcal X/B)$ be an irreducible component of $\Sec(\mathcal X/B)$. 

We say $M$ is a $T$-accumulating component if the following holds.
\begin{enumerate}
\item We have a smooth projective $B$-morphism $\pi : \mathfrak Y \to S \times B$ with a $B$-morphism $f : \mathfrak Y \to \mathcal X$ representing  a family of $B$-morphisms $\mathfrak{Y}_{s} \to \mathcal X$ parametrized by $s \in S$.
\item For every fiber $\mathfrak{Y}_{s} \to B$ over a closed point $s \in S$, the induced map $f_{\eta}|_{\mathfrak{Y}_{s,\eta}}: \mathfrak{Y}_{s,\eta} \to \mathcal{X}_{\eta}$ is an accumulating map, i.e.~it is generically finite to the image and $a(\mathcal{Y}_{s,\eta},-f^{*}K_{\mathcal{X}_{\eta}}) \geq a(\mathcal{X}_{\eta}, -K_{\mathcal{X}_{\eta}})$.
\item There is a component $N \subset \Sec_S(\mathfrak Y/B)$ of the relative space of sections over $S$ (i.e.~a family of sections contained in the fibers $\mathfrak{Y}_{s} \to B$) such that $f$ induces a dominant map $f_* : N \to M$.
\item Let $N_s \subset N$ be the loci parametrizing sections contained in $\mathfrak{Y}_{s}$ for a general $s \in S$. Then we have $\mathrm{codim}_M(f_*(N_s)) \leq T$.
\end{enumerate}
\end{definition}

The second definition comes from the proposed exceptional set in \cite{LST18}.  Inspired by Definition~\ref{definition:breakingthinmap}, we first introduce an auxiliary definition:

\begin{definition}
Let $\kappa$ be a field of characteristic $0$.
Let $X$ be a smooth Fano variety defined over $\kappa$. 
Let $f : Y \to X$ be a thin map from a smooth projective variety $Y$ defined over $\kappa$. 

We say $f$ is an $a$-cover if $f$ is dominant and we have $a(Y, -f^*K_X) = a(X, -K_X)$.

We say $f$ is an exceptional map if $f$ satisfies one of the following conditions:
\begin{itemize}
\item $f$ is a non-dominant accumulating map;
\item $f$ is an $a$-cover and $\kappa(-f^*K_X + K_Y) >0$;
\item $f$ is an $a$-cover with $\kappa(-f^*K_X + K_Y) = 0$ that is a geometrically non-Galois, or;
\item $f$ is an $a$-cover $\kappa(-f^*K_X + K_Y) = 0$ that is geometrically Galois and face-contracting in the sense of \cite[Definition 4.26]{LST18}.
\end{itemize}
(Note that there is no restriction on the $b$-invariant in this definition.  If we impose a condition on the $b$-invariant we recover the exact analogue of Definition \ref{definition:breakingthinmap}.)
\end{definition}

With these notations, we introduce the following classification of components of $\Sec(\mathcal X/B)$:

\begin{definition}
\label{definition:exceptional}
Let $\kappa$ be an algebraically closed field of characteristic $0$ and let $\pi : \mathcal X \to B$ be a Fano fibration over a smooth projective curve $B$ defined over $\kappa$. Let $M \subset \Sec(\mathcal X/B)$ be a component. 

We say $M$ is an exceptional component if there exists a $B$-morphism $f : \mathcal Y \to \mathcal X$ from a smooth projective $B$-variety $\mathcal Y$ and a dominant component $N \subset \Sec(\mathcal Y/B)$ on $\mathcal Y$ such that the base change $f_\eta : \mathcal Y_\eta \to \mathcal X_\eta$ is an exceptional map and $f$ induces a dominant map $f_* : N \to M$.

We say $M$ is a Manin component if $M$ is not an exceptional component.
\end{definition}

We separate out these two definitions because they are useful in different contexts.  Accumulating components (Definition~\ref{definition:accumulating}) are useful for geometry; indeed, this notion implicitly appears in the main theorems recorded in the introduction.  However, it is not sensitive enough for counting problems -- for example, it does not distinguish between components which are ``asymptotically negligible'', components which must be thrown away, and components which must be counted.  Exceptional components (Definition~\ref{definition:exceptional}) are suitable for counting problems but have a weaker link to the geometry.

\begin{remark} \label{rema:sufficientlyfree}
While thus far we have focused on irreducible components of $\Sec(\mathcal{X}/B)$, to obtain the closest analogue to the arithmetic situation one should also remove certain subloci of the Manin components $M \subset \Sec(\mathcal{X}/B)$ and it is important to identify exactly which subloci should be removed. 
Note that if we remove subloci which have low codimension they could feasibly affect the leading constant in the asymptotic growth rate.  Thus they must be handled correctly to obtain the expected value of Peyre's constant.  Here are two possible approaches:
\begin{enumerate}
\item As in \cite{LST18}, whenever $f: \mathcal{Y} \to \mathcal{X}$ induces an exceptional map $f_{\eta}: \mathcal{Y}_{\eta} \to \mathcal{X}_{\eta}$ we can remove all subloci of the form $f_{*}(N)$ for irreducible components $N \subset \Sec(\mathcal{Y}/B)$ (regardless of whether or not $N$ dominates $M$).
\item As in \cite{Peyre17}, we could remove the sublocus of $M$ parametrizing curves that are not ``sufficiently free''.  For example, we could simply remove the non-free locus of $M$.
\end{enumerate}
In fact, Theorem \ref{theo:maintheorem3} suggests that the two approaches are equivalent: if we discount the contributions of the exceptional locus (as in (1)), then the codimension of what is left of the non-free locus increases linearly in the degree and thus will drop out of the asymptotic calculations (as in (2)).  We will not further address this subtle issue here.
\end{remark}

\subsubsection{Classifying non-nef sections}
Suppose $\pi: \mathcal{X} \to B$ is a Fano fibration.  As in the previous section, our first step is to describe the possible numerical classes of sections.  As above, we let $V$ denote the subspace of $N_{1}(\mathcal{X})$ consisting of numerical classes which have intersection $0$ against a general fiber $F$, and $A$ its affine translate of classes with intersection $1$ against $F$.

In contradistinction to Batyrev's heuristic \ref{simp:one}, many Fano fibrations admit non-nef sections.  However, according to the philosophy of Manin's Conjecture we can expect all non-dominant families of sections to be contained in the exceptional set, at least if the degree is sufficiently large.  This issue is addressed in Geometric Manin's Conjecture \ref{guid:one}.

\begin{guid} \label{guid:one}
%\sho{I made changes to this.}
Let $\pi: \mathcal{X} \to B$ be a Fano fibration. Then there exists $T = T(\pi)$ such that the components $M \subset \Sec(\mathcal{X}/B)$ which parametrize non-nef sections satisfy the following properties:
\begin{enumerate}
\item There is a Zariski-closed proper subset $\mathcal{Y} \subsetneq \mathcal{X}$ such that every non-nef section $C$ of $\pi$ is contained in $\mathcal{Y}$.
\item All but finitely many such irreducible components are $T$-accumulating.
\end{enumerate}
\end{guid}

Over a field of characteristic zero, Geometric Manin's Conjecture \ref{guid:one} is established by Theorem \ref{theo:maintheorem2}.(1) and its proof.

\subsubsection{Structure of nef sections}
Our next task is to understanding the possible numerical classes of nef sections, i.e.~$\Nef_{1}(\mathcal{X}) \cap A$.
Unfortunately $\Nef_{1}(\mathcal{X}) \cap A$ can be much more complicated than $\Nef_{1}(\mathcal{X}_{\eta})$.  (For example this intersection need not be polyhedral, see Example \ref{exam:notpolyhedral}.)

Guiding Conjecture \ref{guid:two} shows that we obtain better behavior if we restrict our attention to the convex hull of the $\mathbb{Z}$-classes $\Nef_{1}(\mathcal{X})_{\mathbb{Z}} \cap A$.  In particular the ``infinite'' part of this set does indeed come from $\Nef_{1}(\mathcal{X}_{\eta})$.

\begin{guid} \label{guid:two}
Let $\pi: \mathcal{X} \to B$ be a Fano fibration.  Let $\mathcal{P}$ denote the convex hull in $N_{1}(\mathcal{X})$ of all nef classes in $A$.  Then $\mathcal{P}$ is a rational polyhedron whose recession cone is $\Nef_{1}(\mathcal{X}_{\eta})$.

In particular, the lattice points $\mathcal{P}_{\mathbb{Z}}$ are contained in a finite union of translates of $\Nef_{1}(\mathcal{X}_{\eta})_{\mathbb{Z}}$.
\end{guid}

Over a field of characteristic zero, this is exactly Corollary \ref{coro:polyhedron}.

\subsubsection{Dimension}
The next step is to ``count points'' on these irreducible components.  Of course, we must explain what we mean by ``counting points'' when our ground field is not finite.  If we want a coarse interpretation, we can simply try to measure the dimensions and numbers of the irreducible components $M \subset \Sec(\mathcal{X}/B)$ parametrizing sections with a given numerical class.

Note that Batyrev's heristic \ref{simp:three} is not true -- a family of sections need not have the expected dimension.  However, we can expect all Manin components to have the expected dimension.

\begin{guid} \label{guid:four}
Let $\pi: \mathcal{X} \to B$ be a Fano fibration over an algebraically closed field $\kappa$ of characteristic $0$.  Every Manin component $M$ of $\Sec(\mathcal{X}/B)$ in $\alpha + \mathcal P$ will generically parametrize relatively free sections.  In particular, such components $M$ have the expected dimension. 
\end{guid}

Over $\mathbb C$, Theorem~\ref{theo:maintheorem1} shows that any component of $\Sec(\mathcal X/B)$ only parametrizing non-free sections is an accumulating component.

\subsubsection{Manin components with fixed numerical class}

After identifying the set of possible numerical classes and the dimension of the corresponding moduli spaces, we next need to address how many irreducible components of $\Sec(\mathcal{X}/B)$ represent each numerical class.  There are two issues.  First, which nef numerical classes in $A$ are actually represented by sections?  Recall that for any given section we can construct more by gluing on $\pi$-vertical free curves in a fiber $F$ and smoothing.  If we only care about the asymptotic behavior, the key question is whether every coset of $N_{1}(F)_{\mathbb{Z}}$ inside of $A_{\mathbb{Z}}$ can be represented by a section.  We expect the answer to be yes, and in particular, that any ``sufficiently positive'' nef class in $A$ is represented by a section. 

Second, what is the maximal number of irreducible components of $\Sec(\mathcal{X}/B)$ which represent a fixed nef numerical class?  This number can be unbounded if we allow exceptional components, so  we restrict our attention to the Manin components.  Again, we expect that when the class of the section is ``sufficiently positive'' the number of irreducible Manin components of a given class will exhibit systematic behavior.

The following conjecture subsumes both expectations:

\begin{guid} \label{guid:three}
Let $\pi: \mathcal{X} \to B$ be a Fano fibration over an algebraically closed field $\kappa$ of characteristic $0$. Let $\mathcal{P}$ denote the convex hull in $N_{1}(\mathcal{X})$ of all $\Nef_{1}(\mathcal{X}) \cap A$.  There is some ``sufficiently positive'' $\alpha \in \Nef_{1}(\mathcal{X}_{\eta})$ such that every algebraic equivalence class of curves whose numerical class lies in $\alpha + \mathcal{P}$ is represented by a unique Manin component.
\end{guid}

Note that the translate $\alpha +\mathcal{P}$ will ``asymptotically contain 100\%'' of the lattice points in $\mathcal{P}$, so considering such a translate should suffice for the purposes of estimating the asymptotics in Manin's Conjecture.  Guiding Conjecture \ref{guid:three} has been established for del Pezzo surface fibrations in \cite{LT21b}.  Due to the numerous classification results for rational curves on Fano varieties, the conjecture has also been verified for many fibrations of the form $X \times \mathbb{P}^{1} \to \mathbb{P}^{1}$; one notable recent example is \cite{BJ22} which establishes Guiding Conjecture \ref{guid:three} for Fano threefolds. 

\begin{remark}
Note that it is \emph{algebraic} equivalence and not \emph{numerical} equivalence which appears in Geometric Manin's Conjecture \ref{guid:three}.  Thus the conjecture implicitly relies on an understanding of the relationship between these two equivalence notions.  We will briefly discuss this relationship over the ground field $\mathbb{C}$.  Recall that the Griffiths group of a smooth projective variety measures the difference between algebraic and homological equivalence: it is the quotient of the group of homologically trivial $1$-cycles by the group of algebraically trivial $1$-cycles. 

The underlying theoretical framework depends on two questions raised by Voisin.  Suppose $X$ is a smooth projective Fano variety.  Then:
\begin{enumerate}
\item Is the Griffiths group of $1$-cycles on $X$ trivial?  
\item Does the Integral Hodge Conjecture hold for $1$-cycles on $X$?
\end{enumerate}
There are no known counterexample to either question.  Since a Fano fibration over a curve has a similar ``motive'' to a Fano variety, we can optimistically hope that both questions have affirmative answers for Fano fibrations $\mathcal{X}$.  (This is true when $\mathcal{X}$ has dimension $\leq 3$ by \cite{BS83} and \cite{Voisin06} respectively.)

Assuming that both questions above have an affirmative answer for Fano fibrations, we find that the number of algebraic equivalence classes of curves representing a single numerical curve class is the same as $|H_{2}(\mathcal{X},\mathbb{Z})_{tors}|$.  Since $H^{2,0}(\mathcal{X}) = 0$ by \cite[Theorem 10.17]{Voi03}, Poincar\'e duality implies that
\begin{align*}
 |H_{2}(X,\mathbb{Z})_{tors}| = |H^{3}(X,\mathbb{Z})_{tors}| = |\Br(\mathcal{X})|.
\end{align*}
We conclude that (under our assumption) there are $|\Br(\mathcal{X})|$ different algebraic equivalence classes representing a single numerical class.  (Note that in this context it is $\Br(\mathcal{X})$ and not $\Br(\mathcal{X}_{\eta})$ that controls the number of irreducible components; see \cite[Example 8.6]{LT21a}.)
\end{remark}

\begin{remark}
There is another question hidden in Geometric Manin's Conjecture \ref{guid:three}: are the extremal rays of the nef cone of a Fano variety generated by free rational curves?  Only in this situation can we hope to use the gluing structure to obtain the tightest link between the structure of nef classes in $\Nef_{1}(\mathcal{X}_{\eta})$ and the existence of sections of $\pi$.

The existence of free curves representing extremal rays is roughly equivalent to the existence of free rational curves in the smooth locus of a log Fano variety.  This well-known question has been answered for complex surfaces in \cite{KM99} but remains open in general. 
\end{remark}

\subsubsection{Stability}

Finally, we turn to the problem of counting points on the Manin components $M \subset \Sec(\mathcal{X}/B)$.  (As discussed in Remark \ref{rema:sufficientlyfree} we may want to remove additional subloci of $M$, but we will not address this subtlety here.)  There are two ways to make sense of ``point counts'' on $M$ over an arbitrary ground field.

The first option is to study the class of $M$ in the Grothendieck ring $K_{0}(\mathrm{Var}_{k})$.  This option was pioneered in \cite{Bourqui09} and \cite{CLL16} and has seen recent development in the papers, e.g., \cite{Bilu, BB23, Faisant23, FaisantANT}.  Let $\mathcal{M}_{k}$ denote the localization of $K_{0}(\mathrm{Var}_{k})$ along the multiplicative system generated by $\mathbb{L}$ and $1-\mathbb{L}^{a}$ as we vary $a \in \mathbb{N}$.  Loosely speaking, we expect the classes of the irreducible components $M$ in $\mathcal{M}_{k}$ to ``stabilize'' as the degree increases to infinity.   For example, one might hope that the sequence of renormalized classes $\mathbb{L}^{-\dim(M)} [M]$ converges in one of the topologies discussed in \cite{BDH22}; alternatively, one could ask for the convergence of an associated zeta function.  The correct formulation of this principle is currently a topic of ongoing research.

The second option is to study the cohomology groups of the moduli spaces $M$.  The first steps in this direction predate Manin's Conjecture.  In his influential paper \cite{Segal79}, Segal proved that for any Riemann surface $B$ the cohomology groups of the pointed moduli spaces $\Mor_{*}(B,\mathbb{P}^{n})_{d}$ of degree $d$ morphisms stabilize as $d \to \infty$ to the cohomology groups of the pointed topological mapping spaces $\mathrm{Map}_{*}(B,\mathbb{P}^{n})$.  Segal's results were later extended to other Fano varieties besides projective space (e.g., \cite{Kirwan86, Guest95}).  Despite some remarkable progress, a complete picture is still out of reach.

In their work on Hurwitz schemes, Ellenberg and Venkatesh suggested that homological stability can be used as a tool to verify Manin's Conjecture over global function fields.  Suppose that one had a stability result for the \'etale cohomology groups of the moduli space of sections of a Fano fibration in characteristic $p$.  \cite{EVW16, ETW17, EL23}
use the Grothendieck-Lefschetz formula to show that if the ``stable range'' of cohomology increases linearly in the degree (and the Frobenius actions are the same), then number of points on the corresponding irreducible components satisfies a nice asymptotic formula.

In summary, we have a vague principle:

\begin{guid} \label{guid:five}
Manin components of $\Sec(\mathcal{X}/B)$ exhibit motivic or homological stability as the degree increases.
\end{guid}

\section{Numerical properties of sections} \label{sect:cone}

In this section we prove several fundamental properties about the numerical classes of sections of Fano fibrations, culminating in a polyhedrality statement (Theorem \ref{theo:polyhedrality}).  While we do not use Theorem \ref{theo:polyhedrality} in later sections, some of the steps in the proof will be used again later, particularly the Cone Theorem (Theorem \ref{theo:nefconetheorem}) and the existence of good models of Fano fibrations (Lemma \ref{lemm:birationalantiamplemodel}).

\subsection{Numerical equivalence on Fano fibrations}

Suppose $\pi: \mathcal{X} \to B$ is a Fano fibration.  In this section we give some reminders about the basic properties of numerical equivalence and the cone of nef curves on $\mathcal{X}$.

\begin{lemma}
\label{lemm:relativeprops}
Let $\pi: \mathcal{X} \to B$ be a Fano fibration.  Then for every smooth fiber $F$ the space $N^{1}(F)_{\mathbb{R}}$ has the same dimension.  Furthermore, if $L$ is a $\mathbb{Q}$-Cartier divisor on $\mathcal{X}$ then the following are equivalent:
\begin{enumerate}
\item $L|_{F}$ is ample for some smooth Fano fiber $F$.
\item $L|_{F}$ is ample for all smooth Fano fibers $F$.
\item $L|_{\mathcal{X}_{\eta}}$ is ample.
\end{enumerate}
The analogous statement is true for nefness, for bigness, and for pseudo-effectiveness.
\end{lemma}

\begin{proof}
It follows from \cite[IV.3.5.3 Theorem]{Kollar} that every smooth fiber is rationally connected.  Then Hodge theory tells us that $N^{1}(F)_{\mathbb{R}}$ coincides with $H^{2}(F,\mathbb{R})$.  Since the fibers of a smooth proper family are topologically isotrivial, we deduce that this space has constant dimension. The equivalence of the three conditions for ampleness and nefness follows from \cite[Theorem 1]{Wisniewski09}. For the equivalence of the three conditions for bigness and pseudo-effectiveness, see the paragraph before \cite[Theorem 6.8]{dFH11}.
\end{proof}

We have a restriction map $N^{1}(\mathcal{X})_{\mathbb R} \to N^{1}(\mathcal{X}_{\eta})_{\mathbb R}$ which is surjective (since any Weil divisor on $\mathcal{X}_{\eta}$ can be extended to a Weil divisor on $\mathcal{X}$).  Since linear and numerical equivalence coincide for $\mathcal{X}_{\eta}$, the kernel of the restriction map is spanned by $\pi$-vertical divisors. %\sho{How \cite[Theorem 10.19]{Voisin07} implies this statement?}  \brian{In Voisin's notation $X=\mathcal{X}$, $Y = B$, and we are looking at divisors so $k = \dim(\mathcal{X})-1$.  Choose a finite set of divisors which generate the Picard group of $\mathcal{X}_{\eta}$ and take their closures in $\mathcal{X}$.  Let $X'$ denote the union of these divisors.  This choice of $X'$ satisfies (*) for every divisor $Z$, and the conclusion of the theorem gives our desired statement.}\sho{Voisin's book assumes $X'$ to be a subvariety. Does she allow $X'$ to be reducible? Actually using the fact that numerical equivalence and $\mathbb Q$-linear equivalence are equivalent, we can give a direct proof of this fact. Shall we prove it?} \brian{In my perspective this is a well-known statement so we don't need to prove it carefully.  But we can if you like.}\sho{I am fine with no proof. But not so sure whether Voisin's book is a right reference.} \brian{We can use Theorem 1.1 in the paper by Gounelas and Javanpeykar instead.} %By the Invariant Cycles Theorem, $N^{1}(\mathcal{X}_{\eta})$ is isomorphic to the monodromy-invariant part of $N^{1}(F)$ for a general fiber $F$.
Dually, we have an injective map $N_{1}(\mathcal{X}_{\eta})_{\mathbb R} \to N_{1}(\mathcal{X})_{\mathbb R}$.  Henceforth we will not distinguish between $N_{1}(\mathcal{X}_{\eta})$ and its image under this inclusion.  Note that $N_{1}(\mathcal{X}_{\eta})_{\mathbb R}$ is precisely the subspace of $N_{1}(\mathcal{X})_{\mathbb R}$ consisting of classes $\alpha$ that satisfy $E \cdot \alpha = 0$ for every $\pi$-vertical divisor $E$.

\begin{lemma} \label{lemm:relativenef}
Let $\pi: \mathcal{X} \to B$ be a Fano fibration.  Let $V$ be the codimension $1$ subspace of $N_{1}(\mathcal{X})$ consisting of classes with vanishing intersection against a general fiber $F$.  Then $\Nef_{1}(\mathcal{X}) \cap V = \Nef_{1}(\mathcal{X}_{\eta})$.
\end{lemma}

\begin{proof}
The containment $\supset$ is clear.  Conversely, suppose that $\alpha$ is a nef class contained in $V$.  Since $F \cdot \alpha = 0$, it is also true that $E \cdot \alpha = 0$ for every $\pi$-vertical divisor $E$.  Thus $\alpha \in N_{1}(\mathcal{X}_{\eta})$.  Note that the restriction map $N^{1}(\mathcal{X}) \to N^{1}(\mathcal{X}_{\eta})$ induces a surjection from the effective cone of divisors on $\mathcal{X}$ to the effective cone of divisors on $\mathcal{X}_{\eta}$.  Since $\pi$ is a Fano fibration, the pseudo-effective cone and effective cone for $\mathcal{X}_{\eta}$ coincide.  Altogether we obtain a surjection $\Eff^{1}(\mathcal{X}) \to \Eff^{1}(\mathcal{X}_{\eta})$.  
%\sho{Implicitly are we using the fact that $\mathcal X_\eta$ is Fano so that its pseudo-effective cone coincides with the effective cone? If so, we should emphasize that.} \brian{No, I don't think so.  This statement is true for any fibration.  We just take the closure of an effective class on the generic fiber to get an effective class on the whole thing.}\sho{I agree for effective classes, but it is not clear for pseudo-effective classes. Note that the image of a closed set by projection does not need to be closed. We need some polyhedrality of the cone to claim this.} \brian{Ok, I changed this.} 
We conclude that $\alpha$ must be a nef class on $\mathcal{X}_{\eta}$.
\end{proof}

\subsection{The Cone Theorem for nef curves}

We will also need the following version of the Cone Theorem for nef curves.  This theorem was proved by \cite{Araujo10} conditional on the Borisov-Alexeev-Borisov Conjecture which has subsequently been proved in \cite{Birkar21}.

\begin{theorem}[\cite{Araujo10}] \label{theo:conetheorem}
Let $X$ be a normal $\mathbb{Q}$-factorial projective variety and let $\Delta$ be an effective $\mathbb{Q}$-Cartier divisor on $X$ such that $(X,\Delta)$ is $\epsilon$-lc.  There is a constant $\zeta = \zeta(\dim(X),\epsilon)$ such that
\begin{equation*}
\Eff_{1}(\mathcal{X})_{K_{\mathcal{X}}+\Delta \geq 0} + \Nef_{1}(\mathcal{X}) = \Eff_{1}(\mathcal{X})_{K_{\mathcal{X}}+\Delta \geq 0} + \sum_{i} \mathbb{R}_{\geq 0}[C_{i}]
\end{equation*}
where $\{C_{i}\}$ is a countable collection of curves which satisfy $0 < -(K_{X} + \Delta) \cdot C_{i} \leq \zeta$.  If $\Delta$ is big, then the set $\{ C_{i} \}$ is finite.
\end{theorem}

For Fano fibrations, the Cone Theorem for nef curves allows us to isolate the behavior of vertical curves:

\begin{theorem}  \label{theo:nefconetheorem}
Let $\pi: \mathcal{X} \to B$ be a flat morphism from a normal $\mathbb{Q}$-factorial projective variety $\mathcal{X}$ to a smooth projective curve $B$.  Suppose that $\Delta$ is an effective $\mathbb{Q}$-Cartier divisor on $\mathcal{X}$ such that $\Delta$ is $\pi$-relatively big and $(\mathcal{X},\Delta)$ is $\epsilon$-lc.  Let $F$ denote a general fiber of $\pi$.  There is a positive integer $m = m(\dim(\mathcal{X}),\epsilon)$ such that we have an equality
\begin{equation*}
\Eff_{1}(\mathcal{X})_{K_{\mathcal{X}}+\Delta+mF \geq 0} + \Nef_{1}(\mathcal{X}) = \Eff_{1}(\mathcal{X})_{K_{\mathcal{X}}+\Delta+mF \geq 0} + \sum_{i} \mathbb{R}_{\geq 0}[C_{i}]
\end{equation*}
where $\{ C_{i} \}$ is a finite set of $\pi$-vertical moving curves which satisfy $0 < -(K_{X} + \Delta + mF) \cdot C_{i} \leq m$.
\end{theorem}

\begin{proof}
Since $\Delta$ is effective and $\pi$-relatively big, we see that $\Delta + \delta F$ is big for any $\delta > 0$.
By choosing $\delta$ sufficiently small we may ensure that $\delta < 1$ and that $(\mathcal{X},\Delta+\delta F)$ is $\epsilon/2$-lc.
Applying Theorem \ref{theo:conetheorem} there is a constant $\zeta = \zeta(\dim(\mathcal{X}),\epsilon)$ such that
\begin{equation*}
\Eff_{1}(\mathcal{X})_{K_{\mathcal{X}} + \Delta + \delta F \geq 0} + \Nef_{1}(\mathcal{X}) = \Eff_{1}(\mathcal{X})_{K_{\mathcal{X}} + \Delta + \delta F \geq 0} + \sum_{j} [C_{j}]
\end{equation*}
where the $C_{j}$ are a finite set of movable curves satisfying $0 \leq -(K_{\mathcal{X}'} + \Delta + \delta F) \cdot C_{i} \leq \zeta$.  Choose a positive integer $m > \zeta +1$.  Then $(K_{\mathcal{X}} + \Delta + mF) \cdot C_{j} > 0$ for every one of our movable curves $C_{j}$ that dominates $B$ under $\pi$.  Thus we have
\begin{equation*}
\Eff_{1}(\mathcal{X})_{K_{\mathcal{X}} + \Delta + \delta F \geq 0} + \Nef_{1}(\mathcal{X}) = \Eff_{1}(\mathcal{X})_{K_{\mathcal{X}} + \Delta + mF \geq 0} + \sum_{i} \mathbb{R}_{\geq 0}[C_{i}]
\end{equation*}
where now the $C_{i}$ are $\pi$-vertical and still satisfy $0 \leq -(K_{\mathcal{X}'} + \Delta + m F) \cdot C_{i} \leq \zeta < m$. 
\end{proof}

In particular, when $\pi$ is a Fano fibration then Lemma \ref{lemm:relativenef} gives an equality
\begin{equation*}
\Eff_{1}(\mathcal{X})_{K_{\mathcal{X}}+H+mF \geq 0} + \Nef_{1}(\mathcal{X})  = \Eff_{1}(\mathcal{X})_{K_{\mathcal{X}}+H+mF \geq 0} +  \Nef_{1}(\mathcal{X}_{\eta}).
\end{equation*}

\subsection{Classes of nef sections}

We will need more specific information about the classes of nef sections.  
%\brian{Right now there might be some overlap with the next lemma and Proposition \ref{prop:birmodelposL}.  If so, we should remove the latter one and refer here instead...} \brian{Looking again they are not so similar, so I recommend we leave it as it is.}

\begin{lemma} \label{lemm:birationalantiamplemodel}
Let $\pi: \mathcal{X} \to B$ be a Fano fibration.  Then there is a morphism $\pi': \mathcal{X}' \to B$ satisfying the following conditions:
\begin{enumerate}
\item $\mathcal{X}'$ is a normal projective variety with $\mathbb{Q}$-factorial singularities,
\item $\pi'$ is birationally equivalent to $\pi$ and the generic fibers of $\pi$ and $\pi'$ are isomorphic,
\item there is an effective $\mathbb{Q}$-Cartier divisor $D'$ on $\mathcal{X}'$ such that $(\mathcal{X}',D')$ is klt and $D'$ is relatively $\mathbb{Q}$-linearly equivalent to $-K_{\mathcal{X}'/B}$,
\item $-K_{\mathcal{X}'/B}$ is $\pi$-relatively big and nef.
\end{enumerate}
\end{lemma}

\begin{proof}
We construct $\pi'$ in two steps.  First, we claim that there is a morphism $\widetilde{\pi}: \widetilde{\mathcal{X}} \to B$ satisfying conditions (1)-(3). % and
%\begin{itemize}
%\item[(4)'] There is an effective $\mathbb{Q}$-divisor $\widetilde{D}$ that is $\pi$-relatively big such that $(\widetilde{\mathcal{X}},\widetilde{D})$ is klt and $\widetilde{D}$ is $\mathbb{Q}$-linearly equivalent over $B$ to $-K_{\widetilde{\mathcal{X}}/B}$.
%\end{itemize}
To achieve this, choose an effective $\mathbb{Q}$-divisor $D$ on $\mathcal{X}$ such that $D|_{\mathcal{X}_{\eta}} \sim_{\mathbb{Q}} -K_{\mathcal{X}_{\eta}}$, the coefficients of $D$ are in $(0,1)$, and the restriction to the generic fiber has SNC support.  Let $\pi_{\mathcal{W}}: \mathcal{W} \to B$ be a log resolution of $(\mathcal{X},\Supp(D))$ that is an isomorphism along the generic fiber of $\pi$ and let $D_{\mathcal{W}}$ be the strict transform of $D$.
%\sho{We should insist that the generic fiber of $\mathcal W$ is isomorphic to the generic fiber of $\mathcal X$, right?}  \brian{You are right, I added it.} 
We now run the relative $(\mathcal{W},D_{\mathcal{W}})$-MMP over $B$ with scaling.  Since the restriction of $K_{\mathcal{W}} + D_{\mathcal{W}}$ to the generic fiber is numerically trivial, we see that the outcome of the MMP will be a birational model $\phi: \mathcal{W} \dashrightarrow \widetilde{\mathcal{X}}$ such that $K_{\widetilde{\mathcal{X}}} + \phi_{*}D_{\mathcal{W}} \sim_{\mathbb{Q},B} 0$.  For convenience we define $\widetilde{D} := \phi_{*}D_{\mathcal{W}}$.

We next run the MMP one more time for $\widetilde{\mathcal{X}}$.  Choose a rational $\epsilon >0$ sufficiently small so that $(\widetilde{\mathcal{X}},(1+\epsilon)\widetilde{D})$ is klt.  Run the MMP with scaling for this pair.  Since the restriction of $\widetilde{D}$ to the generic fiber of $\pi$ is ample, we see that the outcome of the MMP is a birational map $\psi: \widetilde{\mathcal{X}} \dashrightarrow \mathcal{X}'$ such that $K_{\mathcal{X}'}+(1+\epsilon)\psi_{*}\widetilde{D}$ is relatively big and nef.  %\brian{Ok, I weakened the conclusion here.}  
%\sho{Is this true? Don't we need to take the canonical model to obtain a relatively ample model? A MMP gives only relatively big and nef I believe.} \brian{I am taking a canonical model.  This is a birational map.  (Note that I am also taking a canonical model in the first step, to achieve $\sim 0$.)  Did I make a mistake about the singularities or $\mathbb{Q}$-factoriality?  I will think about it...}\sho{To me the outcome of a MMP is a minimal model, but not the canonical model, so after applying a MMP, we should take the canonical model. To achieve $\sim 0$, a minimal model is enough but to achieve a relatively ample model, one needs to take the canonical model.} \brian{Sorry, of course you are right.  By passing to a canonical model, I think we might lose $\mathbb{Q}$-factoriality.  This is not a big deal, but I will need to generalize a couple statements to the non-$\mathbb{Q}$-factorial setting.  I will do it later.}  \brian{Actually this is more serious than I thought.  The cone theorem for nef curves does not rely on $\mathbb{Q}$-factoriality.  But the birational argument in Theorem 5.7 does.  I will think about it.}  
But by construction $K_{\mathcal{X}'} + (1+\epsilon)\psi_{*}\widetilde{D}$ is relatively $\mathbb{Q}$-linearly equivalent to $-\epsilon K_{\mathcal{X}'}$.  Thus $\mathcal{X}'$ equipped with the divisor $D' := \psi_{*}\widetilde{D}$ satisfies all the desired conditions.
\end{proof}

%\begin{corollary} \label{coro:relativeamplebirationalmodel}
%Let $\pi: \mathcal{X} \to B$ be a Fano fibration.  There is a birational model $\phi: \mathcal{X}^{+} \to \mathcal{X}$ that is an isomorphism on the generic fiber of $\pi$  such that $\mathcal{X}^{+}$ is smooth and there is a $\pi \circ \phi$-vertical effective divisor $E$ such that $-K_{\mathcal{X}^{+}/B}+E$ is ample.
%\end{corollary}

%\begin{proof}
%Let $\mathcal{X}'$ be the birational model of $\mathcal{X}$ and let $D'$ be the effective $\mathbb{Q}$-Cartier divisor constructed in Lemma \ref{lemm:birationalantiamplemodel}.  Choose a rational $\epsilon > 0$ sufficiently small so that $(\mathcal{X}',(1+\epsilon)D')$ is klt.  Since $-\epsilon K_{\mathcal{X}'/B} \sim_{\mathbb{Q}} K_{\mathcal{X}'/B} + (1+\epsilon)D'$ is $\pi'$-relatively big and nef, by taking a relative canonical model we see it is the pullback of a relatively ample divisor under a birational map.

%Let $\mathcal{X}^{+}$ be a smooth projective variety resolving the rational map $\psi: \mathcal{X} \dashrightarrow \mathcal{X}'$.  Since $\psi$ is an isomorphism along the generic fiber of $\pi$, we may ensure that the birational map $\mathcal{X}^{+} \to \mathcal{X}$ is also an isomorphism over the generic fiber of $\pi$.  It is then clear that $\mathcal{X}^{+}$ has the desired property.
%\sho{I think we should explain how to construct $E$ for $\mathcal X^+$.}
%\end{proof}

We are now prepared to prove the main structural results for sections of a Fano fibration.

\begin{theorem} \label{theo:compactnessforaffineslice}
Let $\pi: \mathcal{X} \to B$ be a Fano fibration and let $F$ be a general fiber of $\pi$.  Let $A \subset N_{1}(\mathcal{X})$ be the affine plane consisting of all numerical classes $\alpha$ which satisfy $\alpha \cdot F = 1$.  Then there is a compact subset $\mathcal{T} \subset A$ such that $\Nef_{1}(\mathcal{X}) \cap A$ is contained in $\mathcal{T} + \Nef_{1}(\mathcal{X}_{\eta})$.
\end{theorem}

\begin{proof}
Let $\pi': \mathcal{X}' \to B$ be a birational model of $\pi$ as in Lemma \ref{lemm:birationalantiamplemodel}.  We first prove the analogous statement for $\mathcal{X}'$.  We will let $A' \subset N_{1}(\mathcal{X}')$ denote the affine subspace consisting of numerical classes with intersection $1$ against a fiber.

By relative Wilson's Theorem (see e.g.~\cite[Theorem 3.2]{LT17}) there is an effective divisor $E$ such that $(\mathcal{X}',E)$ is klt and $-(K_{\mathcal{X}'} + E)$ is $\pi$-relatively ample.  Choose a small ample $\mathbb{Q}$-Cartier divisor $H$ so that $(\mathcal{X}',E+H)$ is still klt and $-(K_{\mathcal{X}'} + E + H)$ is still $\pi$-relatively ample.  Applying Theorem \ref{theo:nefconetheorem} there is a constant $m$ such that
\begin{equation*}
\Nef_{1}(\mathcal{X}') + \Eff_{1}(\mathcal{X}')_{K_{\mathcal{X}'} + E + H + mF \geq 0} = \Eff_{1}(\mathcal{X}')_{K_{\mathcal{X}'} + E + H  + mF \geq 0} + \sum_{i} [C_{i}]
\end{equation*}
where the $C_{i}$ are a finite set of $\pi$-vertical movable curves satisfying $0 < -(K_{\mathcal{X}'} + E + H) \cdot C_{i} \leq m$. %\sho{This should be $<0$.} \brian{No, this is correct.} %By possibly increasing $m$ we may suppose $m \geq 2\dim(\mathcal{X})$.  Thus by Theorem \ref{theo:psefconetheorem} we can also write
%\begin{equation*}
%\Eff_{1}(\mathcal{X}') = \Eff_{1}(\mathcal{X}')_{K_{\mathcal{X}'} + H + mF \geq 0} + \sum \mathbb{R}_{\geq 0} [R_{i}]
%\end{equation*}
%where the $R_{i}$ are all $\pi'$-vertical curves.  

Let $\mathcal{T}'$ denote the intersection of the cone $\Eff_{1}(\mathcal{X}')_{K_{\mathcal{X}'} + E + H + mF \geq 0}$ with $A'$. Since $-(K_{\mathcal X'} + E)$ is $\pi$-relatively ample, there is an integer $k$ such that $-(K_{\mathcal X'} + E+H) + kF$ is ample.
Then the intersection of a class in $\mathcal T'$ against $-(K_{\mathcal X'} + E+H) + kF$ is bounded by $k+m$.
%Since $-(K_{\mathcal{X}'/B} + E + H)$ is relatively ample, the cone $\Eff_{1}(\mathcal{X}')_{K_{\mathcal{X}'} + E + H + mF \geq 0}$ only intersects the subspace of classes with vanishing intersection against $F$ at the point $0$.  
This implies that $\mathcal{T}'$ is compact. %\sho{This reasoning looks strange to me. There is a more direct argument as I explained before.}
Furthermore by construction it is clear that $\Nef_{1}(\mathcal{X}') \cap A' \subset \mathcal{T}' + \sum_{i} [C_{i}]$. 

We now pass the desired statement to $\mathcal{X}$ using a birational argument.  Consider a diagram of $B$-morphisms
\begin{equation*}
\xymatrix
   { & \mathcal{W} \ar[dl]_{g} \ar[dr]^{g'} & \\
     \mathcal{X} & & \mathcal{X}'
   }
\end{equation*}
where $\mathcal{W}$ is smooth projective and $g$ and $g'$ are birational maps which are isomorphisms along the generic fiber.  We next prove that $\mathcal{W}$ satisfies the desired statement.  Let $E_{1},\ldots,E_{r}$ denote the divisors contracted by the map $g': \mathcal{W} \to \mathcal{X}'$.  Consider the induced linear map $g'_{*}: N_{1}(\mathcal{W}) \to N_{1}(\mathcal{X}')$.  We let $\mathcal{T}_{\mathcal{W}}$ denote the subset of $g_{*}'^{-1}(\mathcal{T}')$ consisting of those classes $\alpha$ satisfying $0 \leq E_{i} \cdot \alpha \leq 1$ for every $i$.  By construction $\mathcal{T}_{\mathcal{W}}$ is a compact subset of $N_{1}(\mathcal{W})$.

We claim that 
\begin{equation*}
\Nef_{1}(\mathcal{W}) \cap A_{\mathcal{W}} \subset \mathcal{T}_{\mathcal{W}} + \Nef_{1}(\mathcal{W}_{\eta}).
\end{equation*}
Indeed, note that every $\pi$-vertical moving curve $C_{i}$ on $\mathcal{X}'$ is the image of some $\pi$-vertical moving curve $T_{i}$ on $\mathcal{W}$.  Thus for every element of $\beta \in \Nef_{1}(\mathcal{W}) \cap A_{\mathcal W}$ there is some sum $\gamma = \sum t_{i} [T_{i}]$ such that $\rho_{*}(\beta) - \rho_{*}\gamma \in \mathcal{T}'$.  Equivalently, $\beta - \gamma \in g_{*}'^{-1}(\mathcal{T}')$.  But since each $E_{i}$ is $\pi$-vertical, we know that $0 \leq \beta \cdot E_{i} \leq 1$ for every $i$.  Since $\gamma$ has vanishing intersection against every $E_{i}$, we see that in fact $\beta - \gamma \in \mathcal{T}_{\mathcal{W}}$.

Finally, we set $\mathcal{T} = g_{*}\mathcal{T}_{\mathcal{W}}$.  Since the map $g_{*}: \Nef_{1}(\mathcal{W}) \to \Nef_{1}(\mathcal{X})$ is surjective and $g$ is an isomorphism along the generic fibers, we obtain the desired statement for $\mathcal{X}$.
\end{proof}

\begin{theorem} \label{theo:polyhedrality}
Let $\pi: \mathcal{X} \to B$ be a Fano fibration.  Let $A$ denote the affine plane consisting of numerical classes with intersection $1$ against a general fiber $F$.  There is a finite set of nef classes $\beta_{1},\ldots,\beta_{s} \in A \cap N_{1}(X)_{\mathbb{Z}}$
%\sho{Actually we can take $\beta_i$'s in $A \cap \Nef_{1}(\mathcal{X})_{\mathbb{Z}}$, right?} \brian{How is this different than what is written currently?} 
such that every class in $A \cap \Nef_{1}(\mathcal{X})_{\mathbb{Z}}$ has the form $\beta_{i} + \gamma$ for some index $i$ and for some $\gamma \in \Nef_{1}(\mathcal{X}_{\eta}) \cap N_{1}(\mathcal{X})_{\mathbb{Z}}$. 
\end{theorem}

\begin{proof}
Let $A_{\mathbb{Z}}$ denote the intersection $A \cap N_{1}(\mathcal{X})_{\mathbb{Z}}$.  Since $A$ is a translate of a rational subspace of $N_{1}(\mathcal{X})$ and contains a curve class, there is an isomorphism between $A_{\mathbb{Z}}$ and $\mathbb{Z}^{\dim A}$.  Choosing such an identification gives $A$ the structure of a vector space and $A_{\mathbb{Z}}$ the structure of a lattice in this vector space.  Up to translation $\Nef_{1}(\mathcal{X}_{\eta})$ can be identified with a rational polyhedral cone inside the vector space $A$.  By Theorem \ref{theo:compactnessforaffineslice} there is a compact set $\mathcal{T} \subset A$ such that
\begin{equation*}
A \cap \Nef_{1}(\mathcal{X})_{\mathbb{Z}} \subset \mathcal{T} + \Nef_{1}(\mathcal{X}_{\eta}).
\end{equation*}
The desired statement then follows from \cite[Theorem 1.1.(a)]{HW07}.
\end{proof}

\begin{corollary} \label{coro:polyhedron}
Let $\pi: \mathcal{X} \to B$ be a Fano fibration.  Let $\mathcal{P}$ denote the convex hull of all integral nef curve classes $\alpha \in N_{1}(\mathcal{X})_{\mathbb{Z}}$ satisfying $F \cdot \alpha = 1$ for a general fiber $F$ of $\pi$.  Then $\mathcal{P}$ is a rational polyhedron whose recession cone is isomorphic to $\Nef_{1}(\mathcal{X}_{\eta})$.
\end{corollary}

There is an important subtlety: in the setting of Theorem \ref{theo:polyhedrality} and Corollary \ref{coro:polyhedron}, it is not necessarily true that the intersection $A \cap \Nef_{1}(\mathcal{X})$ is a polyhedron.  This is demonstrated by the following example.

\begin{example} \label{exam:notpolyhedral}
Let $\mathcal{X}$ be the blow-up of $\mathbb{P}^{2}$ at nine very general points (i.e.~a set of nine points that is in the complement of a countable union of proper closed subvarieties of $\Sym^{9}(\mathbb{P}^{2})$).  If we let $\mathbb{F}_{1}$ denote the blow-up of $\mathbb{P}^{2}$ at a point and consider the intermediate blow-up $\mathcal{X} \to \mathbb{F}_{1} \to \mathbb{P}^{2}$, the projective bundle map $\mathbb{F}_{1} \to \mathbb{P}^{1}$ also gives $\mathcal{X}$ the structure of a Fano fibration $\pi: \mathcal{X} \to \mathbb{P}^{1}$.

Let $A$ be the affine subspace of curve classes $\alpha$ such that $F \cdot \alpha = 1$ for a general fiber $F$ of $\pi$.  We claim that $A \cap \Nef_{1}(\mathcal{X})$ is not polyhedral.  Indeed, by \cite{Borcea91} $\Nef_{1}(\mathcal{X})$ has countably many extremal rays not contained in $N_{1}(\mathcal{X}_{\eta})$.  By taking the intersection of these rays with $A$ we obtain a countable number of points which are extremal in $A \cap \Nef_{1}(\mathcal{X})$.  %\brian{I think this works.  In this example we can work out the nef cone very explicitly, see e.g.~\url{https://www.math.stonybrook.edu/~cschnell/pdf/notes/minusone.pdf}.  I started writing it down but got tired and frustrated.  I hope it is automatic and if so we should just take the easy way out.}
\end{example}

\section{Grauert-M\"ulich} \label{sect:gm}

For a good fibration $\pi: \mathcal{Z} \to B$ the deformation theory of a section $C$ is controlled by the Harder-Narasimhan filtration of the restriction $T_{\mathcal{Z}/B}|_{C}$.  In this section, we show that (under certain hypotheses) the Harder-Narasimhan filtration of $T_{\mathcal{Z}/B}|_{C}$ is ``approximately'' the restriction of the $[C]$-Harder-Narasimhan filtration of $T_{\mathcal{Z}/B}$.  Due to the similarity to the Grauert-M\"ulich theorem (\cite{GM75}) describing the restriction of semistable bundles to lines in $\mathbb{P}^{n}$, we will refer to such statements as ``Grauert-M\"ulich'' results.
The material in this section is motivated by \cite[Section 3]{PatelRiedlTseng} and by \cite[Chapter II, Section 2]{OSS80}.

\

Suppose that $Z$ is a smooth projective variety and $W$ is a variety parametrizing a family of maps $s: C \to Z$.   Let $\mathcal{E}$ be a torsion-free sheaf on $Z$ and let $\mathcal{F}$ be a term in the relative Harder-Narasimhan filtration of $\mathcal{E}$ pulled back to the universal family over an open subset of $W$. We would like to determine when $\mathcal{F}$ is the pullback of a sheaf $\mathcal{F}_{Z}$ from $Z$.  In this case we can expect $\mathcal{F}_{Z}$ to be a term in the Harder-Narasimhan filtration of $\mathcal{E}$ with respect to the numerical class of the curves $s_{*}C$.  

In Section \ref{sect:descent} we prove a general criterion for determining when a torsion-free sheaf on a variety is isomorphic to a pullback.  We apply this result in Section \ref{sect:slopecomputation} to show that our ability to descend $\mathcal{F}$ to $Z$ is controlled by the comparison between several invariants of the normal sheaf of $s$ and the ``gaps'' in slope between $\mathcal{F}$ and adjacent terms of the relative Harder-Narasimhan filtration.  Finally, in Section \ref{sect:genusanddimbounds} we show that for sections of a good fibration $\pi: \mathcal{Z} \to B$ these invariants of the normal sheaf are bounded by functions of $\dim(\mathcal{Z})$ and $g(B)$.

\subsection{Descending sheaves}  \label{sect:descent}
The first step is to develop a criterion for identifying when a sheaf is pulled back from the base of a morphism.

\begin{lemma} \label{lemm:rigidity}
Suppose we have morphisms $f: U \to V$, $g: U \to G$ satisfying the following properties:
\begin{enumerate}
\item $U,V,G$ are smooth varieties.
\item $f$ is dominant.  
\item Every fiber of $f$ is contracted to a point by $g$.
\end{enumerate}
Then there is some open set $V^{\circ} \subset V$ such that $g|_{f^{-1}(V^{\circ})}$ factors through $f|_{f^{-1}(V^{\circ})}$.
\end{lemma}

\begin{proof}
Consider the induced map $(f,g): U \to V \times G$ and let $\Gamma$ denote the closure of the image.  Note that $\Gamma$ is still irreducible.  For a general point $v \in V$ there is a unique point in $(f,g)(U) \cap \pi_{1}^{-1}(v)$.  Taking closures, we see that the general fiber of $\Gamma \to V$ is set-theoretically a single point.  Since we are in characteristic $0$, by generic smoothness we see that the general fiber is scheme-theoretically a single point.  We deduce that $\Gamma \to V$ is birational.  If we let $V^{\circ}$ denote an open subset where $\Gamma \to V$ is an isomorphism, then the desired statement follows.
\end{proof}

Recall that for a coherent sheaf $\mathcal{F}$ on a variety we denote by $\mathcal{F}_{tors}$ the torsion subsheaf and by $\mathcal{F}_{tf}$ the quotient of $\mathcal{F}$ by its torsion subsheaf.

\begin{lemma}[Descent Lemma] \label{lemm:descent}
Let $\mathcal{U}$ and $Z$ be smooth varieties with a dominant flat morphism $ev: \mathcal{U} \to Z$ such that the general fiber of $ev$ is connected.  Let $\mathcal{E}$ be a locally free sheaf on $Z$.  Suppose that $ev^{*}\mathcal{E} \to \mathcal{Q}$ is a surjection onto a locally free sheaf and let $\mathcal{S}$ denote the kernel.  If
\begin{equation*}
\mathcal{H}om(\mathcal{H}om(\mathcal{Q},\mathcal{S}), (\Omega_{\mathcal{U}/Z})_{tf}) = 0
\end{equation*}
then there is a subsheaf $\mathcal{S}_{Z} \subset \mathcal{E}$ such that $ev^{*}\mathcal{S}_{Z} = \mathcal{S}$ as subsheaves of $ev^{*}\mathcal{E}$.
\end{lemma}

\begin{proof}
Note that the surjection $ev^*\mathcal{E} \to \mathcal{Q}$ corresponds to a map $\phi: \mathcal{U} \to \mathbb{G}(\mathcal{E},k) = \mathbb{G}$, where $k$ is the rank of $\mathcal{Q}$ and $\mathbb{G}(\mathcal{E},k)$ is the relative Grassmannian of rank $k$ quotients. The map $\phi$ is a map of $Z$-schemes.  We first show that after replacing $\mathcal{U}$ by an open subset the map $\phi$ factors through $\ev: \mathcal{U} \to Z$.  The map $\phi$ induces a map
\[ (d\phi)^* : \phi^{*}\Omega_{\mathbb{G} / Z} \to \Omega_{\mathcal{U}/Z}  .\]
Since $\phi^{*}\Omega_{\mathbb{G} / Z} = \mathcal{H}om(\mathcal{Q}, \mathcal{S})$, our assumption implies that $(d\phi)^*$ is the 0 map on the complement of the support of the torsion subsheaf of $\Omega_{\mathcal U/Z}$.  We denote this open subset by $\widetilde{\mathcal U}$.

Fix a general point $z \in Z$ and consider the map of fibers $\widetilde{\mathcal{U}}_{z} \to \mathbb{G}_{z}$.  Using compatibility of cotangent bundles with base change, we see that $\phi|_{\widetilde{\mathcal{U}}_{z}}^{*}\Omega_{\mathbb{G}_{z}} \to \Omega_{\widetilde{\mathcal{U}}_{z}}$ must be $0$.  Dually, the map $T_{\widetilde{\mathcal{U}}_{z}} \to \phi|_{\widetilde{\mathcal{U}}_{z}}^{*}T_{\mathbb{G}_{z}}$ is zero, and thus if we precompose by the locally closed embedding $(\widetilde{\mathcal{U}}_{z,red})^{smooth} \to \widetilde{\mathcal{U}}_{z,red} \to \widetilde{\mathcal{U}}_{z}$ the induced map on tangent sheaves still vanishes.  By generic smoothness, it follows that $\phi$ must contract $(\widetilde{\mathcal{U}}_{z,red})^{smooth}$ to a point; taking closures, it also contracts each fiber $\widetilde{\mathcal{U}}_{z}$ to a point.
 There is an open subset $Z^{\circ} \subset Z$ with preimage $\widetilde{\mathcal{U}}^{\circ} = ev^{-1}(Z^{\circ}) \cap \widetilde{\mathcal U}$ such that $ev|_{\widetilde{\mathcal{U}}^{\circ}}$ has connected fibers. By Lemma \ref{lemm:rigidity}, after possibly shrinking $Z^{\circ}$ and $\widetilde{\mathcal{U}}^{\circ}$ we can ensure that $\phi|_{\widetilde{\mathcal{U}}^{\circ}}$ factors through $ev|_{\widetilde{\mathcal{U}}^{\circ}}$. Hence, $\mathcal{S}|_{\mathcal{U}^{\circ}}$ must be the pullback of a locally free sheaf $\mathcal{R} \subset \mathcal{E}|_{Z^{\circ}}$ on $Z^{\circ}$.

Let $\mathcal{S}_{Z}$ denote the unique torsion-free saturated subsheaf of $\mathcal{E}$ whose restriction to $Z^{\circ}$ is $\mathcal{R}$.  Since $ev$ is flat, $ev^{*}\mathcal{S}_{Z}$ is a torsion-free saturated subsheaf of $ev^*\mathcal{E}$ whose restriction to $\mathcal{U}^{\circ}$ agrees with $\mathcal{S}|_{\mathcal{U}^{\circ}}$.  This implies that $ev^{*}\mathcal{S}_{Z} = \mathcal{S}$.
\end{proof}

\subsection{Slope computations} \label{sect:slopecomputation}

The key to using the descent lemma is to understand homomorphisms into $\Omega_{\mathcal U/Z}$.  When $\mathcal{U}$ is a family of curves mapping to $Z$, we will control the existence of such homomorphisms using slope calculations.  The next step is to show that in this situation the slope of $\Omega_{\mathcal U/Z}$ is controlled by Lazarsfeld bundles.

\begin{definition}
Let $Y$ be a variety and $\mathcal{E}$ be a globally generated vector bundle on $Y$.  The \emph{Lazarsfeld bundle} $M_{\mathcal{E}}$ is the kernel of the evaluation map $\OO_Y\otimes H^0(Y,\mathcal{E})\to \mathcal{E}$.
\end{definition}

Given a morphism $s: C \to Z$, we will denote by $N_{s}$ the normal sheaf of $s$, i.e.~the cokernel of $T_{C} \to s^{*}T_{Z}$.  We will also let $\overline{\mathcal M}_{g,0}(Z)$ denote the Kontsevich moduli stack of maps from genus $g$ curves to $Z$.  

\begin{lemma}
\label{lem-LazMukHoms}
Let $Z$ be a smooth projective variety, let $W$ be a variety equipped with a generically finite morphism $W \to \overline{\mathcal M}_{g,0}(Z)$ and let $p: \cU_{W} \to W$ be the universal family over $W$ equipped with the evaluation map $ev_{W}: \cU_{W} \to Z$.  Suppose that a general fiber of $p$ is smooth and irreducible and that $ev_{W}$ is dominant. 

Let $C$ denote a general fiber of $\cU_{W} \to W$ equipped with the induced morphism $s: C \to Z$. Let $t$ be the length of the torsion part of $N_{s}$, let $\mathcal{G}$ be the subsheaf of $(N_{s})_{tf}$ generated by global sections, and let $V$ be the tangent space to $W$ at $s$. Let $q$ be the dimension of the cokernel of the composition
\begin{equation*}
V \to T_{ \overline{\mathcal M}_{g,0}(Z),s} = H^{0}(C,N_{s})  \to H^{0}(C,(N_{s})_{tf}).
\end{equation*}
Then
\begin{equation*}
\mu^{max}((\Omega_{\cU_{W}/Z}|_{C})_{tf}) \leq (q+1)\mu^{max}(M_{\mathcal{G}}^{\vee}) + t.
\end{equation*}
\end{lemma}

\begin{proof}
Since the conclusion only involves a general curve, we may shrink $W$ and thus assume that $W$ is smooth.  After perhaps shrinking $W$ further we may assume that $p$ is smooth, and thus $\mathcal{U}_{W}$ is a smooth variety.  We denote by $h$ the map $(p,ev_{W}): \mathcal{U}_{W} \to W \times Z$.

Fix $s: C \to Z$ as in the statement of the lemma.  Let $K_{1}$ denote the kernel of $s^{*}\Omega_{Z} \to \Omega_{C}$ and let $K_{2}$ denote the kernel of $h^{*}\Omega_{W \times Z} \to \Omega_{\cU_{W}}$.  Note that $K_{1}$ is isomorphic to the dual of $(N_{s})_{tf}$.  We claim that the following diagram has exact rows and columns:
\begin{equation*}
\xymatrix{
& & 0 \ar[d] & 0 \ar[d] &    \\
& & \mathcal{O}_{C}^{\dim(W)} \ar[r]^{=} \ar[d] & \mathcal{O}_{C}^{\dim(W)} \ar[d] &  \\
0 \ar[r] & K_{2}|_{C} \ar[r] \ar[d] & h^{*}\Omega_{W \times Z}|_{C} \ar[r] \ar[d] &  \Omega_{\cU_{W}}^{im}|_{C} \ar[r] \ar[d] & 0 \\
0 \ar[r] & K_{1} \ar[r] & s^{*}\Omega_{Z} \ar[r] \ar[d] &  \Omega_{C}^{im} \ar[r]  \ar[d] & 0 \\
& &  0 & 0 & 
}
\end{equation*}
where $\Omega_{\cU_{W}}^{im}$ is the image of $h^{*}\Omega_{W \times Z} \to \Omega_{\cU_{W}}$ and $\Omega_{C}^{im}$ is the image of $s^{*}\Omega_{Z} \to \Omega_{C}$.  The bottom row is exact by definition and the middle row is the restriction of an exact sequence to a general fiber of $p$ and thus remains exact.  
The middle column is exact since $C$ is vertical for the map $p: \mathcal{U}_{W} \to W$.  By comparing the rightmost column against the middle it is clear that $\Omega_{\cU_{W}}^{im}|_{C}$ maps surjectively onto $\Omega_{C}^{im}$ and that the kernel contains $p^{*}\Omega_{W}|_{C} \cong \mathcal{O}_{C}^{\dim(W)}$.  On the other hand the kernel must be contained in the kernel of $\Omega_{\cU_{W}}|_{C} \to \Omega_{C}$ which is also isomorphic to $p^{*}\Omega_{W}|_{C}$.  So the rightmost column is exact.  Finally, by the nine lemma we deduce that $K_{2}|_{C} \cong K_{1}$.

Let $K_{3}$ denote the kernel of the map $p^{*}\Omega_{W} \to \Omega_{\cU_{W/Z}}$ and let $\Omega_{\cU_{W/Z}}^{im}$ denote the image of this map.  We can make a further comparison via the following diagram.
\begin{equation*}
\xymatrix{
& & 0 \ar[d] & 0 \ar[d] &    \\
& & s^{*}\Omega_{Z} \ar[r]^{=} \ar[d] & s^{*}\Omega_{Z} \ar[d] &  \\
0 \ar[r] & K_{2}|_{C} \ar[r] \ar[d] & h^{*}\Omega_{W \times Z}|_{C} \ar[r] \ar[d] &  \Omega_{\cU_{W}}^{im}|_{C} \ar[r] \ar[d] & 0 \\
0 \ar[r] & K_{3}|_{C} \ar[r] & p^{*}\Omega_{W}|_{C} \ar[r] \ar[d] &   \Omega_{\cU_{W/Z}}^{im}|_{C}  \ar[r]  \ar[d] & 0 \\
& &  0 & 0 & 
}
\end{equation*}
We claim that the rows and columns are exact.  Since $C$ is general, an exact sequence of torsion-free sheaves on $\mathcal{U}_{W}$ will remain exact upon restriction.  Thus it suffices to show that the rightmost column is exact.  Since the map $ev_{W}: \mathcal{U}_{W} \to Z$ is dominant, the map $ev_{W}^{*}\Omega_{Z} \to \Omega_{\mathcal{U}_{W}}$ is generically injective.  Since $\Omega_{Z}$ is locally free the map must be injective.  Restricting to the curve $C$, we see that $s^{*}\Omega_{Z} \to \Omega_{\mathcal{U}_{W}}|_{C}$ is injective and it is clear that its image is contained in $\Omega_{\cU_{W}}^{im}|_{C}$.  This shows the rightmost column is left-exact.  Furthermore, we see that the composed map $h^{*}\Omega_{W \times Z} \to p^{*}\Omega_{W} \to \Omega_{\cU_{W/Z}}^{im}$ is surjective, showing that the map $\Omega_{\cU_{W}}^{im} \to  \Omega_{\cU_{W/Z}}^{im}$ is also surjective.  
Finally, the exact sequence
\[
0 \to ev_W^*\Omega_Z \to \Omega_{\mathcal U_W} \to  \Omega_{\mathcal U_W/Z} \to 0
\]
implies that
\[
0 \to s^*\Omega_Z \to \Omega_{\mathcal U_W}|_C \to  \Omega_{\mathcal U_W/Z}|_C \to 0
\]
is exact. Thus the rightmost column must be exact at the middle. 
 By the nine-lemma, we conclude that $K_{3}|_{C} \cong K_{2}|_{C} \cong K_{1}$.
 
Recall that $V$ denotes the tangent space to $W$ at $s$.  Let $\zeta$ denote the composition
\begin{equation*}
V \to H^{0}(C,N_{s})  \to H^{0}(C,(N_{s})_{tf}).
\end{equation*}
Then the map $K_{1} \cong K_{3}|_{C} \to p^{*}\Omega_{W}|_{C}$ is the dual of the composition
\begin{equation*}
V \otimes \mathcal{O}_{C} \xrightarrow{\zeta} H^{0}(C,(N_{s})_{tf}) \otimes \mathcal{O}_{C} \to (N_{s})_{tf}.
\end{equation*}
Let $\mathcal{G}$ denote the subsheaf of $(N_{s})_{tf}$ that is generated by its global sections, so we have an exact sequence of locally free sheaves
\begin{equation*}
0 \to M_{\mathcal{G}} \to H^{0}(C,(N_{s})_{tf}) \otimes \mathcal{O}_{C} \to \mathcal{G} \to 0.
\end{equation*}
Since $C$ deforms in a dominant family, $N_{s}$ is generically globally generated and thus the inclusion $\mathcal{G} \to (N_{s})_{tf}$ is generically surjective.  Taking duals, we see that $\mathcal{G}^{\vee}$ is the saturation of $K_{1}$ inside of $H^{0}(C,(N_{s})_{tf})^{\vee} \otimes \mathcal{O}_{C}$.  In other words, if we let $\mathcal{S}$ denote the cokernel of the map $K_{1} \to H^{0}(C,(N_{s})_{tf})^{\vee} \otimes \mathcal{O}_{C}$ then the torsion free part of $\mathcal{S}$ is isomorphic to $M_{\mathcal{G}}^{\vee}$.

 Consider the following diagram of short exact sequences
\begin{equation*}
\xymatrix{
0 \ar[r] & K_{1} \ar[r] \ar[d]^{=} & H^{0}(C,(N_{s})_{tf})^{\vee} \otimes \mathcal{O}_{C}  \ar[r] \ar[d]^{\zeta^{\vee}} & \mathcal{S}  \ar[r]  \ar[d]^{h} & 0 \\
0 \ar[r] & K_{1} \ar[r] & p^{*}\Omega_{W}|_{C} \ar[r]  & \Omega_{\cU_{W}/Z}^{im}|_{C} \ar[r]  & 0
}
\end{equation*}
By the snake lemma we obtain an exact sequence
\begin{align*}
0 \to \mathcal{O}_{C}^{\oplus q} \to \mathcal{S}  \to \Omega_{\cU_{W}/Z}^{im}|_{C}  \to \mathcal{O}_{C}^{\oplus e} \to 0
\end{align*}
where $q,e$ are respectively the dimensions of the cokernel and kernel of $\zeta$.

Recall that $\mathcal{S}_{tf} \cong M_{\mathcal{G}}^{\vee}$ and note that the torsion part of $\mathcal{S}$ injects into the torsion part of $\Omega_{\cU_{W}/Z}^{im}|_{C}$.  We will denote by $\mathcal{R}$ the saturation of $\mathcal{O}_{C}^{\oplus q}$ in $M_{\mathcal{G}}^{\vee}$ so that we obtain an exact sequence
\begin{align*}
0 \to \mathcal{R} \to M_{\mathcal{G}}^{\vee}  \to (\Omega_{\cU_{W}/Z}^{im}|_{C})_{tf}  \to \mathcal{O}_{C}^{\oplus e} \to 0
\end{align*}
Let $\mathcal{F}$ denote the maximal destabilizing subsheaf of $(\Omega_{\cU_{W}/Z}^{im}|_{C})_{tf}$, let $\mathcal{F}' = \mathcal{F} \cap im(M_{\mathcal{G}}^{\vee})$, and let $\mathcal{Q}$ denote the preimage of $\mathcal{F}'$ in $M_{\mathcal{G}}^{\vee}$.  Our next goal is to prove an upper bound on $\mu(\mathcal{F})$.  First suppose that $\mu(\mathcal{F}) > 0$.  Then the image of $\mathcal{F}$ in $\mathcal{O}_{C}^{\oplus e}$ is $0$ and so $\mathcal{F} = \mathcal{F}'$.  This means we have an exact sequence
\begin{equation*}
0 \to \mathcal{R} \to \mathcal{Q} \to \mathcal{F} \to 0
\end{equation*}
and thus 
\begin{align*}
\mu(\mathcal{F}) & = \frac{\deg(\mathcal{Q}) - \deg(\mathcal{R})}{\rk(\mathcal{Q}) - \rk(\mathcal{R})} \\
& \leq \frac{\deg(\mathcal{Q})}{\rk(\mathcal{Q}) - q} \\
& \leq (q+1) \mu(\mathcal{Q})
\end{align*}
where the final line follows from the elementary inequality $\frac{1}{b-a} \leq \frac{a+1}{b}$ when $b-1 \geq a \geq 0$.  Thus in this case we see that $\mu^{max}(\Omega_{\cU_{W}/Z}^{im}|_{C}) \leq (q+1)\mu^{max}(M_{\mathcal{G}}^{\vee})$.  If $\mu(\mathcal{F}) \leq 0$, then the same inequality still holds: the right-hand side is a non-negative number since $M_{\mathcal{G}}^{\vee}$ is globally generated.

Consider the following diagram
\begin{equation*}
\xymatrix{
& 0  & 0  & 0  & \\
0 \ar[r] & \Omega_C^{im} \ar[r] \ar[u] & \Omega_C \ar[r] \ar[u] & \mathcal{T} \ar[r] \ar[u] & 0 \\
0 \ar[r] & s^{*}\Omega_Z \ar[r] \ar[u] & \Omega_{\mathcal{U}_W}|_{C} \ar[r] \ar[u] & \Omega_{\mathcal{U}_W/Z}|_C  \ar[r]  \ar[u] & 0 \\
0 \ar[r] & K_{3}|_{C} \ar[r] \ar[u] & p^{*}\Omega_{W}|_{C} \ar[r] \ar[u] & \Omega_{\cU_{W}/Z}^{im}|_{C} \ar[r] \ar[u] & 0 \\
& 0 \ar[u] & 0 \ar[u] &0 \ar[u]  &
}
\end{equation*}
where $\Omega_C^{im}$ is the image and $\mathcal{T}$ is the cokernel of $s^{*}\Omega_Z \to \Omega_C$.  Every row and column is exact except possibly the rightmost column, thus by nine-lemma we see the rightmost column is also exact.  We have
\begin{equation*}
\mathrm{len}(\mathcal{T}) = \mathrm{len}(\cok(\Omega_{C}^{im} \to \Omega_{C})) = \mathrm{len}(\cok(T_{C} \to T_{C}^{sat})) = t
\end{equation*}
where $T_{C}^{sat}$ is the saturation of $T_{C}$ in $s^{*}T_{Z}$.  It follows that
\begin{align*}
\mu^{max}((\Omega_{\cU_{W}/Z}|_{C})_{tf}) & \leq \mu^{max}((\Omega_{\cU_{W}/Z}^{im}|_{C})_{tf}) + t \\
& \leq (q+1) \mu^{max}(M_{\mathcal{G}}^{\vee}) + t.
\end{align*}
\end{proof}

By combining this analysis with the descent theorem, we obtain:

\begin{theorem} \label{theo:gmtheorem}
Let $Z$ be a smooth projective variety.
Let $W$ be a variety equipped with a generically finite morphism $W \to \overline{\mathcal M}_{g,0}(Z)$ and let $p: \mathcal{U}_{W} \to W$ denote the universal family over $W$ with evaluation map $\ev_{W}: \mathcal{U}_{W} \to Z$.  Assume that a general map parametrized by $W$ has smooth  irreducible domain, that $ev_{W}$ is dominant, that the general fiber of the composition of the normalization map for $\mathcal{U}_{W}$ with $ev_{W}$ is connected, and that a general fiber of $p$ is contained in the locus where $ev_{W}$ is flat. 

Suppose that $\mathcal{E}$ is a torsion-free sheaf on $Z$ that is semistable with respect to a general curve $s: C \to Z$ parametrized by $W$.  Write
\begin{equation*}
0 = \mathcal{F}_{0} \subset \mathcal{F}_{1} \subset \mathcal{F}_{2} \subset \ldots \subset \mathcal{F}_{k} = s^{*}\mathcal{E}
\end{equation*}
for the Harder-Narasimhan filtration of $s^{*}\mathcal{E}$.  Let $t$ be the length of the torsion part of $N_{s}$, let $\mathcal{G}$ be the subsheaf of $(N_{s})_{tf}$ generated by global sections, and let $V$ be the tangent space to $W$ at $s$. Let $q$ be the dimension of the cokernel of the composition
\begin{equation*}
V \to T_{ \overline{\mathcal M}_{g,0}(Z),s} = H^{0}(C,N_{s})  \to H^{0}(C,(N_{s})_{tf}).
\end{equation*}
Then for every index $1 \leq i \leq k-1$ we have
\begin{equation*}
\mu(\FF_{i}/\FF_{i-1})-\mu(\FF_{i+1}/\FF_i) \leq (q+1)\mu^{max}(M_{\mathcal{G}}^{\vee}) + t.
\end{equation*}
\end{theorem}

Note that by the flatness assumption we may ensure that the image of a general map $s$ parametrized by $W$ will avoid any codimension $2$ locus on $Z$.  In particular this implies that $s^{*}\mathcal{E}$ will be a locally free sheaf.

\begin{proof}
Since we are assuming that the general fiber of the composition of the normalization map for $\mathcal{U}_{W}$ with $ev_{W}$ is connected, if we replace $W$ by a smaller open subset then the general fiber of the evaluation map will still have connected fibers.  Thus to prove the statement, we may replace $W$ by a smaller open subset.  (For ease of notation we continue to call the smaller subset $W$ and its $p$-preimage by $\mathcal{U}_{W}$.)  Thus we may suppose that $W$  and $\mathcal{U}_{W}$ are smooth and that $ev_{W}|_{\mathcal{U}_{W}}$ is flat.   After possibly shrinking $W$ further, by \cite[Theorem 2.3.2]{HL97} we may suppose that for every index $1 \leq i \leq k-1$ there is a torsion-free sheaf $\mathcal{S}_{i} \subset ev_{W}^*\mathcal{E}$ obtained from the relative Harder-Narasimhan filtration of $ev_{W}^{*}\mathcal{E}$ over $W$ such that $\mathcal{S}_{i}|_{C} \cong \mathcal{F}_{i}$ for every fiber $C$ over $W$.  Since torsion-free sheaves are locally free on the complement of a codimension $2$ subset, after perhaps shrinking $W$ again we may suppose that each $\mathcal{S}_{i}$ is locally free.

Suppose for a contradiction that there is an index $i$ such that
\begin{equation*}
\mu(\FF_{i}/\FF_{i-1})-\mu(\FF_{i+1}/\FF_i) > (q+1)\mu^{max}(M_{\mathcal{G}}^{\vee}) + t.
\end{equation*}
If there were a non-zero homomorphism $\mathcal{H}om(ev_{W}^{*}\mathcal{E}/\mathcal{S}_i,\mathcal{S}_i) \to (\Omega_{\mathcal{U}_{W}/Z})_{tf}$,
then its restriction to a general fiber $C$ of $p$ would yield a map that is non-zero on the generic point of $C$, and thus would induce a non-zero map $\mathcal{H}om(ev_{W}^{*}\mathcal{E}/\mathcal{S}_i|_{C},\mathcal{S}_i|_{C}) \to ((\Omega_{\mathcal{U}_{W}/Z})|_{C})_{tf}$. But then
\begin{align*}
\mu^{min}(\mathcal{H}om(ev_{W}^{*}\mathcal{E}/\mathcal{S}_i|_{C},\mathcal S_{i}|_{C})) & = \mu^{min}(\mathcal{S}_i|_{C}) - \mu^{max}(ev_{W}^{*}\mathcal{E}/\mathcal{S}_i|_{C}) \\
& > (q+1)\mu^{max}(M_{\mathcal{G}}^{\vee}) + t \\
& \geq \mu^{max}((\Omega_{\mathcal{U}_{W}/Z}|_{C})_{tf})
\end{align*}
where the last inequality follows from Lemma \ref{lem-LazMukHoms}.  We conclude that there is no non-zero homomorphism $\mathcal{H}om(ev_{W}^{*}\mathcal{E}/\mathcal{S}_i,\mathcal{S}_i) \to (\Omega_{\mathcal{U}_{W}/Z})_{tf}$.  By Lemma \ref{lemm:descent}, we see that there is a sheaf $\mathcal{S}_{Z}$ on $Z$ such that $ev_{W}^{*}\mathcal{S}_{Z} = \mathcal{S}_i$.  But such a sheaf would destabilize $\mathcal{E}$: we have
\begin{equation*}
\mu_{[s_{*}C]}(\mathcal{S}_{Z}) = \mu(\mathcal{S}_{i}|_{C}) > \mu(ev_{W}^{*}\mathcal{E}|_{C}) = \mu_{[s_{*}C]}(\mathcal{E})
\end{equation*}
where the equalities follow from the flatness of the evaluation map and the inequality follows from \eqref{eq:mediant}.  This gives a contradiction and we conclude the desired inequalities for every $i$.
\end{proof}

It will be helpful to have a version of the previous theorem that holds for non-semistable sheaves as well.  The next theorem controls the difference between the Harder-Narasimhan filtration of $\mathcal{E}|_{C}$ and the restriction to $C$ of the $[C]$-Harder-Narasimhan filtration of $\mathcal{E}$.  It is a formal consequence of Theorem \ref{theo:gmtheorem}.

\begin{corollary} \label{coro:hnfversion}
Let $Z$ be a smooth projective variety and let $\mathcal{E}$ be a torsion free sheaf on $Z$ of rank $r$.  Let $W$ be a variety equipped with a generically finite morphism $W \to \overline{\mathcal M}_{g,0}(Z)$ and let $p: \mathcal{U}_{W} \to W$ denote the universal family over $W$ with evaluation map $\ev_{W}: \mathcal{U}_{W} \to Z$.  Assume that a general map parametrized by $W$ has smooth irreducible domain, that $ev_{W}$ is dominant, that the general fiber of the composition of the normalization map for $\mathcal{U}_{W}$ with $ev_{W}$ is connected, and that a general fiber of $p$ is contained in the locus where $ev_{W}$ is flat.

Let $C$ denote a general fiber of $\cU_{W} \to W$ equipped with the induced morphism $s: C \to Z$.  Let $t$ be the length of the torsion part of $N_{s}$, let $\mathcal{G}$ be the subsheaf of $(N_{s})_{tf}$ generated by global sections, and let $V$ be the tangent space to $W$ at $s$. Let $q$ be the dimension of the cokernel of the composition
\begin{equation*}
V \to T_{ \overline{\mathcal M}_{g,0}(Z),s} = H^{0}(C,N_{s})  \to H^{0}(C,(N_{s})_{tf}).
\end{equation*}
Let $\Vert - \Vert$ denote the sup norm on $\mathbb{Q}^{\oplus r}$. Then we have
\begin{equation*}
\Vert \SP_{Z,[C]}(\mathcal{E}) - \SP_{C}(s^{*}\mathcal{E}) \Vert \leq \frac{1}{2} \left( (q+1)\mu^{max}(M_{\mathcal{G}}^{\vee}) + t \right)  (\rk(\mathcal{E})-1).
\end{equation*}
\end{corollary}

\begin{proof}
Set $\delta = (q+1)\mu^{max}(M_{\mathcal{G}}^{\vee}) + t$. Given two $r$-tuples of real numbers $a_{\bullet}, b_{\bullet}$, we write $a_{\bullet} \geq b_{\bullet}$ if for every index $i=1,2,\ldots,r$ we have $a_{i} \geq b_{i}$.

Write $0 = \mathcal{F}_{0} \subset \mathcal{F}_{1} \subset \ldots \subset \mathcal{F}_{k} = \mathcal{E}$ for the $[C]$-Harder-Narasimhan filtration of $\mathcal{E}$.  We prove this statement by induction on the length $k$. We start with the base case $k=1$. By Theorem \ref{theo:gmtheorem}  the slopes of successive quotients of the Harder-Narasimhan filtration of $s^{*}\mathcal{E}$ differ by at most $\delta$.  We conclude the desired statement by Lemma \ref{lemm:weightedAverage} (where the pairs $(a_{i},b_{i})$ record the degrees and ranks of the various successive quotients in the Harder-Narasimhan filtration of $s^{*}\mathcal{E}$).

For the induction step, write $0 = \mathcal{G}_{0} \subset \mathcal{G}_{1} \subset \ldots \subset \mathcal{G}_{\ell} = s^{*}\mathcal{E}$ for the Harder-Narasimhan filtration of $s^{*}\mathcal{E}$.  For convenience we define three tuples:
\begin{itemize}
\item $c_{\bullet} = (c_{1},c_{2},\ldots,c_{r})$ to be $\SP_{C}(s^{*}\mathcal{E})$.
\item $g_{\bullet} = (g_{1},g_{2},\ldots,g_{r})$ defined as follows: we start with the tuple $\SP_{C}(s^{*}\mathcal{F}_{k-1})$ (of length $\rk(\mathcal{F}_{k-1})$) and then replace every entry that is below $\mu_{[C]}(\mathcal{E}/\mathcal{F}_{k-1}) + \frac{\delta}{2}(\rk(\mathcal{E})-1)$ by this number.  We then append entries equal to $\mu_{[C]}(\mathcal{E}/\mathcal{F}_{k-1}) + \frac{\delta}{2}(\rk(\mathcal{E})-1)$ to the end so that $g_{\bullet}$ has length equal to $\rk(\mathcal{E})$.
\item $h_{\bullet} = (h_{1},h_{2},\ldots,h_{r})$ defined as follows: we start with the tuple $\SP_{C}(s^{*}(\mathcal{E}/\mathcal{F}_{1}))$ (of length $\rk(\mathcal{E}/\mathcal{F}_{1})$) and then replace every entry that is above $\mu_{[C]}(\mathcal{F}_{1}) - \frac{\delta}{2}(\rk(\mathcal{E})-1)$ by this number.  We then insert entries equal to $\mu_{[C]}(\mathcal{F}_{1}) - \frac{\delta}{2}(\rk(\mathcal{E})-1)$ at the beginning so that $h_{\bullet}$ has length equal to $\rk(\mathcal{E})$.
\end{itemize}
Since $C$ is a general member of a flat family every term $\mathcal{F}_{i}$ is locally free along $C$.  Thus we have an exact sequence
\begin{equation*}
0 \to s^{*}\mathcal{F}_{k-1} \to s^{*}\mathcal{E} \to s^{*}(\mathcal{E}/\mathcal{F}_{k-1}) \to 0.
\end{equation*}
Fix an index $1 \leq j \leq \ell$ and suppose that $\mu^{min}(\mathcal{G}_{j}) > \mu_{[C]}(\mathcal{E}/\mathcal{F}_{k-1}) + \frac{\delta}{2}(\rk(\mathcal{E})-1)$.  By the base case of our result we see that $\mu^{min}(\mathcal{G}_{j}) > \mu^{max}(s^{*}(\mathcal{E}/\mathcal{F}_{k-1}))$ so that $\mathcal{G}_{j} \subset s^{*}\mathcal{F}_{k-1}$.  Thus we have $g_{\bullet} \geq c_{\bullet}$.

Similarly, consider the exact sequence
\begin{equation*}
0 \to s^{*}\mathcal{F}_{1} \to s^{*}\mathcal{E} \to s^{*}(\mathcal{E}/\mathcal{F}_{1}) \to 0.
\end{equation*}
Fix an index $1 \leq j \leq \ell$ and suppose that $\mu^{max}(s^{*}\mathcal{E}/\mathcal{G}_{j}) < \mu_{[C]}(\mathcal{F}_{1}) - \frac{\delta}{2}(\rk(\mathcal{E})-1)$.  By the base case of our result we see that $\mu^{max}(s^{*}\mathcal{E}/\mathcal{G}_{j}) < \mu^{min}(s^{*}\mathcal{F}_{1})$ so that there can be no non-zero homomorphism from $s^{*}\mathcal{F}_{1}$ to $s^{*}\mathcal{E}/\mathcal{G}_{j}$.  We conclude that $s^{*}\mathcal{F}_{1} \subset \mathcal{G}_{j}$ and thus $s^{*}\mathcal{E}/\mathcal{G}_{j}$ is a quotient of $s^{*}(\mathcal{E}/\mathcal{F}_{1})$.  Thus we have $c_{\bullet} \geq h_{\bullet}$.

Altogether, suppose we define
\begin{itemize}
\item $q_{\bullet}^{+}$ to be $\SP_{Z,C}(\mathcal{E}) + (\frac{\delta}{2}(\rk(\mathcal{E})-1), \frac{\delta}{2}(\rk(\mathcal{E})-1), \ldots, \frac{\delta}{2}(\rk(\mathcal{E})-1))$.
\item $q_{\bullet}^{-}$ to be $\SP_{Z,C}(\mathcal{E}) - (\frac{\delta}{2}(\rk(\mathcal{E})-1), \frac{\delta}{2}(\rk(\mathcal{E})-1), \ldots, \frac{\delta}{2}(\rk(\mathcal{E})-1))$.
\end{itemize}
The induction assumption for $\mathcal{F}_{k-1}$ shows that $q^{+}_{\bullet} \geq g_{\bullet}$ and the induction assumption for $\mathcal{E}/\mathcal{F}_{1}$ shows that $h_{\bullet} \geq q^{-}_{\bullet}$.  By the argument above we conclude that
\begin{equation*}
q_{\bullet}^{+} \geq g_{\bullet} \geq c_{\bullet} \geq h_{\bullet} \geq q_{\bullet}^{-}
\end{equation*}
yielding the desired result.
\end{proof}

\begin{lemma}
\label{lemm:weightedAverage}
Let $\{ (a_{i},b_{i}) \}_{i=1}^{r}$ be pairs of integers with $b_{i} > 0$ such that the fractions $\frac{a_{i}}{b_{i}}$ are nonincreasing as $i$ gets larger.  Let $m = \frac{a_{1}+ \ldots + a_{r}}{b_{1} + \ldots + b_{r}}$ denote the mediant of fractions $\frac{a_{i}}{b_{i}}$ and let $\delta$ denote the maximum of the successive differences $\frac{a_{i}}{b_{i}} - \frac{a_{i+1}}{b_{i+1}}$.  Then for every $i$ we have
\begin{equation*}
\left| \frac{a_{i}}{b_{i}} - m \right| \leq \frac{\delta}{2} \left( \sum_{i=1}^{r} b_{i} - 1 \right)
\end{equation*}
\end{lemma}

\begin{proof}
It suffices to prove the inequality for the extremal fractions $\frac{a_{1}}{b_{1}}$, $\frac{a_{r}}{b_{r}}$.  Replacing each $a_{i}$ by $-a_{i}$, we see that the latter case is implied by the former.  So it suffices to prove the statement for $\frac{a_{1}}{b_{1}}$.  
Since $\frac{a_{1}}{b_{1}} - \frac{a_{i}}{b_{i}} \leq (i-1)\delta$, we have
\begin{align*}
\sum_{i=1}^{r} (a_{1}b_{i} - a_{i}b_{1}) & \leq \sum_{i=1}^{r} b_{1} b_{i} \delta (i-1) \\
& \leq b_{1} \delta \left( \sum_{i=1}^{r} b_{i} \left( \sum_{j=1}^{i-1}  b_{j} \right) \right) \\
& = \frac{b_{1} \delta}{2} \left( \sum_{i \neq j}b_{i} b_{j} \right) \\
& = \frac{b_{1} \delta}{2} \left( \left( \sum_{i=1}^{r} b_{i} \right)^{2} - \left( \sum_{i=1}^{r} b_{i}^{2} \right) \right) \\
& \leq \frac{b_{1} \delta}{2} \left( \left( \sum_{i=1}^{r} b_{i} \right)^{2} - \left( \sum_{i=1}^{r} b_{i} \right) \right)
\end{align*}
Dividing by $b_{1} \left( \sum_{i=1}^{r} b_{i} \right)$ gives us the result.
\end{proof}

\subsection{Genus and dimension bounds}  \label{sect:genusanddimbounds}
Note that the Hilbert scheme of sections admits an embedding into the stack $\overline{\mathcal{M}}_{g,0}(\mathcal{Z})$. To apply Corollary \ref{coro:hnfversion} to moduli spaces of sections one needs to be able to bound the quantities $q$, $t$, and $\mu^{max}(M_{\mathcal{G}}^{\vee})$ appearing in the statement.  We will show how to bound these quantities for sections of a good fibration $\pi: \mathcal{Z} \to B$ using the genus of $B$ and the dimension of $\mathcal{Z}$.  Note that since sections are always smooth the quantity $t=0$.

We first discuss the slope of Lazarsfeld bundles.  These can be bounded using only the genus of $B$ using the following result of Butler.

\begin{theorem}[\cite{Butler94}] \label{theo:butler}
Let $\mathcal{E}$ be a globally generated locally free sheaf on a curve $C$ of genus $g$ and let $M_{\mathcal{E}}$ be its Lazarsfeld bundle.
\begin{enumerate}
\item If $\mu^{min}(\mathcal{E}) \geq 2g$ then $\mu^{min}(M_{\mathcal{E}}) \geq -2$.
\item If $\mu^{min}(\mathcal{E}) < 2g$ then $\mu^{min}(M_{\mathcal{E}}) \geq -2g \rk(\mathcal{E}) - 2$.
\end{enumerate}
\end{theorem}

\begin{proof}
Statement (1) follows from \cite[1.3 Corollary]{Butler94} except when $C$ is a rational curve.  When $C$ is rational $M_{\mathcal{E}}$ is a direct sum of $\mathcal{O}(-1)$'s (see for example \cite[Lemma 3.10]{PatelRiedlTseng}) and thus the statement still holds.

For statement (2), first suppose that $\mu^{max}(\mathcal{E}) \geq 2g$. Let $\mathcal{F}$ denote the maximal destabilizing subsheaf of $\mathcal{E}$.  Our degree assumption implies that $\mathcal{F}$ is globally generated and that $h^{1}(C,\mathcal{F}) = 0$.  By the nine lemma we obtain an exact sequence
\begin{equation*}
0 \to M_{\mathcal{F}} \to M_{\mathcal{E}} \to M_{\mathcal{E}/\mathcal{F}} \to 0.
\end{equation*}
This implies that $\mu^{min}(M_{\mathcal{E}}) \geq \min \{ \mu^{min}(M_{\mathcal{F}}), \mu^{min}(M_{\mathcal{E}/\mathcal{F}}) \}$.  By (1) the first quantity is at least $-2$.  Arguing by induction on the rank we reduce to the case when $\mu^{max}(\mathcal{E}) < 2g$.

When $\mu^{max}(\mathcal{E}) < 2g$, \cite[1.5 Proposition]{Butler94} proves a statement stronger than (2) except in two cases.  The first is when $g = 1$. 
In this case, since $M_{\mathcal{E}}$ is a subsheaf of $H^{0}(C,\mathcal{E}) \otimes \mathcal{O}_{C}$ every subsheaf of $M_{\mathcal{E}}$ has non-positive slope. Since $\deg(M_{\mathcal{E}}) = -\deg(\mathcal{E})$, we conclude that every quotient of $M_{\mathcal{E}}$ has degree at least $-\deg(\mathcal{E})$, which implies that it has slope at least $-\deg(\mathcal{E})$. 
Since we are in the case where $\deg(\mathcal{E})/\rk(\mathcal{E}) < 2$ we conclude $\mu^{min}(M_{\mathcal{E}}) \geq -2 \rk(\mathcal{E})$.  The second is when $g \geq 2$ and $\mathcal{E}$ has a trivial summand.  Write $\mathcal{E} = \mathcal{E}' \oplus \mathcal{O}_{C}^{\oplus k}$.  Note that $M_{\mathcal{E}} = M_{\mathcal{E}'}$.  Since we still have $\mu^{max}(\mathcal{E}') < 2g$ we can apply \cite[1.5 Proposition]{Butler94} to $\mathcal{E}'$ to obtain the desired lower bound.  
\end{proof}

Next we discuss the quantity $q$.

\begin{lemma} \label{lemm:sbound}
Let $\pi: \mathcal{Z} \to B$ be a good fibration.  Suppose that $M \subset \Sec(\mathcal{Z}/B)$ is an irreducible component parametrizing a dominant family of sections on $Z$ and let $W = M_{red}$.  For a general section $C$ parametrized by $M$ let $V \subset H^{0}(B,T_{\mathcal{Z}/B}|_{C})$ denote the tangent space to $W$ at $C$.  Then the codimension of $V$ in $H^{0}(B,T_{\mathcal{Z}/B}|_{C})$ is at most $g(B) (\dim(\mathcal{Z})-1)$.
\end{lemma}

\begin{proof}
We have $h^{0}(B,T_{\mathcal{Z}/B}|_{C}) - \dim(V)  \leq h^{0}(B,T_{\mathcal{Z}/B}|_{C}) - \dim(M)$ and Corollary \ref{coro:domfamilyexpdim} shows that this latter quantity is bounded above by $g(B) (\dim(\mathcal{Z})-1)$.
\end{proof}

Putting these results together, we obtain a version of the Grauert-M\"ulich theorem for sections.

\begin{theorem} \label{theo:hnfforsections}
Let $\pi: \mathcal{Z} \to B$ be a good fibration and let $\mathcal{E}$ be a torsion-free sheaf on $\mathcal{Z}$.  Let $M$ be an irreducible component of $\Sec(\mathcal{Z}/B)$ parameterizing a dominant family of sections of $\pi$ and let $p: \mathcal{U}^\nu \to M_{red}$ be the normalization of the universal family over $M_{red}$ with evaluation map $ev: \mathcal{U}^\nu \to \mathcal{Z}$.  Assume that $ev$ has connected fibers.

Let $C$ be a general section parametrized by $M$.  Then:
\begin{enumerate}
\item Suppose there is an open subset $M_{red}^{\circ} \subset M_{red}$ such that if we define $\cU^{\nu, \circ} = p^{-1}M_{red}^{\circ}$ then $ev|_{\cU^{\nu, \circ}}$ is flat.  Then we have
\begin{equation*}
\Vert \SP_{\mathcal{Z},[C]}(\mathcal{E}) - \SP_{C}(\mathcal{E}|_{C}) \Vert \leq (g(B)\dim(\mathcal{Z}) - g(B) + 1)^{2} \rk(\mathcal{E})
\end{equation*}
where $\Vert - \Vert$ denotes the sup norm on $\mathbb{Q}^{\oplus r}$.
\item Suppose that the general curve $C$ is HN-free.  Then we have
\begin{equation*}
\Vert \SP_{\mathcal{Z},[C]}(\mathcal{E}) - \SP_{C}(\mathcal{E}|_{C}) \Vert \leq \rk(\mathcal{E})
\end{equation*}
where $\Vert - \Vert$ denotes the sup norm on $\mathbb{Q}^{\oplus r}$.
\end{enumerate}
\end{theorem}

In fact, the stronger statements with $\rk(\mathcal{E})$ replaced by $\rk(\mathcal{E})-1$ are also true.  However, since we are not aiming for precise bounds in this paper we will use this simpler version in the future.

\begin{proof}
(1)  Theorem \ref{theo:butler} shows that $\mu^{max}(M_{\mathcal{G}}^{\vee}) \leq 2g(B) (\dim(\mathcal{Z})-1)+2$.  By Lemma \ref{lemm:sbound} we have $q \leq g(B) (\dim(\mathcal{Z})-1)$.  We then apply Corollary \ref{coro:hnfversion} with $t=0$.

(2) When $C$ is HN-free then $M$ is generically smooth and $T_{\mathcal{Z}/B}|_{C'}$ is globally generated for a general deformation $C'$ of $C$.  Furthermore, there is an open subset of $M$ over which the evaluation map is smooth and thus flat.  Theorem \ref{theo:butler} shows that $\mu^{max}(M_{\mathcal{G}}^{\vee}) \leq 2$. We then apply Corollary \ref{coro:hnfversion} with $q=t=0$.
\end{proof}

\section{Sections through general points} \label{sect:genpoints}

Suppose that $\pi: \mathcal{Y} \to B$ is a good fibration and $M$ is an irreducible component of $\Sec(\mathcal{Y}/B)$ parametrizing a dominant family of sections.  Let $C$ be a general section parametrized by $M$.  By Proposition~\ref{prop:deffixpoints} we can identify a lower bound on the slopes in the Harder-Narasimhan filtration of $T_{\mathcal{Y}/B}|_{C}$ by computing how many general points of $\mathcal{Y}$ we can impose on the sections parametrized by $M$.

Even when $C$ has very large anticanonical degree, deformations of $C$ do not need to go through many general points of $\mathcal{Y}$.  In this section we construct a dominant family of subvarieties $\mathcal{W} \subset \mathcal{Y}$ such that deformations of $C$ go through many general points in $\mathcal{W}$.   Results of this type were used earlier in \cite{Shen12}, \cite{LT21a}, and \cite{LT21b}.

Here is the idea behind the construction.  Suppose that $M$ parametrizes a family of sections $C$ which have large degree but do not go through many points of $\mathcal{Y}$.  This implies that the Harder-Narasimhan filtration of $T_{\mathcal{Y}/B}|_{C}$ has a large gap in the slopes between two consecutive terms $\mathcal{G}_{k} \subset \mathcal{G}_{k+1}$ for some index $k$.  When $M$ satisfies the conditions of Corollary \ref{coro:hnfversion}, we can deduce that there is a foliation $\mathcal{F}$ on $X$ that restricts to $\mathcal{G}_{k}$.  Appealing to the results developed in the sequence of papers \cite{BM01}, \cite{KSCT07}, \cite{CP19}, the foliation is induced by a rational map $\phi: \mathcal{Y} \dashrightarrow \mathcal{Z}$.  We can expect that there will be many deformations of $C$ in directions tangent to the fibers of $\phi$.  In particular, deformations of $C$ should go through many general points of the main component $\mathcal{W}$ of $\overline{\phi^{-1}(\phi(C))}$.

\subsection{General points and foliations} \label{sect:genpointfoliations}

We will need the following construction describing the relationship between foliations and relative tangent bundles which is adapted from \cite[Remark 19]{KSCT07}.

\begin{construction} \label{cons:foliationrestriction}
Let $\pi: \mathcal{Y} \to B$ be a good fibration.  Suppose that $\mathcal{F}$ is a foliation on $\mathcal{Y}$ that is contained in the relative tangent bundle $T_{\mathcal{Y}/B}$.  Assume that $\mathcal{F}$ is induced by a rational map $\phi: \mathcal{Y} \dashrightarrow \mathcal{Z}$ that has connected fibers.  Note that $\phi$ must be a rational map over $B$ and we may assume that $\mathcal{Z}$ is a projective $B$-variety. 

Suppose that $C$ is a section of $\pi$ that is contained in the regular locus of $\mathcal{F}$ and goes through a general point of $\mathcal{Y}$.  Since $C$ is transversal to $\mathcal{F}$, the Frobenius theorem shows that there is an irreducible analytic submanifold $W \subset \mathcal{Y}$ containing $C$ such that the fibers of $\pi|_{W}$ are smooth analytic manifolds which are open subsets of the leaves of $\phi$ with the property $N_{C/W} \cong \mathcal{F}|_{C}$. 

Let $\mathcal{W}$ denote the main component of $\overline{\phi^{-1}(\phi(C))}$ (that is, the unique irreducible component that dominates $\phi(C)$ under $\phi$) and let $\mathcal{W}'$ denote its normalization.  Using the universal property of normalization, we see that $W$ admits an embedding into $\mathcal{W}'$.  In particular, $\mathcal{W}'$ admits a section $C'$ in its smooth locus that maps to $C$ and has normal bundle $N_{C'/\mathcal{W}'} \cong \mathcal{F}|_{C}$.  Thus we can choose a resolution $\widetilde{\mathcal{W}}$ of $\mathcal{W}$ that admits a section $\widetilde{C}$ which maps to $C$ and satisfies $T_{\widetilde{\mathcal{W}}/B}|_{\widetilde{C}} \cong \mathcal{F}|_{C}$.
\end{construction}

Grauert-M\"ulich only applies when a general curve is contained in the flat locus of the evaluation map.  This can always be achieved after a birational modification as in the following construction:

\begin{construction} \label{cons:flatteningfamilyofcurves}
Let $S$ be a normal projective variety and let $M$ be a variety admitting a morphism $M \to \overline{\mathcal M}_{g,0}(S)$ that is generically finite onto its image.  Let $\mathcal{U}_{M}$ denote the universal family over $M$ and let $\mathcal{U}_{M}^{\nu}$ denote the normalization of $\mathcal{U}_{M}$.  Then $\mathcal{U}_{M}^{\nu}$ is equipped with a map $p: \mathcal{U}_{M}^{\nu} \to M$ and an evaluation map $ev_{M}: \mathcal{U}_{M}^{\nu} \to S$.

Assume that a general fiber of $p$ is a smooth projective curve.  We claim there is a birational map $\phi: S' \to S$ from a smooth variety $S'$ and an open subset $M^{\circ} \subset M$ such that the preimage $\mathcal{U}_{M}^{\nu,\circ} := p^{-1}M^{\circ}$ admits a flat morphism $ev': \mathcal{U}_{M}^{\nu,\circ} \to S'$ satisfying $ev_{M}|_{\mathcal{U}_{M}^{\nu,\circ}} = \phi \circ ev'$.

Indeed, suppose we take a flattening of $ev$, i.e.~a diagram
\begin{equation*}
\xymatrix{
\mathcal{V} \ar[r]^{\widetilde{ev}} \ar[d]_{\psi} & \widetilde{S} \ar[d]^{\psi_{S}} \\
\mathcal{U}_{M}^{\nu} \ar[r]_{ev_{M}} & S
}
\end{equation*}
where $\mathcal{V}$ and $\widetilde{S}$ are varieties, $\psi$ and $\psi_{S}$ are projective birational morphisms, and $\widetilde{ev}$ is flat.  Let $\rho: S' \to \widetilde{S}$ be a resolution of singularities.  Since $\widetilde{ev}$ is flat, $\mathcal{V}' := \mathcal{V} \times_{\widetilde{S}} S'$ is also a variety and the projection map $ev': \mathcal{V}' \to S'$ is still flat.   
The induced map $\psi': \mathcal{V}' \to \mathcal{U}_{M}^{\nu}$ is still birational.  Since $p$ defines a family of curves, there is an open subset $M^{\circ} \subset M$ such that $p^{-1}M^{\circ}$ is disjoint from every $\psi'$-exceptional center.  Then $M^{\circ}$ has the desired properties.
\end{construction}

We are now ready to state the main result of this section.  

\begin{theorem} \label{theo:betterpointsandfoliations}
Let $\pi: \mathcal{Y} \to B$ be a good fibration.  Fix a positive integer $J \geq 2g(B)+3$.  Suppose $M$ is an irreducible component of $\Sec(\mathcal{Y}/B)$ and let $ev: \mathcal{U}^{\nu} \to \mathcal{Y}$ denote the evaluation map for the normalization of the universal family over $M_{red}$.  Assume that $ev$ is dominant.
Let $g: \mathcal{S} \to \mathcal{Y}$ denote the finite part of the Stein factorization of $ev$ and let $N$ denote the family of sections on $\mathcal{S}$ corresponding to general members of $M$.  Let $\rho: \mathcal{S}' \to \mathcal{S}$ be a birational map from a smooth projective variety that flattens the evaluation map for the normalization of the universal family over $N$ as in Construction \ref{cons:flatteningfamilyofcurves}.  Let $C'$ denote the strict transform on $\mathcal{S}'$ of a general section on $\mathcal{S}$ parametrized by $N$.

Suppose that $\mathcal{S}'$ is equipped with a dominant rational map $\psi: \mathcal{S}' \dashrightarrow \mathcal{T}$ over $B$ where $\mathcal{T}$ is a normal projective $B$-variety and $\psi$ has connected fibers.  Let $\mathcal{G}$ denote the foliation induced by $\psi$.   Furthermore assume that
\begin{align*}
\mu^{max}_{[C']}(\mathcal{G}) & \geq  (J + 2g(B) + \gamma - 1) > \mu^{max}_{[C']}(T_{\mathcal{S}'}/\mathcal{G}).
\end{align*}
where we define $\gamma = (g(B)\dim(\mathcal{Y}) - g(B) + 1)^{2}(\dim(\mathcal{Y})-1)$.  Then either:
\begin{enumerate}
\item the deformations of $C'$ in the main component $\mathcal{P}$ of $\overline{\psi^{-1}(\psi(C'))}$ contain at least $J$ general points of $\mathcal{P}$, or
\item there is a dominant rational map $\phi: \mathcal{S}' \dashrightarrow \mathcal{Z}$ over $B$ to a normal projective $B$-variety $\mathcal Z$ such that $\psi$ factors rationally through $\phi$ and the following holds.  Let $C'$ be a general section in our family.  Let $\mathcal{W}$ denote the main component of $\overline{\phi^{-1}(\phi(C'))}$.  Then there is a resolution $\widetilde{\mathcal{W}}$ of $\mathcal{W}$ and a section $\widetilde{C}$ on $\widetilde{\mathcal{W}}$ that maps to $C'$ such that:
\begin{enumerate}
\item The deformations of $\widetilde{C}$ in $\widetilde{\mathcal{W}}$ contain at least $J$ general points of $\widetilde{\mathcal{W}}$.
\item The space of deformations of $C'$ in $\mathcal{W}$ has codimension at most $(\dim(\mathcal{P})-1)(J + 2g(B) + \gamma)$ in the space of deformations of $C'$ in $\mathcal{P}$.
\item Letting $\mathcal{H}$ denote the foliation induced by $\phi$, we have $\mu_{[C']}^{max}(T_{\mathcal{S}'}/\mathcal{H}) < J+2g(B) + \gamma - 1 \leq \mu_{[C']}^{min}(\mathcal{H})$.
\end{enumerate}
\end{enumerate}
\end{theorem}

\begin{proof}
Let us assume that deformations of $C'$ do not go through $J$ general points of $\mathcal{P}$.  Since $C'$ deforms in a flat family on $\mathcal{S'}$, a general section $C'$ in the family will be contained in the smooth locus of $\mathcal{G}$.  Thus, if we take a resolution $\widetilde{\mathcal{P}}$ and consider the strict transform $C^{\dagger}$ then Construction \ref{cons:foliationrestriction} shows that the normal bundle of $C^{\dagger}$ in $\widetilde{\mathcal{P}}$ is isomorphic to $\mathcal{G}|_{C'}$.  By Lemma \ref{lemma:hnfreecurves} our deformation assumption implies that
\begin{equation*} 
\mu^{min}(\mathcal{G}|_{C'}) < J + 2g(B) - 1.
\end{equation*}
Applying Theorem \ref{theo:hnfforsections} we obtain
\begin{align} 
\mu^{min}_{[C']}(\mathcal{G}) & \leq \mu^{min}(\mathcal{G}|_{C'}) + \gamma \nonumber \\
& < J + 2g(B) + \gamma - 1 \label{eq:minslopegbound}
\end{align}
Write the Harder-Narasimhan filtration of $\mathcal{G}$ with respect to $[C']$ as  
\begin{equation*}
0= \FF_0 \subset \FF_1 \subset \dots \subset \FF_k = \mathcal{G}.
\end{equation*}
By assumption $\mu^{max}_{[C']}(\mathcal{G}) \geq J + 2g(B) + \gamma - 1$ so there is some index $i \geq 1$ such that we have $\mu_{[C']}^{min}(\FF_{i}) \geq J + 2g(B) + \gamma - 1$.  Let $i$ be the maximum index for which this inequality holds.  On the one hand, since Equation \eqref{eq:minslopegbound} shows that $\mu^{min}_{[C']}(\mathcal{G}) < J + 2g(B) + \gamma - 1$ we must have $i < k$.  On the other hand, since $i$ was selected to be as large as possible we must have
\begin{equation*}
\mu^{max}_{[C']}(\mathcal{G} / \FF_i) <  J + 2g(B) + \gamma - 1.
\end{equation*}

We claim that $\FF_{i}$ is a foliation on $\mathcal{S}'$.  By Theorem \ref{theo:HNisfoliation} it suffices to check that $\FF_{i}$ is a term in the Harder-Narasimhan filtration of $T_{\mathcal{S}'}$ with respect to $[C']$.
By assumption $\mu_{[C']}^{max}(T_{\mathcal{S}'}/\mathcal{G}) < J+2g(B)+\gamma - 1 \leq \mu_{[C']}^{min}(\mathcal{F}_{i})$ and thus the Harder-Narasimhan filtration of $T_{\mathcal{S}'}$ agrees with the Harder-Narasimhan filtration of $\mathcal{G}$ up to the $i$th entry, proving the claim.

By \cite[Theorem 1.1]{CP19} the foliation $\FF_i$ is induced by a rational map $\phi: \mathcal{S}' \dashrightarrow \mathcal{Z}$ over $B$ that has connected fibers.  Since $i < k$ this rational map is not trivial.  By our flatness assumption a general section $C'$ will be contained in the regular locus of $\FF_{i}$.
Let $\widetilde{\mathcal{W}}$ denote a resolution of the main component of $\overline{\phi^{-1}(\phi(C'))}$ and let $\widetilde{C}$ denote the section chosen as in Construction \ref{cons:foliationrestriction}.  In particular we have
\begin{equation*}
T_{\widetilde{\mathcal{W}}/B}|_{\widetilde{C}} \cong \mathcal{F}_{i}|_{C'}.
\end{equation*}
Theorem \ref{theo:hnfforsections} implies that
\begin{align*}
\mu^{min}(\mathcal{F}_{i}|_{C'}) & \geq \mu_{[C']}^{min}(\mathcal{F}_{i}) - \gamma \\
& \geq J + 2g(B) - 1
\end{align*}
and so by Lemma \ref{lemma:hnfreecurves} we see that $\widetilde{C}$ can go through at least $J$ general points of $\widetilde{\mathcal{W}}$ verifying (a).  To prove (b), let $N_{\mathcal{P}}$ denote the space of deformations of the strict transform of $C'$ in $\widetilde{\mathcal{P}}$ and let $N_{\mathcal{W}}$ denote the space of deformations of $\widetilde{C}$ in $\widetilde{\mathcal W}$.  Appealing to Corollary \ref{coro:domfamilyexpdim} to give upper and lower bounds on the dimension of a dominant family of sections based on the degree of the normal bundle, we see that
\begin{align*}
\dim(N_{\mathcal{P}}) - \dim(N_{\mathcal{W}}) & \leq (c_{1}(\mathcal{G}) \cdot C' + (\dim(\widetilde{\mathcal{P}}) - 1)) - (c_{1}(\FF_{i}) \cdot C' + (\dim(\widetilde{\mathcal{W}}) - 1)(1-g(B))) \\
& = c_{1}(\mathcal{G}/\FF_{i}) \cdot C' + (\dim(\mathcal{P}) - \dim(\mathcal{W})) + g(B) (\dim(\widetilde{\mathcal{W}}) - 1)
\end{align*}
Since $c_{1}(\mathcal{G}/\FF_{i}) \cdot C' \leq \rk(\mathcal{G} / \FF_i) \cdot  \mu^{max}_{[C']}(\mathcal{G} / \FF_i) < (\dim(\mathcal{P}) - \dim(\mathcal{W})) \cdot (J+2g(B)+\gamma-1)$, we can continue the chain of inequalities as follows:
\begin{align*}
\dim(N_{\mathcal{P}}) - \dim(N_{\mathcal{W}})& < (\dim(\mathcal{P}) - \dim(\mathcal{W})) (J + 2g(B) + \gamma) + g(B) (\dim(\widetilde{\mathcal{W}}) - 1) \\
& \leq (\dim(\mathcal{P})-1)(J + 2g(B) + \gamma)
\end{align*}
Since the dimension of the space of sections is birationally invariant, we obtain (b).  Finally, by construction $\mu^{min}_{[C']}(\mathcal{F}_{i}) \geq J+2g(B)+\gamma - 1$.  On the other hand we have already seen that both $\mu^{max}_{[C']}(\mathcal{G}/\mathcal{F}_{i})$ and $\mu^{max}_{[C']}(T_{\mathcal{S}'}/\mathcal{G})$ are strictly less than $J+2g(B)+\gamma - 1$.  This implies (c).
\end{proof}

\section{Twists over function fields of complex curves} \label{sect:twists}

Let $B$ be a smooth projective curve over an algebraically closed field $k$ of characteristic $0$.  
Suppose we have a dominant  generically finite morphism $f_{\eta}: \mathcal{Y}_{\eta} \to \mathcal{X}_{\eta}$ between normal projective $K(B)$-varieties.  In this section we study the set of twists of $f_{\eta}$. Recall that a twist of $f_\eta$ is a generically finite $K(B)$-morphism $f'_\eta : \mathcal Y'_\eta \to \mathcal X_\eta$ such that there is an $\mathcal X_{\overline{\eta}}$-isomorphism between $\mathcal Y_{\overline{\eta}}$ and $\mathcal Y'_{\overline{\eta}}$ (where the subscript $\overline{\eta}$ denotes the base change to $\Spec \overline{K(B)}$).

In Section \ref{sect:hurwitzspace} we discuss the Hurwitz space as described by \cite{Wewers}.  Using this construction, we show in Section \ref{sect:familyoftwists} that the set of twists of a dominant generically finite map $f_{\eta}: \mathcal{Y}_{\eta} \to \mathcal{X}_{\eta}$ can be parametrized by a scheme with countably many components.  We will not construct a universal stack, since there are some steps in the construction which might not be valid in the setting of stacks.  Instead, we will construct a morphism of schemes such that every twist of $f_{\eta}$ is the fiber over some closed point.

The remainder of the section is devoted to analyzing the canonical divisor for twists.  In Section \ref{sect:localtoglobal}, we prove a local-to-global principle (Corollary \ref{coro:localtoglobal}) for the Galois cohomology group parametrizing twists of $f_{\eta}$.  In particular the local invariant gives us a convenient way to identify the places of $K(B)$ where two twists are ``the same'' locally.  In Section \ref{sect:henselslemma} we apply Hensel's Lemma to give a geometric criterion that will guarantee the vanishing of the local invariant.  Finally, in Section \ref{sect:splittingfields} we analyze how the canonical divisor changes as we choose different twists of $f_{\eta}$.  The key point is that its positivity is controlled by the places of $K(B)$ where the local invariant does not vanish.  In particular, we show that if we bound the positivity of the canonical divisor then the parameter space of twists has finite type (Corollary~\ref{coro:boundedintimpliesboundedtwists}).

\subsection{Hurwitz space} \label{sect:hurwitzspace}
The starting point is the following version of the Hurwitz space: 

\begin{theorem}[\cite{Wewers}]
Let $B$ be a smooth projective curve.  Fix a positive integer $r$ and a finite group $G$.  There is a smooth Deligne-Mumford stack $\mathcal{H}(G,r,B)$ parametrizing pairs $(q,\psi)$ where $q: C \to B$ is a Galois morphism from a smooth projective curve $C$ that has $r$ branch points and $\psi$ is an isomorphism $\psi: \Aut(C/B) \to G$.
\end{theorem}

Suppose we fix a finite group $G$ and set $\mathcal{H}(G,B) = \sqcup_{r} \mathcal{H}(G,r,B)$.  Then we can think of $\mathcal{H}(G,B)$ as a parameter space for pairs $(C/B,\psi)$ where $C/B$ is a finite Galois cover and $\psi: \Gal(K(C)/K(B)) \to G$ is an isomorphism.  We denote the universal family over $\mathcal{H}(G,B)$ by $\mathcal{U}(G,B) \to \mathcal{H}(G,B)$.  This means that there is a morphism $\mathcal{U}(G,B) \to \mathcal{H}(G,B) \times B$ which over every point of the form $\Spec (k) \to \mathcal{H}(G,B)$ is the corresponding cover $C \to B$ with an isomorphism $\psi: \Gal(K(C)/K(B)) \to G$.

Since our parameter space includes the data of an isomorphism $\psi: \Aut(C/B) \to G$, the fiber of $\mathcal{G} := G \times \mathcal H(G, B)$ over $(C/B, \psi) \in \mathcal{H}(G,B)$ can be canonically identified with the Galois group $\mathrm{Gal}(K(C)/K(B))$.

\subsection{The space of twists} \label{sect:familyoftwists}

We fix a dominant generically finite map $f_{\eta}: \mathcal{Y}_{\eta} \to \mathcal{X}_{\eta}$ of normal geometrically integral projective $K(B)$-schemes.

Let $G$ be a finite group.
To avoid some stacky issues in our constructions, we will fix an \'etale cover $\mathcal H_G \to \mathcal H(G, B)$ from a scheme $\mathcal{H}_{G}$ whose irreducible components are varieties of finite type over $\mathbb{C}$.  We denote the pullback of the universal family by $\mathcal{U}_{\mathcal H_G} \to \mathcal H_G$.  We will use $(C/B,\psi)$ to denote any closed point of $\mathcal{H}_G$ such that the corresponding fiber of $\mathcal{U}_{\mathcal H_G} \to \mathcal{H}_G$ is the map $C \to B$ equipped with the isomorphism $\psi$. We let $\mathcal{G}_{\mathcal H_G} \to \mathcal{H}_G$ denote the morphism whose fiber over $(C/B, \psi) \in \mathcal{H}_G$ is the corresponding Galois group, i.e., $\mathcal G_{\mathcal H_G} = G \times \mathcal H_G$.
We let $\overline{\mathcal{K}}(\mathcal Y_\eta/\mathcal X_\eta) _{\mathcal H_G} := \overline{K}(\mathcal Y_\eta/\mathcal X_\eta)  \times \mathcal{H}_G$ denote the trivial group scheme over $\mathcal{H}_G$ associated to
\[
\overline{K}(\mathcal Y_\eta/\mathcal X_\eta) = \mathrm{Aut}(\mathcal Y_{\overline{\eta}}/\mathcal X_{\overline{\eta}}).
\]   
We consider the universal family $\mathcal U_{\mathcal H_G} \to \mathcal H_G \times B$ and its base change $\mathcal U^*_{\mathcal H_G} \to \mathcal H_G \times \Spec K(B)$.
Since $\mathcal Y_\eta \times_{\Spec K(B)}\mathcal U^*_{\mathcal H_G}$ is flat and projective over $\mathcal H_G\times \Spec K(B)$ and $\mathcal X_\eta \times_{\Spec K(B)} \mathcal U^*_{\mathcal H_G}$ is projective over $\mathcal H_G\times \Spec K(B)$, by \cite[I.1.10 Theorem]{Kollar} we can define the relative automorphism scheme
\[
\widetilde{\mathcal K}(\mathcal Y_\eta/\mathcal X_\eta) _{\mathcal H_G} := \mathrm{Aut}_{\mathcal H_G \times \Spec K(B)}(\mathcal Y_\eta \times_{\Spec K(B)}\mathcal U^*_{\mathcal H_G}/\mathcal X_\eta \times_{\Spec K(B)} \mathcal U^*_{\mathcal H_G}).
\]
This is a quasi-finite group scheme over $\mathcal H_G \times \Spec K(B)$.

Since $\widetilde{\mathcal K}(\mathcal Y_\eta/\mathcal X_\eta)_{\mathcal H_G}$ can be embedded into $\overline{\mathcal K}(\mathcal Y_\eta/\mathcal X_\eta) _{\mathcal H_G} \times \Spec K(B)$ as a $\mathcal H_G \times \Spec K(B)$-scheme, we conclude that $\widetilde{\mathcal K}(\mathcal Y_\eta/\mathcal X_\eta) _{\mathcal H_G}$ is quasi-affine over $\mathcal H_G \times \Spec K(B)$.
Using the functoriality of the $\mathrm{Aut}_{\mathcal H_G \times \Spec K(B)}$-functor we can construct descent data for the quasi-finite group scheme $\widetilde{\mathcal K}(\mathcal Y_\eta/\mathcal X_\eta) _{\mathcal H_G}\to \mathcal H_G \times \Spec K(B)$ with respect to the map $\mathcal{H}_G \times \Spec K(B) \to \mathcal{H}_G$. Indeed, we denote by $p_1, p_2$ the two projections $\mathcal H_G \times \Spec K(B) \times \Spec K(B) \to \mathcal H_G \times \Spec K(B)$. Then both $p_1^*\widetilde{\mathcal K}(\mathcal Y_\eta/\mathcal X_\eta) _{\mathcal H_G}$ and $p_2^*\widetilde{\mathcal K}(\mathcal Y_\eta/\mathcal X_\eta) _{\mathcal H_G}$ are canonically isomorphic to 
\[
\mathrm{Aut}_{\mathcal H_G \times \Spec K(B)\times \Spec K(B)}(\mathcal Y_\eta \times_{\Spec K(B)}\mathcal U^*_{\mathcal H_G} \times \Spec K(B)/ \mathcal X_\eta \times_{\Spec K(B)} \mathcal U^*_{\mathcal H_G} \times \Spec K(B)).
\] 
This defines the canonical descent data $p_1^*\widetilde{\mathcal K}(\mathcal Y_\eta/\mathcal X_\eta) _{\mathcal H_G} \to p_2^*\widetilde{\mathcal K}(\mathcal Y_\eta/\mathcal X_\eta) _{\mathcal H_G}$. Then it is easy to check that this data satisfies the gluing condition.
By the fpqc descent theory for quasi affine schemes as in \cite[Theorem 4.3.5(ii)]{Poonen} we conclude the existence of a quasi-finite group scheme $\mathcal K(\mathcal Y_\eta/\mathcal X_\eta) _{\mathcal H_G} \to \mathcal H_G$ whose base change to $\mathcal{H}_G \times \Spec K(B)$ is isomorphic to $\widetilde{\mathcal K}(\mathcal Y_\eta/\mathcal X_\eta) _{\mathcal H_G}$.
Note that this is a locally closed subgroup scheme of $\overline{\mathcal K}(\mathcal Y_\eta/\mathcal X_\eta)_{\mathcal H_G}$, so in particular $\mathcal K(\mathcal Y_\eta/\mathcal X_\eta) _{\mathcal H_G}$ is quasi-affine over $\mathcal H_G$.  For $(C/B, \psi) \in \mathcal H_G$ consider the Galois action by conjugation 
\[
\phi_{C/B, \psi}: \Gal(K(C)/K(B)) \times (\mathcal K(\mathcal Y_\eta/\mathcal X_\eta)_{\mathcal H_G})_{(K(C)/K(B), \psi)} \to (\mathcal K(\mathcal Y_\eta/\mathcal X_\eta)_{\mathcal H_G})_{(K(C)/K(B), \psi)}.
\]
This fiberwise Galois action defines a group scheme action
\[
\phi: \mathcal G_{\mathcal H_G} \times_{\mathcal H_G} \mathcal K(\mathcal Y_\eta/\mathcal X_\eta)_{\mathcal H_G} \to \mathcal K(\mathcal Y_\eta/\mathcal X_\eta)_{\mathcal H_G}.
\]
Consider the morphism scheme $\mathrm{Mor}_{\mathcal H_G}(\mathcal G_{\mathcal H_G}, \mathcal K(\mathcal Y_\eta/\mathcal X_\eta)_{\mathcal H_G})$.
We define the space 
\[
\mathcal C^1(\mathcal G_{\mathcal H_G}, \mathcal K(\mathcal Y_\eta/\mathcal X_\eta)_{\mathcal H_G})
\] 
of 1-cocycles as the closed subscheme of $\mathrm{Mor}_{\mathcal H_G}(\mathcal G_{\mathcal H_G}, \mathcal K(\mathcal Y_\eta/\mathcal X_\eta)_{\mathcal H_G})$ consisting of $1$-cocycles $(C/B, \psi, \sigma : G \to (\mathcal K(\mathcal Y_\eta/\mathcal X_\eta)_{\mathcal H_G})_{(K(C)/K(B), \psi)} )$ which satisfy the cocycle condition $\sigma_{st} = \sigma_{s} \phi_{(C/B, \psi)}(s)(\sigma_{t})$.

Next our goal is to construct a family of twists of $(\mathcal Y_\eta/\mathcal X_\eta)$ over $\mathcal C^1(\mathcal G_{\mathcal H_G}, \mathcal K(\mathcal Y_\eta/\mathcal X_\eta)_{\mathcal H_G})$. Define $\mathcal{Y}' = \mathcal Y_\eta \times_{\Spec K(B)} \mathcal U^*_{\mathcal H_G}$.   
Then the fiber of the projection $\mathcal{Y}' \to \mathcal{H}_G$ over $(C/B,\psi)$ is isomorphic to $\mathcal Y_\eta\otimes K(C)$.  We consider
\[
\mathcal Y' \times_{\mathcal H_G} \mathcal C^1(\mathcal G_{\mathcal H_G}, \mathcal K(\mathcal Y_\eta/\mathcal X_\eta)_{\mathcal H_G}) \to \mathcal X_\eta \times_{\Spec K(B)}\mathcal U^*_{\mathcal H_G}.
\]
There is a group scheme action of $\mathcal G_{\mathcal H_G}$ on this fiber product by
\[
(s, (C/B, \psi)) \cdot (y, \sigma, (C/B, \psi)) = (\sigma_{s}\circ (1\otimes s) (y), \sigma, (C/B, \psi)).
\]
We let $\widetilde{\mathcal{Y}}$ denote the quotient of $\mathcal Y' \times_{\mathcal H_G} \mathcal C^1(\mathcal G_{\mathcal H_G}, \mathcal K(\mathcal Y_\eta/\mathcal X_\eta)_{\mathcal H_G})$ by the finite flat group scheme $\mathcal G_{\mathcal H_G}$.  Then $\widetilde{\mathcal{Y}}$ comes equipped with a map
\[
\widetilde{\mathcal Y} \to \mathcal C^1(\mathcal G_{\mathcal H_G}, \mathcal K(\mathcal Y_\eta/\mathcal X_\eta)_{\mathcal H_G}) \times \mathcal X_\eta.
\]
such that the fiber over $(\sigma, (C/B, \psi)) \in \mathcal C^1(\mathcal G_{\mathcal H_G}, \mathcal K(\mathcal Y_\eta/\mathcal X_\eta)_{\mathcal H_G})$ is the map
\[
\mathcal Y_\eta^\sigma \to \mathcal X_\eta,
\]
where $\mathcal Y_\eta^\sigma$ is the quotient of $\mathcal Y_\eta \otimes K(C)$ by the Galois action
\[
\Gal(K(C)/K(B)) \ni s \mapsto \sigma_s \circ 1\otimes s \in \mathrm{Aut}(\mathcal Y_\eta \otimes K(C)/\mathcal X_\eta).
\]
By construction every twist of $\mathcal{Y}_{\eta} \to \mathcal{X}_{\eta}$ is parametrized by the fiber over some point $(\sigma, (C/B, \psi)) \in \mathcal C^1(\mathcal G_{\mathcal H_G}, \mathcal K(\mathcal Y_\eta/\mathcal X_\eta)_{\mathcal H_G})$.  

Note that the scheme $\mathcal C^1(\mathcal G_{\mathcal H_G}, \mathcal K(\mathcal Y_\eta/\mathcal X_\eta)_{\mathcal H_G})$ constructed above need not have finite type over $K(B)$.  However, if we fix certain invariants then the corresponding subscheme will have finite type.

\begin{lemma} \label{lemm:twistboundedconditions}
Fix a smooth projective curve $B$ and positive integers $d,b$. 
Suppose we have a generically finite dominant $K(B)$-morphism $f_{\eta}: \mathcal{Y}_{\eta} \to \mathcal{X}_{\eta}$ where $\mathcal{X}_{\eta}$ and $\mathcal{Y}_{\eta}$ are normal projective varieties.  Let $S$ denote the set of twists $f^{\sigma}_{\eta}$ such that $\mathcal{Y}^{\sigma}_{\eta}$ and $\mathcal{Y}_{\eta}$ become isomorphic after a base change by a Galois extension $K(C)/K(B)$ whose degree is $\leq d$ and whose branch locus consists of at most $b$ points. 

There is a finite type scheme $R$ over $\mathbb C$ and morphisms $\psi: \mathcal{U}_R \to R$, $g: \mathcal{U}_R \to \mathcal{X}_{\eta}$ such that every element $\mathcal{Y}^{\sigma}_{\eta} \in S$ is isomorphic to the fiber of $\psi$ over some closed point $t \in R$ and $f^{\sigma}_{\eta} = g|_{\psi^{-1}t}$.
\end{lemma}

\begin{proof}
There are finitely many isomorphism classes of finite groups $G$ of order $\leq d$.  As we vary $G$ over all such groups and vary over all $r \leq b$, we obtain a finite type Deligne-Mumford stack $\sqcup \mathcal{H}(G,r,B)$ parametrizing extensions $K(C)/K(B)$ and automorphisms $\Aut(C/B) \to G$. Let $\mathcal H_{G, r}$ be the preimage of $\mathcal H(G, r, B)$ via $\mathcal H_G \to \sqcup_{r} \mathcal{H}(G,r,B)$.  Then the space 
\[\sqcup_{G, r}\mathcal C^1(\mathcal G_{\mathcal H_{G, r}}, \mathcal K(\mathcal Y_\eta/\mathcal X_\eta)_{\mathcal H_{G, r}})
\]
 is a finite type scheme over $\mathbb C$ where $\mathcal C^1(\mathcal G_{\mathcal H_{G, r}}, \mathcal K(\mathcal Y_\eta/\mathcal X_\eta)_{\mathcal H_{G, r}})$ is the base change of 
\[
 \mathcal C^1(\mathcal G_{\mathcal H_{G}}, \mathcal K(\mathcal Y_\eta/\mathcal X_\eta)_{\mathcal H_{G}})
\]
via $\mathcal H_{G, r} \to \mathcal H_G$. Thus our assertion follows.
\end{proof}

\subsubsection{The spaces of twists in families}
Here we perform the constructions in the previous section in families.
As before we fix a smooth projective curve $B$ defined over $\mathbb C$.
Let $\mathfrak X \to \mathfrak S \times \Spec K(B), \mathfrak Y \to \mathfrak S \times \Spec K(B)$ be flat families of normal projective $K(B)$-varieties where $\mathfrak S$ is a scheme of finite type over $\mathbb C$. We also assume that we have a $\mathfrak S \times \Spec K(B)$-morphism $f : \mathfrak Y \to \mathfrak X$ which is fiberwise dominant and generically finite.  

We fix a finite group $G$ and take an \'etale open cover $\mathcal H_G \to \mathcal H(G, B)$ where $\mathcal H_G$ is a scheme over $\mathbb C$.  As before we denote the pullback of the universal family by $\mathcal U_{\mathcal H_G} \to \mathcal H_G \times B$ and consider its base change $\mathcal U^*_{\mathcal H_G} \to \mathcal H_G \times \Spec K(B)$. Since $\mathfrak Y \times_{\Spec K(B)}\mathcal U^*_{\mathcal H_G}$ is projective and flat over $\mathfrak S \times \mathcal H \times \Spec K(B)$ and $\mathfrak X \times_{\Spec K(B)} \mathcal U^*_{\mathcal H_G}$ is projective over $\mathfrak S \times \mathcal H_G \times \Spec K(B)$, we can define the relative automorphism group
\[
\widetilde{\mathcal K}(\mathfrak Y/\mathfrak X)_{\mathfrak S \times \mathcal H_G} = \mathrm{Aut}_{\mathfrak S \times \mathcal H_G \times \Spec K(B)}(\mathfrak Y \times_{\Spec K(B)}\mathcal U^*_{\mathcal H_G}/\mathfrak X \times_{\Spec K(B)} \mathcal U^*_{\mathcal H_G}).
\]
This is a quasi-finite group scheme over $\mathfrak S \times \mathcal H_G \times \Spec K(B)$.
Since the above relative automorphism group is also separated over $\mathfrak S \times \mathcal H_G \times \Spec K(B)$,
$\widetilde{\mathcal K}(\mathfrak Y/\mathfrak X)_{\mathfrak S \times \mathcal H_G}$ is quasi-affine over $\mathfrak S \times \mathcal H_G \times \Spec K(B)$.
 Using fpqc descent theory, $\widetilde{\mathcal K}(\mathfrak Y/\mathfrak X)_{\mathfrak S \times \mathcal H_G} \to \mathfrak S \times \mathcal H_G \times \Spec K(B)$ descends to $\mathcal K(\mathfrak Y/\mathfrak X)_{\mathfrak S \times \mathcal H_G} \to \mathfrak S \times \mathcal H_G$.

Let $\mathcal G_{\mathfrak S \times \mathcal H_G} = G \times \mathfrak S \times \mathcal H_G$ and consider the natural conjugation group action
\[
\phi : \mathcal G_{\mathfrak S \times \mathcal H_G} \times_{\mathfrak S \times \mathcal H_G} \mathcal K(\mathfrak Y/\mathfrak X)_{\mathfrak S \times\mathcal H_G} \to \mathcal K(\mathfrak Y/\mathfrak X)_{\mathfrak S \times \mathcal H_G}.
\]
Consider the morphism scheme $\mathrm{Mor}_{\mathfrak S \times \mathcal H_G} (\mathcal G_{\mathfrak S \times \mathcal H_G}, \mathcal K(\mathfrak Y/\mathfrak X)_{\mathfrak S \times \mathcal H_G})$ and the closed subscheme consisting of the space of $1$-cycles $\mathcal C^1(\mathcal G_{\mathfrak S \times \mathcal H_G}, \mathcal K(\mathfrak Y/\mathfrak X)_{\mathfrak S \times \mathcal H_G})$.
Define $\mathfrak Y' = \mathfrak Y \times_{\Spec K(B)} \mathcal U^*_{\mathcal H_G}$ as a scheme over $\mathfrak S \times \mathcal H_G$.
Again we have a natural group action of $\mathcal G_{\mathfrak S \times \mathcal H_G}$ on $\mathfrak Y' \times_{\mathfrak S \times \mathcal H_G} \mathcal C^1(\mathcal G_{\mathfrak S \times \mathcal H_G}, \mathcal K(\mathfrak Y/\mathfrak X)_{\mathfrak S \times \mathcal H_G})$.
We define $\widetilde{\mathfrak Y}$ to be the quotient of this group action.  It comes equipped with a morphism $\widetilde{\mathfrak Y} \to \mathcal C^1(\mathcal G_{\mathfrak S \times \mathcal H_G}, \mathcal K(\mathfrak Y/\mathfrak X)_{\mathfrak S \times \mathcal H_G}) \times_{\mathfrak S} \mathfrak X$ realizing $$\mathcal C^1(\mathcal G_{\mathfrak S \times \mathcal H_G}, \mathcal K(\mathfrak Y/\mathfrak X)_{\mathfrak S \times \mathcal H_G})$$ as the parameter space of twists of the maps $f_{s, \eta} : \mathfrak Y_{s, \eta} \to \mathfrak X_{s, \eta}$ for closed points $s \in \mathfrak S$.

Regarding this family we have the following boundedness statement:

\begin{lemma} \label{lemm:twistboundedconditions2}
Fix a smooth projective curve $B$ and positive integers $d,b$. Let $p : \mathfrak X \to \mathfrak S \times \Spec K(B), q : \mathfrak Y \to \mathfrak S \times \Spec K(B)$ be flat families of normal projective $K(B)$-varieties where $\mathfrak S$ is a scheme of finite type over $\mathbb C$. We also assume that we have a $\mathfrak S \times \Spec K(B)$-morphism $f : \mathfrak Y \to \mathfrak X$ which is fiberwise dominant and generically finite.

Let $A$ denote the set of twists $f^{\sigma}_{\eta,s}: \mathfrak{Y}_{\eta,s}^{\sigma} \to \mathcal{X}_{\eta}$ where $s$ is a closed point of $\mathfrak S$ and $\mathfrak{Y}^{\sigma}_{\eta,s}$ and $\mathfrak{Y}_{\eta,s}$ become isomorphic after a base change by a Galois extension $K(C)/K(B)$ whose degree is $\leq d$ and whose branch locus consists of at most $b$ points.

There is a finite type scheme $\mathfrak R$ over $\mathfrak S$ and morphisms $\psi: \mathfrak{U}_{\mathfrak R} \to \mathfrak R$, $g: \mathfrak{U}_{\mathfrak R} \to \mathfrak R \times_{\mathfrak S} \mathfrak{X}$ such that every element $ f^{\sigma}_{\eta,s}: \mathcal{Y}^{\sigma}_{\eta,s} \to \mathcal X_{\eta, s} \in A$ is isomorphic to the fiber of $\psi$ over some closed point $t \in \mathfrak R$ and $f^{\sigma}_{\eta,s}= g|_{\psi^{-1}t}$.
\end{lemma}

\begin{proof}
Lemma~\ref{lemm:twistboundedconditions} proves this statement for a single morphism; the proof for a finite-type family of morphisms uses exactly the same argument in a relative setting.
\end{proof}

\subsubsection{Functoriality}
\label{subsec:functoriality}
Let $\mathfrak X \to \mathfrak S \times \Spec K(B),  \mathfrak Y \to \mathfrak S \times \Spec K(B)$ be flat families of normal projective $K(B)$-varieties where $\mathfrak S$ is a smooth scheme of finite type over $\mathbb C$. We also assume that we have a $\mathfrak S \times \Spec K(B)$-morphism $f : \mathfrak Y \to \mathfrak X$ which is fiberwise dominant and generically finite. We further assume that we have flat families $\mathfrak W \to \mathfrak S \times \Spec K(B)$, $\mathfrak T \to \mathfrak S \times \Spec K(B)$ of projective varieties such that $\mathfrak T_s$ is normal for any $s \in \mathfrak S$ and we have a commutative diagram
\begin{equation*}
\xymatrix{ {\mathfrak Y} \ar[r]^{f} \ar[d]_{r} &  {\mathfrak X}\ar[d]_{p} \\
\mathfrak T \ar[r]_{t} & \mathfrak W}
\end{equation*}
over $\mathfrak S \times \Spec K(B)$ where $p, r$ are dominant with connected fibers and $t$ is dominant, finite, and fiberwise Galois over $\mathfrak S$, i.e., each fiber $\mathfrak T_s \to \mathfrak W_s$ is a finite Galois cover.   We also assume that the Stein factorization of $\mathfrak Y \to \mathfrak X \to \mathfrak W$ is given by $\mathfrak Y \to \mathfrak T$. Since a relative automorphism induces a relative automorphism of Stein factorizations, we obtain a homomorphism
\[
\mathrm{Aut}_{\mathfrak S \times \mathcal H_G}(\mathfrak Y \times_{\Spec K(B)} \mathcal U^*_{\mathcal H_G}/\mathfrak X \times_{\Spec K(B)}\mathcal U^*_{\mathcal H_G}) \to \mathrm{Aut}_{\mathfrak S \times \mathcal H_G}(\mathfrak T \times_{\Spec K(B)} \mathcal U^*_{\mathcal H_G}/\mathfrak W \times_{\Spec K(B)}\mathcal U^*_{\mathcal H_G}),
\]
and this induces a morphism 
\[
\mathcal C^1(\mathcal G_{\mathfrak S \times \mathcal H_G}, \mathcal K(\mathfrak Y/\mathfrak X)_{\mathfrak S \times \mathcal H_G}) \to \mathcal C^1(\mathcal G_{\mathfrak S \times \mathcal H_G}, \mathcal K(\mathfrak T/\mathfrak W)_{\mathfrak S \times \mathcal H_G}).
\]

\subsection{Local-to-global principle} \label{sect:localtoglobal}

Let $K(B)$ be the function field of a smooth projective curve $B$. Let $f : \mathcal Y_\eta \to \mathcal X_\eta$ be a dominant generically finite morphism between normal projective varieties $\mathcal Y_\eta$ and $\mathcal X_\eta$ defined over $K(B)$.

We fix a place $\overline{\nu}$ of $\overline{K(B)}$ over a place $\nu$ of $K(B)$. This specifies for every finite cover $C \to B$ a point $p_{\overline{\nu}, C}$ on $C$ such that for every factoring $C \xrightarrow{g} C' \to B$ we have $g(p_{\overline{\nu}, C}) = p_{\overline{\nu}, C'}$.  Consider the decomposition group
\[
D_{\overline{\nu}} = \{ \sigma \in \mathrm{Gal}(K(B)) \, | \, \sigma (\overline{\nu}) =\overline{\nu} \}.
\]
This is isomorphic to
\[
\mathrm{Gal}(K(B)_\nu) \cong \varprojlim (\mathbb Z/N\mathbb Z)^{\times}.
\]
Note that if we have two places $\overline{\nu}, \overline{\nu}'$ corresponding to the same place $\nu$ on $K(B)$, then $D_{\overline{\nu}}$ and $D_{\overline{\nu}'}$ are conjugate to each other in $\mathrm{Gal}(K(B))$.
Recall that the Galois group acts on $\mathrm{Aut}(\mathcal Y_{\overline{\eta}}/\mathcal X_{\overline{\eta}})$ by conjugation, and in this way one can consider Galois cohomology
\[
H^1(\mathrm{Gal}(K(B)), \mathrm{Aut}(\mathcal Y_{\overline{\eta}}/\mathcal X_{\overline{\eta}}))
\]
The injection $\mathrm{Gal}(K(B)_\nu) \cong D_{\overline{\nu}} \subset \Gal(K(B))$ induces a map on Galois cohomology
\[
H^1(\mathrm{Gal}(K(B)), \mathrm{Aut}(\mathcal Y_{\overline{\eta}}/\mathcal X_{\overline{\eta}})) \to H^1(\mathrm{Gal}(K(B)_\nu), \mathrm{Aut}(\mathcal Y_{\overline{\eta}}/\mathcal X_{\overline{\eta}}))
\]
which we denote by $\mathrm{inv}_{\nu}$.  (Although the choice of isomorphism $\mathrm{Gal}(K(B)_\nu) \cong D_{\overline{\nu}}$ depends on the choice of $\overline{\nu}$, the induced map of Galois cohomology only depends on the place $\nu$ up to an isomorphism of the pointed set, justifying our mild abuse of notation.) 

For any twist $[\sigma]$ of $\mathcal Y_\eta/\mathcal X_\eta$ the local invariant $\mathrm{inv}_{\nu}([\sigma])$ vanishes for all but finitely many places $\nu$ of $B$.  Thus we obtain a map
\[
H^1(\mathrm{Gal}(K(B)), \mathrm{Aut}(\mathcal Y_{\overline{\eta}}/\mathcal X_{\overline{\eta}})) \to \bigoplus_{\nu \in B} H^1(\mathrm{Gal}(K(B)_\nu), \mathrm{Aut}(\mathcal Y_{\overline{\eta}}/\mathcal X_{\overline{\eta}}))
\]
and we would like to use this map to establish a local-to-global principle.  Note that this map does not need to be injective; for example, there can be twists of $\mathcal Y_\eta/\mathcal X_\eta$ which are trivialized by an \'etale cover of $B$.  However, we will show that the fibers of this map are finite, which is good enough for our purposes.

\begin{lemma} \label{lemm:boundingdeganddisc}
Let $f : \mathcal Y_\eta \to \mathcal X_\eta$ be a dominant generically finite morphism between normal projective varieties over $K(B)$.
Then there exists a positive integer $d = d(\mathcal Y_\eta/\mathcal X_\eta)$ and a fixed finite subset $\mathcal{P} \subset B$ (depending only on $\mathcal Y_\eta/\mathcal X_\eta$) such that the following holds.  Suppose that 
$$
[\sigma] \in H^1(\mathrm{Gal}(K(B)), \mathrm{Aut}(\mathcal Y_{\overline{\eta}}/\mathcal X_{\overline{\eta}})),
$$ 
is a cohomology class and let $\mathcal{Q} \subset B$ denote the finite set of places $\nu \in B$ with $\mathrm{inv}_{\nu}([\sigma]) \neq 0$.  Then there exists a Galois cover $C \to B$ of degree at most $d$ and whose branch locus is contained in $\mathcal{P} \cup \mathcal{Q}$ such that $\mathcal Y_\eta^\sigma \to \mathcal X_\eta$ splits over $K(C)$. 
\end{lemma}

\begin{proof}
Let $K(B')/K(B)$ be a fixed Galois extension so that
\[
\mathrm{Aut}(\mathcal Y_{\overline{\eta}}/\mathcal X_{\overline{\eta}}) = \mathrm{Aut}(\mathcal Y_\eta\otimes K(B')/\mathcal X_\eta\otimes K(B')).
\]
We define $\mathcal{P}$ to be the set of branch points for $B' \to B$.

We can restrict our cocycle $\sigma: \Gal(K(B)) \to \mathrm{Aut}(\mathcal Y_{\overline{\eta}}/\mathcal X_{\overline{\eta}})$ to the subgroup $\Gal(K(B'))$ to get a cocycle $\sigma'$.  
Then $\sigma' : \Gal(K(B')) \to \mathrm{Aut}(\mathcal Y_\eta\otimes K(B')/\mathcal X_\eta\otimes K(B'))$ is a honest homomorphism because the Galois action of $\Gal(K(B'))$ on $\mathrm{Aut}(\mathcal Y_\eta\otimes K(B')/\mathcal X_\eta\otimes K(B'))$ is trivial. The kernel is an open subgroup and thus defines a Galois cover $C$ over $B'$.  Then the induced cocycle $[\tau] \in H^1(\Gal(K(C)), \mathrm{Aut}(\mathcal Y_\eta\otimes K(C)/\mathcal X_\eta \otimes K(C)))$ is trivial. Thus $\mathcal Y_\eta^\sigma \to \mathcal X_\eta$ splits over $K(C)$. Now note that the degree of $K(C)/K(B')$ is bounded by the order of $\mathrm{Aut}(\mathcal Y_{\overline{\eta}}/\mathcal X_{\overline{\eta}})$ and the degree of $K(B')/K(B)$ only depends on $\mathcal Y_\eta/\mathcal X_\eta$.
Furthermore for any place $\nu \in B$ with $\mathrm{inv}_\nu ([\sigma]) = 0$ and $\nu \not\in \mathcal P$, we have $D_{\overline{\nu}} \subset \mathrm{Gal}(K(C))$. Thus $C/B'$ cannot be ramified over any place $\nu \in B \backslash (\mathcal{Q}\cup \mathcal P)$ and our assertion follows.
\end{proof}

\begin{corollary} \label{coro:localtoglobal}
Let $f : \mathcal Y_\eta \to \mathcal X_\eta$ be a dominant generically finite morphism between normal projective varieties over $K(B)$. The fibers of the local invariant map
\[
H^1(\mathrm{Gal}(K(B)), \mathrm{Aut}(\mathcal Y_{\overline{\eta}}/\mathcal X_{\overline{\eta}})) \to \bigoplus_{\nu \in B} H^1(\mathrm{Gal}(K(B)_\nu), \mathrm{Aut}(\mathcal Y_{\overline{\eta}}/\mathcal X_{\overline{\eta}}))
\]
are finite.
\end{corollary}

\begin{proof}
Suppose that $\xi$ is a element of the direct sum and let $\mathcal{Q} \subset B$ denote the finite set of indices for which the entries of $\xi$ are non-zero.  By Lemma \ref{lemm:boundingdeganddisc}, there is a fixed integer $d$ and a fixed finite set $\mathcal{P} \subset B$ such that any twist that lies in the fiber over $\xi$ is split by a Galois cover $C \to B$ of degree at most $d$ and whose branch locus is contained in $\mathcal{P} \cup \mathcal{Q}$.  There are only finitely many such maps $C \to B$, and for each such map there are only a finite set of twists trivialized by the map.
\end{proof}

\subsection{Hensel's lemma} \label{sect:henselslemma}
Let $B$ be a smooth projective curve. Let $\pi : \mathcal X \to B$ be a good fibration and $f : \mathcal Y \to \mathcal X$ be a dominant finite morphism from a normal projective variety such that $\mathcal Y_\eta$ is geometrically integral.  In this setting we have the equality
\[
\mathrm{Bir}(\mathcal Y_{\overline{\eta}}/\mathcal X_{\overline{\eta}}) = \mathrm{Aut}(\mathcal Y_{\overline{\eta}}/\mathcal X_{\overline{\eta}}).
\]

For each twist $[\sigma] \in H^1(\mathrm{Gal}(K(B)), \mathrm{Aut}(\mathcal Y_{\overline{\eta}}/\mathcal X_{\overline{\eta}}))$ of $\mathcal{Y}_{\eta}/\mathcal{X}_{\eta}$, we can construct an integral model $\mathcal{Y}^{\sigma} \to B$ in the following way.  Suppose that $K(C)/K(B)$ is a Galois extension such that $\mathcal{Y}_{\eta}/\mathcal{X}_{\eta}$ and $\mathcal{Y}^{\sigma}_{\eta}/\mathcal{X}_{\eta}$ are isomorphic after base change to $K(C)$.  Then the cohomology class $[\sigma]$ is represented by a cocycle $\sigma: \Gal(K(C)/K(B)) \to \mathrm{Aut}(\mathcal Y_\eta\otimes K(C)/\mathcal X_\eta \otimes K(C))$.  Let $\widetilde{\mathcal Y}_C$ be the normalization of $\mathcal Y\times_B C$.  Then $\sigma$ defines a homomorphism from $\Gal(K(C)/K(B))$ to the birational automorphism group of $\widetilde{\mathcal{Y}}_{C, \eta}$ over $\mathcal{X}_\eta$, or equivalently, to the birational automorphism group of $\widetilde{\mathcal{Y}}_{C}$ over $\mathcal{X}$.  By construction a birational automorphism of $\widetilde{\mathcal{Y}}_{C}$ over $\mathcal{X}$ is actually an automorphism, so that $\Gal(K(C)/K(B))$ acts on $\widetilde{\mathcal{Y}}_{C}$ via $\sigma$.  We let $\mathcal{Y}^{\sigma}$ denote the quotient.  Note that $\mathcal{Y}^{\sigma}$ is normal and comes equipped with a finite $B$-morphism $f^\sigma: \mathcal Y^\sigma \to \mathcal X$. 

Here we prove the following birational version of Hensel's lemma.

\begin{lemma}  \label{lemma:birationalHensel}
Let $\pi : \mathcal X \to B$ be a good fibration and let $\mathcal Y$ be a normal projective variety equipped with a dominant morphism $\mathcal Y \to B$ such that  $\mathcal Y_\eta$ is geometrically integral. Suppose that $f : \mathcal Y \to \mathcal X$ is a dominant finite $B$-morphism. 
 Fix a place $\nu \in B$ and assume that $\mathcal X_\nu$ is smooth. We also assume that $\mathrm{Aut}_B(\mathcal Y/\mathcal X) \to B$ is flat at $\nu \in B$ and $\mathcal Y_\nu$ is reduced and normal.
 
Suppose that $\mathcal{Y}^{\sigma}$ is an integral model of a twist of $f_{\eta}$ as constructed above and that there is a birational $\mathcal X_\nu$-map $h_{\nu}: \mathcal Y_\nu \dashrightarrow \mathcal Y^\sigma_\nu$. Then there is an $\mathcal X\otimes K(B)_\nu$-isomorphism between $\mathcal Y\otimes_\eta K(B)_\nu$ and $\mathcal Y^\sigma \otimes_\eta K(B)_\nu$. In particular $\mathrm{inv}_\nu(\sigma) = 0$.
\end{lemma}

\begin{proof}
Since $\mathcal Y_\nu$ is reduced and normal, the automorphism group $\mathrm{Bir}(\mathcal Y_\nu/\mathcal X_\nu) = \mathrm{Aut}(\mathcal Y_\nu/\mathcal X_\nu)$ is a reduced finite group.
The flatness of $\mathrm{Aut}_B(\mathcal Y/\mathcal X) \to B$ at $\nu$ implies that the lengths of $\mathrm{Aut}(\mathcal Y_{\overline{\eta}}/\mathcal X_{\overline{\eta}})$ and $\mathrm{Aut}(\mathcal Y_\nu/\mathcal X_\nu)$ are equal.

Also note that since $\mathcal Y_\nu$ is normal and finite over $\mathcal X_\nu$ and $\mathcal Y^\sigma_\nu$ is also finite over $\mathcal X_\nu$, our birational map $h_\nu$ extends to a birational morphism $\mathcal Y_\nu \to \mathcal Y^\sigma_\nu$.

Let us consider the relative $\mathcal X$-birational morphism scheme over $\Spec(\widehat{\mathcal O}_{B, \nu})$:
\[
\mathcal B= \mathrm{BirMor}_{\Spec(\widehat{\mathcal O}_{B, \nu})}(\mathcal Y\times_B \Spec(\widehat{\mathcal O}_{B, \nu}), \mathcal Y^\sigma \times_B\Spec(\widehat{\mathcal O}_{B, \nu}))
\]
equipped with a morphism $\mathcal{B}  \to \Spec(\widehat{\mathcal O}_{B, \nu})$.
(This scheme can be constructed as an open subscheme of the relative Hilbert scheme parametrizing graphs in $(\mathcal Y\times_{\mathcal X} \mathcal Y^\sigma)\times_B \Spec(\widehat{\mathcal O}_{B, \nu})$.)

Note that $\mathcal B$ is an $\mathrm{Aut}_{\Spec(\widehat{\mathcal O}_{B, \nu})}(\mathcal Y \times_B \Spec(\widehat{\mathcal O}_{B, \nu})/\mathcal X\times_B \Spec(\widehat{\mathcal O}_{B, \nu}))$-torsor. 
Using the fact that $\mathrm{Aut}_B(\mathcal Y/\mathcal X) \to B$ is flat at $\nu \in B$ we conclude that
the above relative birational morphism scheme is also flat at $\nu \in B$.

Since we have a birational morphism of fibers $\mathcal Y_\nu \to \mathcal Y^\sigma_\nu$, Hensel's lemma (see, e.g., \cite[Th\'eor\`em 18.5.17]{EGAIV}) implies that we have an $\mathcal X\times \Spec(\widehat{\mathcal O}_{B, \nu})$-birational morphism from $\mathcal Y \times_B \Spec(\widehat{\mathcal O}_{B, \nu})$ to $\mathcal Y^\sigma\times_B \Spec(\widehat{\mathcal O}_{B, \nu})$.
Since $\mathcal Y\otimes_\eta K(B)_\nu$ and $\mathcal Y^\sigma \otimes_\eta K(B)_\nu$ are normal and $ \mathcal Y\otimes_\eta K(B)_\nu \to \mathcal X\otimes K(B)_\nu$, $\mathcal Y^\sigma\otimes_\eta K(B)_\nu \to \mathcal X\otimes K(B)_\nu$ are finite, this birational morphism induces an isomorphism of the generic fibers.
\end{proof}

\subsection{Splitting fields and ramification} \label{sect:splittingfields}

Suppose $\pi: \mathcal{Y} \to B$ is an algebraic fiber space where $\mathcal{Y}$ is a normal projective variety and $B$ is a smooth projective curve.  Let $B' \to B$ be a dominant finite morphism of smooth projective curves.  We will use the term ``normalized base change'' to refer to the normalization of $\mathcal{Y} \times_{B} B'$ equipped with the structure morphism to $B'$.  Note that the normalized base change $\mathcal{Y}'$ admits a dominant finite morphism $\mathcal{Y}' \to \mathcal{Y}$.

\begin{lemma} \label{lemm:genreducedpreserved}
Let $\mathcal{Y}$ be a normal projective variety equipped with a surjective morphism $\pi: \mathcal{Y} \to B$ with connected fibers.  Suppose $B' \to B$ is a dominant finite morphism of smooth projective curves and that $\pi': \mathcal{Y}' \to B'$ is the normalized base change.  Fix a closed point $t \in B$ and let $t' \in B'$ be any point mapping to it. If $\mathcal{Y}_{t}$ is generically reduced, then $\mathcal{Y}'_{t'}$ is also generically reduced and the morphism $\mathcal{Y}'_{t'} \to \mathcal{Y}_{t}$ is birational.
\end{lemma}

\begin{proof}
Since $\mathcal{Y}$ is normal, the open set $U \subset \mathcal{Y}_{t}$ consisting of points that lie in the smooth locus of $\mathcal Y$ and the smooth locus of $\mathcal{Y}_{t}$ is dense in $\mathcal{Y}_{t}$.  It is clear that the fiber of $\mathcal{Y} \times_{B} B'$ over $t'$ is generically smooth, hence generically reduced.  Furthermore, the preimage of $U$ in $\mathcal{Y} \times_{B} B'$ will be contained in the smooth locus of $\mathcal{Y} \times_{B} B'$.  Thus the normalization map restricts to a birational morphism on this fiber.
\end{proof}

\begin{corollary}
\label{corollary:fiberwisebirational}
Let $\pi: \mathcal{Y} \to B$ and $\pi^{\sigma}: \mathcal{Y}^{\sigma} \to B$ be morphisms from normal projective varieties with connected fibers. Suppose that $f: \mathcal{Y} \to \mathcal{X}$ and $f^{\sigma}: \mathcal{Y}^{\sigma} \to \mathcal{X}$ are dominant finite $B$-morphisms whose generic fibers are twists of each other over $\mathcal X$.  Suppose further that $t \in B$ is a closed point such that the fibers $\mathcal{Y}_{t}, \mathcal{Y}^{\sigma}_{t}$ are generically reduced.  Then the maps $f_{t}, f^{\sigma}_{t}$ are birationally equivalent.
\end{corollary}

\begin{proof}
Choose a dominant finite morphism $B' \to B$ such that the normalized base changes $f': \mathcal{Y}' \to \mathcal{X}'$, $f'^{\sigma}:\mathcal{Y}'^{\sigma} \to \mathcal{X}'$ are birationally equivalent.  Since $f', f'^{\sigma}$ are finite, they are equal to their own Stein factorizations and thus $f'$ and $f'^\sigma$ are isomorphic to each other.
In particular the maps $f'_{t'}: \mathcal{Y}'_{t'} \to \mathcal{X}_{t}$ and $f'^{\sigma}_{t'}: \mathcal{Y}'^{\sigma}_{t'} \to \mathcal{X}_{t}$ are isomorphic to each other.  But by Lemma \ref{lemm:genreducedpreserved} these are birationally equivalent to $f_{t}$ and $f^{\sigma}_{t}$ respectively.
\end{proof}

We next analyze how the canonical divisor changes upon normalized base change.  This is well-understood, e.g.,~in the context of semistable reduction.

\begin{definition}
Suppose $h: \mathcal{Y}' \to \mathcal{Y}$ is a finite morphism of normal projective varieties.  Let $\mathcal{U} \subset \mathcal{Y}$ and $\mathcal{U}' \subset \mathcal{Y}'$ denote the smooth loci and set $\mathcal{V} = \mathcal{U}' \cap h^{-1}(\mathcal{U})$.  Note that the complement of $\mathcal{V}$ in $\mathcal{Y}'$ has codimension $\geq 2$.  The Riemann-Hurwitz formula gives us a distinguished effective representative $E$ in the linear equivalence class of $K_{\mathcal{V}/\mathcal{U}}$.  We define the relative canonical divisor $K_{\mathcal{Y}'/\mathcal{Y}}$ to be the effective Weil divisor obtained by taking the closure of $E$.
\end{definition}

\begin{lemma} \label{lemm:canonicalandnbc}
Let $\mathcal{Y}$ be a normal projective variety equipped with a surjective morphism $\pi: \mathcal{Y} \to B$ with connected fibers.  Suppose $g: B' \to B$ is a dominant finite morphism of curves and consider the normalized base change
\begin{equation*}
\xymatrix{
\mathcal{Y}' \ar[r]^-{\phi} \ar[dr]_{\pi'} & \mathcal{Y} \times_{B} B' \ar[r]^-{h} \ar[d] & \mathcal{Y} \ar[d]^{\pi} \\
& B' \ar[r]^{g} & B
}
\end{equation*}
We define $R$ to be the $\pi'$-pullback of the ramification divisor of $g$ and for any point $t' \in B'$ we let $R_{t'}$ denote the intersection of $R$ with the preimage of $t'$.
\begin{enumerate}
\item We have $K_{\mathcal{Y}'/\mathcal{Y}} \leq R$.
\item If $t \in B$ is a closed point such that the fiber $\mathcal{Y}_{t}$ is generically reduced then for any $t' \in B'$ mapping to $t$ we have $R_{t'} \leq K_{\mathcal{Y}'/\mathcal{Y}}$.
\end{enumerate}
\end{lemma}

\begin{proof}
(1) First note that the support of $K_{\mathcal{Y}'/\mathcal{Y}}$ is contained in the support of $R$.  Indeed, the map $h$ is \'etale away from  $\phi(\Supp(R))$ and thus $\phi$ is an isomorphism over this locus.

Suppose $t' \in B'$ is a ramification point for $g$ with index $e$. Let $t \in B$ be the image of $t'$. We choose local coordinates $s'$ and $s$ at $t'$ and $t$ respectively so that $g$ is defined locally by $s = s'^e$.
Let $T'$ be an irreducible component of the fiber $\pi'^{-1}(t')$ and let $q'$ denote its multiplicity in its component. Let $u'$ be a generic local equation of $T'$ so that generically the map $\mathcal Y' \to B'$ is given by by $s' = u'^{q'}$. Let $T \subset \mathcal U$ be the image of $T'$ and let $u$ denote a local equation of $T$ at a generic point. We denote the multiplicity of $T$ in $\pi^{-1}(t)$ by $q$ so that $\mathcal Y$ is generically defined by $s = u^q$. The coefficient of $T'$ in $R$ is given by $(e-1)q'$. Also $u^q = s = s'^e = u'^{eq'}$ so the coefficient of $T'$ in $K_{\mathcal Y'/\mathcal Y}$ is given by
\begin{equation} \label{eq:multiplicityincanonical}
e\frac{q'}{q} - 1.
\end{equation}
Since $e > 1$, we have $(e-1)q' \geq e\frac{q'}{q} - 1$ when $q > 1$. When $q = 1$, $q' = 1$ by Lemma \ref{lemm:genreducedpreserved} and so the inequality still holds. Thus our assertion follows.

(2) As explained above when $q = 1$ we have $q' = 1$ as well by Lemma \ref{lemm:genreducedpreserved}.  Thus we have our assertion.
\end{proof}

\begin{proposition}
\label{prop:ramificationdivisor}
Let $\mathcal{X} \to B$ be a good fibration.  
Let $\mathcal{Y},\mathcal{Y}^{\sigma}$ be normal projective varieties which admit surjective morphisms $\pi: \mathcal{Y} \to B$ and $\pi^{\sigma}: \mathcal{Y}^{\sigma} \to B$ with connected fibers.  Suppose there are dominant finite $B$-morphisms $f: \mathcal{Y} \to \mathcal{X}$ and $f^{\sigma}: \mathcal{Y}^{\sigma} \to \mathcal{X}$ whose generic fibers are twists of each other over $K(B)$.  Choose a finite Galois morphism $g: B' \to B$ such that the normalized base changes of $\mathcal{Y}$ and $\mathcal{Y}^{\sigma}$ over $B'$ are isomorphic.  
We let $\widehat{\mathcal{Y}}$ denote this abstract variety; it is equipped with finite morphisms $\rho_{1}: \widehat{\mathcal{Y}} \to \mathcal{Y}$ and $\rho_{2}: \widehat{\mathcal{Y}} \to \mathcal{Y}^{\sigma}$.
We denote the degree of $B' \to B$ by $d$.

There is a Weil divisor $E$ on $\widehat{\mathcal Y}$, which we write as $E = E^{+} - E^{-}$ where $E^{+},E^{-}$ are effective with no common divisor in their support, such that
\begin{itemize}
\item $K_{\widehat{\mathcal{Y}}/\mathcal{Y}} - K_{\widehat{\mathcal{Y}}/\mathcal{Y}^{\sigma}} \geq E$;
\item we have $E^{+} \geq \sum_{t'} \widehat{\mathcal{Y}}_{t'}$ as $t' \in B'$ varies over closed points whose image $t \in B$ satisfies that $\mathcal{Y}_{t}$ is normal but the fiber $\mathcal{Y}^{\sigma}_{t}$ is not generically reduced;
\item we have $E^{-} \leq d \cdot \sum_{t'} \widehat{\mathcal{Y}}_{t'}$ as $t' \in B'$ varies over closed points whose image $t \in B$ satisfies that $\mathcal{Y}_{t}$ is not normal.  
\end{itemize}
\end{proposition}

\begin{proof}
We compare these divisors along each fiber separately. Let $t' \in B'$ be a closed point and let $t \in B$ be its image.  If the fiber $\mathcal{Y}_{t}$ is not normal, then it follows from Lemma~\ref{lemm:canonicalandnbc} that 
\[
(K_{\widehat{\mathcal{Y}}/\mathcal{Y}})_{t'} - (K_{\widehat{\mathcal{Y}}/\mathcal{Y}^{\sigma}})_{t'} \geq -(e-1)\widehat{\mathcal Y}_{t'},
\]
where $e$ is the ramification index of $t'$.  In particular this difference is $\geq -d\widehat{\mathcal Y}_{t'}$.
If $\mathcal Y_t$ is normal, then in particular $\mathcal Y_t$ is irreducible and reduced.  Thus the fiber $\widehat{\mathcal{Y}}_{t'}$ is also irreducible and reduced, and we conclude that $\mathcal Y^\sigma_t$ is irreducible.  Let us denote the multiplicity of $\mathcal Y^\sigma_t$ by $q$. 
It follows from Lemma~\ref{lemm:canonicalandnbc} that 
\[
(K_{\widehat{\mathcal{Y}}/\mathcal{Y}})_{t'} - (K_{\widehat{\mathcal{Y}}/\mathcal{Y}^{\sigma}})_{t'} \geq 0
\]
and equality holds when $q=1$.  The two inequalities above together prove the upper bound on $E^{-}$. To prove the lower bound on $E^{+}$, we analyze those fibers such that $\mathcal{Y}_{t}$ is normal and which satisfy $q >1$. Since $\widehat{\mathcal{Y}}_{t'}$ is reduced we must have $e > 1$.  Since the multiplicities of $\mathcal{Y}_{t}$ and $\widehat{\mathcal{Y}}_{t'}$ are $1$, Equation \eqref{eq:multiplicityincanonical} in Lemma~\ref{lemm:canonicalandnbc} shows that
\[
(K_{\widehat{\mathcal{Y}}/\mathcal{Y}})_{t'} = (e-1)\widehat{\mathcal Y}_{t'}, \quad  (K_{\widehat{\mathcal{Y}}/\mathcal{Y}^{\sigma}})_{t'} = \left(\tfrac{e}{q}-1 \right)\widehat{\mathcal Y}_{t'}
\]
so that we conclude
\[
(K_{\widehat{\mathcal{Y}}/\mathcal{Y}})_{t'} - (K_{\widehat{\mathcal{Y}}/\mathcal{Y}^{\sigma}})_{t'} = e\left(1-\tfrac{1}{q} \right)\widehat{\mathcal Y}_{t'} \geq \widehat{\mathcal Y}_{t'}.
\]
Thus our assertion follows.
\end{proof}

\begin{proposition} \label{prop:curveintandrambound}
Let $\mathcal{X} \to B$ be a good fibration.  
Let $\mathcal{Y},\mathcal{Y}^{\sigma}$ be normal projective varieties which admit surjective morphisms $\pi: \mathcal{Y} \to B$ and $\pi^{\sigma}: \mathcal{Y}^{\sigma} \to B$ with connected fibers.  Suppose there are dominant finite $B$-morphisms $f: \mathcal{Y} \to \mathcal{X}$ and $f^{\sigma}: \mathcal{Y}^{\sigma} \to \mathcal{X}$ whose generic fibers are twists of each other over $K(B)$.  Choose a finite Galois morphism $g: B' \to B$ such that the normalized base changes of $\mathcal{Y}$ and $\mathcal{Y}^{\sigma}$ over $B'$ are isomorphic.  We let $\widehat{\mathcal{Y}}$ denote this abstract variety; it is equipped with finite morphisms $\rho_{1}: \widehat{\mathcal{Y}} \to \mathcal{Y}$ and $\rho_{2}: \widehat{\mathcal{Y}} \to \mathcal{Y}^{\sigma}$.  We denote the degree of $B' \to B$ by $d$.

Let ${\mathcal Y^\sigma}'$ be a smooth birational model of $\mathcal Y^\sigma$ equipped with a birational morphism $\beta: {\mathcal Y^\sigma}' \to \mathcal Y^\sigma$.
Assume that there exists a section $C$ on ${\mathcal Y^\sigma}'$ and a constant $R > 0$ such that $C$ corresponds to a rational point on the smooth locus of ${\mathcal Y^\sigma_\eta}$ and 
\[
(K_{{\mathcal{Y}^\sigma}'/B} - \beta^*(f^\sigma)^{*}K_{\mathcal X/B}) \cdot C \leq R.
\]
Let $r$ denote the number of closed points $t \in B$ such that the fiber $\mathcal{Y}_{t}$ is normal but the fiber $\mathcal{Y}^{\sigma}_{t}$ is not generically reduced. Then we have $r \leq dR$.
\end{proposition}

\begin{proof}
We choose smooth models $\mathcal Y',  \widehat{\mathcal Y}'$ of $\mathcal Y, \widehat{\mathcal Y}$ respectively such that there are birational morphisms $\alpha : \mathcal Y' \to \mathcal Y$,  $\gamma : \widehat{\mathcal Y}' \to \widehat{\mathcal Y}$ and generically finite morphisms $\widetilde{\rho}_1 : \widehat{\mathcal Y}' \to \mathcal Y'$, $\widetilde{\rho}_2 : \widehat{\mathcal Y}' \to {\mathcal Y^\sigma}'$ which are birationally equivalent to $\rho_1, \rho_2$ respectively.  We may ensure that $\gamma^{-1}$ is well-defined along the smooth locus of $\mathcal{Y}^{\sigma}_{\eta} \otimes K(B')$.  Our intersection bound implies that
\[
K_{{\mathcal Y^\sigma}'/\mathcal X} \cdot C \leq R.
\]
Since by assumption $C$ is not contained in the $\rho_{2}$-image of the $\gamma$-exceptional centers, there is a section $C'$  of $\widehat{\mathcal Y}'/B'$ such that $(\widetilde{\rho}_2)_* C' = dC$.
Then we have
\[
\widetilde{\rho}_2^*K_{{\mathcal Y^\sigma}'/\mathcal X} \cdot C' \leq dR.
\]
Note that
\begin{align*}
\widetilde{\rho}_2^*K_{{\mathcal Y^\sigma}'/\mathcal X} &= K_{\widehat{\mathcal Y}'/\mathcal X} - K_{\widehat{\mathcal Y}'/{\mathcal Y^\sigma}'}\\
& = \widetilde{\rho}_1^*K_{\mathcal Y'/\mathcal X} + K_{\widehat{\mathcal Y}'/\mathcal Y'} - K_{\widehat{\mathcal Y}'/{\mathcal Y^\sigma}'} \geq K_{\widehat{\mathcal Y}'/\mathcal Y'} - K_{\widehat{\mathcal Y}'/{\mathcal Y^\sigma}'}.
\end{align*}
Let $E$ be the divisor on $\widehat{\mathcal{Y}}$ defined by Proposition~\ref{prop:ramificationdivisor}, so $E^{+}$ is at least as effective as the sum of the $r$ fibers of $\widehat{\mathcal Y}$ corresponding to the $r$ closed points $t \in B$ such that the fiber $\mathcal{Y}_{t}$ is normal but the fiber $\mathcal{Y}^{\sigma}_{t}$ is not generically reduced.  Taking strict transforms, we see that $K_{\widehat{\mathcal Y}'/\mathcal Y'} - K_{\widehat{\mathcal Y}'/{\mathcal Y^\sigma}'}$ is at least as effective as the sum of the strict transforms of these $r$ fibers of $\widehat{\mathcal Y}$.
Furthermore every exceptional divisor of $\beta: {\mathcal Y^\sigma}' \to \mathcal Y^\sigma$ is contracted by $f^\sigma\circ \beta : {\mathcal Y^\sigma}' \to \mathcal X$ as well, and thus appears with positive coefficient in the ramification divisor $K_{{\mathcal Y^\sigma}'/\mathcal X}$.  We conclude that the support of the effective divisor $\widetilde{\rho}_2^*K_{{\mathcal Y^\sigma}'/\mathcal X}$ contains the $r$ reduced fibers over these points.  Since our section $C'$ must meet each fiber in a component of multiplicity one, we conclude that
\[
r \leq \widetilde{\rho}_2^*K_{{\mathcal Y^\sigma}'/\mathcal X} \cdot C' \leq dR.
\]
\end{proof}

\begin{corollary} \label{coro:boundedintimpliesboundedtwists}
Let $\mathcal{X} \to B$ be a good fibration.  
Let $\mathcal{Y},\mathcal{Y}^{\sigma}$ be normal projective varieties which admit surjective morphisms $\pi: \mathcal{Y} \to B$ and $\pi^{\sigma}: \mathcal{Y}^{\sigma} \to B$ with connected fibers.  Suppose there are dominant finite $B$-morphisms $f: \mathcal{Y} \to \mathcal{X}$ and $f^{\sigma}: \mathcal{Y}^{\sigma} \to \mathcal{X}$ whose generic fibers are twists of each other over $K(B)$.  

Let $\widetilde{\mathcal Y}^\sigma$ be a smooth birational model of $\mathcal Y^\sigma$ equipped with a birational morphism $\beta: \widetilde{\mathcal Y}^\sigma \to \mathcal Y^\sigma$.
Assume that there exists a section $C$ on $\widetilde{\mathcal Y}^\sigma$ and a constant $R > 0$ such that $C$ corresponds to a rational point on the smooth locus of ${\mathcal Y^\sigma_\eta}$ and 
\[
(K_{\widetilde{\mathcal{Y}}^\sigma/B} - \beta^*(f^\sigma)^{*} K_{\mathcal X/B}) \cdot C \leq R
\]
Then there exist constants $d = d(\mathcal Y/\mathcal X)$ and $n = n(\mathcal Y/\mathcal X, R)$ such that there exists a finite Galois morphism $B'\to B$ of degree at most $d$ with at most $n$ branch points such that the normalized base changes of $\mathcal Y/B$ and $\mathcal Y^\sigma/B$ by $B' \to B$ become $\mathcal X \times_B B'$-isomorphic.

In particular, the set of such twists is a bounded family.
\end{corollary}

\begin{proof}
It follows from Lemma~\ref{lemm:boundingdeganddisc} that there exist $d = d(\mathcal Y/\mathcal X)$ and a finite Galois morphism $\widetilde{B} \to B$ of degree at most $d$ such that the normalizations of $\mathcal Y \times_B \widetilde{B}$ and $\mathcal Y^\sigma \times_B \widetilde{B}$ are isomorphic.

Let $s = s(\mathcal Y/\mathcal X)$ be the number of $t \in B$ such that $\mathcal Y_t$ is not normal or $\mathrm{Aut}_B(\mathcal Y/\mathcal X) \to B$ is not flat at $t \in B$.  Let $r$ be the number of $t \in B$ such that $\mathcal Y_t$ is normal but $\mathcal Y^\sigma_t$ is not generically reduced.  By Proposition \ref{prop:curveintandrambound} we have $r \leq dR$.

For $t \in B$ such that $\mathcal Y_t$ is normal, $\mathrm{Aut}_B(\mathcal Y/\mathcal X) \to B$ is flat at $t \in B$ and $\mathcal Y^\sigma_t$ is generically reduced, it follows from Corollary~\ref{corollary:fiberwisebirational} and Lemma~\ref{lemma:birationalHensel} that $\mathrm{inv}_t(\sigma) = 0$. Thus the number of $t \in B$ such that $\mathrm{inv}_t(\sigma)\neq 0$ is bounded above by $s + dR$. Thus our first assertion follows from Lemma~\ref{lemm:boundingdeganddisc}.

The final statement then follows from Lemma \ref{lemm:twistboundedconditions}.
\end{proof}

\section{Fujita invariant and sections} \label{sec:fujinv}

Suppose that $\pi: \mathcal{X} \to B$ is a good fibration and $L$ is a generically relatively big and semiample Cartier divisor on $\mathcal{X}$.  In this section the goal is to classify the generically finite $B$-morphisms $f: \mathcal{Y} \to \mathcal{X}$ such that $\mathcal{Y}$ carries a family of sections $N$ with the property that $f_{*}N$ has small codimension in an irreducible component of $\Sec(\mathcal{X}/B)$.  Assuming the sections have large $L$-degree but small degree against $f^{*}(K_{\mathcal{X}/B} + a(\mathcal{X}_{\eta},L|_{\mathcal{X}_{\eta}})L)$, we show that the Fujita invariant of $\mathcal{Y}_{\eta}$ must be at least as large as the Fujita invariant of $\mathcal{X}_{\eta}$.  This puts a strong constraint on the set of morphisms $f$ which have this property.

After addressing some preliminaries in Section \ref{sect:pivertical} and Section \ref{sect:relvsabs}, we show the fundamental result discussed above in Section \ref{sect:fujinvongenfiber}.  When working with the Fujita invariant it is often helpful to know that the pair $(\mathcal{Y}_{\eta},f^{*}L|_{\mathcal{Y}_{\eta}})$ is adjoint rigid; in Section \ref{sect:adjrigid} we show that if we additionally assume that the sections parametrized by $N$ go through many general points of $\mathcal{Y}$ then we can also guarantee adjoint rigidity.

\subsection{Modifying by $\pi$-vertical divisors} \label{sect:pivertical}

Let $\pi: \mathcal{X} \to B$ be a good fibration and let $L$ be a generically relatively big and semiample Cartier divisor on $\mathcal{X}$.  We know that $L|_{\mathcal{X}_{\eta}}$ is $\mathbb{Q}$-linearly equivalent to a divisor which has smooth support.  The following proposition discusses how to reframe this property as a global statement by adding $\pi$-vertical divisors.

\begin{proposition} \label{prop:generalsurjectiontofiber}
Let $\pi: \mathcal{X} \to B$ be a good fibration, let $L$ be a generically relatively big and semiample Cartier divisor on $\mathcal{X}$, and let $a$ be a positive rational number.  Let $b > a$ be a positive integer such that $bL|_{\mathcal{X}_{\eta}}$ defines a basepoint free linear series.  There is some effective $\pi$-vertical $\mathbb{Q}$-Cartier divisor $E$ on $\mathcal{X}$ such that the following property holds.

Suppose $\psi: \mathcal{Y} \to B$ is a good fibration and $f: \mathcal{Y} \to \mathcal{X}$ is a $B$-morphism that is generically finite onto its image.  Then there is an effective $\mathbb{Q}$-Cartier divisor $D$ on $\mathcal{Y}$ that is $\mathbb{Q}$-linearly equivalent to $f^{*}(aL+E)$ such that $D|_{\mathcal{Y}_{\eta}}$ has smooth irreducible support and coefficient $\frac{a}{b}$.  In particular $(\mathcal{Y}_{\eta},D|_{\mathcal{Y}_{\eta}})$ is a terminal pair.
\end{proposition}

\begin{proof}
Let $T_{1},\ldots,T_{r}$ be a $K(B)$-basis for $|bL|_{\mathcal{X}_{\eta}}|$.  Note that $\cap_{i} T_{i} = \emptyset$.

We denote by $\overline{T}_{i}$ the closure of $T_{i}$ in $\mathcal{X}$.  There is some effective $\pi$-vertical $\mathbb{Q}$-Cartier divisor $\widehat{E}$ such that for every $i$ there is an effective $\pi$-vertical divisor $F_{i}$ satisfying $\overline{T}_{i} + F_{i} \sim b(L + \widehat{E})$.

Let $f: \mathcal{Y} \to \mathcal{X}$ be a morphism as in the statement.  By construction we have $\cap_{i} f^{*}(\overline{T}_{i} + F_{i})$ does not intersect $\mathcal{Y}_{\eta}$.  Thus $f^{*}(b(L+\widehat{E}))$ is linearly equivalent to a divisor $\widehat{D}$ whose restriction to $\mathcal{Y}_{\eta}$ is smooth and irreducible.  Then $D = \frac{a}{b}\widehat{D}$ and $E = a\widehat{E}$ have the desired properties.
\end{proof}

We will use the following definition to capture the effect of the extra divisor $E$.
 
\begin{definition}
\label{defi: invariant_tau}
Let $\pi: \mathcal{X} \to B$ be a good fibration.  Suppose that $E$ is an effective $\pi$-vertical $\mathbb{Q}$-Cartier divisor.  Define
\begin{equation*}
\tau(\pi,E) = \sup_{\textrm{sections }C} E \cdot C.
\end{equation*}
Note that this supremum is achieved by some section $C$ since the intersection number is bounded above by the sum of the coefficients of $E$ and is contained in $\frac{1}{r}\mathbb{Z}$ where $r$ is the least common multiple of the denominators of the coefficients of $E$.
\end{definition}

\subsection{Relative versus absolute positivity} \label{sect:relvsabs}
We will also need a couple results comparing relative and absolute positivity for a fibration over a curve.

\begin{lemma} \label{lemm:relvsabsmmp}
Let $\pi: \mathcal{Z} \to B$ be a good fibration.  Suppose that $D$ is an effective $\mathbb{Q}$-Cartier divisor on $\mathcal{Z}$ such that $(\mathcal{Z},D)$ is a terminal pair.

\begin{enumerate}
\item Suppose that $g(B) \geq 1$.   Suppose that $\rho: \mathcal{Z} \dashrightarrow \tilde{\mathcal{Z}}$ is a rational map obtained by running the $(K_{\mathcal{Z}}+D)$-MMP.  Then $\rho$ is also a run of the relative $(K_{\mathcal{Z}}+D)$-MMP over $B$.
\item Suppose that $g(B) = 0$.  There is a constant $m = m(\dim(\mathcal{Z}))$ such that the following holds.  Fix a general fiber $F$ of $\pi$ and suppose that $\rho: \mathcal{Z} \dashrightarrow \tilde{\mathcal{Z}}$ is a birational morphism obtained by running the $(K_{\mathcal{Z}}+D+m F)$-MMP.  Then $\rho$ is also a run of the relative $(K_{\mathcal{Z}}+D+mF)$-MMP over $B$.
\end{enumerate}

In particular, in case (1) (resp.~case (2)) if $K_{\mathcal{Z}}+D$ (resp.~$K_{\mathcal{Z}}+D+mF$) is not pseudo-effective then its restriction to $\mathcal{Z}_{\eta}$ is also not pseudo-effective.
\end{lemma}

\begin{proof}
(1) By \cite{Kawamata91} each step of the $(K_{\mathcal{Z}} + D)$-MMP contracts an extremal ray that is spanned by a rational curve.  This rational curve must be vertical with respect to $\pi$ because $B$ has genus $\geq 1$.

(2) Since $(\mathcal{Z},D)$ is $\frac{1}{2}$-lc, by Theorem \ref{theo:nefconetheorem} there is an integer $m = m(\dim(\mathcal{Z}))$ such that
\begin{equation*}
\Nef_{1}(\mathcal{Z}) + \Eff_{1}(\mathcal{Z})_{K_{\mathcal{Z}} + D + mF \geq 0} = \Eff_{1}(\mathcal{Z})_{K_{\mathcal{Z}} + D + mF \geq 0} + \sum_{j} [C_{j}]
\end{equation*}
where the $C_{j}$ are $\pi$-vertical moving curves.  In particular, any contraction of a $(K_{\mathcal{Z}} + D + mF)$-negative extremal ray must define a relative contraction over $B$.  Furthermore the analogous equality of cones holds for any birational model of $\mathcal{Z}$ obtained by running the MMP (since the $\frac{1}{2}$-lc condition is preserved).  Thus we see that every step of the $(K_{\mathcal{Z}}+D+mF)$-MMP is actually a step of the relative MMP over $B$.

To see the final statement, suppose we are in case (1) and $K_{\mathcal{Z}}+D$ is not pseudo-effective.  Then we can run the $(K_{\mathcal{Z}} + D)$-MMP with scaling of an ample divisor and the outcome will be a Mori fibration.  But then this Mori fibration must be a relative fibration over $B$, so that $(K_{\mathcal{Z}}+D)$ is not relatively pseudo-effective over $B$.  The same argument applies in case (2).
\end{proof}

Our next result shows how to turn intersection inequalities into Fujita invariant inequalities.

\begin{lemma} \label{lemm:terminalsectiontofiber}
Let $\pi: \mathcal{Z} \to B$ be a good fibration.  Suppose that $D$ is an effective $\mathbb{Q}$-Cartier divisor on $\mathcal{Z}$ such that $(\mathcal{Z}_{\eta},D|_{\mathcal{Z}_{\eta}})$ is a terminal pair and $D|_{\mathcal{Z}_{\eta}}$ is big and nef.

\begin{enumerate}
\item Suppose that $g(B) \geq 1$.   If there is a dominant family of HN-free sections $C$ on $\mathcal{Z}$ which satisfy $-(K_{\mathcal{Z}} + D) \cdot C > 0$ then $a(\mathcal{Z}_{\eta},D|_{\mathcal{Z}_{\eta}}) > 1$.
\item Suppose that $g(B) = 0$.  There is a constant $\Xi = \Xi(\dim(\mathcal{Z}))$ such that the following holds.  If there is a dominant family of HN-free sections $C$ on $\mathcal{Z}$ which satisfy $-(K_{\mathcal{Z}} + D) \cdot C > \Xi$ then $a(\mathcal{Z}_{\eta},D|_{\mathcal{Z}_{\eta}}) > 1$.
\end{enumerate}
\end{lemma}

\begin{proof}
Let $\phi: \mathcal{Z}' \to \mathcal{Z}$ be a log resolution of $(\mathcal{Z},D)$ and let $D'$ be the strict transform of the $\pi$-horizontal components of $D$.  After perhaps taking a further blow-up, we may assume that two irreducible components of $D'$ intersect if and only if their restrictions to the generic fiber intersect.  Along the central fiber we can write $K_{\mathcal{Z}'_{\eta}} + D'|_{\mathcal{Z}'_{\eta}} = \phi^{*}(K_{\mathcal{Z}_{\eta}} + D|_{\mathcal{Z}_{\eta}}) + E_{\eta}$ where $E_{\eta}$ is an effective $\phi$-exceptional divisor.  We conclude that $K_{\mathcal{Z}'_{\eta}} + D'|_{\mathcal{Z}'_{\eta}}$ is pseudo-effective if and only if $K_{\mathcal{Z}_{\eta}} + D|_{\mathcal{Z}_{\eta}}$ is pseudo-effective.

We claim that the pair $(\mathcal{Z}',D')$ has terminal singularities.  Since $\Supp(D')$ is an SNC divisor, by \cite[3.11 Lemma]{Kollar97} the pair $(\mathcal{Z}',D')$ will be terminal if and only if when we write $D' = \sum_{i} d_{i}D'_{i}$ in terms of irreducible components we have
\begin{equation*}
\min_{i} \{ 1 - d_{i} \} > 0 \qquad \textrm{and} \qquad \min_{i,j | D_{i} \cap D_{j} \neq \emptyset} \{ 1-d_{i} - d_{j} \} > 0.
\end{equation*}
Recall that by construction two irreducible components of $D'$ intersect if and only if their restrictions to the generic fiber intersect.  Thus this computation can be done on the generic fiber, where the desired inequalities follow from the fact that $(\mathcal{Z}_\eta,D|_{\mathcal Z_\eta})$ is terminal. 

Let $C'$ be the strict transform of a general deformation of $C$.  Since $C$ is HN-free, by Lemma \ref{lemm:HNavoidscodim2} we can assume that $C'$ avoids any codimension $2$ locus in $\mathcal{Z}$ and thus $C'$ has vanishing intersection against every $\phi$-exceptional divisor.  We have
\begin{align*}
(K_{\mathcal{Z'}} + D') \cdot C' & \leq (K_{\mathcal{Z'}} + \phi^{*}D) \cdot C' \\
& = \phi^{*}(K_{\mathcal{Z}} + D) \cdot C'
\end{align*}

(1) We are in the case $g(B) \geq 1$ and
\begin{equation*}
(K_{\mathcal{Z'}} + D') \cdot C' \leq \phi^{*}(K_{\mathcal{Z}} + D) \cdot C' < 0
\end{equation*}
Since $C'$ is a movable curve, we see that $K_{\mathcal{Z'}} + D'$ is not pseudo-effective.  By Lemma \ref{lemm:relvsabsmmp} we see that $(K_{\mathcal{Z}'} + D')|_{\mathcal{Z}'_{\eta}}$ is also not pseudo-effective.  As demonstrated above this means that $(K_{\mathcal{Z}} + D)|_{\mathcal{Z}_{\eta}}$ also fails to be pseudo-effective, showing that $a(\mathcal{Z}_{\eta},D|_{\mathcal{Z}_{\eta}}) > 1$.

(2) We are in the case $g(B) = 0$.  Let $m = m(\dim(\mathcal{Z}))$ be the constant from Lemma \ref{lemm:relvsabsmmp}.(2) and set $\Xi = m+1$.  We have
\begin{equation*}
(K_{\mathcal{Z'}} + D') \cdot C' \leq \phi^{*}(K_{\mathcal{Z}} + D) \cdot C' < -\Xi
\end{equation*}
and thus $(K_{\mathcal{Z}'}+D'+mF) \cdot C' < 0$.  Since $C'$ is a movable curve, we see that that $K_{\mathcal{Z}'} + D' + mF$ is not pseudo-effective.  By Lemma \ref{lemm:relvsabsmmp} we see that $(K_{\mathcal{Z}'} + D')|_{\mathcal{Z}'_{\eta}}$ is also not pseudo-effective.  As demonstrated above this means that $(K_{\mathcal{Z}} + D)|_{\mathcal{Z}_{\eta}}$ also fails to be pseudo-effective, showing that $a(\mathcal{Z}_{\eta},D|_{\mathcal{Z}_{\eta}}) > 1$.
\end{proof}

\subsection{Fujita invariant along the generic fiber} \label{sect:fujinvongenfiber}
In this section we show that the Fujita invariant along the generic fiber controls the expected dimension for families of sections.
We will use the following easy lemma many times.

\begin{lemma} \label{lemm:easyintcalc}
Let $\pi: \mathcal{X} \to B$ be a good fibration and let $L$ be a generically relatively big and semiample Cartier divisor.  Assume that $\mathcal{X}_{\eta}$ is geometrically uniruled.  Fix a positive rational number $a_{rel}$ and define $a = a_{rel} a(\mathcal{X}_{\eta},L|_{\mathcal{X}_{\eta}})$.  Fix a rational number $\beta$.  Fix a positive integer $T$.

Suppose that $\psi: \mathcal{Y} \to B$ is a good fibration equipped with a $B$-morphism $f: \mathcal{Y} \to \mathcal{X}$ that is generically finite onto its image.   Suppose that $N$ is an irreducible component of $\Sec(\mathcal{Y}/B)$ parametrizing a dominant family of sections $C$ on $\mathcal{Y}$ which satisfy $f^{*}(K_{\mathcal{X}/B} + a(\mathcal{X}_{\eta},L|_{\mathcal{X}_{\eta}})L) \cdot C \leq \beta$.  Finally, suppose that
\begin{equation} \label{eq:dimension}
\dim(N) \geq a_{rel}(-f^{*}K_{\mathcal{X}/B} \cdot C + (\dim(\mathcal{X})-1)(1-g(B))) - T.
\end{equation}
Then
\begin{equation*}
(K_{\mathcal{Y}} + af^{*}L) \cdot C \leq a_{rel} \beta + T + a_{rel} (\dim(\mathcal{X}) - 1) (g(B)-1) +(\dim(\mathcal{X})-1) + 2g(B)-2.
\end{equation*}
\end{lemma}

\begin{proof}
Let $C$ denote a general section parametrized by $N$.  By Corollary \ref{coro:domfamilyexpdim}
\begin{equation*}
\dim(N) \leq -K_{\mathcal{Y}/B} \cdot C + (\dim(\mathcal{Y})-1).
\end{equation*}
Combining this equality with Equation \eqref{eq:dimension} and rearranging we get
\begin{align}
(K_{\mathcal{Y}/B}- a_{rel}f^*K_{\mathcal{X}/B}) \cdot C & \leq T + a_{rel} (\dim(\mathcal{X}) - 1) (g(B)-1) + (\dim(\mathcal{Y})-1) \nonumber \\
& \leq T + a_{rel} (\dim(\mathcal{X}) - 1) (g(B)-1) + (\dim(\mathcal{X})-1) \label{eq:stupidinequality}
\end{align}
Adding in the fact that $f^{*}(K_{\mathcal{X}/B} + a(\mathcal{X}_{\eta},L|_{\mathcal{X}_{\eta}})L) \cdot C \leq \beta$, we see that
\begin{equation} \label{eq:alternateform}
(K_{\mathcal{Y}/B} + af^*L) \cdot C \leq a_{rel} \beta + T + a_{rel} (\dim(\mathcal{X}) - 1) (g(B)-1) + (\dim(\mathcal{X})-1).
\end{equation}
or equivalently
\begin{equation*}
(K_{\mathcal{Y}} + af^*L) \cdot C \leq a_{rel} \beta + T + a_{rel} (\dim(\mathcal{X}) - 1) (g(B)-1) + (\dim(\mathcal{X})-1) + 2g(B)-2.
\end{equation*}
\end{proof}

We can now prove our basic result for controlling the Fujita invariant using pathological families of sections.

\begin{theorem} \label{theo:generalainvsections}
Let $\pi: \mathcal{X} \to B$ be a good fibration and let $L$ be a generically relatively big and semiample Cartier divisor on $\mathcal{X}$.  Assume that $\mathcal{X}_{\eta}$ is geometrically uniruled.  Fix a positive rational number $a_{rel}$ and set $a = a_{rel} a(\mathcal{X}_{\eta},L|_{\mathcal{X}_{\eta}})$.  Fix a rational number $\beta$.  Fix a positive integer $T$.  Fix a positive integer $b>a$ such that $bL|_{\mathcal{X}_{\eta}}$ defines a basepoint free linear series.  Use $b$ to construct an effective $\pi$-vertical $\mathbb{Q}$-Cartier divisor $E$ satisfying the conclusion of Proposition \ref{prop:generalsurjectiontofiber} with respect to $aL$.  
There is some constant $\xi = \xi(\dim(\mathcal{X}), g(B), \tau(\pi,E), a_{rel}, a, T, \beta,b)$ with the following property.

Suppose that $\psi: \mathcal{Y} \to B$ is a good fibration equipped with a $B$-morphism $f: \mathcal{Y} \to \mathcal{X}$ that is generically finite onto its image.  Suppose that $N$ is an irreducible component of $\Sec(\mathcal{Y}/B)$ parametrizing a dominant family of sections $C$ on $\mathcal{Y}$ which satisfy $f^{*}L \cdot C \geq \xi$ and $f^{*}(K_{\mathcal{X}/B} + a(\mathcal{X}_{\eta},L|_{\mathcal{X}_{\eta}}) L) \cdot C \leq \beta$.  Finally, suppose that
\begin{equation} \label{eq:eh}
\dim(N) \geq a_{rel} (-f^{*}K_{\mathcal{X}/B} \cdot C + (\dim(\mathcal{X})-1)(1-g(B))) - T.
\end{equation}
Then
\begin{equation*}
a(\mathcal{Y}_{\eta},f^{*}L|_{\mathcal{Y}_{\eta}}) \geq a.
\end{equation*}
\end{theorem}

\begin{proof}
We first prove the statement when the general section $C$ parametrized by $N$ is HN-free on $\mathcal{Y}$.
%(
By Theorem \ref{theo:Dicerbo} there is a rational number $\epsilon > 0$ depending only on $a$ and $\dim(\mathcal{X})$ such that no smooth variety of dimension $\leq \dim \mathcal{X}-1$ has Fujita invariant in the range $[(1-\epsilon)a,a)$ with respect to any big and nef Cartier divisor. Define $\Xi$ as:
%]
\begin{itemize}
\item $\Xi = 0$, if $g(B) \geq 1$.
\item $\Xi$ is the supremum of the constants obtained by applying Lemma \ref{lemm:terminalsectiontofiber} to all dimensions $\leq \dim(\mathcal{X})$, if $g(B) = 0$.
\end{itemize}
We define $\xi_{HN} (\dim(\mathcal{X}), g(B), \tau(\pi,E), a_{rel}, a, T, \beta,b)$ to be
\begin{align*}
\frac{1}{a \epsilon} \left((1-\epsilon) \tau(\pi,E) + a_{rel} \beta + T + a_{rel} (\dim(\mathcal{X}) - 1) (g(B)-1) + (\dim(\mathcal{X})-1) + 2g(B)-2 + \Xi \right) + 1
\end{align*}
and assume that our sections $C$ satisfy $f^{*}L \cdot C \geq \xi_{HN}$.

Let $C$ denote a general section parametrized by $N$.  Since we are assuming $C$ moves in a dominant family on $\mathcal{Y}$ Lemma \ref{lemm:easyintcalc} shows that 
\begin{equation*}
(K_{\mathcal{Y}} + af^*L) \cdot C \leq  a_{rel} \beta + T + a_{rel} (\dim(\mathcal{X}) - 1) (g(B)-1) + (\dim(\mathcal{X})-1) + 2g(B)-2.
\end{equation*}
%Adding in $E$, we get
%\begin{equation*}
%(K_{\mathcal{Y}} + af^*L + f^{*}E) \cdot C \leq a_{rel} \beta + T + a_{rel} (\dim(\mathcal{X}) - 1) (g(B)-1) + (\dim(\mathcal{X})-1) + 2g(B)-2 + \tau(\pi,E).
%\end{equation*}
When the section $C$ has degree $\geq \xi_{HN}$ then $a\epsilon f^{*}L \cdot C \geq a\epsilon \xi_{HN}$ and so this inequality simplifies to
\begin{equation*}
(K_{\mathcal{Y}} + a(1-\epsilon)f^*L) \cdot C  + (1-\epsilon) \tau(\pi,E) + \Xi + a\epsilon \leq 0.
\end{equation*}
Using the fact that $a\epsilon > 0$ and that $(1-\epsilon) \tau(\pi,E) \geq (1-\epsilon) f^{*}E \cdot C$ we get
\begin{equation*}
(K_{\mathcal{Y}} + a(1-\epsilon) f^{*}L + (1-\epsilon) f^{*}E) \cdot C < -\Xi.
\end{equation*}
Since $E$ satisfies the conclusion of Proposition \ref{prop:generalsurjectiontofiber} with respect to $aL$ the pullback $f^{*}(aL + E)$ is $\mathbb{Q}$-linearly equivalent to an effective $\mathbb{Q}$-divisor $D$ such that $(\mathcal{Y}_{\eta},D|_{\mathcal{Y}_{\eta}})$ has terminal singularities.  Of course $(\mathcal{Y}_{\eta},(1-\epsilon)D|_{\mathcal{Y}_{\eta}})$ also has terminal singularities.  Applying Lemma \ref{lemm:terminalsectiontofiber}, we deduce that $a(\mathcal{Y}_{\eta},f^{*}L|_{\mathcal{Y}_{\eta}}) \geq (1-\epsilon)a$.  By construction this implies that $a(\mathcal{Y}_{\eta},f^{*}L|_{\mathcal{Y}_{\eta}}) \geq a$.

Next we prove the statement when $C$ is not HN-free on $\mathcal{Y}$.  Our strategy is to reduce to the HN-free case.  Define 
\begin{equation}
\label{equation:generalainvsections}
\xi
= \sup \left\{  \begin{array}{c} \xi_{HN}(\dim(\mathcal{X}),g(B),\tau(\pi,E),a_{rel}, a, T+(\dim(\mathcal{X})-1)(4g(B) + 3 + \gamma),\beta,b), \\ \frac{1}{a} ((\dim(\mathcal{X})-1)(5g(B)+ 3 + \gamma) + a_{rel}\beta + T + a_{rel} (\dim(\mathcal{X}) - 1) (g(B)-1))  \end{array} \right\}
\end{equation}
where $\gamma = (g(B)\dim(\mathcal{X}) - g(B) + 1)^{2}(\dim(\mathcal{X})-1)$.
Using the second term as a lower bound on $\xi$ and appealing to Equation \eqref{eq:alternateform} in Lemma \ref{lemm:easyintcalc}, we have 
\begin{align*}
-K_{\mathcal{Y}/B} \cdot C & \geq af^{*}L \cdot C - a_{rel}\beta - T - a_{rel} (\dim(\mathcal{X}) - 1) (g(B)-1) - (\dim(\mathcal{X})-1) \\
& \geq (\dim(\mathcal{X})-1)(5g(B)+2 + \gamma) \\
& \geq (\dim(\mathcal{Y})-1)(5g(B)+2 + \gamma_{\mathcal{Y}})
\end{align*}
where $\gamma_{\mathcal{Y}}= (g(B)\dim(\mathcal{Y}) - g(B) + 1)^{2}(\dim(\mathcal{Y})-1)$.  Let $\mathcal{S}'$ be a smooth birational model of the finite part of the Stein factorization of the evaluation map for the normalization of the universal family over $N$ equipped with the map $\pi': \mathcal{S}' \to B$.  Since the dimension of $N$ is the same as the dimension of the corresponding family of sections $C'$ on $\mathcal{S}'$, Corollary \ref{coro:domfamilyexpdim} shows that
\begin{equation*}
-K_{\mathcal{S}'/B} \cdot C' \geq -K_{\mathcal{Y}/B} \cdot C - g(B) (\dim(\mathcal{S}')-1).
\end{equation*}
In particular we see that $\mu^{max}_{[C']}(T_{\mathcal{S}'/B}) \geq \mu_{[C']}(T_{\mathcal{S}'/B}) \geq (4g(B)+2 + \gamma_{\mathcal{Y}})$. On the other hand it is clear that $(4g(B)+2 + \gamma_{\mathcal{Y}}) > 2 \geq \mu^{max}_{[C']}(\pi'^{*}T_{B})$.  These are exactly the slope inequalities needed to apply Theorem \ref{theo:betterpointsandfoliations} to $\mathcal{Y}$ with $J=2g(B)+3$, $\mathcal T = B$, and $\mathcal G = T_{\mathcal S'/B}$.  Note that since our sections are assumed not to be HN-free and $J > 2g(B)$, Theorem \ref{theo:betterpointsandfoliations}.(1) guarantees that the rational map we obtain is non-trivial.

Consider the dominant family of subvarieties $\mathcal{W}$ on $\mathcal{Y}$ obtained by taking images of the subvarieties  constructed on $\mathcal{S}'$ by Theorem \ref{theo:betterpointsandfoliations}.  These subvarieties $\mathcal{W}$ have the following properties.  First, if we take the strict transform of $C$ in a resolution $\widetilde{\mathcal{W}}$ of $\mathcal{W}$ then deformations go through $\geq 2g(B)+3$ general points of $\widetilde{\mathcal{W}}$.  In particular by Proposition \ref{prop:generalimplieshnfree} the strict transform of a general $C$ in $\widetilde{\mathcal{W}}$ is HN-free.  Second, the codimension in $N$ of the space of sections on $\widetilde{\mathcal{W}}$ can only increase by at most $(\dim(\mathcal{Y})-1)(4g(B) + 3 + \gamma_{\mathcal{Y}}) \leq (\dim(\mathcal{X})-1)(4g(B) + 3 + \gamma)$.
Applying the HN-free version of the desired statement to $\widetilde{\mathcal{W}}$ with the constant $T_{new} = T + (\dim(\mathcal{X})-1)(4g(B) + 3 + \gamma)$, we see that
\begin{equation*}
a(\mathcal{W}_{\eta},f^{*}L|_{\mathcal{W}_{\eta}}) \geq a.
\end{equation*}
Since such $\mathcal{W}$ move in a dominant family on $\mathcal{Y}$, \cite[Lemma 4.8]{LST18} shows that the generic Fujita invariant of $\mathcal{Y}$ is at least as large as that of $\mathcal{W}$ so that
\begin{equation*}
a(\mathcal{Y}_{\eta},f^{*}L|_{\mathcal{Y}_{\eta}}) \geq a.
\end{equation*}
\end{proof}

\begin{remark}
Let $M$ denote the component of $\Sec(\mathcal{X}/B)$ containing the pushforward of the sections parametrized by $N$.  Then the right hand side of Equation \eqref{eq:eh} is $a_{rel} \cdot \mathrm{expdim}(M) - T$.  Since the expected dimension is a lower bound on $\dim(M)$, we can replace Equation \eqref{eq:eh} by the stronger assumption
\begin{equation*}
\dim(N) \geq  a_{rel} \cdot \dim(M) - T.
\end{equation*}
In particular, when $a_{rel}=1$ then $T$ should be thought of as the codimension of $N$ in $M$.  The same remark holds for later theorems as well.
\end{remark}

\subsection{Adjoint rigidity} \label{sect:adjrigid}

Our next goal is to establish a strengthening of Theorem \ref{theo:generalainvsections} that allows us to conclude adjoint rigidity at the cost of increasing the constants.

\begin{lemma} \label{lemm:boundinggeneralpoints}
Let $\pi: \mathcal{Z} \to B$ be a good fibration and let $M$ be an irreducible component of $\Sec(\mathcal{Z}/B)$ parametrizing sections $C$.  Suppose that $H$ is a Cartier divisor on $\mathcal{Z}$ satisfying
\begin{equation*}
H \cdot C+1 <  h^{0}(\mathcal{Z},\mathcal{O}_{\mathcal{Z}}(H)).
\end{equation*}
Then the sections parametrized by $M$ go through at most $H \cdot C+1$ general points of $\mathcal{Z}$.
\end{lemma}

\begin{proof}
Set $Q = H \cdot C + 1$.  Since general points impose codimension $1$ conditions on the linear series $|H|$ we see that for any set of $Q$ points in $\mathcal{Z}$ there is a (possibly reducible) divisor $D \in |H|$ containing all $Q$ points.

Suppose for a contradiction that the sections parametrized by $M$ can go through $Q+1$ general points.  This means that the space of sections through $Q$ general points of $\mathcal{Z}$ forms a dominant family.  In particular, if we fix $Q$ general points and a divisor $D \in |H|$ containing those points, then we can find a section $C$ parametrized by $M$ that contains all the points but is not contained in $\Supp(D)$.  Thus $D \cdot C \geq Q > H \cdot C$, yielding a contradiction.
\end{proof}

\begin{lemma} \label{lemm:h0growthbound}
Let $Z$ be a smooth projective variety of dimension $n$ and let $H$ be a Cartier divisor on $Z$ such that $|H|$ defines a birational morphism.  Then for any non-negative integer $m$ we have
\begin{equation*}
h^{0}(Z,\mathcal{O}_{Z}(mH)) \geq \binom{n+m}{n}
\end{equation*}
\end{lemma}

\begin{proof}
The map $|H|$ defines a morphism $g: Z \to \mathbb{P}^{N}$ for some $N \geq n$ such that $\mathcal{O}_{Z}(H) = g^{*}\mathcal{O}(1)$.  By composing with a generic projection, we obtain a morphism $h: Z \to \mathbb{P}^{n}$ such that $\mathcal{O}_{Z}(H) = h^{*}\mathcal{O}(1)$.  Thus we have
\begin{equation*}
h^{0}(Z,\mathcal{O}_{Z}(mH)) \geq h^{0}(\mathbb{P}^{n},\mathcal{O}(m)) = \binom{n+m}{n}.
\end{equation*}
\end{proof}

We can now prove the criterion for adjoint rigidity.

\begin{theorem} \label{theo:adjointrigidcriterion}
Let $\pi: \mathcal{X} \to B$ be a good fibration and let $L$ be a generically relatively big and semiample Cartier divisor on $\mathcal{X}$.  Assume that $\mathcal{X}_{\eta}$ is geometrically uniruled.  Fix a positive rational number $a_{rel}$ and set $a = a_{rel} a(\mathcal{X}_{\eta},L|_{\mathcal{X}_{\eta}})$.  Fix a rational number $\beta$.  Fix a positive integer $T$.  Fix a positive integer $b>a$ such that $bL|_{\mathcal{X}_{\eta}}$ defines a basepoint free linear series.  Use $b$ to construct an effective $\pi$-vertical $\mathbb{Q}$-Cartier divisor $E$ satisfying the conclusion of Proposition \ref{prop:generalsurjectiontofiber} with respect to $aL$.  
There is some constant $\Gamma = \Gamma(\dim(\mathcal{X}), g(B), \tau(\pi,E), a_{rel}, a, T, \beta,b)$ with the following property.

Suppose that $\psi: \mathcal{Y} \to B$ is a good fibration equipped with a $B$-morphism $f: \mathcal{Y} \to \mathcal{X}$ that is generically finite onto its image.   Suppose that $N$ is an irreducible component of $\Sec(\mathcal{Y}/B)$ parametrizing a dominant family of sections $C$ on $\mathcal{Y}$ which satisfy $f^{*}(K_{\mathcal{X}/B} + a(\mathcal{X}_{\eta},L|_{\mathcal{X}_{\eta}})L) \cdot C \leq \beta$.  Suppose that
\begin{equation*}
\dim(N) \geq a_{rel}(-f^{*}K_{\mathcal{X}/B} \cdot C + (\dim(\mathcal{X})-1)(1-g(B))) - T.
\end{equation*}
Suppose that
\begin{equation*}
a(\mathcal{Y}_{\eta},f^{*}L|_{\mathcal{Y}_{\eta}}) = a.
\end{equation*}
Then either:
\begin{enumerate}
\item $(\mathcal{Y}_{\eta},f^{*}L|_{\mathcal{Y}_{\eta}})$ is adjoint rigid, or
\item deformations of $C$ go through at most $\Gamma$ general points of $\mathcal{Y}$.
\end{enumerate}
\end{theorem}

\begin{proof}
Assume that $(\mathcal{Y}_{\eta},f^{*}L|_{\mathcal{Y}_{\eta}})$ is not adjoint rigid.  We may assume that the general section $C$ is HN-free in $\mathcal{Y}$, since otherwise by Proposition \ref{prop:generalimplieshnfree} the sections parametrized by $N$ can go through at most $2g(B)$ general points of $\mathcal{Y}$.  Applying Lemma \ref{lemm:easyintcalc} we see that
\begin{equation*}
(K_{\mathcal{Y}} + af^{*}L) \cdot C \leq a_{rel} \beta + T + a_{rel} (\dim(\mathcal{X}) - 1) (g(B)-1) + (\dim(\mathcal{X})-1) + 2g(B)-2.
\end{equation*}
Adding in $E$ and rearranging slightly, we obtain
\begin{equation*}
(K_{\mathcal{Y}} + a f^{*}L + f^{*}E) \cdot C \leq a_{rel} \beta + T + a_{rel} (\dim(\mathcal{X}) - 1) (g(B)-1) + \tau(\pi,E) + (\dim(\mathcal{X})-1) + 2g(B) - 2.
\end{equation*}
We denote the right hand side of this equation by $R = R(\dim(\mathcal{X}),g(B),\tau(\pi,E),a_{rel},a,T,\beta,b)$.

Since $E$ satisfies the conclusion of Proposition \ref{prop:generalsurjectiontofiber} with respect to $aL$, the pullback $f^{*}(aL + E)$ is $\mathbb{Q}$-linearly equivalent to an effective $\mathbb{Q}$-divisor $D$ such that $D|_{\mathcal{Y}_{\eta}}$ has smooth irreducible support and coefficient $\frac{a}{b}$. 
Let $\phi: \mathcal{Y}' \to \mathcal{Y}$ be a log resolution of $(\mathcal{Y},D)$ and let $D'$ denote the strict transform of the $\pi$-horizontal components of $D$.  We may ensure that $\phi$ is an isomorphism on an open neighborhood of $\mathcal{Y}_{\eta}$.  In particular $D'$ is still generically relatively big and nef and is irreducible with coefficient $\frac{a}{b}$, so $(\mathcal{Y}',D')$  has terminal singularities.   Since we are assuming the sections are HN-free, the strict transform $C'$ of a general deformation of $C$ avoids any $\phi$-exceptional divisor and thus satisfies
\begin{align*}
(K_{\mathcal{Y'}} + D') \cdot C' & \leq (K_{\mathcal{Y'}} + \phi^{*}D) \cdot C' \\
& = \phi^{*}(K_{\mathcal{Y}} + D) \cdot C' \\
&  \leq R
\end{align*}
Let $F'$ denote a general fiber of $\mathcal{Y}' \to B$.  By Lemma \ref{lemm:relvsabsmmp}, there is some integer $m$ only depending on $\dim(\mathcal{X})$ such that the $(K_{\mathcal{Y}'} + D' + mF')$-MMP is the same as a relative MMP over $B$.  By assumption on the Fujita invariants $K_{\mathcal{Y}'} + D'+mF'$ is on the boundary of the relative pseudo-effective cone over $B$.  Thus the result of the MMP will be a relative Iitaka fibration $\psi: \mathcal{Y}' \dashrightarrow \mathcal{Z}$ for this divisor.  Furthermore since we are assuming $(\mathcal{Y}_{\eta},f^{*}L|_{\mathcal{Y}_{\eta}})$ is not adjoint rigid we know that $\dim(\mathcal{Z}) \geq 2$.

Since $D'$ is relatively big over $B$, it is also relatively big over $\mathcal{Z}$ in the sense of \cite[Definition 2.3]{HX15}.  By \cite[Theorem 1.4]{HX15}, there is a positive integer $k$ only depending on $\dim(\mathcal{X})$ and $\frac{a}{b}$ such that $|k(K_{\mathcal{Y}'} + D' + mF')|$ defines a rational map birational to the Iitaka fibration.  Possibly replacing $k$ by $bk$ we may assume that $k(K_{\mathcal{Y}'} + D' + mF')$ is Cartier.  Applying the canonical bundle formula as in \cite[Section 4]{FM00}, we obtain the following:

\begin{claim}
There is a birational model $\rho: \mathcal{W} \to \mathcal{Y}'$, a morphism $\psi_{W}: \mathcal{W} \to \mathcal{Z}_{W}$ birationally equivalent to $\psi$, and a DCC set $\Lambda \subset [0,1]$ that only depends on $\dim(\mathcal{X})$ and $a/b$ such that:
\begin{itemize}
\item $\mathcal{W}$ and $\mathcal{Z}_{W}$ are smooth,
\item there is an effective $\mathbb{Q}$-Cartier divisor $B_{W}$ and a nef $\mathbb{Q}$-Cartier divisor $M_{W}$ such that $(\mathcal{Z}_{W},B_{W})$ is klt, $B_{W}$ has coefficients in $\Lambda$, and $kM_{W}$ is Cartier,
\item $K_{\mathcal{Z}_{W}} + B_{W} + M_{W}$ is big, and
\item for every integer $p$ divisible by $k$ we have that
\begin{equation*}
h^{0}(\mathcal{Y}',\mathcal{O}_{\mathcal{Y}'}(p(K_{\mathcal{Y}'} + D' + mF'))) = h^{0}(\mathcal{Z}_{W},\mathcal{O}_{\mathcal{Z}_{W}}(\lfloor p(K_{\mathcal{Z}_{W}} + B_{W} + M_{W})\rfloor))
\end{equation*}
and the linear series $| \lfloor p(K_{\mathcal{Z}_{W}} + B_{W} + M_{W}) \rfloor |$ defines a birational map.  
\end{itemize}
\end{claim}

\begin{proof}[Proof of claim:]
We refer to \cite[Theorem 4.2]{FL20} for a modern statement of the canonical bundle formula of \cite{FM00}.  Note that a resolution of the map $\mathcal{Y}' \dashrightarrow \mathcal{Z}$ will be an example of a parabolic map as denoted by $X \to Y$ in the notation of \cite[Theorem 4.2]{FL20}.  The canonical bundle formula yields a birationally equivalent morphism (denoted by $f': X' \to Y'$ in \cite[Theorem 4.2]{FL20}) which in our notation is $\mathcal{W} \to \mathcal{Z}_{W}$.  We suppress the divisors $B^{+},B^{-}$ as they are not necessary for our purposes.  The second bullet point of the claim is \cite[Theorem 4.2.(v)]{FL20} except for the claims about the coefficients of $B_{W}$ and $M_{W}$ which follow from \cite[Claim 3.2.(ii)-(iii)]{HX15}.  The fourth bullet point of the claim is explained in the last paragraph of the proof of \cite[Theorem 1.4]{HX15}.  (Although \cite{HX15} assumes that $(K_{\mathcal{Y}'} + D' + mF')$ is relatively semiample, this is valid in our situation.)  Finally, the third bullet point of the claim follows from the fact that $\mathcal{Y}' \dashrightarrow \mathcal{Z}$ is birationally equivalent to the Iitaka fibration for $(K_{\mathcal{Y}'} + D' + mF')$; using the identification of sections in the fourth bullet point, this is only possible if $K_{\mathcal{Z}_{W}} + B_{W} + M_{W}$ is big.
\end{proof}

We next claim that there is an integer $Q = Q(\dim(\mathcal{X}),g(B),\tau(\pi,E),a_{rel},a,T,\beta,b)$ such that
\begin{equation*}
h^{0}(\mathcal{Z}_{W},\mathcal{O}_{\mathcal{Z}_{W}}(\lfloor Q(K_{\mathcal{Z}_{W}} + B_{W} + M_{W}) \rfloor)) > Q (R +m) + 1
\end{equation*}
Indeed, consider the birational map $\mathcal{Z}_{W} \dashrightarrow V$ defined by $| \lfloor k(K_{\mathcal{Z}_{W}} + B_{W} + M_{W}) \rfloor |$.  Let $\mu: \mathcal{Z}' \to \mathcal{Z}_{W}$ be a smooth model resolving the map and let $H$ denote the basepoint free part of $\mu^{*}| \lfloor k(K_{\mathcal{Z}_{W}} + B_{W} + M_{W}) \rfloor|$.  There are only finitely many possible values of $\dim(\mathcal{Z}_{W})$ which satisfy $\dim(\mathcal{X}) \geq \dim(\mathcal{Z}_{W}) \geq 2$, and thus Lemma \ref{lemm:h0growthbound} gives a quadratic lower bound (that depends only on $\dim(\mathcal{X})$) on the growth rate of sections of multiples of $H$.  In particular there is a constant $Q' = Q'(\dim(\mathcal{X}),g(B),\tau(\pi,E),a_{rel},T,\beta,b)$ such that
\begin{equation*}
h^{0}(\mathcal{Z}',\mathcal{O}_{\mathcal{Z}'}(Q'H)) > (Q')k  (R+ m) + 1.
\end{equation*}
Then we have
\begin{align*}
h^{0}(\mathcal{Z}',\mathcal{O}_{\mathcal{Z}'}(Q'H)) & \leq h^{0}(\mathcal{Z}',\mathcal{O}_{\mathcal{Z}'}(Q'\mu^{*}(\lfloor k(K_{\mathcal{Z}_{W}} + B_{W} + M_{W}) \rfloor ))) \\
& \leq h^{0}(\mathcal{Z}_{W},\mathcal{O}_{\mathcal{Z}_{W}}( \lfloor Q'k(K_{\mathcal{Z}_{W}} + B_{W} + M_{W}) \rfloor )) 
\end{align*}
finishing the proof of the claim with $Q = Q'k$.  Using the comparison of spaces of sections above, we conclude that also
\begin{align*}
h^{0}(\mathcal{Y}',\mathcal{O}_{\mathcal{Y}'}(Q(K_{\mathcal{Y}'} + D' + mF'))) & > Q  (R+ m) + 1 \\
& \geq  Q(K_{\mathcal{Y}'} + D' + mF') \cdot C' + 1
\end{align*}
By applying Lemma \ref{lemm:boundinggeneralpoints} to the divisor $Q(K_{\mathcal{Y}'} + D' + mF')$ we obtain an upper bound
\begin{equation*}
\Gamma(\dim(\mathcal{X}), g(B), \tau(\pi,E), a_{rel} ,a,T, \beta,b) = Q(R + m) + 1
\end{equation*}
on the number of general points that can be contained in deformations of the sections $C'$ on $\mathcal{Y}'$.  But this also implies an upper bound $\Gamma$ on the number of general points that can be contained in deformations of the sections $C$ on $\mathcal{Y}$.
\end{proof}

Suppose that $\mathcal{Y}$ carries a family of sections which have large $L$-degree.  Although we cannot necessarily use Theorem \ref{theo:adjointrigidcriterion} to show that $(\mathcal{Y}_{\eta},f^{*}L|_{\mathcal{Y}_{\eta}})$ is adjoint rigid, by combining with the results of Section \ref{sect:genpoints} we can at least find a covering family of subvarieties of $\mathcal{Y}$ whose generic fibers are adjoint rigid.

\begin{corollary} \label{coro:adjointrigidcodim}
Let $\pi: \mathcal{X} \to B$ be a good fibration and let $L$ be a generically relatively big and semiample Cartier divisor on $\mathcal{X}$.  Assume that $\mathcal{X}_{\eta}$ is geometrically uniruled.  Fix a positive rational number $a_{rel}$ and set $a = a_{rel} a(\mathcal{X}_{\eta},L|_{\mathcal{X}_{\eta}})$.  Fix a rational number $\beta$.  Fix a positive integer $T$.  Fix a positive integer $b>a$ such that $bL|_{\mathcal{X}_{\eta}}$ defines a basepoint free linear series.  Use $b$ to construct an effective $\pi$-vertical $\mathbb{Q}$-Cartier divisor $E$ satisfying the conclusion of Proposition \ref{prop:generalsurjectiontofiber} with respect to $aL$.   
There are constants $\xi^{+} = \xi^{+}(\dim(\mathcal{X}), g(B), \tau(\pi,E), a_{rel}, a, T, \beta,b)$ and $T^{+} = T^{+}(\dim(\mathcal{X}), g(B), \tau(\pi,E), a_{rel}, a, T, \beta,b)$ with the following property.

Suppose that $\psi: \mathcal{Y} \to B$ is a good fibration equipped with a $B$-morphism $f: \mathcal{Y} \to \mathcal{X}$ that is generically finite onto its image. Suppose that $N$ is an irreducible component of $\Sec(\mathcal{Y}/B)$ parametrizing a dominant family of sections $C$ on $\mathcal{Y}$ which satisfy $f^{*}L \cdot C \geq \xi^{+}$ and $f^{*}(K_{\mathcal{X}/B} + a(\mathcal{X}_{\eta},L|_{\mathcal{X}_{\eta}})L) \cdot C \leq \beta$.  Suppose that
\begin{equation*}
\dim(N) \geq a_{rel}(-f^{*}K_{\mathcal{X}/B} \cdot C + (\dim(\mathcal{X})-1)(1-g(B))) - T.
\end{equation*}
Suppose that
\begin{equation*}
a(\mathcal{Y}_{\eta},f^{*}L|_{\mathcal{Y}_{\eta}}) = a.
\end{equation*} 
Let $g: \mathcal{S} \to \mathcal{Y}$ denote the finite part of the Stein factorization of the evaluation map for the normalization of the universal family over $N$.  Then there is a dominant rational $B$-map $\phi: \mathcal{S} \dashrightarrow \mathcal{T}$ to  a normal projective $B$-variety such that the following holds.  For a general section $C^{\dagger}$ on $\mathcal{S}$ parametrized by $N$ let $\mathcal{W}$ denote the main component of the closure of $\phi^{-1}(\phi(C^{\dagger}))$.
Then:
\begin{enumerate}
\item We have $a(\mathcal{W}_{\eta},g^{*}f^{*}L|_{\mathcal{W}_{\eta}}) = a$ and the pair $(\mathcal{W}_{\eta},g^{*}f^{*}L|_{\mathcal{W}_{\eta}})$ is adjoint rigid,
\item $\mathcal{W}$ is swept out by the sections parametrized by a sublocus $N_{\mathcal{W}} \subset N$ whose closure has codimension $\leq T^{+}$ in $N$, and
\item  there is a resolution of $\mathcal{W}$ such that the strict transform of a general section in $N_{\mathcal{W}}$ to the resolution goes through $\geq 2g(B)+1$ general points and is HN-free.
\end{enumerate}
\end{corollary}

\begin{proof}
Define $d = \dim(\mathcal{Y})$ and set $\gamma = (dg(B) - g(B) + 1)^{2}(d-1)$.  We also define constants $T_{k}$ and $\Gamma_{k}$ for $2 \leq k \leq d$ as follows.  We first set $T_{d} = 0$ and
\begin{equation*}
\Gamma_{d} = \sup \{ 2g(B)+3, \Gamma(\dim(\mathcal{X}), g(B), \tau(\pi,E), a_{rel},a,T,\beta,b) + 1 \}
\end{equation*}
where $\Gamma$ is the constant defined in Theorem \ref{theo:adjointrigidcriterion}.  Then for $2 \leq k < d$ we define via a downward induction  
\begin{equation*}
T_{k} = k(\Gamma_{k+1} + 2g(B) + \gamma) + T_{k+1}
\end{equation*}
and
\begin{equation*}
\Gamma_{k} = \sup \{ 2g(B)+3, \Gamma_{k+1}, \Gamma(\dim(\mathcal{X}), g(B), \tau(\pi,E), a_{rel},a,T + T_{k},\beta,b) + 1 \}.
\end{equation*}
Finally, we set $T^{+} = T + \sup_{k=2,\ldots,d} T_{k}$, $\Gamma^{+} = \sup_{k=2, \ldots, d} \Gamma_{k}$ and $\xi^{+}$ to be the maximum of the constant $\xi(\dim(\mathcal{X}),g(B),\tau(\pi,E),a_{rel}, a, T^{+}, \beta,b)$ as in Theorem \ref{theo:generalainvsections} and of
\begin{equation*}
\frac{1}{a} (\dim(\mathcal{X})(\Gamma^{+} + 2g(B) + \gamma+1) + a_{rel}\beta + T^{+} +  (a_{rel} + 1)( \dim(\mathcal{X}) - 1) (g(B)-1)).
\end{equation*}

Recall that $\mathcal{S}$ denotes the finite part of the Stein factorization of the evaluation map for the normalization of the universal family over $N$.  We let $\mathcal{S}'$ denote a smooth birational model of $\mathcal{S}$ that flattens the family of sections on $\mathcal{S}$ as in Construction \ref{cons:flatteningfamilyofcurves}.  We denote the strict transform of a general section in our family on $\mathcal{S}'$ by $C'$ and denote the family of deformations of $C'$ by $N'$.  Appealing to Corollary \ref{coro:domfamilyexpdim} we have
\begin{align*}
-K_{\mathcal{S}'/B} \cdot C' + (\dim(\mathcal{S}')-1) & \geq \dim(N') \\
& = \dim(N) \\
& \geq -K_{\mathcal{Y}/B} \cdot C + (\dim(\mathcal{Y})-1)(1-g(B)) \\
& \geq  af^{*}L \cdot C - a_{rel} \beta - T^+ -  (\dim(\mathcal{X})-1) \\
& \qquad -(a_{rel}\dim(\mathcal{X}) + \dim(\mathcal{Y}) - a_{rel} - 1) (g(B)-1) 
\end{align*}
where the last line follows from Equation \eqref{eq:alternateform} of Lemma \ref{lemm:easyintcalc}.
Combining with the bound $f^{*}L \cdot C \geq \xi^{+}$, we conclude that 
\begin{equation} \label{eq:sprimebound}
-K_{\mathcal{S}'/B} \cdot C' \geq \dim(\mathcal{S}') (\Gamma^{+} + 2g(B) + \gamma - 1) + 2.
\end{equation}
We next inductively define foliations $\mathcal{G}_{d},\mathcal{G}_{d-1},\ldots$ on $\mathcal{S}'$ by repeatedly applying Theorem \ref{theo:betterpointsandfoliations} to $\mathcal{Y}$ using the constants $\Gamma_{d},\Gamma_{d-1},\ldots$.  We will also denote by $\psi_{i}$ the rational map on $\mathcal{S}'$ induced by the foliation $\mathcal{G}_{i}$.  We will inductively verify the inequalities 
\begin{align*}
\mu_{[C]}^{max}(\mathcal{G}_{i}) & \geq \Gamma_{i} + 2g(B) + \gamma - 1 \\
\mu_{[C]}^{max}(T_{\mathcal{S}'}/\mathcal{G}_{i}) & < \Gamma_{i} + 2g(B) + \gamma - 1
\end{align*}
which show the requirements necessary to inductively apply Theorem \ref{theo:betterpointsandfoliations}.

For the base case we set $\mathcal{G}_{d} = T_{\mathcal{S}'/B}$.  By construction we have $\Gamma_{d} > 2 \geq \mu_{[C]}^{max}(\pi^{*}T_{B})$ and Equation \eqref{eq:sprimebound} shows that 
\begin{equation*}
\mu_{[C]}^{max}(T_{\mathcal{S}'/B}) \geq \mu_{[C]}(T_{\mathcal{S}'/B}) \geq \Gamma^+ + 2g(B) + \gamma - 1 \geq \Gamma_{d} + 2g(B) + \gamma - 1.
\end{equation*}
Thus we have verified the two necessary inequalities in the base case.

Now suppose inductively that we apply Theorem \ref{theo:betterpointsandfoliations} to $\mathcal{Y}$ for the foliation $\mathcal{G}_{i}$ on $\mathcal{S}'$ and the constant $\Gamma_{i}$.  There are two possible outcomes.  The first possibility is that if we set $\mathcal{P}_{i}$ to be the main component of $\overline{\psi_{i}^{-1}(\psi_{i}(C'))}$ for a general section $C'$ in our family then the deformations of $C'$ go through at least $\Gamma_{i}$ general points of $\mathcal{P}_{i}$.  In this case we stop the inductive process.  The second possibility is that we obtain a new foliation $\mathcal{G}_{i-1}$ and a new rational map $\psi_{i-1}$.
Note that Theorem \ref{theo:betterpointsandfoliations}.(c) shows that 
\begin{equation*}
\mu_{[C]}^{max}(T_{\mathcal{S}'}/\mathcal{G}_{i-1}) < \Gamma_{i} + 2g(B) + \gamma - 1 \leq \Gamma_{i-1} + 2g(B) + \gamma - 1.
\end{equation*}
On the other hand, letting $r$ denote the rank of $\mathcal G_{i-1}$ we have
\begin{align*}
\mu^{max}_{[C]}(\mathcal{G}_{i-1}) \geq \mu_{[C]}(\mathcal{G}_{i-1}) & = \frac{c_{1}(T_{\mathcal{S}'/B}) \cdot C - c_{1}(T_{\mathcal{S}'/B}/\mathcal{G}_{i-1}) \cdot C}{r} \\
& = \frac{c_{1}(T_{\mathcal{S}'/B}) \cdot C - c_{1}(T_{\mathcal{S}'}/\mathcal{G}_{i-1}) \cdot C + 2g(B) - 2}{r} \\
& \geq \frac{\dim(\mathcal{S}') (\Gamma^{+} + 2g(B) + \gamma - 1) - (\dim(\mathcal{S}')-r)(\Gamma_{i} + 2g(B) + \gamma - 1)}{r} \\
& \geq \Gamma^{+} + 2g(B) + \gamma - 1 \\
& \geq \Gamma_{i-1} + 2g(B) + \gamma - 1
\end{align*}
where the third line is a consequence of Equation \eqref{eq:sprimebound}.
Thus we have verified the necessary inequalities for continuing the inductive process.

 Since the dimension of $\overline{\psi_{i}^{-1}(\psi_{i}(C'))}$ is always at least $2$, this process stops after at most $d-2$ steps with either
\begin{itemize}
\item a foliation $\mathcal{G}_{k}$ such that if we set $\mathcal{P}_{k}$ to be the main component of $\overline{\psi_{k}^{-1}(\psi_{k}(C'))}$ for a general section $C'$ in our family then the deformations of $C'$ go through at least $\Gamma_{k}$ general points of $\mathcal{P}_{k}$, or
\item a rank $1$ foliation $\mathcal{G}_{k}$ such that if we set $\mathcal{P}_{k}$ to be the main component of $\overline{\psi_{k}^{-1}(\psi_{k}(C'))}$ for a general section $C'$ in our family then the deformations of $C'$ go through at least $\Gamma_{k+1}$ general points of $\mathcal{P}_{k}$.
\end{itemize}
Then the foliation $\mathcal{G}_{k}$ induces a rational map $\mathcal{S}' \dashrightarrow \mathcal{T}$, and hence also a rational map $\mathcal{S} \dashrightarrow \mathcal{T}$.  We prove that in either case this map has the desired properties.  Note that the subvarieties $\mathcal{P}_{k}$ of $\mathcal{S}'$ are birational to the subvarieties $\mathcal{W}$ in the statement of the theorem, and it suffices to prove that $\mathcal{P}_{k}$ has the desired properties.

By applying Theorem \ref{theo:betterpointsandfoliations}.(2).(b) inductively, we see that the codimension of the space of deformations of $C'$ in $\mathcal{P}_{k}$ inside of $N$ is at most 
\begin{equation*}
\sum_{j=k}^{d} (j-1)(\Gamma_{j} + 2g(B) + \gamma)
\end{equation*}
verifying (2).  Note that by construction the deformations of $C'$ in $\mathcal{P}_{k}$ go through either $\Gamma_{k}$ or $\Gamma_{k+1}$ general points of $\mathcal{P}_{k}$.  Both quantities are at least $2g(B)+1$.  This property is preserved by passing to the strict transform, and curves through this many general points must be HN-free by Proposition \ref{prop:generalimplieshnfree}, proving (3).  Using the lower bound $\xi \leq \xi^{+}$, we can apply Theorem \ref{theo:generalainvsections} to $\mathcal{P}_{k}$ equipped with the family of deformations of $C'$ to see that
\begin{equation*}
a(\mathcal{P}_{k,\eta},f^{*}L|_{\mathcal{P}_{k,\eta}})  \geq a
\end{equation*}
But since the deformations of $\mathcal{P}_{k,\eta}$ form a dominant family of subvarieties on $\mathcal{Y}_{\eta}$ the equality must be achieved.  To prove (1), it only remains to verify the adjoint rigidity.  If $\mathcal{G}_{k}$ has rank $1$, $\mathcal{P}_{k}$ is a $\mathbb{P}^{1}$-fibration over $B$ and thus is automatically adjoint rigid.  If $\mathcal{G}_{k}$ has rank $>1$, then deformations of $C'$ on a resolution  $\widetilde{\mathcal{P}}_{k}$ of $\mathcal{P}_{k}$ go through at least $\Gamma_{k}$ general points.   We then apply Theorem \ref{theo:adjointrigidcriterion} on $\widetilde{\mathcal{P}}_{k}$ to determine adjoint rigidity.  This proves (3).
\end{proof}

\section{Boundedness statements} \label{sect:boundedness}

We now turn to proving boundedness statements for the set of morphisms $f: \mathcal{Y} \to \mathcal{X}$ such that $\mathcal{Y}$ carries a family of sections which is ``large'' on $\mathcal{X}$.  In Section \ref{sect:mainboundedresults} we prove several technical statements which combine \cite{birkar22} with our work on twists. In Section \ref{sect:boundedconsequences} we state and prove our main boundedness result, Theorem~\ref{theo:combinedbounded}.

\subsection{Boundedness} \label{sect:mainboundedresults}

Our first statement appeals to the recent results of \cite{birkar22} to prove birational boundedness when the generic fiber is adjoint rigid.

\begin{theorem}  \label{theo:adjointrigidbounded}
Let $\pi: \mathcal{X} \to B$ be a good fibration and let $L$ be a generically relatively big and semiample Cartier divisor on $\mathcal{X}$.  Assume that $\mathcal{X}_{\eta}$ is geometrically uniruled.  Fix a positive rational number $a_{rel}$ and set $a = a_{rel} a(\mathcal{X}_{\eta},L|_{\mathcal{X}_{\eta}})$.  Fix a rational number $\beta$.  Fix a positive integer $T$.  Fix a positive integer $b>a$ such that $bL|_{\mathcal{X}_{\eta}}$ defines a basepoint free linear series.  Use $b$ to construct an effective $\pi$-vertical $\mathbb{Q}$-Cartier divisor $E$ satisfying the conclusion of Proposition \ref{prop:generalsurjectiontofiber} with respect to $aL$.  There is some constant $\xi = \xi(\dim(\mathcal{X}), g(B), \tau(\pi,E), a_{rel}, a, T, \beta,b)$ with the following property.

Suppose that $\psi: \mathcal{Y} \to B$ is a good fibration equipped with a $B$-morphism $f: \mathcal{Y} \to \mathcal{X}$ that is generically finite onto its image.  Suppose that $a(\mathcal{Y}_{\eta},f^{*}L|_{\mathcal{Y}_{\eta}}) = a$ and that $(\mathcal{Y}_{\eta},f^{*}L|_{\mathcal{Y}_{\eta}})$ is adjoint rigid.  Suppose that $N$ is an irreducible component of $\Sec(\mathcal{Y}/B)$ parametrizing a dominant family of HN-free sections $C$ on $\mathcal{Y}$ which satisfy $f^{*}L \cdot C \geq \xi$ and $f^{*}(K_{\mathcal{X}/B} + a(\mathcal{X}_{\eta},L|_{\mathcal{X}_{\eta}})L) \cdot C \leq \beta$.  Finally, suppose that
\begin{equation}
\dim(N) \geq a_{rel}(-f^{*}K_{\mathcal{X}/B} \cdot C + (\dim(\mathcal{X})-1)(1-g(B))) - T.
\end{equation}
Then:
\begin{enumerate}
\item The set of such projective varieties $\mathcal{Y}$ is birationally bounded.
\item Suppose that $L$ is big and semiample.  Choose a positive integer $b>a$ such that $|bL|$ is basepoint free.  Then there is a constant $\daleth = \daleth(\dim(\mathcal{X}), g(B), a_{rel}, a, T, \beta,b)$ such that $\vol(f^{*}L) \leq \daleth$.
\end{enumerate}
\end{theorem}

\begin{proof}
Lemma \ref{lemm:easyintcalc} shows that
\begin{equation*}
(K_{\mathcal{Y}} + af^{*}L) \cdot C \leq a_{rel} \beta + T + a_{rel} (\dim(\mathcal{X}) - 1) (g(B)-1) + (\dim(\mathcal{X})-1) + 2g(B)-2.
\end{equation*}
Proposition \ref{prop:generalsurjectiontofiber} shows that $f^{*}(aL + E)$ is $\mathbb{Q}$-linearly equivalent to an effective $\mathbb{Q}$-Cartier divisor $D$ such that $(\mathcal{Y}_{\eta},D|_{\mathcal{Y}_{\eta}})$ is a terminal pair.  Let $\phi: \mathcal{Y}' \to \mathcal{Y}$ be a log resolution of this pair and let $D'$ denote the strict transform of the $\pi$-horizontal components of $D$.  Since $D'$ is irreducible we see that $(\mathcal{Y}',D')$ is a terminal pair.

Since we are assuming the sections are HN-free, the strict transform $C'$ of a general deformation of $C$ avoids any $\phi$-exceptional divisor and thus satisfies
\begin{align*}
(K_{\mathcal{Y'}} + D') \cdot C' & \leq (K_{\mathcal{Y'}} + \phi^{*}D) \cdot C' \\
& = \phi^{*}(K_{\mathcal{Y}} + D) \cdot C' \\
&  \leq a_{rel} \beta + T +  \tau(\pi,E) + a_{rel} (\dim(\mathcal{X}) - 1) (g(B)-1) \\
& \qquad \qquad + (\dim(\mathcal{X})-1) + 2g(B)-2.
\end{align*}
Furthermore $C'$ is HN-free on $\mathcal{Y}'$.

Run the relative MMP for $K_{\mathcal{Y}'} + D'$ over $B$.  Due to our adjoint rigidity assumption on the generic fiber, the result will be a birational model $\rho: \mathcal{Y}' \dashrightarrow \widetilde{\mathcal{Y}}$ where $K_{\widetilde{\mathcal{Y}}} + \rho_{*}D'$ is relatively $\mathbb{Q}$-linearly equivalent to $0$.  We denote by $\widetilde{\psi}$ the structural map $\widetilde{\psi}: \widetilde{\mathcal{Y}} \to B$.  Write $K_{\widetilde{\mathcal{Y}}} + \rho_{*}D' \sim_{\mathbb{Q}} \widetilde{\psi}^{*}P$.  Since the map $\rho$ is a $(K_{\mathcal{Y}'} + D')$-negative birational contraction,
\begin{align*}
\widetilde{\psi}^{*}P \cdot \rho_{*}C' & = \rho^{*}(K_{\widetilde{\mathcal{Y}}} + \rho_{*}D') \cdot C' \\
& \leq (K_{\mathcal{Y}'} + D') \cdot C' \\
& \leq a_{rel} \beta + T +  \tau(\pi,E) + a_{rel} (\dim(\mathcal{X}) - 1) (g(B)-1) + (\dim(\mathcal{X})-1) + 2g(B)-2.
\end{align*}
Then \cite[Theorem 1.3]{birkar22} applies with $Z=B$ and $A = \lceil a_{rel} \beta + T +  \tau(\pi,E) + a_{rel} (\dim(\mathcal{X}) - 1) (g(B)-1) + (\dim(\mathcal{X})-1) + 2g(B)- 1 \rceil p$ for a point $p \in B$, showing that the set of minimal models $(\widetilde{\mathcal{Y}},\rho_{*}D')$ is log bounded.

Now suppose that $L$ is big and semiample. Note that the effective divisor $E=0$ satisfies the conclusion of Proposition \ref{prop:generalsurjectiontofiber} with respect to $aL$, and we make this choice for $E$.  Recall that the divisor $D$ is constructed by applying Proposition \ref{prop:generalsurjectiontofiber}.  However, in our setting we can choose $D \sim_{\mathbb{Q}} aL$ where $\Supp(D)$ is smooth and irreducible and has coefficient $\frac{a}{b}$.  In particular we can take $\mathcal{Y}' = \mathcal{Y}$ and $D' = D$.

Repeating the argument above, we run a relative MMP $\rho: \mathcal{Y} \dashrightarrow \widetilde{\mathcal Y}$ and see that the resulting pairs $(\widetilde{\mathcal Y},\rho_{*}D)$ are log bounded.  By construction $\rho_{*}D$ is irreducible with coefficient $a/b$.  This implies that there is some constant $\daleth$ such that
\begin{equation*}
\vol(f^{*}L) = \frac{\vol(D)}{a^{\dim \mathcal Y}} \leq \frac{\vol(\rho_{*}D)}{a^{\dim \mathcal Y}} \leq \daleth.
\end{equation*}
\end{proof}

\begin{remark}
In the setting of Theorem \ref{theo:adjointrigidbounded}.(1) the variety $\mathcal{Y}$ can be replaced by a higher birational model and $N$ can be replaced by the strict transform family of curves without affecting the hypotheses.  Thus birational boundedness is the best one can hope for.
\end{remark}

Construction \ref{cons:rigidsubvarieties} and Theorem \ref{theo:ainvboundedandtwists} construct certain families of varieties over $K(B)$.  We next modify these constructions to apply to integral models.

\begin{construction} \label{cons:allsubvarieties}
Let $\pi: \mathcal{X} \to B$ be a good fibration and let $L$ be a big and semiample Cartier divisor on $\mathcal{X}$.  Assume that $\mathcal{X}_{\eta}$ is geometrically uniruled.  Set $a = a(\mathcal{X}_{\eta},L|_{\mathcal{X}_{\eta}})$. 

By applying Construction \ref{cons:rigidsubvarieties} to $\mathcal{X}_{\eta}$ we obtain a proper closed subset $\mathcal V_{\eta} \subset \mathcal{X}_{\eta}$ and a finite collection of families $p_{i,\eta} : \mathcal{U}_{i,\eta} \to \mathcal{W}_{i,\eta}$ whose smooth fibers are birational to closed subvarieties of $\mathcal{X}_{\eta}$. Let $p_{i}: \mathcal{U}_{i} \to \mathcal{W}_{i}$ denote any smooth integral model such that the structural morphism $\mathcal{U}_{i,\eta} \to \mathcal{X}_{\eta}$ extends to $\mathcal{U}_{i}$ and let $\mathcal V$ denote the closure of $\mathcal V_\eta$.  The subvarieties parametrized by $p_{i,\eta}$ correspond to the $K(B)$-points of $\mathcal{W}_{i,\eta}$, or equivalently, to sections of $\mathcal{W}_{i}$ over $B$.  Let $\mathfrak{W}_i$ denote $\Sec(\mathcal{W}_i/B)$. We let $\mathfrak W = \sqcup_i \mathfrak W_i$.
We first shrink $\mathfrak{W}$ so that the generic point of every section parametrized by $\mathfrak{W}$ is contained in the open locus over which $\sqcup_{i} p_{i}$ is smooth. We enlarge $\mathcal V$ by adding the images in $\mathcal{X}$ of the loci where the maps $p_i$ fail to be smooth. Consider the universal family $\mathfrak W_i \times B \to \mathfrak W_i$ with the evaluation map $\mathfrak W_i \times B \to  \mathcal{W}_i$.  By taking a base change of $p_{i}: \ \mathcal{U}_{i} \to  \mathcal{W}_{i}$ over this morphism, we obtain a morphism which we denote by $\mathfrak{Z}_i \to \mathfrak{W}_i \times B$. We let $\mathfrak Z = \sqcup_i \mathfrak Z_i$.

Note that $\mathfrak{W}$ is a countable union of quasiprojective schemes and $\mathfrak{Z} \to \mathfrak{W} \times B$ is a finite type morphism such that for every closed point $w \in \mathfrak{W}$ the $B$-scheme $\mathfrak{Z}_{w} \to \{w\} \times B$ has the property that $\mathfrak{Z}_{w,\eta}$ is isomorphic to a fiber of $p_{i,\eta}$ over a $K(B)$-point of $\mathcal{W}_{i,\eta}$. 
By repeatedly stratifying $\mathfrak{W}$ into locally closed subsets and taking resolutions of components of $\mathfrak{Z}$, we may also ensure that the fibers of $\mathfrak{Z} \to \mathfrak{W}$ are smooth.  We denote the evaluation map by $\iota : \mathfrak Z \to \mathcal X$.
\end{construction}

\begin{construction} \label{cons:allfamilies}

Let $\sqcup_{G} \mathcal H(G, B)$ be the Hurwitz stack.  Fix an \'etale covering $\sqcup_G \mathcal H_G \to \sqcup_G \mathcal H(G, B)$ from a scheme.

Let $\pi: \mathcal{X} \to B$ be a good fibration and let $L$ be a big and semiample Cartier divisor on $\mathcal{X}$.  Assume that $\mathcal{X}_{\eta}$ is geometrically uniruled.  Set $a = a(\mathcal{X}_{\eta},L|_{\mathcal{X}_{\eta}})$.  Let $\mathfrak{Z} \to \mathfrak{W} \times B$ be the morphism constructed in Construction \ref{cons:allsubvarieties}.

By Theorem \ref{theo:ainvboundedandtwists} there is a finite set of smooth projective $K(B)$-varieties $\mathcal{Y}_{i,j, \eta}$ equipped with dominant generically finite morphisms $h_{i, j, \eta}: \mathcal{Y}_{i, j,\eta} \to \mathcal U_{i, \eta}$ and a closed set $\mathcal R_{\eta} \subset \mathcal{X}_{\eta}$ that has the following property.  Suppose that $g: \mathcal{Y}_{\eta} \to \mathcal{X}_{\eta}$ is a generically finite morphism from a geometrically integral smooth projective variety such that $g(\mathcal{Y}_{\eta})$ is not contained in $\mathcal R_{\eta}$.  Suppose furthermore that $a(\mathcal{Y}_{\eta},g^{*}L|_{\mathcal{Y}_{\eta}}) = a$ and that $(\mathcal{Y}_{\eta},g^{*}L|_{\mathcal{Y}_{\eta}})$ is adjoint rigid.  Since $\mathcal Y_{\eta}$ is geometrically rationally connected by \cite[Theorem 4.5]{LTT18}, \cite[Theorem 1.1]{GHS03} and \cite[Theorem 12]{HT06} show that $\mathcal{Y}_{\eta}$ carries a dense set of rational points.  Thus Theorem \ref{theo:ainvboundedandtwists} shows that there are indices $i,j$ such that the map $g$ factors rationally through a twist $h^{\sigma}_{i, j, \eta}$ and $\mathcal{Y}_{\eta}$ maps birationally to a fiber of a morphism $\mathcal{Y}^{\sigma}_{i,j,\eta} \to \mathcal{T}^{\sigma}_{i,j,\eta}$. 

Let $\mathcal V_\eta$ be the union of $\mathcal R_\eta$ with the generic fiber of the closed set from Construction~\ref{cons:allsubvarieties}.
Then we enlarge $\mathcal V_\eta$ by adding $s_i(\mathcal B_{i, j, \eta})$ where $\mathcal B_{i, j, \eta}$ is the union of the irreducible components of the branch locus of $h_{i, j, \eta}$.  We further enlarge $\mathcal{V}_\eta$ by adding the Zariski closure of the union of the images of the fibers of $r_{i, j, \eta}$ which fail to be smooth, fail to have the same $a$-invariant as $\mathcal Y_{i, j}$, or fail to be adjoint rigid. 
If $\mathcal{Y}_{i,j, \eta}$ is a component such that some twist of $\mathcal{Y}_{i,j,\eta}$ admits a $K(B)$-rational point mapping to $\mathcal X_\eta \setminus \mathcal V_\eta$, then we replace $\mathcal{Y}_{i,j, \eta}$ by this twist.  If $\mathcal{Y}_{i,j, \eta}$ is a component such that no twist of $\mathcal{Y}_{i,j,\eta}$ admits a $K(B)$-rational point mapping to $\mathcal X_\eta \setminus \mathcal V_\eta$, then we remove $\mathcal{Y}_{i,j, \eta}$ from our set.

Set $\mathfrak{D}_{i, j} = \sqcup_{G} \mathcal C^1(\mathcal G_{\mathcal H_G}, \mathcal K(\mathcal Y_{i, j, \eta}/\mathcal U_{i, \eta}) _{\mathcal H_G})$. 
By construction $\mathfrak{D}_{i, j}$ is a countable union of finite type schemes over $\mathbb{C}$.  As described in Section \ref{sect:familyoftwists}, there is a morphism 
\[
h_{i, j, \eta}: \mathfrak Y_{i, j, \eta} \to \mathfrak{D}_{i, j} \times \mathcal U_{i, \eta}
\]
which parametrizes twists of $h_{i, j, \eta}: \mathcal{Y}_{i, j, \eta} \to \mathcal{U}_{i, \eta}$.  After perhaps replacing $\mathfrak{D}_{i, j}$ by a stratification into locally closed subsets, we can construct integral models in families to obtain a map $\mathfrak Y_{i, j} \to \mathfrak{D}_{i, j}\times \mathcal{U}_{i} \to \mathfrak D_{i, j} \times \mathcal X$ whose composition we denote by $f_{i, j}$.
After perhaps again replacing $\mathfrak{D}_{i, j}$ by a stratification by locally closed subsets and taking resolutions of irreducible components of $\mathfrak{Y}_{i, j}$, we may ensure that for every closed point $d \in \mathfrak{D}_{i, j}$ the fiber $\mathfrak{Y}_{i, j, d}$ is a good fibration equipped with a $B$-morphism $f_{i, j, d}: \mathfrak{Y}_{i, j, d} \to \mathcal{X}$.  By construction every twist of $h_{i,j,\eta}$ has an integral model $h_{i, j}^\sigma: \mathcal{Y}^{\sigma}_{i, j} \to \mathcal{U}_i$ that is a member of our family.   After again replacing $\mathfrak D_{i, j}$ by a stratification into locally closed subsets, we may ensure that the Stein factorization of the composition $\mathfrak{Y}_{i, j} \to \mathfrak D_{i, j} \times\mathcal{U}_i \to \mathfrak D_{i, j} \times \mathcal W_i$ induces for every fiber over a closed point in $\mathfrak{D}_{i, j}$ the Stein factorization of $\mathcal{Y}_{i,j}^{\sigma} \to \mathcal{W}_{i}$.  Denote the Stein factorization of $\mathfrak{Y}_{i, j} \to \mathfrak{D}_{i, j} \times \mathcal{W}_{i}$ by $r_{i, j} : \mathfrak Y_{i, j} \to \mathfrak T_{i, j}$.  Then $r_{i, j} : \mathfrak Y_{i, j} \to \mathfrak T_{i, j}$ and $t_{i, j} : \mathfrak T_{i, j}\to \mathcal W_i$ define a family of Stein factorizations $r^{\sigma}_{i, j}: \mathcal{Y}^{\sigma}_{i, j} \to \mathcal T^\sigma_{i, j}$, $t_{i, j}^\sigma : \mathcal T^\sigma_{i, j} \to \mathcal W_i$. 
Due to the functoriality in Section~\ref{subsec:functoriality}, we see that $\mathfrak T_{i, j, \eta} \to \mathfrak D_{i, j} \times \mathcal W_{i, \eta}$ parametrizes the family of twists of $\mathcal T_{i, j, \eta} \to \mathcal{W}_{i, \eta}$ which are induced by twists of $\mathcal{Y}_{i,j, \eta} \to \mathcal{U}_{i, \eta}$.

Let $\mathcal V$ be the closure of $\mathcal V_\eta$.
We further enlarge $\mathcal{V}$ by adding the Zariski closure of the union of the images of the fibers of $r_{i, j}$ which fail to be smooth, fail to have the same $a$-invariant as $\mathcal Y_{i, j}$, or fail to be adjoint rigid.  We let $\mathcal B_{i, j}$ be the closure of $\mathcal B_{i, j, \eta}$ and define $\mathcal B : = \cup_{i, j} \mathcal B_{i, j}$. We also define $\mathcal B_{i, j}' \subset \mathcal W_i$ as the union of components of the branch locus of $t_{i, j}$ which dominate $B$ and the closures of the images of loci where fibers of $r_{i, j}$ fail to be smooth, fail to have the same $a$-invariant as $\mathcal Y_{i, j}$, or fail to be adjoint rigid. We define $\mathcal B' : = \cup_{i, j} \mathcal B_{i, j}'$.

Recall that we assume that $\mathcal Y_{i, j, \eta}$ admits a $K(B)$-rational point $y$ mapping to $\mathcal X_\eta \setminus \mathcal V_\eta$.
Let $\widetilde{h}_{i, j, \eta} : \widetilde{\mathcal Y}_{i, j, \eta} \to \mathcal U_{i, \eta}$ be a geometric Galois closure of $h_{i, j, \eta} : \mathcal Y_{i, j, \eta} \to \mathcal U_{i, \eta}$ such that $\mathrm{Bir}(\widetilde{\mathcal Y}_{i, j, \overline{\eta}}/\mathcal U_{i, \overline{\eta}}) = \mathrm{Aut}(\widetilde{\mathcal Y}_{i, j, \overline{\eta}}/\mathcal U_{i, \overline{\eta}})$ and $\widetilde{\mathcal Y}_{i, j, \eta}$ admits a $K(B)$-rational point $\widetilde{y}$ mapping to $y$. 
Let $\mathcal Y_{i, j} \xrightarrow{\widehat{h}_{i,j}} \mathcal P_{i, j} \xrightarrow{\ell_{i, j}} \mathcal U_i$ be the cover corresponding to the normalizer of $\mathrm{Aut}(\widetilde{\mathcal Y}_{i, j, \overline{\eta}}/\mathcal Y_{i, j, \overline{\eta}})$ in $\mathrm{Aut}(\widetilde{\mathcal Y}_{i, j, \overline{\eta}}/\mathcal U_{i, \overline{\eta}})$ such that $\mathcal P_{i, j}$ is normal and $\ell_{i, j}$ is finite. Note that every twist $\mathcal Y^\sigma_{i, j} \to \mathcal U_{i}$ factors through $\ell_{i, j}: \mathcal P_{i, j} \to \mathcal U_{i}$.  By taking the Stein factorization, we have a commutative diagram
\begin{equation*}
\xymatrix{ {\mathcal Y}_{i, j} \ar[r]^{\widehat{h}_{i, j}} \ar[d]_{r_{i, j}} & \mathcal P_{i, j} \ar[d]^{b_{i, j}} \ar[r]^{\ell_{i, j}}&  {\mathcal U}_{i} \ar[d]_{p_{i}} \\
\mathcal T_{i, j} \ar[r]_{c_{i, j}} & \mathcal S_{i, j} \ar[r]_{a_{i, j}} & \mathcal W_{i}}
\end{equation*}
where $\mathcal S_{i, j}$ is projective and normal, $b_{i, j}$ has connected fibers, and $a_{i, j}$ is finite.

We also let $\mathfrak M_{i, j} = \sqcup_G \mathcal C^1(\mathcal G_{\mathcal H_G}, \mathcal K(\mathcal T_{i, j, \eta}/\mathcal W_{i, \eta})_{\mathcal H_G})$ denote the parameter space for all twists of $\mathcal{T}_{i,j,\eta} \to \mathcal{W}_{i,\eta}$.  We denote the universal family over $\mathfrak M_{i, j}$ by $\mathfrak T_{i, j}' \to \mathfrak M_{i, j} \times \mathcal W_i$. Then by the functoriality established in Section~\ref{subsec:functoriality} we have a natural isomorphism
\[
\mathfrak T_{i, j} \to \mathfrak T'_{i, j} \times_{\mathfrak M_{i, j}} \mathfrak D_{i, j}.
\]
We set $\mathfrak Y = \sqcup_{i, j}\mathfrak Y_{i, j}$, $\mathfrak T = \sqcup_{i, j} \mathfrak T_{i, j}$, $\mathfrak T' = \sqcup_{i, j} \mathfrak T_{i, j}'$, $\mathfrak D = \sqcup_{i, j} \mathfrak D_{i, j}$, and $\mathfrak M= \sqcup_{i, j} \mathfrak M_{i, j}$ with morphisms $\mathfrak Y \to \mathfrak D \times \mathcal X$, $\mathfrak T \to \mathfrak D$ and $\mathfrak T' \to \mathfrak M$.

For each $\mathcal{Y}^{\sigma}_{i,j} \to \mathcal{T}_{i,j}^{\sigma}$, the varieties described by Theorem \ref{theo:ainvboundedandtwists}.(3).(b) are parametrized by the $K(B)$-points of $\mathcal T_{i, j,\eta}^\sigma$, or equivalently, by the closed points of $\Sec(\mathcal{T}^{\sigma}_{i,j}/B)$.  
We consider the relative space of sections $\mathfrak{S} = \Sec_{\mathfrak D} (\mathfrak T/B) = \Sec_{\mathfrak M}(\mathfrak T'/B) \times_{\mathfrak M} \mathfrak D$.  We first shrink $\mathfrak{S}$ so that the generic point of every section parametrized by $\mathfrak{S}$ is contained in the locus in $\mathfrak{T}$ over which $\mathfrak{Y} \to \mathfrak{T}$ is smooth.  
Then by taking a base change of $\mathfrak{Y} \to \mathfrak{T}$ over the evaluation map $\mathfrak{S} \times B \to \mathfrak{T}$, we obtain a morphism $\mathfrak{F}' \to \mathfrak{S} \times B$ whose fibers are closed subvarieties of the various $\mathcal{Y}^{\sigma}_{i,j}$. Note that we have a morphism $\mathfrak S \to \Sec(\sqcup_i \mathcal W_i/B)$, and we replace $\mathfrak S$ by the fiber product $\mathfrak S \times_{\Sec(\sqcup_i \mathcal W_i/B)} \mathfrak W$ so that we have a compatible morphism $\mathfrak S \to \mathfrak W$. By repeatedly stratifying $\mathfrak{S}$ into locally closed subsets, throwing away components whose intersection with a fiber $\mathfrak Y_s$ does not dominate $B$, taking resolutions of irreducible components of $\mathfrak{F}'$, and taking Stein factorizations, we obtain a commuting diagram
\begin{equation*}
\xymatrix{ \mathfrak{F} \ar[r] \ar[d] &  \mathfrak{Z} \ar[d] \\
\mathfrak{S} \times B \ar[r] & \mathfrak{W} \times B
}
\end{equation*}
where for every closed point $s \in \mathfrak{S}$ the fiber $\mathfrak F_s$ is a normal projective $B$-variety such that $\mathfrak F_s \to B$ has connected fibers and $\mathfrak F_s \to \mathfrak Z_w$ is a finite morphism where $w$ denotes the image of $s$ in $\mathfrak{W}$. After taking a stratification of $\mathfrak S$, we may assume that $\mathfrak F \to \mathfrak S$ is a flat family.  Altogether, we have constructed a family $\mathfrak{F} \to \mathfrak{S} \times B$ whose base is a countable union of finite type schemes and a morphism $g: \mathfrak{F} \to \mathfrak S \times \mathcal{X}$ such that
\begin{enumerate}
\item for every closed point $s \in \mathfrak{S}$ the fiber $\mathfrak{F}_{s}$ is a normal projective $B$-variety such that $\mathfrak F_s \to B$ has connected fibers;
\item for every closed point $s \in \mathfrak{S}$ the map $g_{s}: \mathfrak{F}_{s} \to \mathcal{X}$ is a $B$-morphism that is generically finite onto its image and the corresponding morphism $\mathfrak F_s \to \mathfrak Z_w$ is a finite morphism;
\item for every closed point $s \in \mathfrak{S}$ we have $a(\mathfrak{F}_{s,\eta},g_{s}^{*}L|_{\mathfrak{F}_{s,\eta}}) = a$ and $(\mathfrak{F}_{s,\eta},g_{s}^{*}L|_{\mathfrak{F}_{s,\eta}})$ is adjoint rigid, 
\item if $\mathcal{Y}$ is a good fibration over $B$ and $f: \mathcal{Y} \to \mathcal{X}$ is a generically finite $B$-morphism such that  $a(\mathcal{Y}_{\eta},f^{*}L|_{\mathcal{Y}_{\eta}}) = a$ and $(\mathcal{Y}_{\eta},f^{*}L|_{\mathcal{Y}_{\eta}})$ is adjoint rigid, either the map $f$ is birationally equivalent to $g_{s}$ for some closed point $s$ in our family or $f(\mathcal{Y}_{\eta}) \subset \mathcal V_{\eta}$.
\end{enumerate}
We also have a family $\mathfrak Y \to \mathfrak D \times \mathcal X$ parametrizing integral models $h^\sigma_{i, j}:\mathcal Y^\sigma_{i, j} \to \mathcal U_i$ of twists $h^\sigma_{i, j, \eta}: \mathcal Y^\sigma_{i, j, \eta} \to \mathcal U_{i, \eta}$. Note that all such twists $h^\sigma_{i, j}:\mathcal Y^\sigma_{i, j} \to \mathcal U_i$ factor through $\widehat{h}_{i, j} : \mathcal Y^\sigma_{i, j} \to \mathcal P_{i, j}$ and that $\widehat{h}_{i, j} : \mathcal Y^\sigma_{i, j} \to \mathcal P_{i, j}$ is Galois.
\end{construction}

We will also need two additional lemmas.

\begin{lemma} \label{lemm:deforminghnfreetonearbyfibers}
Suppose that $\mathfrak{Y} \to \mathfrak{S}\times B$ is a family of good fibrations over $B$ with $\mathfrak{S}$ irreducible.  Suppose that for some closed point $s \in \mathfrak{S}$ we have an HN-free section $C$ of $\mathfrak{Y}_{s}/B$.  Then the deformations of $C$ form a dominant family on $\mathfrak{Y}$.
\end{lemma}

\begin{proof}
Let $M$ denote the space of deformations of $C$ in $\mathfrak{Y}$ and for a closed point $s' \in \mathfrak{S}$ let $M_{s'}$ denote the sublocus parametrizing spaces of sections of $\mathfrak{Y}_{s'}/B$.  Since $H^{1}(C,T_{\mathfrak{Y}_{s}/B}|_{C}) = 0$, \cite[Theorem I.2.15.(2)]{Kollar} shows that
\begin{equation*}
\dim(M) \geq \dim(M_{s}) + \dim(\mathfrak{S}).
\end{equation*}
Since $C$ is an HN-free section in $\mathfrak{Y}_{s}$, by replacing $C$ by a general deformation we may ensure that it avoids any codimension $2$ locus of $\mathfrak{Y}_{s}$.  We conclude that $T_{\mathfrak{Y}/\mathfrak{S} \times B}|_{C}$ is locally free, and thus the restriction of this sheaf to a general deformation of $C$ is also locally free.
As the universal family over $\Sec(\mathfrak{Y}/B)$ is smooth, Lemma \ref{lemm:minslopelower} shows that the minimal slope of a quotient of $T_{\mathfrak{Y}/\mathfrak{S}\times B}|_{C'}$ is a lower semicontinuous function as we vary $C'$ in $\Sec(\mathfrak{Y}/B)$.  Thus there is an open subset of $M$ parametrizing sections which are HN-free in their fiber. 
For a general section $C'$ parametrized by $M$ denote by $s'$ the closed point of $\mathfrak{S}$ parametrizing the good fibration $\mathfrak{Y}_{s'} \to B$ containing $C'$.  As explained above $C'$ is HN-free in $\mathfrak{Y}_{s'}$, and in particular for such points $s'$ we have $\dim(M_{s'}) = \dim(M_{s})$.  Thus the dimension computation shows that sections parametrized by $M$ must dominate $\mathfrak{Y}$.
\end{proof}

The next lemma is well-known, but we include a proof for completeness.

\begin{lemma}
\label{lemma:uniquenessoftwists}
Let $k$ be a field of characteristic $0$ and let $f : Y \to X$ be a dominant generically finite morphism between normal projective varieties defined over $k$. Assume that
\[
\mathrm{Bir}(\overline{Y}/\overline{X}) = \mathrm{Aut}(\overline{Y}/\overline{X}).
\]
Let $X^\circ \subset X$ be a Zariski open subset such that $f|_{f^{-1}(X^\circ)} : f^{-1}(X^\circ) \to X^{\circ}$ is \'etale.
Let $f^\sigma : Y^\sigma \to X$ be a twist of $f$ over $X$ and suppose that there are $k$-rational points $p \in Y(k)$ and $p^{\sigma} \in Y^{\sigma}(k)$ which define the same geometric point on $f^{-1}(X^\circ)_{\overline{k}}$. Then $f$ and $f^\sigma$ are isomorphic as $X$-schemes.
\end{lemma}

\begin{proof}
Let $x$ be a common rational point on $f^{-1}(X^{\circ})$ and $(f^\sigma)^{-1}(X^{\circ})$. We prove that 
$\pi_1^{\text{\'et}}(f^{-1}(X^{\circ}), x)$ and $\pi_1^{\text{\'et}}((f^\sigma)^{-1}(X^{\circ}), x)$ define the same subgroup in $\pi^{\text{\'et}}_1(X^{\circ}, f(x))$. Our assertion follows from this fact. We have the following the homotopy exact sequence:
\[
1 \to \pi^{\text{\'et}}_1(\overline{X^{\circ}}, f(x)) \to \pi^{\text{\'et}}_1(X^{\circ}, f(x)) \to \mathrm{Gal}(k) \to 1,
\]
and moreover this exact sequence admits a section $s$ as $f(x)$ is a $k$-rational point. Since we have
\[
\pi_1^{\text{\'et}}(f^{-1}(X^{\circ}), x) \cap \pi^{\text{\'et}}_1(\overline{X^{\circ}}, f(x)) = \pi_1^{\text{\'et}}(\overline{f^{-1}(X^{\circ})}, x) = \pi_1^{\text{\'et}}((f^\sigma)^{-1}(X^{\circ}), x) \cap \pi^{\text{\'et}}_1(\overline{X^{\circ}}, f(x))
\]
and $x$ induces sections of $\pi_1^{\text{\'et}}(f^{-1}(X^{\circ}), x)$ and $\pi_1^{\text{\'et}}((f^\sigma)^{-1}(X^{\circ}), x)$ which are compatible with $s$, we conclude that
\[
\pi_1^{\text{\'et}}(f^{-1}(X^{\circ}), x) =  \pi_1^{\text{\'et}}(\overline{f^{-1}(X^{\circ})}, x) \cdot s(\mathrm{Gal}(k)) = \pi_1^{\text{\'et}}((f^\sigma)^{-1}(X^{\circ}), x).
\]
Thus our assertion follows.
\end{proof}

We are now ready to prove our  main boundedness theorems.  For these results we will be in the situation $a_{rel} = 1$.

\begin{theorem} \label{theo:positivebounded}
Let $\pi: \mathcal{X} \to B$ be a good fibration and let $L$ be a big and semiample Cartier divisor.  Assume that $\mathcal{X}_{\eta}$ is geometrically uniruled.  Set $a = a(\mathcal{X}_{\eta},L|_{\mathcal{X}_{\eta}})$.  Fix a rational number $\beta$.  Fix a positive integer $T$.  Fix a positive integer $b > a$ such that $bL$ defines a basepoint free linear series.  There is:
\begin{itemize}
\item a constant $\xi^{\dagger} = \xi^{\dagger}(\dim(\mathcal{X}), g(B), a, T, \beta,b)$,
\item a closed subset $\mathcal{V} \subset \mathcal{X}$, and
\item a bounded family of smooth projective varieties $q: \widehat{\mathfrak{F}} \to \widehat{\mathfrak{S}}$ equipped with $\widehat{\mathfrak{S}}$-morphisms $p: \widehat{\mathfrak{F}} \to \widehat{\mathfrak{S}} \times B$ and $g: \widehat{\mathfrak{F}} \to \widehat{\mathfrak{S}} \times \mathcal{X}$
\end{itemize}
which have the following properties:
\begin{enumerate}
\item For every closed point $s \in \widehat{\mathfrak{S}}$, $\widehat{\mathfrak F}_s \to B$ is a good fibration. 
\item For every closed point $s \in \widehat{\mathfrak{S}}$ the morphism $g_{s}: \widehat{\mathfrak{F}}_{s} \to \mathcal{X}$ is a $B$-morphism that is  generically finite onto its image.  
\item For every irreducible component $\widehat{\mathfrak{F}}_{i}$ of $\widehat{\mathfrak{F}}$ the composition of $g|_{\widehat{\mathfrak{F}}_{i}}: \widehat{\mathfrak{F}}_{i} \to \widehat{\mathfrak{S}} \times \mathcal{X}$ with the projection $\widehat{\mathfrak{S}} \times \mathcal X \to \mathcal X$ is dominant.  
\item For every closed point $s \in \widehat{\mathfrak{S}}$ we have $a(\widehat{\mathfrak{F}}_{s,\eta},g_{s}^{*}L|_{\widehat{\mathfrak{F}}_{s,\eta}}) = a(\mathcal{X}_{\eta},L|_{\mathcal{X}_{\eta}})$ and $(\widehat{\mathfrak{F}}_{s,\eta},g_{s}^{*}L|_{\widehat{\mathfrak{F}}_{s,\eta}})$ is adjoint rigid.
\item Suppose that $\psi: \mathcal{Y} \to B$ is a good fibration equipped with a $B$-morphism $f: \mathcal{Y} \to \mathcal{X}$ that is generically finite onto its image and satisfies $a(\mathcal{Y}_{\eta},f^{*}L|_{\mathcal{Y}_{\eta}}) \geq a$.  Suppose that $N$ is an irreducible component of $\Sec(\mathcal{Y}/B)$ parametrizing a dominant family of sections $C$ on $\mathcal{Y}$ which satisfy $f^{*}L \cdot C \geq \xi$ and $f^{*}(K_{\mathcal{X}/B} + aL) \cdot C \leq \beta$.  Let $M \subset \Sec(\mathcal{X}/B)$ be the irreducible component containing the pushforward of the sections parametrized by $N$.  Finally, suppose that
\begin{equation*}
\dim(N) \geq \dim(M) - T.
\end{equation*}

For a general section $C$ parametrized by $N$, either:
\begin{itemize}
\item $C$ is contained in $\mathcal{V}$, or
\item there is an irreducible component $\widehat{\mathfrak{F}}_{i}$ of $\widehat{\mathfrak{F}}$ and an irreducible component $N'$ of $\Sec(\widehat{\mathfrak{F}}_{i}/B)$ parametrizing a dominant family of sections on $\widehat{\mathfrak{F}}_{i}$ such that $f(C)$ is the image of a section $C'$ parametrized by $N'$ and if $\widehat{\mathfrak{F}}_{i,s}$ denotes the fiber containing $C'$ then the strict transform of $C'$ in a resolution of $\widehat{\mathfrak{F}}_{i,s}$ is HN-free.
\end{itemize}
\end{enumerate}
\end{theorem}

We will divide the proof into five steps.  Step 1 is devoted to some preliminary work.  In Step 2, we construct a bounded family of varieties $\mathfrak{P} \to \mathfrak{Q}$ such that every $\mathcal{Y}$ as in Theorem \ref{theo:positivebounded}.(5) is a birational to a twist of a fiber over a closed point of $\mathfrak{Q}$.  In Step 3, we bound the invariants from Corollary \ref{coro:boundedintimpliesboundedtwists} for the varieties in our family $\mathfrak{P} \to \mathfrak{Q}$.  Then in Step 4 we use this corollary to construct a bounded family of twists $\widehat{\mathfrak F} \to \widehat{\mathfrak S}$ which carry sections of low $L$-degree.  Finally, in Step 5 we verify that $\widehat{\mathfrak F} \to \widehat{\mathfrak S}$ has the desired properties.

\begin{proof}
\textbf{Step 1:} Let $m$ be the maximum of the degrees of the morphisms $\mathcal Y_{i, j} \to \mathcal U_i$ from Construction~\ref{cons:allfamilies} and set $d = (m + 1)!$.  Let $\sqcup_{G} \mathcal H(G, B)$ be the Hurwitz stack and fix an \'etale covering $\sqcup_{G, |G| \leq d} \mathcal H_G \to \sqcup_{G, |G| \leq d} \mathcal H(G, B)$ by a scheme. We will work over this base for the entire proof.

Note that the divisor $E = 0$ satisfies the condition of Proposition \ref{prop:generalsurjectiontofiber}.  Define $\xi = \xi(\dim(\mathcal{X}), g(B), 0, 1, a, T, \beta,b)$ as in Theorem \ref{theo:generalainvsections}.  We then choose 
$$
\xi^{+} = \xi^{+}(\dim(\mathcal{X}), g(B), 0, 1, a, T, \beta, b) \quad \text{ and } \quad T^{+} = T^{+}(\dim(\mathcal{X}), g(B), 0, 1, a, T, \beta,b)
$$ 
as in Corollary \ref{coro:adjointrigidcodim}.  Define $\daleth = \daleth(\dim(\mathcal{X}), g(B), 1, a, T, \beta,b)$ as in Theorem \ref{theo:adjointrigidbounded}.  Finally we define $\xi^{\dagger} = \sup \left\{ \xi, \xi^{+} \right\}$.

Since $L$ is big and semiample, there is a closed subvariety $\mathcal{V}_{1} \subset \mathcal{X}$ such that the family of subvarieties of $\mathcal{X}$ that are not contained in $\mathcal{V}_{1}$ and have $L$-degree $\leq \daleth$ is bounded. By \cite[Theorem 4.18.(2)]{LST18} there is a closed sublocus $\mathcal{V}_{2,\eta} \subset \mathcal{X}_{\eta}$ that contains all subvarieties with larger generic Fujita invariant and we let $\mathcal{V}_{2}$ denote its closure. Let $\mathcal V_3$ be the exceptional closed set from Construction~\ref{cons:allfamilies}. We start by setting $\mathcal{V}$ to be the union of $\mathcal{V}_{1}$, $\mathcal{V}_{2}$, and $\mathcal V_3$; we will later enlarge it.

Let $\mathfrak{Z}_i \to \mathfrak{W}_i\times B$ be the families in Construction \ref{cons:allsubvarieties}.  Then there is a finite-type subscheme $\mathfrak{R}_i \subset \mathfrak{W}_i$ parametrizing those varieties whose images in $\mathcal{X}$ have $L$-degree $\leq \daleth$ and are not contained in $\mathcal V$. Indeed, this follows from the fact that being the integral main component of the pullback of a section in $\mathfrak{W}_i$ is an open condition so that perhaps after taking a finer stratification of $\mathfrak{W}_i$ there is a bijective map from $\mathfrak{W}_i$ to a Zariski open subset of a component of the Hilbert scheme parametrizing integral main components of pullbacks of sections. We denote by $\mathfrak{Z}_{i, \mathfrak R_i} \to \mathfrak{R}_i$ the universal family over $\mathfrak{R}_i$. Set $\mathfrak R = \sqcup_i \mathfrak R_i$.

Let $\mathfrak{Y} \to \mathfrak{T} \to \mathfrak{D}$, $\mathfrak T' \to \mathfrak M$, $\mathfrak{F} \to \mathfrak{S}$ and $g: \mathfrak{F} \to \mathcal{X} \times B$ be defined as in Construction \ref{cons:allfamilies}. 
For any closed point $s \in \mathfrak{S}$ the map $g_{s}: \mathfrak{F}_{s} \to \mathcal{X}$ has image that is birationally equivalent to a fiber $\mathcal{Z}_{s}$ of $\mathfrak{Z} \to \mathfrak{W}$.   By \cite[Lemma 4.7]{LST18} we have 
\begin{equation*}
a(\mathfrak{F}_{s,\eta},g_s^{*}L|_{\mathfrak{F}_{s,\eta}}) \leq a(\mathcal{Z}_{s,\eta},\iota_{s}^*L|_{\mathcal{Z}_{s,\eta}})
\end{equation*}
and if equality is achieved then \cite[Lemma 4.9]{LST18} shows that $(\mathcal{Z}_{s,\eta},\iota_{s}^*L|_{\mathcal{Z}_{s,\eta}})$ is adjoint rigid.  If we shrink $\mathfrak{S}$ to remove all maps $g_{s}$ whose image lies in $\mathcal{V}$, then we will always have equality of Fujita invariants.  We let $\mathfrak{S}'$ denote the sublocus of $\mathfrak{S}$ consisting of maps $g_{s}$ whose image is a member of our fixed bounded family $\mathfrak{Z}_{\mathfrak R} \to \mathfrak{R}$ and denote by $\mathfrak{F}' \to \mathfrak{S}'$ the corresponding family.

\textbf{Step 2:}  We next claim that there is a morphism $\mathfrak{Q}\to \mathfrak S'\subset \mathfrak S$ such that $\mathfrak Q$ is of finite type over $\mathbb C$ and for every map $g_{s}$ parametrized by $\mathfrak{S}'$ the map $g_{s,\eta}$ is a twist of the generic fiber of a map parametrized by $\mathfrak{Q}$.  Indeed, recall from Theorem \ref{theo:ainvboundedandtwists} that we have a finite number of smooth projective varieties $\mathcal{Y}_{i,j,\eta}$ equipped with morphisms
\begin{equation*}
\xymatrix{ {\mathcal Y}_{i, j,\eta} \ar[r]^{h_{i, j,\eta}} \ar[d]_{r_{i, j,\eta}} &  {\mathcal U}_{i,\eta} \ar[d]_{p_{i,\eta}} \\
\mathcal T_{i, j,\eta} \ar[r]_{t_{i, j,\eta}} & \mathcal W_{i,\eta}}
\end{equation*}
where $t_{i,j,\eta}$ is Galois.

It follows from our construction that $\mathfrak R_i$ is contained in finitely many irreducible components of $\Sec(\mathcal W_i/B)$. 
Let $\Sec(\mathcal W_i/B, \mathcal B')$ denote the space of sections not contained in $\mathcal B'$ and define $\Sec(\mathcal S_{i, j}/B, a_{i,j}^{-1}(\mathcal B'))$ analogously.
Then $a_{i, j, *}: \Sec(\mathcal S_{i, j}/B, a_{i,j}^{-1}(\mathcal B')) \to  \Sec(\mathcal W_i/B, \mathcal B')$ is of finite type, so the fiber product
\[
\widehat{\mathfrak R}_{i, j}:= \Sec(\mathcal S_{i, j}/B, a_{i,j}^{-1}(\mathcal B'))\times_{\Sec(\mathcal W_i/B, \mathcal B')} \mathfrak R_i
\]
is of finite type over $\mathbb C$.
Moreover for any $C \in \widehat{\mathfrak R}_{i, j}$, the fiber $b_{i, j}^{-1}(C_\eta)$ is geometrically rationally connected so $\widehat{\mathfrak R}_{i, j}$ is in the image of $\Sec(\mathcal P_{i, j}/B, \ell_{i,j}^{-1}(\mathcal B)) \to \Sec( \mathcal S_{i, j}/B, a_{i,j}^{-1}(\mathcal B'))$.
Thus there is a finite disjoint union of locally closed subschemes of finite type $\widetilde{\mathfrak R}_{i, j} \subset \Sec(\mathcal P_{i, j}/B, \ell_{i,j}^{-1}(\mathcal B))$ with a surjective morphism $\widetilde{\mathfrak R}_{i, j} \to  \widehat{\mathfrak R}_{i, j}$. 
We denote the base change of $\mathfrak Z_{i, \mathfrak R_{i}} \to \mathfrak R_{i}$ over $\widetilde{\mathfrak R}_{i, j} \to  \widehat{\mathfrak R}_{i, j} \to \mathfrak R_i$ by $\mathfrak Z_{\widetilde{\mathfrak R}_{i, j}} \to \widetilde{\mathfrak R}_{i, j}$.

Recall that we are working over $\sqcup_{G, |G| \leq d} \mathcal H_G$. We claim that every twist of $\mathcal Y_{i, j, \eta}/ \mathcal U_{i, \eta}$ splits over an extension $K(B')/K(B)$ of degree $\leq d$.  Indeed, if we denote $\mathrm{Aut}( \mathcal Y_{i, j, \overline{\eta}}/\mathcal U_{i, \overline{\eta}})$ by $G$, then Lemma~\ref{lemm:boundingdeganddisc} shows that one may use a Galois base change of degree $\leq |G| \cdot \# \mathrm{Aut}(G)$. In particular $d = (m+1)!$ gives an upper bound on this degree.

Since each rational point of $\mathcal{W}_{i,\eta}$ not contained in $\mathcal{B'}$ will lift to a unique twist of $\mathcal{T}_{i,j,\eta}$ and the number of $1$-cycles representing the same Galois cohomology class is at most $m$, the pushforward morphism
\[
t_{i, j, *} : \Sec_{\mathfrak M_{i, j}}(\mathfrak T'_{i, j}/B, t_{i, j}^{-1}(\mathcal B')) \to \Sec(\mathcal W_i/B, \mathcal B') \times (\sqcup_{G, |G|\leq d} \mathcal H_G)
\]
is a quasi-finite morphism onto its image of degree at most $m^2$.

Now for each $C \in \Sec(\mathcal P_{i, j}/B, \ell_{i, j}^{-1}(\mathcal B))$, $\widehat{h}_{i, j}^{-1}(C)$ decomposes into a union of curves which are Galois conjugate to each other over $B$ where $\widehat{h}_{i, j} : \mathcal Y_{i, j} \to \mathcal P_{i, j}$ is the morphism constructed in Construction~\ref{cons:allfamilies}. Then after taking a stratification of $\widetilde{\mathfrak R}_{i, j}$ by locally closed subsets and replacing $\widetilde{\mathfrak R}_{i, j}$ by an \'etale cover, the universal property of the Hurwitz stack yields a morphism
\[
\psi_{i, j}: \widetilde{\mathfrak R}_{i, j}   \to \sqcup_{G, |G|\leq d} \mathcal H(G, B)
\]
that sends a section $C$ of $\mathcal{P}_{i,j}/B$ to the cover $C' \to B$ obtained by normalizing an irreducible component of $\widehat{h}_{i, j}^{-1}(C)$. We denote by $\mathfrak R_{i, j}$ the fiber product
\[
\widetilde{\mathfrak R}_{i, j}\times_{\sqcup_{G, |G|\leq d} \mathcal H(G, B)} (\sqcup_{G, |G|\leq d} \mathcal H_G)
\]
which is of finite type over $\mathbb C$.
Thus using the morphism $\mathfrak R_{i, j} \to \Sec(\mathcal W_i/B, \mathcal B_{i, j}) \times (\sqcup_{G, |G|\leq d} \mathcal H_G)$ we define the scheme 
\[
\mathfrak Q'_{i, j} = \Sec_{\mathfrak M_{i, j}}(\mathfrak T'_{i, j}/B, t_{i, j}^{-1}(\mathcal B_{i, j}))\times_{\Sec(\mathcal W_i/B, \mathcal B_{i, j}) \times (\sqcup_{G, |G|\leq d} \mathcal H_G)} \mathfrak R_{i, j}
\]
which is a finite type scheme over $\mathbb C$. Note that $\mathfrak{Q}'_{i,j}$ parametrizes the sections of $\mathcal T_{i,j}^{\sigma}/B$ which map to $\mathfrak{R}_{i}$ such that the twist $\mathcal T_{i,j}^{\sigma}$ is trivialized by a base change $C' \to B$ coming from $\widehat{h}_{i,j}^{-1}(C)$ as constructed above.
Let $\mathfrak Q' = \sqcup_{i, j} \mathfrak Q'_{i, j}$.
Then we have a morphism $\mathfrak Q' \to \Sec_{\mathfrak M}(\mathfrak T'/B)$ and $\mathfrak S' \to \Sec_{\mathfrak D}(\mathfrak T/B) \to \Sec_{\mathfrak M}(\mathfrak T'/B)$. 

We set $\mathfrak Q = \mathfrak Q' \times_{\Sec_{\mathfrak M}(\mathfrak T'/B)} \mathfrak S'$.
Since $\mathfrak S' \to \Sec_{\mathfrak M}(\mathfrak T'/B)$ is of finite type over each $\mathcal H_G$, $\mathfrak Q$ is a scheme of finite type over $\mathbb C$. 
Then we denote the base change of $\mathfrak Y \to \mathfrak D$ over $\mathfrak Q \to \mathfrak D$ by $\widehat{\mathfrak Y} \to \mathfrak Q$ and we denote the base change of $\mathfrak F \to \mathfrak S$ over $\mathfrak Q \to \mathfrak S' \hookrightarrow \mathfrak S$ by $\mathfrak P \to \mathfrak Q$.

We still must verify that $\mathfrak{P} \to \mathfrak{Q}$ satisfies the claimed property.
Let $s \in \mathfrak S'$  and consider the corresponding $g_s :  \mathfrak F_s' \to \mathfrak Z_r \to \mathcal X$ with $r \in \mathfrak R$. By Theorem~\ref{theo:ainvboundedandtwists}, $g_{s, \eta} : \mathfrak F_{s, \eta}' \to \mathcal X_\eta$ is birationally equivalent to the map to $\mathcal{X}_{\eta}$ from a fiber of a twist $\mathcal Y^\sigma_{i, j, \eta} \to \mathcal T^\sigma_{i, j,\eta}$ for some $i, j$. Then since $\mathcal Y^{\sigma}_{i, j}$ factors through $\mathcal P_{i, j}$, by the construction we find a point $\widehat{r} \in {\mathfrak R}_{i, j}$ mapping to $r$. This point $\widehat{r}$ specifies a point $(B'/B) \in \sqcup_{G, |G|\leq d} \mathcal H_G$. On the other hand one can find a twist $\mathcal Y^{\tau}_{i, j, \eta} \to \mathcal T^\sigma_{i, j,\eta}$ such that the preimage $(\widehat{h}^\tau_{i, j, \eta})^{-1}(C_\eta)$ of the section $C \in \Sec(\mathcal P_{i, j}/B)$ corresponding to $\widehat{r}$ consists entirely of $K(B)$-rational points. Such a twist will be trivialized by the base change $B'\to B$. This means that $r$ is in the image of the map $\mathfrak Q = \mathfrak Q' \times_{\Sec_{\mathfrak M}(\mathfrak T'/B)} \mathfrak S' \to \Sec(\sqcup_i \mathcal W_i/B, \mathcal B')$.  Thus there is a point $q \in \mathfrak Q$ mapping to $r$ such that $\mathfrak{F}'_{s,\eta} \to \mathcal{X}_{\eta}$ is a twist of $\mathfrak P_{q, \eta} \to \mathcal X_\eta$ for the fiber $\mathfrak{P}_{q}$ over $q$.  This finishes the verification of the desired property.

\textbf{Step 3:} Note that the degree of $h^\sigma_{q,\eta}: \mathfrak{P}^{\sigma}_{q,\eta} \to \mathcal{Z}_{q,\eta}$ is bounded by the maximum of the degrees of $h_{i, j} : \mathcal Y_{i, j} \to \mathcal U_i$. In particular the size of $\mathrm{Aut}( \mathfrak{P}^{\sigma}_{q,\eta} / \mathcal{Z}_{q,\eta})$ is uniformly bounded by the integer $m$.  It follows from Lemma~\ref{lemm:boundingdeganddisc} that for every closed point $q \in \mathfrak{Q}$ the map $h^\sigma_{q,\eta}: \mathfrak{P}^{\sigma}_{q,\eta} \to \mathcal{Z}_{q,\eta}$ becomes isomorphic to $h_{q,\eta}: \mathfrak{P}_{q,\eta} \to \mathcal{Z}_{q,\eta}$ after a Galois base change $\widetilde{B} \to B$ of degree $\leq d$.

Next we define an integer $t$ by using the family $\mathfrak{P} \to \mathfrak{Q}$.  
Since normality is a constructible property in proper families (\cite[Th\'eor\`eme 12.2.4]{EGAIV3}), by Noetherian induction there is a positive integer $t_{1}$ that bounds the number of non-normal fibers of $\mathfrak{P}_{q} \to B$ as we vary over all $q \in \mathfrak{Q}$.  Since the relative automorphism scheme $\Aut_{B}(\mathfrak{P}_{q}/\mathcal{Z}_{q})$ is quasifinite over $\mathfrak{Q} \times B$ and since flatness is a constructible property, as we vary over all closed points $q \in \mathfrak{Q}$ there is a positive integer $t_{2}$ that bounds the number of places in $B$ where the restriction of $\Aut_{B}(\mathfrak{P}_{q}/\mathcal{Z}_{q})$ to $\{q \} \times B$ is not flat.  We set $t = t_{1}+t_{2}$.

\textbf{Step 4:} Lemma \ref{lemm:twistboundedconditions2} and Corollary \ref{coro:boundedintimpliesboundedtwists}  show that as we vary the closed point $q \in \mathfrak{Q}$ the set of twists of $h_{q}: \mathfrak{P}_{q} \to \mathcal{Z}_{q}$ which are trivialized by a base change $B' \to B$ of degree at most $d$ and with at most $t+d(T+T^{+})$ branch points is parametrized by a bounded family.  We denote by $\widetilde{\mathfrak{F}}\to \widetilde{\mathfrak{S}}$ the bounded subfamily of $\mathfrak{F}' \to \mathfrak{S}'$ parametrizing maps $g_{s}: \mathfrak{F}_{s} \to \mathcal{Z}_{s}$ satisfying these properties. 
After taking smooth resolutions and stratifying the base, we obtain $\widehat{\mathfrak F} \to \widehat{\mathfrak S}$ such that each fiber is a good fibration over $B$.
We then shrink $\widehat{\mathfrak{S}}$ by removing all irreducible components $\mathfrak{S}_{j}$ such that the corresponding family $\widehat{\mathfrak{F}}_{j}$ fails to dominate $\mathcal{X}$ and we enlarge $\mathcal{V}$ by taking the union with the closures of the images of these families.

\textbf{Step 5:}  We are now ready to verify the desired properties of $\widehat{\mathfrak{F}} \to \widehat{\mathfrak{S}}$.
 Properties (1)-(4)  follow from the construction, and we only need to check (5).  Suppose $f: \mathcal{Y} \to \mathcal{X}$ is as in (5).  Applying Theorem \ref{theo:generalainvsections} with $a_{rel} = 1$ we see that $a(\mathcal{Y}_{\eta},f^{*}L|_{\mathcal{Y}_{\eta}}) \geq a$.  If we have a strict inequality then $f(\mathcal{Y}) \subset \mathcal{V}$ and so the sections on $\mathcal{Y}$ are accounted for by $\mathcal{V}$.  From now on we assume that $f(\mathcal{Y}) \not \subset \mathcal{V}$ which implies that $a(\mathcal{Y}_{\eta},f^{*}L|_{\mathcal{Y}_{\eta}}) = a$.

We apply Corollary \ref{coro:adjointrigidcodim} to construct subvarieties on the Stein factorization of the evaluation map over $\mathcal{Y}$ and then take images in $\mathcal{Y}$ to obtain a dominant family of subvarieties $\mathcal{F} \subset \mathcal{Y}$ satisfying:
\begin{itemize}
\item the codimension in $N$ of the space of deformations of $C$ in $\mathcal{F}$ is at most $T^{+}$,
\item the strict transform of $C$ to a resolution of $\mathcal{F}$ is HN-free,
\item$(\mathcal{F}_{\eta},f^{*}L|_{\mathcal{F}_{\eta}})$ is adjoint rigid.
\end{itemize}
We will show that the conclusion of (5) holds for the sections on the general subvariety $\mathcal{F}$ in our family.

Consider a general subvariety $\mathcal{F}$ in our family and set $\mathcal{Z} = f(\mathcal{F})$.  Since it is not possible for $\mathcal{Z}_{\eta}$ to have larger $a$-invariant (as it is not contained in $\mathcal V_\eta$), we have $a(\mathcal{Z}_{\eta},L|_{\mathcal{Z}_{\eta}}) = a$ and thus by \cite[Lemma 4.9]{LST18} $(\mathcal{Z}_{\eta},L|_{\mathcal{Z}_{\eta}})$ is adjoint rigid.  Theorem \ref{theo:adjointrigidbounded}.(2) shows that $\mathcal{Z}$ is birationally equivalent to a smooth $\mathcal{Z}'$ that is parametrized by the bounded family $\mathfrak{Z}_{\mathfrak{R}} \to \mathfrak{R}$.  In particular the map $\mu: \mathcal{F}_{\eta} \to \mathcal{Z}_{\eta}$ is birationally equivalent to a twist of $h_{q}: \mathfrak{P}_{q, \eta} \to \mathcal{Z}_{q, \eta}$ for some closed point $q \in \mathfrak{Q}$.

Choose a morphism $\mu': \mathcal{F}' \to \mathcal{Z}'$ birationally equivalent to $\mu$ where $\mathcal{F}'$ is smooth.  Let $N_{\mathcal{F}}$ denote the moduli space of deformations of the strict transforms $C'$ of $C$ in $\mathcal{F}'$ and let $M_{\mathcal{Z}}$ denote the moduli space of deformations of the image in $\mathcal{Z}'$.  Then
\begin{align*}
\dim(M_{\mathcal{Z}}) - \dim(N_{\mathcal{F}}) \leq \dim(M) - (\dim(N) - T^{+})
\leq T + T^{+}
\end{align*}
so that $N_{\mathcal{F}}$ has codimension at most $T^{+} + T$ in $M_{\mathcal{Z}}$.  Then we have
\begin{align*}
-K_{\mathcal{F}'/B} \cdot C' + (\dim(\mathcal{F}')-1)(1-g(B)) & = \dim(N_{\mathcal{F}}) \\
& \geq \dim(M_{\mathcal{Z}}) - T^{+} - T \\
& \geq -K_{\mathcal{Z}'/B} \cdot \mu'_{*}C' + (\dim(\mathcal{Z}')-1)(1-g(B)) - T^{+} - T
\end{align*}
which rearranges to
\begin{equation*}
(K_{\mathcal{F}'/B} - \mu^{*}K_{\mathcal{Z}'/B}) \cdot C' \leq T^{+} + T.
\end{equation*}
This intersection bound and Corollary~\ref{coro:boundedintimpliesboundedtwists} imply that $\mu': \mathcal{F}' \to \mathcal{Z}'$ is birationally equivalent to a twist of $h_{q}$ that is trivialized by a base change $B' \to B$ that has degree at most $d$ and has at most $t + d(T+T^{+})$ branch points.  Thus $\mu'$ is birationally equivalent to one of the maps $h_{s}: \widetilde{\mathfrak{F}}_{s} \to \mathcal{Z}_{s}$ parametrized by our bounded family $\widetilde{\mathfrak{F}} \to \widetilde{\mathfrak{S}}$.

Consider the strict transform of our family of sections in the fiber $\widehat{\mathfrak{F}}_{s}$.  Since these sections go through at least $2g(B)+1$ general points of $\widehat{\mathfrak{F}}_{s}$, they are HN-free in this fiber.  Lemma \ref{lemm:deforminghnfreetonearbyfibers} shows that the sections deform to give a dominant family on the entire irreducible component $\widehat{\mathfrak{F}}_{i}$ containing $\widehat{\mathfrak{F}}_{s}$.  Since by construction every irreducible component of $\widehat{\mathfrak{F}}$ dominates $\mathcal{X}$, we deduce that the family of sections gives a dominant family on $\mathcal{X}$.  Furthermore, we see that the general section parametrized by $N_{\mathcal{F}}$ is in the image of the map $\Sec(\widehat{\mathfrak{F}}_{i}/B) \to \Sec(\mathcal{X}/B)$.  Thus the same property is true for $N$, proving (5).
\end{proof}

Our next boundedness statement is closer in spirit to the results of \cite{LST18}: instead of using a bounded family $\widehat{\mathfrak{F}} \to \widehat{\mathfrak{S}}$ such that the fibers $\widehat{\mathfrak{F}}_{s,\eta}$ are adjoint rigid, one can instead use a bounded family $\widetilde{\mathfrak{Y}} \to \widetilde{\mathfrak{S}}$ such that the fibers $\widetilde{\mathfrak{Y}}_{s,\eta}$ are twists of the finite set of universal families constructed in Theorem \ref{theo:ainvboundedandtwists}.

\begin{theorem} \label{theo:bounded_bigandsemiample}
Let $\pi: \mathcal{X} \to B$ be a good fibration and let $L$ be a big and semiample Cartier divisor.  Assume that $\mathcal{X}_{\eta}$ is geometrically uniruled.  Set $a = a(\mathcal{X}_{\eta},L|_{\mathcal{X}_{\eta}})$.  Fix a constant $\beta$. Fix a positive integer $b>a$ such that $bL$ defines a basepoint free linear series.

There is a constant $\xi^{\dagger} = \xi^{\dagger}(\dim(\mathcal{X}), g(B), a, \beta,b)$, a proper closed subset $\mathcal{V} \subset \mathcal{X}$, and a bounded family $\widetilde{\mathfrak Y} \to \widetilde{\mathfrak S} \times B$ of good fibrations equipped with a $\widetilde{\mathfrak{S}} \times B$-morphism $\widetilde{f} : \widetilde{\mathfrak Y} \to \widetilde{\mathfrak S} \times \mathcal X$ such that:
\begin{enumerate}
\item for every closed point $s \in \widetilde{\mathfrak{S}}$ the map $\widetilde{f}_{s}$ is dominant and generically finite but not birational;
\item for every closed point $s \in \widetilde{\mathfrak{S}}$ we have $a(\widetilde{\mathfrak{Y}}_{s,\eta}, -\widetilde{f}_{s}^*K_{\mathcal X/B}|_{\widetilde{\mathfrak{Y}}_{s,\eta}}) = a(\mathcal{X}_{\eta},-K_{\mathcal{X}/B}|_{\mathcal X_\eta})$;
\item as we vary over all closed points $s \in \widetilde{\mathfrak{S}}$ the set of birational equivalence classes of the maps $\{ \widetilde{f}_{s,\overline{\eta}}: \widetilde{\mathfrak{Y}}_{s,\overline{\eta}} \to \mathcal{X}_{\overline{\eta}} \}$ obtained by base changing to $\Spec(\overline{K(B)})$ is finite;
\item if $M \subset \Sec(\mathcal{X}/B)$ is an irreducible component  that generically parametrizes non-HN-free sections $C$ with $L \cdot C \geq \xi^\dagger$ and $(K_{\mathcal{X}/B} + a(\mathcal{X}_{\eta},L|_{\mathcal{X}_{\eta}})L) \cdot C \leq \beta$ then a general section $C$ parametrized by $M$ satisfies either $C \subset \mathcal{V}$ or $C \in \widetilde{f}_{*}(\Sec(\widetilde{\mathfrak{Y}}_{s}/B))$ for some closed point $s \in \widetilde{\mathfrak{S}}$.
\end{enumerate}
\end{theorem}

\begin{proof}
We begin by making exactly the same constructions and definitions as in the proof of Theorem \ref{theo:positivebounded}; we continue from the end of this proof.  Additionally we set $T=0$.

Let $\widetilde{f} : \widetilde{\mathfrak{Y}} \to \widetilde{\mathfrak{S}} \times \mathcal X$ be the family of twists of $h_{i, j} : \mathcal Y_{i, j} \to \mathcal U_i$ which becomes isomorphic to a member of $\widehat{\mathfrak Y} \to \mathfrak Q$ by a finite base change $B'\to B$ of degree $\leq d$ and with at most $t +dT^+$ branch points.  By Lemma \ref{lemm:twistboundedconditions2} $\widetilde{\mathfrak{S}}$ has finite type over $\mathbb{C}$.

Properties (1), (2), (3) follow from the construction and we only need to verify (4). Suppose $M \subset \Sec(\mathcal{X}/B)$ is an irreducible component that generically parametrizes non-relatively free sections $C$ with $-K_{\mathcal{X}/B} \cdot C \geq \xi$ and $(K_{\mathcal{X}/B} + a(\mathcal{X}_{\eta},L|_{\mathcal{X}_{\eta}})L) \cdot C \leq \beta$.  We may assume that $M$ generically parametrizes sections which are not contained in $\mathcal{V}$. Then $M$ parametrizes a dominant family of sections due to Theorem \ref{theo:positivebounded}.(5). Let $\mathcal Y \to \mathcal X$ be the finite part of the Stein factorization for the evaluation map for $M$.

Applying Corollary \ref{coro:adjointrigidcodim}, we find a dominant family of subvarieties $\mathcal{F} \subset \mathcal{Y}$ satisfying:
\begin{itemize}
\item the codimension in $N$ of the space of deformations of $C$ in $\mathcal{F}$ is at most $T^{+}$,
\item the strict transform of $C$ is HN-free in a resolution of $\mathcal{F}$,
\item$(\mathcal{F}_{\eta},f^{*}L|_{\mathcal{F}_{\eta}})$ satisfies $a(\mathcal X_\eta, L) = a(\mathcal{F}_{\eta},f^{*}L|_{\mathcal{F}_{\eta}})$ and is adjoint rigid.
\end{itemize}
Using the universal property described in Theorem~\ref{theo:ainvboundedandtwists}, there exists a twist $\mathcal Y_{i, j}^\sigma \to \mathcal T_{i, j}^\sigma$ over $\mathcal U_i$ such that $\mathcal F$ is birational to the main component $\mathcal F'_C$ of the preimage of a section $C$ under the map $\mathcal Y_{i, j}^\sigma \to \mathcal T_{i, j}^\sigma$. Then note that $\mathcal Y_{i, j}^\sigma \to \mathcal U_i$ factors through $\mathcal P_{i, j} \to \mathcal U_i$.

We claim that there is some closed point $q \in \mathfrak{Q}$ such that there is an $\mathcal X_{\overline{\eta}}$-isomorphism between $\mathcal Y_{i, j, \overline{\eta}}^\sigma$ and $\widehat{\mathfrak Y}_{q,\overline{\eta}} $ which maps $\mathcal F_{C, \overline{\eta}} $ to the image $\mathfrak P'_{q,\overline{\eta}}$ of $\mathfrak P_{q,\overline{\eta}}$ under the map $\mathfrak P_{q,\overline{\eta}} \to \widehat{\mathfrak Y}_{q,\overline{\eta}}$.
Indeed, by the defining property of $\mathfrak{P} \to \mathfrak{Q}$ we know that $\mathcal F_{C, \eta}$ is a twist of $\mathfrak P_{q',\eta}$ for some $q' \in \mathfrak Q$.  The point $q'$ specifies a twist $\mathcal Y^\tau_{i, j} \to \mathcal T^\tau_{i, j}$ and a point $q''$ on $\Sec(\mathcal T_{i, j}^\tau/B, t_{i, j}^{-1}(\mathcal B'))$.  We let $p$ denote the rational point on $\Sec(\mathcal W_{i}/B, \mathcal B')$ obtained by taking the image of $q''$ under
\begin{equation*}
\Sec(\mathcal T_{i, j}^\tau/B, t_{i, j}^{-1}(\mathcal B)) \to \Sec(\mathcal W_{i}/B, \mathcal B').
\end{equation*}
Let $C_p$ denote the section of $\mathcal{W}_{i} \to B$ corresponding to $p$.
Since $t^\tau_{i, j}: \mathcal{T}^{\tau}_{i,j} \to \mathcal{W}_{i}$ is Galois, every geometric point in $(t^{\tau}_{i, j})^{-1}(C_{p, \eta})$ is a $K(B)$-rational point on $\mathcal T_{i, j, \eta}^\tau$.  Note that the geometric fiber corresponding to $\mathcal{F}_{C,\eta}$ will lie over one of these points; we replace $q''$ by this point. Moreover by construction the image $\mathfrak Z_r$ of $\mathcal F_C$ has $L$-degree $\leq \daleth$. In particular this point will lift to $q'''\in \mathfrak{Q}'$. Thus we can define a point $q = (q''', s') \in \mathfrak Q' \times_{\Sec_{\mathfrak M}(\mathfrak T'/B)} \mathfrak S' = \mathfrak Q$ where $s' \in \mathfrak S'$ is a point corresponding to $(q'', d')$ with $d'$ is the image of $q'$ via $\mathfrak Q \to \mathfrak D$.  This point $q$ maps to $p$ and  $\mathfrak P_{q,\overline{\eta}}$ is birational to the same geometric fiber of $\widehat{\mathfrak Y}_{q,\overline{\eta}}$ as $\mathcal F_{C, \overline{\eta}}$ in $\mathcal Y_{i, j, \eta}^\sigma$.

By the construction of the families $\widehat{\mathfrak{F}} \to \widehat{\mathfrak{S}}$, $\mathcal F_{C,\eta}$ and $\mathfrak P'_{q,\eta}$ are trivialized by a base change $B' \to B$ of degree $\leq d$ with the number of branch points $\leq t + dT^+$. By Lemma~\ref{lemma:uniquenessoftwists} $\mathcal Y_{i, j, \eta}^\sigma$ and $\widehat{\mathfrak Y}_{q,\eta}$ are trivialized by the same base change. Thus our assertion follows.
\end{proof}

\subsection{General statements} \label{sect:boundedconsequences}

The following proposition will allow us to remove the global positivity assumption in Theorem \ref{theo:bounded_bigandsemiample}.

\begin{proposition} \label{prop:birmodelposL}
Let $\pi: \mathcal{X} \to B$ be a good fibration and let $L$ be a generically relatively ample $\mathbb{Q}$-Cartier divisor.  There is a birational model $\phi: \mathcal{X}^{+} \to \mathcal{X}$ that restricts to an isomorphism of generic fibers over $B$ such that $\mathcal{X}^{+}$ is smooth and there is a $\pi \circ \phi$-vertical effective $\mathbb{Q}$-Cartier divisor $G$ such that $\phi^{*}L+G$ is a big and semiample $\mathbb{Q}$-Cartier divisor. 
\end{proposition}

\begin{proof}
Choose a positive integer $p$ such that $A = pL - K_{\mathcal{X}}$ is generically relatively ample.  Choose $E$ as in Proposition \ref{prop:generalsurjectiontofiber} applied to $A$.  Thus there is an effective $\mathbb{Q}$-divisor $D \sim_{\mathbb{Q}} A + E$ such that $D|_{\mathcal{X}_{\eta}}$ has SNC support and has positive coefficients $<1$.  Let $\psi: \mathcal{X}' \to \mathcal{X}$ denote a log resolution and let $D'$ denote the strict transform of the $\pi$-horizontal components of $D$.  Note that $\psi$ is an isomorphism over $\mathcal{X}_{\eta}$ and so $D'$ and $\psi^{*}D$ only differ by a $\pi$-vertical divisor.  Thus we can choose some positive integer $m$ such that $D' + mF - \psi^{*}D$ is $\mathbb{Q}$-linearly equivalent to an effective $\mathbb{Q}$-Cartier divisor, where $F$ denotes a general fiber of $\mathcal{X}' \to B$.

By passing to a relative canonical model, we obtain a birational map $\rho: \mathcal{X}' \dashrightarrow \widehat{\mathcal{X}}$ such that $\rho_{*}(K_{\mathcal{X}'} + D' + mF)$ is relatively ample.  Note that $\rho$ is an isomorphism along $\mathcal{X}'_{\eta}$ since $(K_{\mathcal{X}'} + D' + mF)|_{\mathcal{X}'_{\eta}}$ was already ample.  By increasing $m$, we can ensure that $\rho_{*}(K_{\mathcal{X}'} + D' + mF)$ is ample.  Let $\mathcal{X}^{+}$ denote a birational model admitting morphisms to $\mathcal{X}$ and to $\widehat{\mathcal{X}}$.  Then the difference between the pullback of $\frac{1}{p}\rho_{*}(K_{\mathcal{X}'} + D' + mF)$ to $\mathcal{X}^{+}$ and the pullback of $L$ to $\mathcal{X}^{+}$ is $\mathbb{Q}$-linearly equivalent to a $\pi$-vertical $\mathbb{Q}$-Cartier divisor $G'$.  Since we may add any fiber of $\mathcal{X}^{+} \to B$ to the pullback of $\frac{1}{p}\rho_{*}(K_{\mathcal{X}'} + D' + mF)$ without affecting semiampleness, we may eliminate the negative part of $G'$ to obtain the desired effective $\pi$-vertical $\mathbb{Q}$-Cartier divisor $G$.
\end{proof}

Putting everything together, we obtain the following variant of Theorem \ref{theo:positivebounded}.  (One can easily develop an analogous variant of Theorem \ref{theo:bounded_bigandsemiample} using a similar argument.)

\begin{theorem} \label{theo:combinedbounded}
Let $\pi: \mathcal{X} \to B$ be a good fibration and let $L$ be a generically relatively ample $\mathbb{Q}$-Cartier divisor.  Fix a constant $\beta$.
\begin{enumerate}
\item There is a proper closed subset $\mathcal{R} \subsetneq \mathcal{X}$ such that if $M \subset \Sec(\mathcal{X}/B)$ is an  irreducible component parametrizing a non-dominant family of sections with $(K_{\mathcal{X}/B} + a(\mathcal{X}_{\eta},L|_{\mathcal{X}_{\eta}})L) \cdot C \leq \beta$ then the sections parametrized by $M$ are contained in $\mathcal{R}$.
\item There is a constant $\xi$, a proper closed subset $\mathcal{V} \subset \mathcal{X}$, and a bounded family of smooth projective $B$-varieties $\mathcal{Y}$ equipped with $B$-morphisms $f: \mathcal{Y} \to \mathcal{X}$ satisfying:
\begin{enumerate}
\item $\dim(\mathcal{Y}) < \dim(\mathcal{X})$ and $f$ is generically finite onto its image; 
\item $a(\mathcal Y_\eta, -f^*L|_{\mathcal{Y}_{\eta}}) = a(\mathcal{X}_{\eta},L|_{\mathcal{X}_{\eta}})$ and the Iitaka dimension of $K_{\mathcal{Y}_{\eta}} + f^*L|_{\mathcal{Y}_{\eta}}$ is $0$;
\item if $M \subset \Sec(\mathcal{X}/B)$ is an irreducible component  that generically parametrizes non-HN-free sections $C$ with $L \cdot C \geq \xi$  and $(K_{\mathcal{X}/B} + a(\mathcal{X}_{\eta},L|_{\mathcal{X}_{\eta}})L) \cdot C \leq \beta$ then for a general section $C$ parametrized by $M$ we have either
\begin{enumerate}
\item $C \subset \mathcal{V}$, or
\item for some $f: \mathcal{Y} \to \mathcal{X}$ in our family there is a HN-free section $C'$ of $\mathcal{Y}/B$ such that $C = f(C')$.
\end{enumerate}
\end{enumerate}
\end{enumerate}
\end{theorem}

\begin{proof}
Let $\phi: \mathcal{X}^{+} \to \mathcal{X}$ be a birational morphism and $G$ be a $\pi$-vertical effective $\mathbb{Q}$-Cartier divisor as in Proposition \ref{prop:birmodelposL}.  Choose a positive integer $k$ such that $k\phi^{*}(L+G)$ is Cartier.  We define $L^{+} = k(\phi^{*}L + G)$ and $\beta^{+} = \beta + \tau(\pi \circ \phi,K_{\mathcal{X}^{+}/\mathcal{X}} + ka(\mathcal{X}^{+}_{\eta},L^{+}|_{\mathcal{X}^{+}_{\eta}})G)$.

Choose a $b > a(\mathcal{X}^{+}_{\eta},L^{+}|_{\mathcal{X}^{+}_{\eta}})$ such that $bL^{+}|_{\mathcal{X}^{+}_{\eta}}$ defines a basepoint free linear series.  
We then apply the constructions of the proof of Theorem \ref{theo:positivebounded} to $(\mathcal{X}^{+},L^{+})$ with our chosen constants and with $T = 0$ to obtain a constant $\xi^{\dagger}$, a closed subset $\mathcal{V}^{+} \subset \mathcal{X}^{+}$, and a bounded family of normal projective varieties $q: \widehat{\mathfrak{F}} \to \widehat{\mathfrak{S}}$ equipped with $\widehat{\mathfrak{S}}$-morphisms $p: \widehat{\mathfrak{F}} \to \widehat{\mathfrak{S}} \times B$ and $g: \widehat{\mathfrak{F}} \to \widehat{\mathfrak{S}} \times \mathcal{X}^{+}$.  We define $\mathcal{V}$ to be the union of $\phi(\mathcal{V}^{+})$ with the locus where $\phi^{-1}$ is not defined.  We define $\xi = \frac{1}{k}\xi^{\dagger}$.

Suppose $M \subset \Sec(\mathcal{X}/B)$ parametrizes a family of sections satisfying $L \cdot C \geq \xi$ and $(K_{\mathcal{X}/B} + a(\mathcal{X}_{\eta},L|_{\mathcal{X}_{\eta}})L) \cdot C \leq \beta$.  If the locus swept out by the curves parametrized by $M$ meets the locus where $\phi^{-1}$ is defined, by taking strict transforms we obtain a family of sections $C^{+}$ on $\mathcal{X}^{+}$.  These sections satisfy $L^{+} \cdot C^{+} \geq \xi^{\dagger}$.  Furthermore since $a(\mathcal{X}^{+}_{\eta},L^{+}|_{\mathcal{X}^{+}_{\eta}}) = \frac{1}{k} a(\mathcal{X}_{\eta},L|_{\mathcal{X}_{\eta}})$ we have
\begin{equation*}
(K_{\mathcal{X}^{+}/B} + a(\mathcal{X}^{+}_{\eta},L^{+}|_{\mathcal{X}^{+}_{\eta}})L^{+}) \cdot C^{+} \leq \beta^{+}.
\end{equation*}

First we prove the statement (1). We define $\mathcal{R}$ to be the union of $\mathcal{V}$ with the images of the (finitely many) non-dominant families of sections satisfying $L \cdot C < \xi$ and $(K_{\mathcal{X}/B} + a(\mathcal{X}_{\eta},L|_{\mathcal{X}_{\eta}})L) \cdot C \leq \beta$.  If we have a non-dominant family of sections $C$ such that $L \cdot C \geq \xi$ then it follows from Theorem \ref{theo:positivebounded}.(5) applied to $(\mathcal{X}^{+},L^{+})$ that the sections will be contained in $\mathcal{V}$ and thus in $\mathcal{R}$.  Altogether we see that $\mathcal{R}$ has the desired property.

Next we prove (2).  By Theorem \ref{theo:positivebounded} the bounded family $\widehat{\mathfrak{F}} \to \widehat{\mathfrak{S}}$ equipped with the composition $\widehat{\mathfrak{F}} \to \widehat{\mathfrak{S}} \times \mathcal{X}^{+} \to  \widehat{\mathfrak{S}} \times \mathcal{X}$ satisfies all the properties except possibly Theorem \ref{theo:combinedbounded}.(2).(c).  If $M$ parametrizes a non-dominant family of sections, then as explained above the sections are contained in $\mathcal{V}$.  On the other hand, if $M$ parametrizes a dominant family of non-HN-free sections, then the inclusion $T_{\mathcal{X}^{+}} \to \phi^{*}T_{\mathcal{X}}$ is still injective upon restriction to a general section $C^{+}$.  Thus the family of strict transforms $C^{+}$ is a dominant family of non-HN-free curves on $\mathcal{X}^{+}$.
Theorem \ref{theo:positivebounded} shows that the general section parametrized by $M$ will be the pushforward of an HN-free section on some fiber $\widehat{\mathfrak{F}}_{s}$.
\end{proof}

\section{Fano fibrations} \label{sect:fanofib}

In this section we apply previous results to Fano fibrations.

\subsection{The $\Upsilon$-invariant}

\begin{definition}
\label{defi: invariant_upsilon}
Let $\pi: \mathcal{X} \to B$ be a Fano fibration.  Fix a positive rational number $a$.  By \cite[Theorem 0.2]{KMM92} there is a positive integer $b = b(\dim(\mathcal{X}),a)$ such that $|-bK_{\mathcal X_\eta}|$ is very ample and $b>a$.  Define $\Upsilon_{a}(\pi)$ to be the minimal value of $\tau(\pi,E)$ as we vary over all effective $\pi$-vertical $\mathbb{Q}$-Cartier divisors $E$ constructed as in Proposition \ref{prop:generalsurjectiontofiber} with respect to our choice of $b$ and $a$.  (Note that there is a divisor $E$ achieving this infimum since if $a=\frac{p}{q}$ then each $\tau(\pi,E)$ lies in $\frac{1}{bq}\mathbb{Z}$.)

We also define $\Upsilon(\pi) = \Upsilon_{1}(\pi)$.  
\end{definition}

\begin{remark}
The invariant $\Upsilon(\pi)$ measures the ``failure'' of $\pi$ to be a trivial fibration.
For example, if $\mathcal{X} = X \times B$ for some Fano variety $X$ then we have $\Upsilon(\pi) = 0$. 
\end{remark}

Applying the results of Section \ref{sec:fujinv} we obtain:

\begin{theorem} \label{theo:ainvariantandsections}
Let $\pi: \mathcal{X} \to B$ be a Fano fibration.  Fix a positive rational number $a_{rel}$.  Fix a positive integer $T$.  There is some constant $\xi = \xi(\dim(\mathcal{X}), g(B), \Upsilon_{a_{rel}}(\pi), a_{rel}, T)$ with the following property.

Suppose that $\psi: \mathcal{Y} \to B$ is a good fibration equipped with a $B$-morphism $f: \mathcal{Y} \to \mathcal{X}$ that is generically finite onto its image.   Suppose that $N$ is an irreducible component of $\Sec(\mathcal{Y}/B)$ parametrizing a dominant family of sections $C$ on $\mathcal{Y}$ which satisfy $-f^{*}K_{\mathcal{X}/B} \cdot C \geq \xi$.  Finally, suppose that
\begin{equation*}
\dim(N) \geq a_{rel}(-f^{*}K_{\mathcal{X}/B} \cdot C + (\dim(\mathcal{X})-1)(1-g(B))) - T.
\end{equation*}
Then
\begin{equation*}
a(\mathcal{Y}_{\eta},-f^{*}K_{\mathcal{X}/B}|_{\mathcal{Y}_{\eta}}) \geq a_{rel}.
\end{equation*}
\end{theorem}

\begin{proof}
Set $L = -K_{\mathcal{X}/B}$.  By \cite[Theorem 0.2]{KMM92} there is a positive integer $b > a_{rel}$ depending only on $\dim(\mathcal{X})$ and $a_{rel}$ such that $|-bK_{\mathcal X_\eta}|$ is very ample.  We apply Theorem \ref{theo:generalainvsections} with $\beta = 0$, with $a=a_{rel}$, with our choice of $b$, and with an effective $\pi$-vertical Cartier divisor $E$ as in Definition \ref{defi: invariant_upsilon} that achieves the bound $\tau(\pi,E) = \Upsilon_{a_{rel}}(\pi)$.

The explicit bound (\ref{equation:generalainvsections}) for $\xi$ is a max of two terms.  Most of the summands of the second term also appear in the first with a coefficient of $1/\epsilon$  where $\epsilon$ is a rational number chosen so that no smooth projective variety of dimension $\leq \dim(\mathcal{X})-1$ has a Fujita invariant in the range %(
$[1-\epsilon,1)$ %]
with respect to any big and nef Cartier divisor.  Since $\epsilon < 1$, each of these summands in the second term is no more than the corresponding summand in the first term.  However, there is one more summand in the second term which has a factor of $(\dim(\mathcal{X})-1)(5g(B) + 3 + \gamma)$ where the corresponding term in the first summand has the smaller factor $(\dim(\mathcal{X})-1)(4g(B) + 3 + \gamma)$ in the first.  In all, a term-by-term comparison shows that the max defining the explicit bound is no more than
\begin{align} \label{eq:fanofibrationbound}
\xi = \frac{1}{a_{rel} \epsilon} ((1-\epsilon) \Upsilon_{a_{rel}}(\pi) +&  T + (\dim(\mathcal{X})-1)(5g(B) + 3 + \gamma)) \\
& + \frac{1}{\epsilon} ((\dim(\mathcal{X}) - 1) (g(B)-1) + 2g(B)-2 + \Xi) + 1 \nonumber
\end{align}
where $\gamma = (g(B)\dim(\mathcal{X}) - g(B) + 1)^{2}(\dim(\mathcal{X})-1)$, and
\begin{itemize}
\item $\Xi = 0$, if $g(B) \geq 1$.
\item $\Xi$ is the supremum of the constants obtained by applying Lemma \ref{lemm:terminalsectiontofiber} to all dimensions $\leq \dim(\mathcal{X})$, if $g(B) = 0$.
\end{itemize}
Theorem \ref{theo:generalainvsections} immediately implies the desired conclusion.
\end{proof}

\begin{remark} \label{rema:exceptionalrelativelyfree}
The exceptional set in Geometric Manin's Conjecture as described in \cite{LST18} can include families of relatively free sections as well as families of non-relatively free sections.  For example, sometimes we must discount the contributions of irreducible components $M \subset \Sec(\mathcal{X}/B)$ which parametrize relatively free sections when the evaluation map for the universal family over $M$ has disconnected fibers. Let $f: \mathcal{Y} \to \mathcal{X}$ denote the finite part of the Stein factorization of the evaluation map.  Theorem \ref{theo:ainvariantandsections} shows that $a(\mathcal{Y}_{\eta},-f^{*}K_{\mathcal{X}/B}|_{\mathcal{Y}_{\eta}}) = a(\mathcal{X}_{\eta},-K_{\mathcal{X}/B}|_{\mathcal{X}_{\eta}})$ so that such sections can be accounted for by the exceptional set of \cite{LST18}.
\end{remark}

\subsection{Proofs of main results}

We now prove the theorems stated in the introduction (except for Theorem~\ref{theorem:arithmeticapp} which is postponed to Section \ref{sect:application}).

\begin{proof}[Proof of Theorem \ref{theo:maintheorem1}:]
(1) Let $\mathcal{Y}'$ be a resolution of $\mathcal{Y}$ and let $N$ parametrize the strict transforms on $\mathcal{Y}'$ of the general sections on $\mathcal{Y}$ parametrized by $M$.  Since the Fujita invariant is birationally invariant, the desired statement follows from Theorem \ref{theo:ainvariantandsections} applied to $\mathcal{Y}'$ and $N$ with $a_{rel}=1$ and $T=0$.

(2) To see the equality of Fujita invariants for $\mathcal{Y}$, we let $\mathcal{Y}'$ denote a resolution of singularities of $\mathcal{Y}$ and let $N$ denote the family of sections on $\mathcal{Y}'$ such that the pushforward map $f_* : N \to \Sec(\mathcal X/B)$ defines a dominant morphism $f_* : N \to M$.   The inequality $a(\mathcal{Y}_{\eta},-f^{*}K_{\mathcal{X}/B}) \geq 1$ follows from Theorem \ref{theo:ainvariantandsections} applied to $\mathcal{Y}'$, $N$ with $a_{rel}=1$ and $T=0$ and the birational invariance of the Fujita invariant.  The reverse inequality follows from \cite[Lemma 4.7]{LST18} showing that the Fujita invariant cannot increase when pulling back a big and nef divisor under a generically finite and dominant morphism.

We next construct the rational map $\phi$.  By \cite[Theorem 0.2]{KMM92} there is a positive integer $b > a_{rel}$ depending only on $\dim(\mathcal{X})$ and $a_{rel}$ such that $|-bK_{\mathcal X_\eta}|$ is very ample.  Define $\xi^{+}$ as in Corollary \ref{coro:adjointrigidcodim} applied to $\mathcal{Y}'$, $N$, $L = -K_{\mathcal{X}/B}$, $a_{rel} = 1$, $\beta = 0$, $T=0$, with our choice of $b$, and with an effective $\pi$-vertical Cartier divisor $E$ as in Definition \ref{defi: invariant_upsilon} that achieves the bound $\tau(\pi,E) = \Upsilon(\pi)$.  Then Corollary \ref{coro:adjointrigidcodim} constructs a rational map $\phi: \mathcal{Y} \dashrightarrow \mathcal{Z}$ over $B$ that has all the desired properties.  (Since $\mathcal{Y}$ was constructed as the finite part of the normalization of the evaluation map over $M$, the evaluation map for $N$ has connected fibers, and thus the variety $\mathcal{S}$ constructed in Corollary \ref{coro:adjointrigidcodim} can be taken to be $\mathcal{Y}$ itself.)  The only thing left to check is that $\dim(\mathcal{Z}) \geq 2$, or in other words, that the rational map $\phi$ is not trivial.  Note that we have an inclusion $T_{\mathcal{Y}'/B} \to f^{*}T_{\mathcal{X}/B}$.  Since $C'$ deforms in a dominant family on $\mathcal{Y}'$, this map remains injective upon restriction to a general $C'$ and we conclude that $C'$ is not an HN-free section on $\mathcal{Y}'$.  But then the map $\mathcal{Y}' \to B$ does not satisfy Corollary \ref{coro:adjointrigidcodim}.(3), showing that $\phi$ must be non-trivial.
\end{proof}

\begin{proof}[Proof of Theorem \ref{theo:maintheorem2}:]
This follows from Theorem \ref{theo:combinedbounded} applied with $L = -K_{\mathcal{X}/B}$ and $\beta = 0$.
\end{proof}

We will deduce Theorem \ref{theo:maintheorem3} as a consequence of a slightly more precise theorem that takes dominant maps into account:

\begin{theorem} \label{theo:nonconnectedfiberscase}
Let $\pi: \mathcal{X} \to B$ be a Fano fibration.  There is a linear function $R(d)$ whose leading coefficient is a positive number depending only on $\dim(\mathcal{X})$ such that the following property holds.

Suppose that $M$ is an irreducible component of $\Sec(\mathcal{X}/B)$ parametrizing a family of sections $C$ which satisfy $-K_{\mathcal{X}/B} \cdot C = d$.  Let $W \subset M$ be a locally closed subvariety.  Let $\mathcal{U}_{W}^{\nu}$ denote the normalization of the universal family over $W$ equipped with the evaluation map $ev_{W}: \mathcal{U}_{W}^{\nu} \to \mathcal{X}$.

Then either:
\begin{enumerate}
\item the codimension of $W$ in $M$ is at least $\lfloor R(d) \rfloor$,
\item $ev_{W}$ is dominant with connected fibers, or
\item $ev_{W}$ factors through a generically finite non-birational morphism $f: \mathcal{Y} \to \mathcal{X}$ that satisfies
\begin{equation*}
a(\mathcal{Y}_{\eta},-f^{*}K_{\mathcal{X}/B}|_{\mathcal{Y}_{\eta}}) \geq a(\mathcal{X}_{\eta},-K_{\mathcal{X}/B}|_{\mathcal X_\eta}).
\end{equation*}
%\sho{It should be equal, right?} \brian{You are right.}
\end{enumerate}
\end{theorem}

We emphasize that case (3) of Theorem \ref{theo:nonconnectedfiberscase} must be included: an accumulating morphism can contribute families of bounded codimension in arbitrarily large degrees.

\begin{proof}
Fix a positive integer $T$.  Define the constant $\xi(T)$ (depending also on $\dim(\mathcal{X})$, $g(B)$, and $\Upsilon(\pi)$) as in Theorem \ref{theo:ainvariantandsections} using the constants $a_{rel}=1$, an effective $\pi$-exceptional divisor $E$ as in Definition \ref{defi: invariant_upsilon} such that $\tau(\pi,E) = \Upsilon(\pi)$, and the chosen value of $T$.  The explicit description of $\xi(T)$ in Equation \eqref{eq:fanofibrationbound} shows that $\xi(T)$ is linear in $T$ with leading coefficient $1/\epsilon$ where $\epsilon = \epsilon(\dim(\mathcal{X}))$ is a positive rational number such that no smooth projective variety of dimension $\leq \dim(\mathcal{X})-1$ has Fujita invariant in  %(
$[1-\epsilon,1)$ %]
with respect to a big and nef Cartier divisor.  By inverting the linear function $\xi(T)$ we obtain a linear function $R$.

Suppose that $ev_{W}$ does not satisfy (2).  If $ev_{W}$ is not dominant, let $f: \mathcal{Y} \to \mathcal{X}$ denote the inclusion of the image of $ev_{W}$.  If $ev_{W}$ is dominant with disconnected fibers, let $f: \mathcal{Y} \to \mathcal{X}$ denote the finite part of the Stein factorization of $ev$.  For a fixed value of $T$, if $d \geq \xi(T)$ and the codimension of $W$ is smaller than $T$ then Theorem \ref{theo:ainvariantandsections} shows that 
\begin{equation*}
a(\mathcal{Y}_{\eta},-f^{*}K_{\mathcal{X}/B}|_{\mathcal{Y}_{\eta}}) \geq a(\mathcal{X}_{\eta},-K_{\mathcal{X}/B}|_{\mathcal X_\eta}).
\end{equation*}
Equivalently, either the codimension of $W$ is at least $\lfloor R(d) \rfloor$ or we have the desired inequality of Fujita invariants.
\end{proof}

\begin{proof}[Proof of Theorem \ref{theo:maintheorem3}] 
Since by assumption the sections parametrized by $N$ do not dominate $\mathcal{X}$, case (2) in Theorem \ref{theo:nonconnectedfiberscase} does not hold.  If case (1) does not hold, then $ev_{N}$ factors through a generically finite non-dominant morphism $f: \mathcal{Y}' \to \mathcal{X}$ satisfying
\begin{equation*}
a(\mathcal{Y}'_{\eta},-f^{*}K_{\mathcal{X}/B}|_{\mathcal{Y}'_{\eta}}) \geq a(\mathcal{X}_{\eta},-K_{\mathcal{X}/B}|_{\mathcal X_\eta}).
\end{equation*}
Set $\mathcal{Y} = f(\mathcal{Y}')$.  Since the Fujita invariant of $\mathcal{Y}_{\eta}$ is at least as large as the Fujita invariant of $\mathcal{Y}'_{\eta}$, we deduce the desired statement.
\end{proof}

\begin{proof}[Proof of Theorem \ref{theo:introgm}:]
This is the special case of Theorem \ref{theo:hnfforsections}.(1) when $\mathcal{E}$ is $[C]$-semistable.
\end{proof}

\begin{proof}[Proof of Theorem \ref{theo:introainvariant}:]
This follows from Theorem \ref{theo:ainvariantandsections} applied to an irreducible component $M \subset \Sec(\mathcal{X}/B)$ using $a_{rel} = a$ and the inequality
\begin{equation*}
\dim(M) \geq -K_{\mathcal{X}/B} \cdot C + (\dim(\mathcal{X})-1)(1-g(B)).
\end{equation*}

\end{proof}

\begin{proof}[Proof of Theorem \ref{theo:introtwists}:]
Let $C$ be a general section parametrized by $\widetilde{N}$.  By Corollary \ref{coro:domfamilyexpdim} we have
\begin{align*}
-K_{\widetilde{\mathcal{Y}}/B} \cdot C + (\dim(\widetilde{\mathcal{Y}})-1) & \geq \dim(\widetilde{N}) \\
& \geq \dim(M) - T \\
& \geq -\widetilde{f}^{*}K_{\mathcal{X}/B} \cdot C + (\dim(\mathcal{X})-1)(1-g(B)) - T
\end{align*}
Rearranging we see that
\begin{equation*}
(K_{\widetilde{\mathcal{Y}}/B} -\widetilde{f}^{*}K_{\mathcal{X}/B}) \cdot C \leq T + g(B)(\dim(\mathcal{X})-1).
\end{equation*}
We conclude the desired statement by Corollary \ref{coro:boundedintimpliesboundedtwists}.
\end{proof}

\section{Examples} \label{sect:examples}

Our first example illustrates how Theorem \ref{theo:maintheorem1} can be used in practice to understand sections.

\begin{example}[Cubic hypersurface fibrations] \label{exam:cubichyp}
Suppose that $\pi: \mathcal{X} \to B$ is a Fano fibration whose general fiber is a smooth cubic hypersurface of dimension $n \geq 4$.  We will analyze the irreducible components of $\Sec(\mathcal{X}/B)$ parametrizing non-relatively free sections of large degree.  (In the special case when $X$ is a smooth cubic hypersurface and $B = \mathbb{P}^{1}$, \cite{CS09} proves a stronger statement for $X \times \mathbb{P}^{1}$ by classifying all the irreducible components of $\Mor(\mathbb{P}^{1},X)$.)

One can classify the generically finite morphisms $f: \mathcal{Y}_{\eta} \to \mathcal{X}_{\eta}$ such that $a(\mathcal{Y}_{\eta},-f^{*}K_{\mathcal{X}_{\eta}}) \geq 1$ by combining \cite[1.3 Proposition]{Horing10} with the techniques of \cite[Theorem 11.1]{LTRMS}.  For the sake of the reader we reproduce the argument here.  We claim that:
\begin{itemize}
\item If $n \geq 5$ then there are no non-birational generically finite morphisms $f: \mathcal{Y}_{\eta} \to \mathcal{X}_{\eta}$ with $a(\mathcal{Y}_{\eta},-f^{*}K_{\mathcal{X}_{\eta}}) \geq 1$.
\item If $n=4$ then there are no non-birational generically finite morphisms $f: \mathcal{Y}_{\eta} \to \mathcal{X}_{\eta}$ with $a(\mathcal{Y}_{\eta},-f^{*}K_{\mathcal{X}_{\eta}}) \geq 1$ unless $\mathcal{X}_{\eta}$ contains a plane.  When $\mathcal{X}_{\eta}$ contains a plane, the only possibility is that $f$ is the composition of a birational map $\phi: \mathcal{Y}_{\eta} \to \mathbb{P}^{2}_{\eta}$ and the inclusion of a plane $\mathbb{P}^{2}_{\eta} \subset \mathcal{X}_{\eta}$.
\end{itemize}
First suppose that $\mathcal{Y}_{\eta} \subsetneq \mathcal{X}_{\eta}$ is a proper subvariety with $a(\mathcal{Y}_{\eta},-f^{*}K_{\mathcal{X}_{\eta}}) \geq 1$.  Letting $H$ denote the hyperplane class on $\mathbb{P}^{n+1}_{\eta}$, we have $a(\mathcal{X}_{\eta},H) = n-1$.  By \cite[Theorem 3.16]{LTRMS} the Fujita invariant with respect to a big and nef Cartier divisor is always at most one more than the dimension.  We conclude that $\dim(\mathcal{Y}_{\eta}) \geq n-2$.  

If $\dim(\mathcal{Y}_{\eta}) = n-1$, then by \cite[1.3 Proposition]{Horing10} we see that $(\mathcal{Y}_{\eta},H|_{\mathcal{Y}_{\eta}})$ is birationally equivalent to either $(\mathbb{P}_{\eta}^{n-1},\mathcal{O}(1))$, $(Q_{\eta},\mathcal{O}(1)|_{Q_{\eta}})$, or a $\mathbb{P}_{\eta}^{n-2}$-bundle such that $H$ has degree $1$ along the fibers.  In the first and second cases, $\mathcal{Y}$ must either be a linear subvariety or a degree $2$ subvariety in $\mathcal{X}_{\eta}$ with codimension $1$, an impossibility by the Lefschetz hyperplane theorem.  In the third case, the $\mathbb{P}_{\eta}^{n-2}$-fibers must define a $1$-dimensional family of codimension $2$ linear subvarieties of $\mathcal{X}_{\eta}$, again an impossibility.

If $\dim(\mathcal{Y}_{\eta}) = n-2$, then by \cite[1.3 Proposition]{Horing10} we see that $(\mathcal{Y}_{\eta},H|_{\mathcal{Y}_{\eta}})$ is birationally equivalent to $(\mathbb{P}^{n-2},\mathcal{O}(1))$.  In particular $H^{n-2} \cdot \mathcal{Y}_{\eta} = 1$, showing that $\mathcal{Y}_{\eta}$ must be a codimension $2$ linear subvariety contained in $\mathcal{X}_{\eta}$.  This is impossible if $n \geq 5$ and accounts for the unique exception if $n = 4$.

Next suppose that $f: \mathcal{Y}_{\eta} \to \mathcal{X}_{\eta}$ is a generically finite morphism.  \cite[Lemma 4.7]{LST18} shows that the Fujita invariant of the image of $f$ is at least as large as the Fujita invariant of $\mathcal{Y}_{\eta}$.  In particular, if $f$ is not dominant then $n=4$ and the image of $f$ must be $\mathbb{P}_{\eta}^{2}$.  Furthermore, any generically finite dominant non-birational morphism $f: \mathcal{Y}_{\eta} \to \mathbb{P}^{2}_{\eta}$ satisfies $a(\mathcal{Y}_{\eta},-f^{*}K_{\mathcal{X}_{\eta}}) < a(\mathbb{P}^{2}_{\eta},-f^{*}K_{\mathcal{X}_{\eta}}) = 1$ by \cite[Corollary 9.2]{LT17}.

It only remains to classify the dominant generically finite morphisms $f: \mathcal{Y}_{\eta} \to \mathcal{X}_{\eta}$ such that $a(\mathcal{Y}_{\eta},-f^{*}K_{\mathcal{X}_{\eta}}) \geq 1$.  In this case \cite[Lemma 4.7]{LST18} shows we must have Fujita invariant exactly $1$.  We separate into two cases based on the Iitaka dimension of $K_{\mathcal{Y}_{\eta}} - f^{*}K_{\mathcal{X}_{\eta}}$.  If the Iitaka dimension is $>0$, then by \cite[Lemma 4.9]{LST18} we see that $\mathcal{X}_{\eta}$ is covered by a dominant family of proper subvarieties with Fujita invariant $1$, an impossibility by our earlier discussion.  If the Iitaka dimension is $=0$, by \cite[Corollary 2.8]{Sengupta21} any irreducible component $B_{\eta}$ of the branch locus of $f$ will satisfy $a(B_{\eta},-K_{\mathcal{X}_{\eta}}) > 1$.  However, we have already verified that $\mathcal{X}_{\eta}$ has no codimension $1$ subvarieties with this property.  By purity of the branch locus, we deduce that $f$ must be \'etale, and thus birational.

Returning to the discussion of sections, let $M$ be a component of $\Sec(\mathcal{X}/B)$ of sufficiently large anticanonical degree.  Then Theorem \ref{theo:maintheorem1} shows:
\begin{enumerate}
\item If $n \geq 5$ then $M$ will generically parametrize relatively free sections and the evaluation map for its universal family will have connected fibers as in Remark \ref{rema:exceptionalrelativelyfree}.
\item If $n=4$, then $M$ can only fail to generically parametrize relatively free sections if it parametrizes a family of sections whose intersection with $\mathcal{X}_{\eta}$ is contained in a plane.
\end{enumerate}
This finishes the classification of irreducible components parametrizing non-free curves of large degree.
\end{example}

Our second example addresses the non-generically-globally-generated locus.  Let $\pi: \mathcal{X} \to B$ be a Fano fibration and let $M$ be an irreducible component of $\Sec(\mathcal{X}/B)$ of large degree.  Theorem \ref{theo:maintheorem3} shows that the codimension in $M$ of the non-generically-globally-generated locus will grow linearly in degree except possibly when the sections sweep out a subvariety with large Fujita invariant.  The following example demonstrates that it is possible for the non-generically-globally-generated locus to have constant codimension.

\begin{example} \label{exam:largenonfreelocus}
Let $X$ be a smooth cubic threefold.  Suppose $M$ is a component of $\Mor(\mathbb{P}^{1},X)$ parametrizing maps of anticanonical degree $\geq 3$.  We will see that $M$ admits a (possibly reducible) codimension $1$ sublocus parametrizing multiple covers of non-free lines.  In particular, the non-free locus in $M$ will always have codimension $1$.

\cite{CS09} shows that for any degree $d \geq 2$ the moduli stack $\overline{\mathcal{M}}_{0,0}(X,d)$ has two irreducible components: a component $M_{d}$ that generically parametrizes irreducible free curves and a component $R_{d}$ that parametrizes degree $d$ covers of lines.  We will be interested in the intersection $T_{d}$ of these two components.

Inside of the parameter space of lines on $X$ the sublocus parametrizing non-free lines has codimension $1$.  Thus the locus $Q_{d} \subset R_{d}$ parametrizing multiple covers of non-free lines also has codimension $1$.  Note that $T_{d}$ will be contained in $Q_{d}$.  On the other hand, since all components of the moduli stack $\overline{\mathcal{M}}_{0,0}(X,d)$ have the expected dimension $\overline{\mathcal{M}}_{0,0}(X,d)$ has only LCI singularities.  Thus $T_{d}$ must have codimension $1$.  Altogether, we see that $T_{d}$ consists of a (non-empty) union of irreducible components of $Q_{d}$.

Since the general stable map parametrized by $T_{d}$ has irreducible domain, we see that every irreducible component of $\Mor(\mathbb{P}^{1},X)$ of degree $\geq 3$ will have a codimension $1$ sublocus parametrizing non-free morphisms consisting of  multiple covers of non-free lines.

This result illustrates Theorem \ref{theo:maintheorem3} applied to the projection $\pi: X \times \mathbb{P}^{1} \to \mathbb{P}^{1}$.  Let $Y \subset X$ denote a subvariety swept out by the curves parametrized by an irreducible component of $T_{d}$.  Then $Y$ is also swept out by a one-parameter family of non-free lines; in particular we have $a(Y,-K_{X}|_{Y}) = 1$.  Passing to the relative situation, we see that any codimension $1$ locus of $M_{d}$ parametrizing sections whose normal bundle is not generically globally generated will sweep out a subvariety $Y \times \mathbb{P}^{1}$ which has generic Fujita invariant $\geq 1$.
\end{example}

\section{An arithmetic application}
\label{sect:application}

In this section we prove Theorem~\ref{theorem:arithmeticapp}.
We freely use the notations set up in Section~\ref{subsec:arithmeticapp}.

First of all, consider the open subscheme of the relative Hilbert scheme that parametrizes sections of anticanonical height $\leq d$ which are not contained in $\mathcal R$.  There exists a finite set of places $S_d \supset S$ such that this open subset is flat over $\Spec \mathfrak o_{F, S_d}$. Let $\psi: \mathcal H _d\to \Spec \mathfrak o_{F, S_d}$ be the closure of this set inside the relative Hilbert scheme equipped with the reduced structure.  Then $\psi: \mathcal H_d \to \Spec \mathfrak o_{F, S_d}$ is projective and flat.
 
 We denote the generic fiber of $\psi$ by $H_d$. By Theorem \ref{theo:maintheorem2}.(1) every irreducible component of $H_d$ parametrizes a dominant family of sections.  Corollary~\ref{coro:domfamilyexpdim} shows that the dimension of such a component is bounded by $d + \dim X_\eta$.  Define 
 \[
C_{d} := \sum_{i = 0}^{2(d+\dim X_\eta)} h_{\text{sing}}^{i}(H_d^{an}, \mathbb Q).
 \]
 Since the $\ell$-adic sheaf $R^i\psi_*\underline{\mathbb Q}_\ell$ is constructible in the pro-\'etale topology by \cite[Lemma 6.7.2]{BS15}, 
 there exists a pro-\'etale open 
 $i : U \to \Spec \mathfrak o_{F, S_d}$ such that $i^{-1}R^i\psi_*\underline{\mathbb Q}_\ell$ is a constant sheaf in the pro-\'etale topology.
In particular  there exists a finite set of places $S'_d \supset S_d$ such that we have
 \[
 h^{i}_{\text{\'et}}(H_{d, \overline{v}}, \mathbb Q_\ell) = h^i_{\text{sing}}(H_{d}^{an}, \mathbb Q)
 \]
 for all $i$ and $v\not\in S'_d$ where $H_{d, \overline{v}}$ is the base change of $H_{d, v}$ to the algebraic closure. By applying the Grothendieck-Lefschetz trace formula and a version of the Weil conjectures for singular projective varieties (\cite{Deligne}), we conclude that
 \[
 N(X_v\setminus R_v, -K_{X_v/B_v}, d) \leq \#H_{d, v}(k_v) \leq C_dq_v^{d + \dim X_\eta} 
 \]
 Thus assuming $d\epsilon > \dim X_\eta$, we can conclude that
 \[
 \frac{ N(X_v\setminus R_v, -K_{X_v/B_v}, d)}{q_v^{d(1+\epsilon)}} \to 0
 \]
 as $v \to \infty$.

 %\nocite{*}
\bibliographystyle{alpha}
\bibliography{nonfreecurves}

\end{document}